## UNIVERSITÉ PARIS 7 – DENIS-DIDEROT
U.F.R. de Mathématiques

# Thèse

Présentée pour obtenir le diplôme de
Docteur de l'Université Paris 7
Spécialité : Mathématiques

Présentée et soutenue publiquement le vendredi 9 décembre 2005 par

## François Brunault

---

**Étude de la valeur en $s = 2$ de la fonction $L$ d'une courbe elliptique**

---

**Directeur :**
Prof. Loïc Merel

**Rapporteurs :**
Prof. Christopher Deninger
Prof. Guido Kings

**Jury :**
Prof. Guido Kings
Prof. Loïc Merel
Prof. Jan Nekovář
Prof. Joseph Oesterlé
Prof. Jörg Wildeshaus
Prof. Don Zagier

# Remerciements

Je tiens en premier lieu à exprimer ma reconnaissance à Loïc Merel pour avoir accepté de diriger cette thèse. Lors du stage de D.E.A., il m'a proposé l'étude d'un article riche, dense et intéressant [32], qui m'a ouvert de nombreuses voies de recherche. J'ai largement bénéficié de la grande culture mathématique de Loïc Merel et de sa passion pour l'enseignement. J'ai également pu compter sur son soutien constant et sa grande disponibilité. Je le remercie profondément.

Je suis très reconnaissant envers Guido Kings pour m'avoir accueilli en décembre dernier à l'Université de Ratisbonne. Les nombreuses discussions que j'ai eues avec lui m'ont permis d'envisager les polylogarithmes de manière plus conceptuelle et donneront lieu, je l'espère, à des travaux intéressants. L'ambiance chaleureuse au sein de l'équipe Géométrie arithmétique a sans aucun doute également contribué à l'aboutissement de ce travail. Regensburg est effectivement une belle ville ! Je garderai un excellent souvenir de ce séjour pré-doctoral.

Je remercie Christopher Deninger et Guido Kings d'avoir accepté le lourd travail de rapporter ma thèse. Je remercie tout aussi chaleureusement Jan Nekovář, Joseph Oesterlé, Jörg Wildeshaus et Don Zagier de me faire l'honneur de faire partie de mon jury. Je suis très touché que Jörg Wildeshaus puisse finalement venir à ma soutenance.

Durant ces quatre années de thèse, de nombreux professeurs ont pris le temps de répondre à mes questions mathématiques : je tiens à remercier en particulier Yves André, Dominique Bernardi, Marie-José Bertin, Daniel Bertrand, Pierre Colmez, Marc Hindry, Bruno Kahn, Odile Lecacheux, Vincent Maillot, Jan Nekovář, Joseph Oesterlé et Jörg Wildeshaus. J'exprime toute ma gratitude envers Don Zagier pour avoir répondu à mes questions lors de l'Arbeitstagung de Bonn en 2003 et pour m'avoir communiqué un article qui n'était pas encore publié [27].

Mes remerciements vont également à Jean-Benoît Bost, responsable du D.E.A. de mathématiques pures de l'Université d'Orsay en 2000/2001. Son cours d'introduction à la géométrie algébrique et son article d'introduction aux surfaces de Riemann [15] ont été essentiels dans ma formation mathématique. Plus largement, je remercie tous mes professeurs de mathématiques pour leur enseignement.

À l'Institut de Mathématiques de Jussieu, sur le site de Chevaleret, j'ai bénéficié de conditions de travail excellentes pour préparer cette thèse, grâce à une allocation couplée du Ministère de l'Éducation nationale, de l'Enseignement supérieur et de la Recherche. Je remercie également le LAGA de l'Université Paris 13, pour m'avoir accueilli pendant un an en tant qu'A.T.E.R. Enseigner là-bas fut pour moi une expérience à la fois riche et agréable.

J'adresse un grand merci aux différents secrétariats pour leur dévouement et leur gentillesse, notamment à Mme Monchanin, Mme Orion, Mme Prosper-Cojande et Mme Wasse ; à Mme Larive et aux bibliothécaires ; à Joël Marchand et à toute l'équipe informatique.

Je souhaite remercier tous les thésards (actuels et anciens) pour leur soutien et la bonne ambiance sur le plateau 7C : Aïcha, Aïni, Anna, Aurélien, Benoît, Cécile, Charlotte, Christophe, Fabien, François M., Gwendal, Joël, Julien M., Julien P., Lionel, Luc, Luca, Magda, Marusia, Michel, Miguel, Najib, Nicolas, Olivier, Pierre, Samy. . . Merci tout spécialement à David Blottière. Je remercie aussi l'équipe de foot de Chevaleret : Amadeo, Cecilia, Claire, Ernesto, Esther, Gio-vanni, Maria, Maria-Paula, Paulo, c'est un plaisir de jouer avec vous ! Je remercie enfin Mostafa pour les nombreuses discussions, les proverbes et le café !

Je n'oublie pas non plus mes camarades de la promotion X97, dont l'amitié m'est chère : en particulier Vianney, Jean-Dominique, Sélim, Thibaud, Guillaume et Matthieu.

J'exprime ma haute estime pour Robert Pagès, Directeur de recherche au C.N.R.S., et je le remercie chaleureusement pour les échanges fructueux.

Je remercie tous mes amis pour leur soutien et la réalité de leur amitié. Je tiens particulièrement à mentionner les montreuillois, Laurent, Virginie et Céline, pour leur bonne humeur quotidienne et les discussions sur l'histoire, le chant et le théâtre ! Merci à tous ceux qui sont venus à la soutenance, et merci à tous les autres !

Pour finir, je remercie toute ma famille, présente en nombre à la soutenance, mes frères, et tout particulièrement mes parents, qui ont toujours été présents à mes côtés. C'est sans compter que je les remercie, et je leur dédie ce travail.

# Table des matières

# Introduction

Il est devenu classique d'étudier les objets de nature arithmétique en leur associant une fonction $L$, définie au moyen d'un produit eulérien encodant leurs propriétés "locales". On espère alors relier la fonction $L$ à des quantités "globales", c'est-à-dire géométriques. L'exemple le plus simple, et le premier historiquement, est donné par la *fonction zêta de Riemann* [9], qui permet une étude fine des nombres premiers. Considérons maintenant une courbe elliptique $E$ définie sur $\mathbf{Q}$, donnée par une équation de Weierstraß

$$E : y^2 = x^3 + ax + b \qquad (a, b \in \mathbf{Z}).$$

Hasse et Weil ont défini la *fonction L associée à E* de la manière suivante. Pour tout nombre premier $p$, notons $E_p$ la réduction modulo $p$ d'une équation de Weierstraß minimale en $p$ pour $E$ [70, §VII.1]. Posons $a_p = p + 1 - \operatorname{Card} E_p(\mathbf{F}_p)$ et

$$L_p(E, X) = \begin{cases} (1 - a_p X + p X^2)^{-1} & \text{si } E_p \text{ est une courbe elliptique sur } \mathbf{F}_p, \\ (1 - a_p X)^{-1} & \text{sinon.} \end{cases}$$

La fonction $L$ associée à $E$ est alors définie par

$$L(E, s) = \prod_{p \text{ premier}} L_p(E, p^{-s}) \qquad \left(\Re(s) > \frac{3}{2}\right).$$

La convergence du produit infini ci-dessus pour $\Re(s) > \frac{3}{2}$ découle de la majoration $|a_p| \leq 2\sqrt{p}$, démontrée par Hasse dans les années 1930. La fonction $L(E, s)$ s'écrit sous forme de série de Dirichlet

$$L(E, s) = \sum_{n=1}^{\infty} \frac{a_n}{n^s} \qquad \left(\Re(s) > \frac{3}{2}\right).$$

Nous étudions dans cette thèse la valeur spéciale $L(E, 2)$. Pour tout caractère de Dirichlet $\chi : (\frac{\mathbf{Z}}{m\mathbf{Z}})^* \to \mathbf{C}^*$, la *série L de E tordue par $\chi$* est définie par

$$L(E, \chi, s) = \sum_{n=1}^{\infty} \frac{a_n \chi(n)}{n^s} \qquad \left(\Re(s) > \frac{3}{2}\right),$$

où par convention $\chi(n) = 0$ lorsque $(n, m) > 1$. Nous savons maintenant, grâce à un théorème profond dû à Breuil, Conrad, Diamond, Taylor et Wiles [78, 74, 17], que les fonctions $L(E, s)$ et $L(E, \chi, s)$ admettent des prolongements holomorphes au plan complexe. Notons $N$ le conducteur de $E$ [71, §IV.10-11], et $w(E)$ l'opposé du signe de l'équation fonctionnelle satisfaite par $L(E, s)$.

**0.1. Détermination de $L(E,2)$ par les $L(E,\chi,1)$, avec $\chi$ caractère de Dirichlet de niveau divisant $N$**

Soit $\mathcal{H} = \{z \in \mathbf{C},\ \Im(z) > 0\}$ le demi-plan de Poincaré, muni de l'action homographique du groupe $\mathrm{SL}_2(\mathbf{Z})$. Pour tout couple $(u,v) \in (\frac{\mathbf{Z}}{N\mathbf{Z}})^2$, posons

$$E^*_{u,v}(z) = \lim_{\substack{s \to 1 \\ \Re(s) > 1}} \Big( \sum_{\substack{m \equiv u \ (N) \\ n \equiv v \ (N) \\ (m,n) \neq (0,0)}} \frac{\Im(z)^s}{|mz+n|^{2s}} - \frac{\pi}{N^2(s-1)} \Big) \qquad (z \in \mathcal{H}). \tag{1}$$

Pour tous $a,b \in \frac{\mathbf{Z}}{N\mathbf{Z}}$, définissons une forme différentielle $\eta(a,b)$ sur $\mathcal{H}$ par

$$\eta(a,b) = E^*_{0,a} \cdot (\partial - \overline{\partial}) E^*_{0,b} - E^*_{0,b} \cdot (\partial - \overline{\partial}) E^*_{0,a}. \tag{2}$$

Pour $\chi : (\frac{\mathbf{Z}}{N\mathbf{Z}})^* \to \mathbf{C}^*$ caractère de Dirichlet modulo $N$, posons

$$\eta_\chi = \sum_{a \in (\frac{\mathbf{Z}}{N\mathbf{Z}})^*} \sum_{b \in (\frac{\mathbf{Z}}{N\mathbf{Z}})^*} \chi(a)\, \overline{\chi}(b)\, \eta(a,b). \tag{3}$$

La *somme de Gauß* de $\chi$ est définie par $\tau(\chi) = \sum_{a=0}^{N-1} \chi(a)\, e^{2\pi i a/N}$. Pour $z_0, z_1 \in \mathcal{H}$, on note $\int_{z_0}^{z_1}$ l'intégrale le long d'une géodésique de $\mathcal{H}$ reliant $z_0$ à $z_1$. Pour $v \in \mathbf{Z}$, posons $g_v = \begin{pmatrix} 0 & -1 \\ 1 & v \end{pmatrix} \in \mathrm{SL}_2(\mathbf{Z})$. Enfin, notons $\rho = e^{\pi i/3} \in \mathcal{H}$.

**Théorème 1.** *Supposons $N = p$ premier. Pour tout caractère de Dirichlet $\chi$ modulo $p$, pair et non trivial, nous avons la formule*

$$L(E,2) L(E,\chi,1) = \frac{p\, w(E)\, \tau(\chi)}{8\pi i (p-1)} \sum_{\chi'} c_{\chi,\chi'} L(E,\chi',1) \tag{4}$$

*où la somme est étendue aux caractères $\chi'$ modulo $p$, pairs et non triviaux, et les coefficients $c_{\chi,\chi'}$ sont donnés par*

$$c_{\chi,\chi'} = \tau(\overline{\chi'}) \sum_{v=1}^{p-1} \chi'(v) \int_{g_v \rho}^{g_v \rho^2} \eta_\chi. \tag{5}$$

*Remarques.*   1. La démonstration du théorème 1 s'inspire de la méthode de Rankin-Selberg, déjà utilisée par Beĭlinson. Elle utilise donc la modularité des courbes elliptiques sur $\mathbf{Q}$.
  2. Les coefficients $c_{\chi,\chi'}$ ne dépendent que de $\chi$, $\chi'$ et $p$ (et pas de $E$).
  3. L'hypothèse $N$ premier n'est pas essentielle. Pour $N$ quelconque, on peut obtenir une formule analogue à (4) en utilisant le théorème A de l'appendice (p. 144).
  4. De même, on peut remplacer $E$ par une forme parabolique primitive de poids 2.
  5. Nous démontrons qu'il existe un caractère $\chi$ modulo $p$, pair et non trivial, tel que $L(E,\chi,1)$ est non nul (cf. démonstration du lemme 99). La valeur spéciale $L(E,2)$ se déduit donc des $L(E,\chi,1)$, où $\chi$ parcourt les caractères pairs et non triviaux modulo $p$.

Voici une méthode alternative permettant de déduire $L(E,2)$ de la formule précédente. Considérons la série de Dirichlet

$$L(E \otimes E, s) = \sum_{n=1}^{\infty} \frac{a_n^2}{n^s} \qquad (\Re(s) > 2).$$

Elle admet un prolongement méromorphe au plan complexe [19] dont le résidu en $s=2$ est non nul (cf. appendice §5, p. 154).

**Théorème 2.** *Supposons $N = p$ premier. Nous avons la formule*

$$L(E,2) = \frac{p^3\,w(E)}{8(p+1)(p-1)^3\,\pi^2} \frac{\sum_{\chi,\chi'} \lambda_{\chi,\chi'} L(E,\chi,1) L(E,\chi',1)}{\operatorname{Res}_{s=2} L(E\otimes E, s)} \qquad (6)$$

*où la somme est étendue aux caractères $\chi$ (resp. $\chi'$) modulo $p$, pairs et non triviaux (resp. impairs), et les coefficients $\lambda_{\chi,\chi'}$ sont donnés par*

$$\lambda_{\chi,\chi'} = \sum_{\chi''} \frac{\tau(\chi'')}{\tau(\chi'\chi'')} c_{\chi'',\chi} \qquad (7)$$

*la dernière somme portant sur les caractères $\chi''$ pairs non triviaux modulo $p$, les nombres $c_{\chi'',\chi}$ étant donnés par (5).*

*Remarque.* D'après le théorème D de l'appendice (p. 154), le résidu en $s=2$ de $L(E\otimes E, s)$ est combinaison bilinéaire des $L(E,\chi,1)L(E,\chi',1)$. La valeur spéciale $L(E,2)$ se déduit donc mécaniquement de $w(E)$ et des $L(E,\chi,1)$, où $\chi$ parcourt les caractères de Dirichlet non triviaux modulo $p$. Explicitement

$$L(E,2) = \frac{p^2\,i\,w(E)}{8(p-1)\,\pi^2} \frac{\sum_{\chi,\chi'} \lambda_{\chi,\chi'} L(E,\chi,1) L(E,\chi',1)}{\sum_{\chi,\chi'} \tau(\overline{\chi\chi'}) L(E,\chi,1) L(E,\chi',1)} \qquad (8)$$

où les sommes portent sur les caractères $\chi$ (resp. $\chi'$) modulo $p$, pairs et non triviaux (resp. impairs).

## 0.2. Généralisation au cadre modulaire

Le théorème 1 admet un analogue naturel dans le cadre modulaire. Soit donc $f$ une forme parabolique *primitive* (i. e. propre pour l'algèbre de Hecke, nouvelle et normalisée) de poids 2 pour $\Gamma_1(N)$, de caractère $\psi$.

La *transformée de Fourier* d'un caractère de Dirichlet $\chi : (\frac{\mathbf{Z}}{N\mathbf{Z}})^* \to \mathbf{C}^*$ est définie par $\widehat{\chi}(b) = \sum_{a=0}^{N-1} \chi(a) e^{-2\pi i ab/N}$, $b \in \frac{\mathbf{Z}}{N\mathbf{Z}}$. Posons

$$\eta(1, \widehat{\chi}) = \sum_{b \in \frac{\mathbf{Z}}{N\mathbf{Z}}} \widehat{\chi}(b)\,\eta(1,b). \qquad (9)$$

Notons $E_N$ l'ensemble des éléments d'ordre $N$ du groupe additif $(\frac{\mathbf{Z}}{N\mathbf{Z}})^2$. Pour tout $x \in E_N$, choisissons une matrice $g_x = \begin{pmatrix} a & b \\ c & d \end{pmatrix} \in \mathrm{SL}_2(\mathbf{Z})$ telle que $(c,d) \in x$. Les formes différentielles (2) étant invariantes sous l'action du groupe $\Gamma_1(N)$, on vérifie que l'intégrale $\int_{g_x\rho}^{g_x\rho^2} \eta(1, \widehat{\chi})$ ne dépend que de $x$. Pour tout $x \in E_N$, le *symbole de Manin* $\xi_f(x)$ est défini par

$$\xi_f(x) = -i \int_{g_x 0}^{g_x \infty} f(z)dz. \qquad (10)$$

Pour $x = (u,v) \in E_N$, posons $x^c = (-u,v)$ et $\xi_f^\pm(x) = \frac{1}{2}\bigl(\xi_f(x) \pm \xi_f(x^c)\bigr)$. Le théorème 1 se déduit du théorème 3 suivant.

**Théorème 3.** *Pour tout caractère de Dirichlet $\chi$ modulo $N$, pair et distinct de $\overline{\psi}$, nous avons*

$$L(f,2)L(f,\chi,1) = \frac{N\,i}{4} \sum_{x \in E_N} \Big(\int_{g_x\rho}^{g_x\rho^2} \eta(1,\widehat{\chi})\Big)\xi_f^+(x). \tag{11}$$

*Remarque.* Le 1-cycle $\sum_{x \in E_N} \xi_f^+(x)\,\{g_x\rho, g_x\rho^2\}$ est fermé (lemme 111).

Appliquons le théorème 3 à la forme primitive $f$ associée à $E$. Nous montrons qu'il existe un caractère $\chi$ modulo $N$, pair et non trivial, tel que $L(f,\chi,1)$ est non nul (cf. démonstration du lemme 99). Nous en tirons le résultat suivant.

**Corollaire.** *Pour toute courbe elliptique $E$ sur $\mathbf{Q}$, de conducteur $N$, la valeur spéciale $L(E,2)$ est combinaison linéaire à coefficients rationnels des quantités $\frac{1}{\pi i}\int_{g_x\rho}^{g_x\rho^2} \eta(1,b)$, où $x$ et $b$ parcourent respectivement $E_N$ et $\frac{\mathbf{Z}}{N\mathbf{Z}}$.*

### 0.3. La théorie de Beĭlinson

Rappelons les résultats de Beĭlinson [5, 6]. Notons $X_1(N)$ la courbe modulaire sur $\mathbf{Q}$ associée au sous-groupe de congruence

$$\Gamma_1(N) = \Big\{\begin{pmatrix} a & b \\ c & d \end{pmatrix} \in \mathrm{SL}_2(\mathbf{Z}),\ \begin{pmatrix} a & b \\ c & d \end{pmatrix} \equiv \begin{pmatrix} 1 & * \\ 0 & 1 \end{pmatrix} \pmod{N}\Big\}.$$

Soit $f$ une forme parabolique primitive de poids 2 pour $\Gamma_1(N)$. Beĭlinson donne une formule pour $L(f,2)L(f,\chi,1)$ où $\chi$ est un caractère de Dirichlet pair quelconque, en termes du régulateur d'un élément du groupe de $K$-théorie algébrique $K_2^{(2)}(X_1(N'))$, où $N'$ est un entier divisible par $N$ (essentiellement $N' = N \cdot m^2$ si $\chi$ est de niveau $m$). Il choisit $\chi$ de telle sorte que $L(f,\chi,1) \neq 0$. Lorsque $f$ est associée à $E$, il suit une formule pour $L(E,2)L(E,\chi,1)$ en termes du régulateur de l'élément de $K_2^{(2)}(E)$ obtenu grâce à la trace $K_2^{(2)}(X_1(N')) \to K_2^{(2)}(X_1(N)) \to K_2^{(2)}(E)$.

La méthode de Beĭlinson souffre de deux imprécisions :
– le choix de du caractère $\chi$, et donc de l'entier $N'$, n'est pas précisé ;
– la paramétrisation modulaire $X_1(N) \to E$ n'est pas explicite. Il en va donc de même de la trace $K_2^{(2)}(X_1(N)) \to K_2^{(2)}(E)$.

Nous reprenons la méthode de Beĭlinson et montrons qu'il est possible de choisir $N' = N$ et $\chi$ parmi les caractères de niveau divisant $N$. Notons $\mathbf{Q}(X_1(N))$ le corps des fonctions de $X_1(N)$. La localisation (3.86) identifie le groupe de $K$-théorie algébrique $K_2(X_1(N)) \otimes \mathbf{Q}$ à un sous-espace du groupe de $K$-théorie de Milnor $K_2(\mathbf{Q}(X_1(N))) \otimes \mathbf{Q}$. Par définition, le *régulateur de Beĭlinson* est la restriction à ce sous-espace de l'application

$$r_N : K_2(\mathbf{Q}(X_1(N))) \otimes \mathbf{Q} \to \mathrm{Hom}_{\mathbf{C}}(S_2(\Gamma_1(N)), \mathbf{C}) \tag{12}$$
$$\{u,v\} \mapsto \Big(f \mapsto \int_{X_1(N)(\mathbf{C})} \log|u| \cdot \omega_f \wedge \overline{\partial}\log|v|\Big),$$

où nous avons noté $\omega_f = 2\pi i f(z)dz$. Le régulateur s'étend par $\mathbf{C}$-linéarité à $K_2(X_1(N)) \otimes \mathbf{C}$. Pour tout caractère de Dirichlet $\chi$ modulo $N$, pair et non trivial, nous définissons une unité modulaire $u_\chi \in \mathcal{O}^*(Y_1(N)) \otimes \mathbf{C}$ vérifiant

$$\log|u_\chi| = \frac{1}{\pi} \sum_{a \in (\frac{\mathbf{Z}}{N\mathbf{Z}})^*} \chi(a)\, E_{0,a}^*. \tag{13}$$

**Proposition.** *Pour tous caractères $\chi$, $\chi'$ modulo $N$, pairs et non triviaux, l'élément $\{u_\chi, u_{\chi'}\}$ appartient à $K_2(X_1(N)) \otimes \mathbf{C} \subset K_2(\mathbf{Q}(X_1(N))) \otimes \mathbf{C}$.*

En explicitant la méthode de Beĭlinson, nous sommes parvenus au résultat suivant. Notons $\varphi$ la fonction indicatrice d'Euler.

**Théorème 4.** *Soit $f$ une forme parabolique primitive de poids $2$ pour $\Gamma_1(N)$, de caractère $\psi$. Pour tout caractère de Dirichlet $\chi$ modulo $N$, pair, distinct de $\overline{\psi}$ et primitif, nous avons*

$$L(f,2)L(f,\chi,1) = \frac{N\pi i}{\varphi(N)}\tau(\chi)\left\langle r_N\big(\{u_{\psi\chi}, u_{\overline{\chi}}\}\big), f\right\rangle. \tag{14}$$

*Remarques.*   1. L'hypothèse $\chi$ primitif n'est pas essentielle. On a une formule analogue sans cette hypothèse. De même, on peut prendre $\chi$ de niveau $m$ divisant $N$ (avec $m \neq 1$).

2. D'après les résultats de l'appendice (corollaire 2, p. 145), il existe un caractère pair $\chi$ modulo $N$ tel que $L(f,\chi,1) \neq 0$.

### 0.4. Une question de Schappacher et Scholl

Notons $\{\mathcal{O}^*(Y_1(N)), \mathcal{O}^*(Y_1(N))\}$ le sous-groupe de $K_2(\mathbf{Q}(X_1(N)))$ engendré par les symboles $\{u,v\}$, $u,v \in \mathcal{O}^*(Y_1(N))$, et posons

$$K_N = \{\mathcal{O}^*(Y_1(N)), \mathcal{O}^*(Y_1(N))\} \cap K_2(X_1(N)) \otimes \mathbf{Q}. \tag{15}$$

La question suivante a été soulevée par Schappacher et Scholl [62, 1.1.3] : l'espace vectoriel réel $V_N$ de l'application régulateur (12) est-il engendré par $r_N(K_N)$ ? Pour le sous-groupe de congruence $\Gamma_1(p)$, $p$ premier, nous montrons le théorème suivant.

**Théorème 5.** *Pour tout nombre premier $p$, l'espace vectoriel réel $V_p$ est engendré par $r_p(K_p)$.*

*Remarque.* L'analogue de ce résultat pour $\Gamma_0(p)$ est faux [62, 1.1.3 (i)].

### 0.5. Rappels sur la conjecture de Zagier pour $L(E,2)$

La valeur spéciale $L(E,2)$, où $E$ une courbe elliptique définie sur $\mathbf{Q}$, a fait l'objet de travaux de Bloch [12] antérieurs à Beĭlinson. Bloch a défini le *dilogarithme elliptique* $D_E$ par la formule très élégante suivante

$$D_E([x]) = \sum_{n=-\infty}^{\infty} D(xq^n) \qquad ([x] \in E(\mathbf{C}) \cong \mathbf{C}^*/q^{\mathbf{Z}}), \tag{16}$$

avec $q \in \mathbf{C}^*$, $|q| < 1$, et où $D$ est la fonction de Bloch-Wigner, version univaluée de la fonction dilogarithme $\mathrm{Li}_2$. On obtient ainsi une fonction continue $D_E : E(\mathbf{C}) \to \mathbf{R}$, bien définie au signe près[1]. Dans le cas où $E$ est à multiplication complexe, Bloch a montré comment exprimer $L(E,2)$ comme combinaison linéaire de valeurs de $D_E$. La généralisation de cet énoncé à toute courbe elliptique est connu sous le nom de *conjecture de Zagier* pour $L(E,2)$ [32, 76, 81]. Les résultats de Bloch et Beĭlinson entraînent l'existence d'un nombre rationnel $c \in \mathbf{Q}^*$, de points $P_1, \ldots, P_k \in E(\overline{\mathbf{Q}}) \subset E(\mathbf{C})$ et d'entiers $n_1, \ldots, n_k \in \mathbf{Z}$ tels que

$$L(E,2) = c\,\pi \sum_{j=1}^{k} n_j\, D_E(P_j). \tag{17}$$

---

[1] Ce signe dépend du choix d'une orientation de $E(\mathbf{R})$.

Grâce à une étude approfondie du groupe $K_2(E)$, Goncharov et Levin [32] ont montré que le diviseur $\sum_{j=1}^{k} n_j \, [P_j]$ de (17) satisfait trois propriétés simples. Fait remarquable, ces conditions suffisent conjecturalement à assurer la proportionnalité du dilogarithme elliptique avec $L(E, 2)$. Nous renvoyons également à [60, 81, 48].

Nombre d'identités numériques ont été découvertes par Bloch et Grayson [13], Grayson et Schappacher (voir la thèse de Rolshausen [59]), Cohen et Zagier (non publié), et sont d'ailleurs à l'origine de la formulation de la conjecture. Il existe même une famille infinie à un paramètre de telles identités [59, 60]. Cependant, la démonstration de ces identités est inaccessible dans la plupart des cas.

### 0.6. Une généralisation du dilogarithme elliptique aux jacobiennes

La complexité de la paramétrisation modulaire $X_1(N) \to E$ suggère de reformuler la conjecture de Zagier au niveau de la jacobienne $J_1(N)$ de $X_1(N)$. Cela nécessite de définir une fonction appropriée sur $J_1(N)$, jouant le rôle du dilogarithme elliptique $D_E$. Nous définissons une telle fonction pour toute jacobienne.

Soient $X$ une surface de Riemann compacte de genre $g \geq 1$, et $J$ la jacobienne de $X$. Soit $\Theta \subset J$ un diviseur thêta. Considérons la fonction $\log\|\theta\| : J - \Theta \to \mathbf{R}$ définie dans [15, p. 196]. Notons $\omega_\Theta \in \Omega^{1,1}(J)$ la forme différentielle duale de Poincaré de $\Theta$. Nous définissons

$$R_J : J \to \mathrm{Hom}_{\mathbf{C}}(\Omega^{1,0}(J), \mathbf{C}) \tag{18}$$
$$x \mapsto \left(\omega \mapsto \frac{1}{(g-1)!} \int_{u \in J} \log\|\theta(u)\| \cdot \overline{\partial}_u \log\|\theta(u-x)\| \wedge \omega_{\Theta,u}^{g-1} \wedge \omega_u \right).$$

La fonction $R_J$ ne dépend pas des choix de $\Theta$ et $\theta$.

Notons $R_X : X \times X \to \mathrm{Hom}_{\mathbf{C}}(\Omega^{1,0}(X), \mathbf{C})$ la fonction définie par Goncharov à l'aide de la fonction de Green de $X$ (nous rappelons la définition de $R_X$ dans le premier chapitre). Pour tous points $x, y \in X$, notons $x - y \in J$ la classe du diviseur $[x] - [y]$ de degré 0 sur $X$.

**Théorème 6.** *Il existe une fonction continue $\Phi_X : X \to \mathrm{Hom}_{\mathbf{C}}(\Omega^{1,0}(X), \mathbf{C})$, telle que nous ayons, via l'identification canonique $\Omega^{1,0}(X) \cong \Omega^{1,0}(J)$*

$$R_J(x - y) = R_X(x, y) + \Phi_X(x) - \Phi_X(y) \qquad (x, y \in X). \tag{19}$$

Appliquons ce résultat au cas où $X = X_1(N)(\mathbf{C})$ et $J = J_1(N)(\mathbf{C})$. L'*algèbre de Hecke* $\mathbf{T} \subset \mathrm{End}_{\mathbf{C}}(S_2(\Gamma_1(N)))$ est le sous-anneau engendré par les opérateurs de Hecke $T_n$, $n \geq 1$ et les opérateurs diamants $\langle d \rangle$, $d \in (\frac{\mathbf{Z}}{N\mathbf{Z}})^*$. Posons

$$L(\mathbf{T}, s) = \sum_{n=1}^{\infty} T_n \otimes \frac{1}{n^s} \qquad L(\mathbf{T}, \chi, s) = \sum_{n=1}^{\infty} T_n \otimes \frac{\chi(n)}{n^s} \qquad \left(\chi : (\frac{\mathbf{Z}}{m\mathbf{Z}})^* \to \mathbf{C}^*\right).$$

Ces fonctions, à valeurs dans $\mathbf{T} \otimes \mathbf{C}$, admettent un prolongement holomorphe au plan complexe. Pour $\psi$ caractère de Dirichlet (pair) modulo $N$, notons $\mathbf{T}^\psi$ la composante $\psi$-isotypique de $\mathbf{T} \otimes \mathbf{C}$ pour l'action des opérateurs diamants. Notons $L(\mathbf{T}^\psi, s)$ et $L(\mathbf{T}^\psi, \chi, s)$ les projections des séries précédentes sur $\mathbf{T}^\psi$.

Nous avons un isomorphisme $\mathbf{T} \otimes \mathbf{C} \cong \mathrm{Hom}_{\mathbf{C}}(S_2(\Gamma_1(N)), \mathbf{C})$ donné par $T \mapsto (f \mapsto a_1(Tf))$. Nous pouvons donc considérer $R_{J_1(N)} = R_{J_1(N)(\mathbf{C})}$ comme une fonction à valeurs dans $\mathbf{T} \otimes \mathbf{C}$. Pour tout $\lambda \in (\frac{\mathbf{Z}}{N\mathbf{Z}})^*$, notons $P_\lambda = \langle \lambda \rangle \infty \in X_1(N)(\mathbf{C})$, où $\infty$ est la pointe infinie de $X_1(N)(\mathbf{C})$.

Les théorèmes 4 et 6 entraînent le résultat suivant.

**Théorème 7.** *Soit $\psi$ un caractère de Dirichlet pair modulo $N$. Pour tout caractère de Dirichlet $\chi$ modulo $N$, pair, primitif et distinct de $\overline{\psi}$, nous avons*

$$L(\mathbf{T}^\psi, 2)\, L(\mathbf{T}^\psi, \chi, 1) = C_{\psi,\chi} \sum_{\lambda,\mu \in (\frac{\mathbf{Z}}{N\mathbf{Z}})^*/\pm 1} \overline{\psi\chi}(\lambda)\, \chi(\mu)\, R_{J_1(N)}(P_\lambda - P_\mu) \qquad (20)$$

*où la constante $C_{\psi,\chi}$ est donnée par*

$$C_{\psi,\chi} = \frac{N\pi i}{\varphi(N)} \cdot \frac{L(\psi\chi, 2)\, \tau(\chi)\, L(\overline{\chi}, 2)}{\pi^4}. \qquad (21)$$

*Remarque.* En particulier, le produit $L(\mathbf{T}^\psi, 2)L(\mathbf{T}^\psi, \chi, 1)$ est combinaison linéaire (explicite) de valeurs de la fonction $R_{J_1(N)}$ en des points $\mathbf{Q}$-rationnels du sous-groupe cuspidal de $J_1(N)$.

En appliquant le théorème 7 au cas $N = 11$, nous obtenons le résultat suivant.

**Théorème 8.** *Soit $E$ la courbe elliptique donnée par l'équation $y^2 + y = x^3 - x^2$, et $P = (0,0)$, point d'ordre $5$ de $E(\mathbf{Q})$. Orientons $E(\mathbf{R})$ dans le sens des $y$ croissants. Soit $\chi$ un caractère de Dirichlet modulo $11$, pair et non trivial. Posons $\zeta = \chi(3) \in \mu_5$. Nous avons la formule*

$$L(E, 2) = \frac{2^2 \cdot 5 \cdot \pi}{11^2} \cdot \frac{1 + 3(\zeta + \overline{\zeta})}{\zeta - \overline{\zeta}} \sum_{a \in \mathbf{Z}/5\mathbf{Z}} \zeta^a D_E(aP). \qquad (22)$$

Nous déduisons du théorème 8 la preuve d'une identité conjecturée par Bloch et Grayson [13].

**Corollaire.** *En conservant les hypothèses du théorème 8, nous avons*

$$L(E, 2) = \frac{10}{11} \cdot \pi \cdot D_E(P) \quad et \quad D_E(2P) = \frac{3}{2} D_E(P). \qquad (23)$$

La seconde des identités (23), appelée *relation exotique*, a été démontrée récemment par Bertin [10].

Le succès de la méthode dans ce cas repose sur deux faits :
– la courbe elliptique $E$ n'est autre que $X_1(11)$ ;
– le sous-groupe de $E(\mathbf{Q})$ engendré par le point $P$ est constitué de pointes.

*Remarque.* La méthode semble également s'appliquer aux courbes elliptiques $X_1(14)$ et $X_1(15)$. Néanmoins, comment montrer la relation conjecturée par Bloch et Grayson pour la courbe $X_0(11)$, isogène à $X_1(11)$ ?

### 0.7. Plan

L'exposition de cette thèse suit la démarche suivante. Le premier chapitre est essentiellement un chapitre de rappels sur la fonction dilogarithme définie par Goncharov dans le cas d'une surface de Riemann compacte. Après en avoir rappelé la définition (section 1.1), nous montrons qu'elle généralise le dilogarithme elliptique défini par Bloch (section 1.2). Nous étudions les propriétés différentielles de cette fonction (section 1.3) ainsi que son comportement vis-à-vis des morphismes finis (section 1.4).

Le deuxième chapitre se place au niveau d'une jacobienne $J$. Nous commençons par construire des éléments dans le $K_2$ de $J$ à partir de points de torsion de la courbe (section 2.1). Nous définissons ensuite une fonction $R_J$ sur $J$ (section 2.2). Dans la section 2.3, nous montrons le théorème 6 qui relie $R_J$ à la fonction dilogarithme du premier chapitre. Nous étudions également le comportement de $R_J$ vis-à-vis de certains automorphismes (section 2.4). En vue de justifier l'introduction de cette fonction, nous calculons le régulateur de Beĭlinson des éléments construits à la section 2.1, en termes de $R_J$ (sections 2.5 et 2.6).

Le chapitre final est consacré à l'explicitation de tous ces objets dans le cas modulaire, en tirant notamment profit de la méthode utilisée par Beĭlinson. Nous rappelons la définition et les propriétés des séries d'Eisenstein en suivant l'exposition de Siegel (section 3.1). Nous exprimons ensuite une intégrale de type Rankin-Selberg en termes de valeurs spéciales de fonctions $L$ (section 3.2). Ce calcul est à la base de tous les résultats présentés dans cette introduction, théorème 6 excepté. Après des rappels sur les unités modulaires (section 3.3), nous présentons une version explicite du théorème de Beĭlinson (section 3.4). Nous montrons en particulier les théorèmes 4, 5 et 7. Le calcul de l'intégrale de Rankin-Selberg admet également le théorème 3 comme conséquence (section 3.5). Nous spécialisons ensuite aux courbes elliptiques (section 3.6), et montrons les théorèmes 1 et 2. Nous traitons en détail l'exemple des courbes modulaires $X_1(11)$ et $X_1(13)$, de genres respectifs 1 et 2 (sections 3.7 et 3.8). Nous montrons en particulier le théorème 8. Enfin, nous indiquons les applications de ces résultats vers la mesure de Mahler, en prenant l'exemple de la courbe $X_1(11)$ (section 3.9).

# Chapitre 1

# Dilogarithme sur une surface de Riemann compacte

Soit $X$ une surface de Riemann compacte connexe non vide, de genre $g \geq 1$. La définition d'un dilogarithme explicite pour $X$ est un problème important, lié en particulier aux conjectures de Beĭlinson sur les valeurs spéciales de fonctions $L$ associées aux courbes projectives lisses sur $\mathbf{Q}$. Une telle définition est particulièrement intéressante car elle permet également d'envisager la formulation d'un analogue de la conjecture de Zagier pour ces mêmes valeurs spéciales.

La fonction $R_X$ introduite dans ce chapitre a été définie par Goncharov [31, Def. 9.1, p. 390], qui la note $G_{1,2}$. D'autre part, cette même fonction a été utilisée implicitement par Deninger et Wingberg [24, §1]. Le point de vue que nous adoptons ici n'est donc pas essentiellement nouveau. En revanche, la propriété différentielle (1.69) et le comportement vis-à-vis des morphismes finis (1.118) n'avaient à notre connaissance pas encore été écrits.

Notons $\Omega^{1,0}(X)$ l'espace vectoriel complexe des 1-formes différentielles holomorphes sur $X$. Pour toute forme différentielle $\omega \in \Omega^{1,0}(X)$, nous allons construire une fonction continue

$$R_\omega : X \times X \to \mathbf{C}, \tag{1.1}$$

de classe $\mathcal{C}^\infty$ hors de la diagonale de $X \times X$. La fonction $R_\omega$ dépend linéairement de $\omega$, donnant lieu à une fonction $R_X$, définie sur $X \times X$ et à valeurs dans le dual de $\Omega^{1,0}(X)$. Nous montrons que la fonction $R_X$ permet d'expliciter une application régulateur

$$r_X : K_2(\mathbf{C}(X)) \to \mathrm{Hom}_{\mathbf{C}}(\Omega^{1,0}(X), \mathbf{C}), \tag{1.2}$$

où $\mathbf{C}(X)$ désigne le corps des fonctions méromorphes sur $X$. Dans la section 1.2, nous nous plaçons dans le cas où $X$ est une courbe elliptique et montrons que la fonction $R_X$ permet de retrouver le dilogarithme elliptique défini par Bloch (proposition 26). Nous retournons ensuite au cas général. Dans la section 1.3, nous montrons que la fonction $R_\omega$ satisfait une propriété différentielle par rapport à chacune de ses variables (théorème 29), et que cette propriété différentielle la caractérise (théorème 33). Enfin, dans la section 1.4, nous étudions le comportement de la fonction $R_X$ vis-à-vis des morphismes finis.

## 1.1 Définition de la fonction $R_X$, selon Goncharov

Une *forme volume* sur $X$ est une 2-forme différentielle réelle de classe $\mathcal{C}^\infty$ sur $X$, partout non nulle et d'intégrale 1 [41, II, §1], [40, p. 329].

Il existe une forme volume canonique sur $X$. Nous en donnons maintenant la définition. Considérons le produit scalaire hermitien

$$(\alpha, \beta) := i \int_X \alpha \wedge \overline{\beta} \qquad (\alpha, \beta \in \Omega^{1,0}(X)) \tag{1.3}$$

sur l'espace vectoriel complexe $\Omega^{1,0}(X)$. Soit $(\omega_j)_{1 \leq j \leq g}$ une base orthonormale pour le produit scalaire (1.3). La forme différentielle

$$\mathrm{vol}_X = \frac{i}{g} \sum_{j=1}^g \omega_j \wedge \overline{\omega_j} \tag{1.4}$$

est une forme volume sur $X$ qui ne dépend pas du choix de la base orthonormale $(\omega_j)_{1 \leq j \leq g}$. Cette forme volume provient naturellement d'un diviseur thêta sur la jacobienne de $X$.

Nous notons $\mathrm{Div}(X)$ le groupe des diviseurs sur $X$, c'est-à-dire le groupe abélien libre de base $X$. Pour tout diviseur $l \in \mathrm{Div}(X)$, nous noterons $(l) \subset X$ son *support* et $\deg l \in \mathbf{Z}$ son *degré*, i. e.

$$(l) = \{x \in X, \ \mathrm{ord}_x(l) \neq 0\}, \tag{1.5}$$

$$\deg l = \sum_{x \in (l)} \mathrm{ord}_x(l). \tag{1.6}$$

Nous rappelons maintenant la définition de la *fonction de Green* associée à $X$. Cette fonction classique joue le rôle de hauteur archimédienne en théorie d'Arakelov. Nous renvoyons à [41, Chap. II] pour les détails ou démonstrations omis ici. Notons

$$\Delta_X = \{(x, x) \mid x \in X\} \subset X \times X \tag{1.7}$$

la diagonale de $X \times X$.

**Proposition 9** (Arakelov [2]). *Il existe une unique fonction*

$$G_X : X \times X - \Delta_X \to \mathbf{R}, \tag{1.8}$$

*appelée* fonction de Green associée à $X$, *de classe $\mathcal{C}^\infty$ et vérifiant les trois conditions suivantes.*

1. *Pour tout $x \in X$, nous avons*

$$\partial_y \overline{\partial}_y G_X(x, y) = \pi i \, \mathrm{vol}_X \qquad (y \in X - \{x\}). \tag{1.9}$$

2. *Pour tout $x \in X$ et pour toute coordonnée locale holomorphe $z(y)$ au point $x$ vérifiant $z(x) = 0$, la fonction*

$$y \mapsto G_X(x, y) - \log|z(y)|, \tag{1.10}$$

*définie sur un voisinage épointé de $x$ dans $X$, s'étend en une fonction de classe $\mathcal{C}^\infty$ sur un voisinage de $x$ dans $X$.*

3. *Pour tout $x \in X$, nous avons*

$$\int_{y \in X} G_X(x, y) \cdot \mathrm{vol}_X = 0. \tag{1.11}$$

*Démonstration.* L'existence de $G_X$ est démontrée par Arakelov dans [2, §1–2]. Coleman en a également donné une preuve [41, II, §4]. L'unicité de $G_X$ résulte quant à elle facilement des propriétés 1, 2 et 3. □

Nous mentionnons maintenant, sans démonstration, quelques propriétés supplémentaires de la fonction de Green.

4. La fonction $G_X$ est symétrique : nous avons

$$G_X(x, y) = G_X(y, x) \qquad (x, y \in X, \ x \neq y). \tag{1.12}$$

5. La fonction $G_X$ a une *singularité logarithmique* le long de $\Delta_X$. Cela signifie que pour tout $x_0 \in X$ et pour toute coordonnée locale holomorphe $z(x)$ au point $x_0$, la fonction

$$(x, y) \mapsto G_X(x, y) - \log|z(x) - z(y)| \tag{1.13}$$

s'étend en une fonction de classe $\mathcal{C}^\infty$ sur un voisinage de $(x_0, x_0)$ dans $X \times X$.

6. Pour toute fonction méromorphe $f \in \mathbf{C}(X)^*$, il existe une constante $C_f \in \mathbf{R}$ telle que

$$\log|f(y)| = C_f + \sum_{x \in (f)} \mathrm{ord}_x(f) \cdot G_X(x, y) \qquad (y \in X - (f)), \tag{1.14}$$

où $(f)$ désigne l'ensemble des zéros et des pôles de $f$, et $\mathrm{ord}_x(f) \in \mathbf{Z}$ désigne l'ordre d'annulation de $f$ en $x \in X$. Nous avons en outre

$$C_f = \int_X \log|f| \cdot \mathrm{vol}_X. \tag{1.15}$$

7. En utilisant le langage des courants, nous pouvons condenser les propriétés 1 et 2 de la fonction $G_X$ en une seule équation : pour tout $x \in X$, nous avons

$$\frac{1}{\pi i} \partial \overline{\partial} G_X(x, \cdot) = \mathrm{vol}_X - \delta_x, \tag{1.16}$$

où $\delta_x$ désigne le courant d'évaluation en $x$.

8. De manière encore plus intrinsèque, la fonction $G_X$ est l'unique distribution sur $X \times X$ satisfaisant l'équation différentielle

$$\frac{1}{\pi i} \partial \overline{\partial} G_X = p^* \mathrm{vol}_X + q^* \mathrm{vol}_X - \Omega_X - \delta_{\Delta_X} \tag{1.17}$$

et normalisée par la condition

$$\int_{X \times X} G_X \cdot p^* \mathrm{vol}_X \wedge q^* \mathrm{vol}_X = 0, \tag{1.18}$$

où les applications $p, q : X \times X \to X$ sont les deux projections naturelles, la forme différentielle $\Omega_X \in \Omega^{1,1}(X \times X)$ est définie par

$$\Omega_X = i \sum_{j=1}^g (p^* \omega_j \wedge q^* \overline{\omega_j} + q^* \omega_j \wedge p^* \overline{\omega_j}) \tag{1.19}$$

pour toute base orthonormale $(\omega_j)_{1 \leq j \leq g}$ de $\Omega^{1,0}(X)$ muni du produit scalaire (1.3), et $\delta_{\Delta_X}$ désigne le courant d'intégration le long de $\Delta_X$.

*Remarque* 10. Dans la proposition 9, il est possible de remplacer $\mathrm{vol}_X$ par une forme volume quelconque $\mathrm{vol}'_X$ sur $X$. La fonction $G'_X$ ainsi obtenue est liée à $G_X$. Nous renvoyons à l'article d'Arakelov [2, Prop. 3.2, p. 1177] pour une version précise de ce lien.

Nous allons maintenant définir la fonction $R_\omega$. Choisissons donc une forme différentielle $\omega \in \Omega^{1,0}(X)$.

**Définition 11.** *Soit $R_\omega$ la fonction définie par*

$$R_\omega : X \times X \to \mathbf{C}$$
$$(x,y) \mapsto \int_{z \in X} G_X(x,z) \cdot \omega_z \wedge \overline{\partial}_z G_X(y,z). \tag{1.20}$$

La convergence absolue de l'intégrale (1.20) résulte de la propriété 2 (p. 16) appliquée aux fonctions $G_X(x,\cdot)$ et $G_X(y,\cdot)$.

**Proposition 12.** *La fonction $R_\omega$ est antisymétrique :*

$$R_\omega(y,x) = -R_\omega(x,y) \qquad (x,y \in X). \tag{1.21}$$

*Démonstration.* Soient $x, y \in X$. Nous avons

$$\begin{aligned} R_\omega(x,y) + R_\omega(y,x) &= \int_{z \in X} \omega_z \wedge \big(G_X(x,z)\overline{\partial}_z G_X(y,z) + G_X(y,z)\overline{\partial}_z G_X(x,z)\big) \\ &= \int_{z \in X} \omega_z \wedge \overline{\partial}_z \big(G_X(x,z) G_X(y,z)\big) \\ &= \int_{z \in X} d\big(G_X(x,z) G_X(y,z) \omega_z\big) \\ &= 0, \end{aligned}$$

d'après la formule de Stokes et car la forme différentielle $G_X(x,z) G_X(y,z) \omega_z$ croît au plus comme le carré du logarithme en $z = x$ et $z = y$. □

**Corollaire 13.** *La fonction $R_\omega$ est nulle sur $\Delta_X$.*

La définition (1.20) de la fonction $R_\omega$ est visiblement linéaire en $\omega$. Il est utile de rassembler toutes ces informations en introduisant une fonction $R_X$ à valeurs dans le dual de $\Omega^{1,0}(X)$.

**Définition 14.** *Soit $R_X$ la fonction définie par*

$$R_X : X \times X \to \mathrm{Hom}_{\mathbf{C}}(\Omega^{1,0}(X), \mathbf{C})$$
$$(x,y) \mapsto \big(\omega \mapsto R_\omega(x,y)\big). \tag{1.22}$$

La définition (1.20) de la fonction $R_\omega$ n'est peut-être pas très parlante. Elle trouve sa justification et son intérêt dans l'explicitation d'une application régulateur

$$r_X : K_2(\mathbf{C}(X)) \to \mathrm{Hom}_{\mathbf{C}}(\Omega^{1,0}(X), \mathbf{C}). \tag{1.23}$$

Rappelons la définition du groupe de $K$-théorie de Milnor $K_2(F)$ associé à un corps $F$.

*Notation.* Soit $F$ un corps. Nous notons $K_2(F)$ le groupe abélien défini par

$$K_2(F) = \frac{F^* \otimes F^*}{\langle x \otimes (1-x) \mid x \in F - \{0,1\} \rangle}. \tag{1.24}$$

Pour tous $x, y \in F^*$, nous appelons *symbole de Milnor associé à $x$ et $y$*, et nous notons $\{x, y\}$, la classe de $x \otimes y$ dans $K_2(F)$.

**Lemme 15.** *Pour toute fonction méromorphe $f \in \mathbf{C}(X) - \{0, 1\}$ et toute forme différentielle holomorphe $\omega \in \Omega^{1,0}(X)$, nous avons*

$$\int_X \log|f| \cdot \omega \wedge \overline{\partial} \log|1-f| = 0. \tag{1.25}$$

*Démonstration.* Notons $D$ la fonction de Bloch-Wigner [79]. Un calcul explicite donne

$$\log|f| \cdot \omega \wedge \overline{\partial} \log|1-f| = \frac{i}{2} d\Big( (D \circ f + i \log|f| \log|1-f|) \omega \Big). \tag{1.26}$$

Le lemme découle donc de la formule de Stokes. □

**Définition 16.** *L'application régulateur $r_X$ associée à $X$ est définie par*

$$\begin{aligned} r_X : K_2(\mathbf{C}(X)) &\to \mathrm{Hom}_{\mathbf{C}}(\Omega^{1,0}(X), \mathbf{C}) \\ \{f, g\} &\mapsto \Big( \omega \mapsto \int_X \log|f| \cdot \omega \wedge \overline{\partial} \log|g| \Big). \end{aligned} \tag{1.27}$$

Cette application est bien définie d'après (1.25). Nous pouvons maintenant faire le lien entre l'application régulateur $r_X$ et la fonction $R_X$.

**Proposition 17.** *Pour toutes fonctions méromorphes $f, g \in \mathbf{C}(X)^*$, nous avons l'égalité*

$$r_X(\{f, g\}) = \sum_{x \in (f)} \sum_{y \in (g)} \mathrm{ord}_x(f) \mathrm{ord}_y(g) R_X(x, y), \tag{1.28}$$

*où nous utilisons les notations de (1.14). Autrement dit, nous avons pour tout $\omega \in \Omega^{1,0}(X)$*

$$\int_X \log|f| \cdot \omega \wedge \overline{\partial} \log|g| = \sum_{x \in (f)} \sum_{y \in (g)} \mathrm{ord}_x(f) \mathrm{ord}_y(g) R_\omega(x, y). \tag{1.29}$$

*Démonstration.* Il suffit de démontrer (1.29). Pour cela, nous utilisons l'égalité (1.14) appliquée à $f$ et $g$, ainsi que le calcul

$$\int_X C_f \cdot \omega \wedge \overline{\partial} \log|g| = -C_f \int_X d\left( \log|g| \cdot \omega \right) = 0$$

utilisant la formule de Stokes. La définition de $R_\omega$ permet de conclure. □

*Remarques.*   1. L'application régulateur $r_X$ définie en (1.27) est essentiellement le régulateur de Bloch et Beĭlinson [5]. Supposons que $X$ est l'ensemble des points complexes d'une courbe projective lisse $X_\mathbf{Q}$ définie sur $\mathbf{Q}$, et notons $\mathbf{Q}(X)$ le corps des fonctions rationnelles de $X_\mathbf{Q}$. Considérons le régulateur de Beĭlinson défini sur le groupe de $K$-théorie de Quillen $K_2^{(2)}(X_\mathbf{Q})$. La localisation en $K$-théorie algébrique induit une inclusion $K_2^{(2)}(X_\mathbf{Q}) \hookrightarrow K_2(\mathbf{Q}(X)) \otimes \mathbf{Q}$, d'où après extension des scalaires un morphisme

$$\eta_X : K_2^{(2)}(X_{\mathbf{Q}}) \to K_2(\mathbf{C}(X)) \otimes \mathbf{Q}.$$

La composition de $\eta_X$ et $r_X \otimes \mathbf{Q}$ donne, à un facteur près, le régulateur de Beĭlinson.

2. L'antisymétrie de la fonction $R_\omega$ se traduit algébriquement par l'égalité bien connue dans $K_2(\mathbf{C}(X))$

$$\{g, f\} = -\{f, g\} \qquad (f, g \in \mathbf{C}(X)^*). \tag{1.30}$$

3. La proposition 1.23 jointe à la relation de Steinberg $\{f, 1 - f\} = 0$ entraîne la relation suivante pour la fonction $R_X$

$$\sum_{x \in (f)} \sum_{y \in (1-f)} \operatorname{ord}_x(f) \operatorname{ord}_y(1-f) R_X(x, y) = 0, \tag{1.31}$$

valable pour toute fonction méromorphe $f \in \mathbf{C}(X) - \{0, 1\}$.

4. Le cas où $X$ est une courbe modulaire nous intéressera particulièrement au chapitre 3. Soit $\mathcal{H}$ le demi-plan de Poincaré et $\Gamma$ un sous-groupe discret de $\mathrm{PSL}_2(\mathbf{R})$, de covolume fini dans $\mathcal{H}$ [15, p. 92]. Notons $Y = \Gamma \backslash \mathcal{H}$ et $X$ la complétion de $Y$, qui est une surface de Riemann compacte, connexe. La forme différentielle définie sur $\mathcal{H}$

$$d\mu = \frac{dx \wedge dy}{y^2} = -2i \cdot \partial\overline{\partial} \log y \qquad (x + iy \in \mathcal{H}) \tag{1.32}$$

induit une 2-forme réelle sur $Y$, partout non nulle et d'intégrale finie. Peut-on définir la fonction $R_X$ en utilisant la forme $d\mu$ plutôt que $\mathrm{vol}_X$ ? Dans cette direction, Gross a construit une fonction de Green sur $Y$ utilisant la forme $d\mu$ [33, §9]. Le passage de $Y$ à $X$ est rendu délicat par le fait que $d\mu$ ne s'étend pas à $X$.

## 1.2   Cas d'une courbe elliptique, selon Bloch

Dans cette section seulement, nous faisons l'hypothèse que $X$ est une courbe elliptique (sur $\mathbf{C}$), c'est-à-dire $g = 1$ et $X$ est muni d'un point distingué $0 \in X$. Nous noterons $X = E$, $G_X = G_E$ et $R_X = R_E$. Nous allons voir que la fonction $R_E$ permet de retrouver le dilogarithme elliptique sur $E$, défini par Bloch. Nous remarquons enfin que l'intérêt de la fonction $R_E$ réside dans son caractère intrinsèque.

**Proposition 18.** *Les fonctions $G_E$ et $R_E$ sont invariantes par translation : pour tout $a \in E$, nous avons les égalités*

$$G_E(x + a, y + a) = G_E(x, y) \qquad (x, y \in E, \ x \neq y) ; \tag{1.33}$$
$$R_E(x + a, y + a) = R_E(x, y) \qquad (x, y \in E). \tag{1.34}$$

*Démonstration.* Fixons $x \in E$ et $a \in E$. Par définition de $G_E$, nous avons

$$\partial_y \overline{\partial}_y G_E(x + a, y) = \pi i (\mathrm{vol}_E - \delta_{x+a}).$$

Notons

$$\begin{aligned} t_a : E &\to E \\ y &\mapsto y + a \end{aligned}$$

la translation par $a$ dans $E$. Nous avons $t_a^* \operatorname{vol}_E = \operatorname{vol}_E$ et $t_a^* \delta_{x+a} = \delta_x$, d'où nous déduisons

$$\partial_y \overline{\partial}_y G_E(x+a, y+a) = \pi i (\operatorname{vol}_E - \delta_x) = \partial_y \overline{\partial}_y G_E(x,y).$$

Puisque toute distribution harmonique sur $E$ est constante, nous en déduisons l'existence d'une constante $C_{x,a} \in \mathbf{R}$ telle que

$$G_E(x+a, y+a) = G_E(x,y) + C_{x,a} \qquad (y \in E, y \neq x). \tag{1.35}$$

En multipliant chacun des membres de (1.35) par $\operatorname{vol}_E$, en intégrant sur $E$ et en utilisant la propriété 3 de $G_E$ (p. 16), nous obtenons $C_{x,a} = 0$, d'où (1.33). Soient maintenant $x, y, a \in E$. Nous avons

$$\begin{aligned}
R_\omega(x+a, y+a) &= \int_{u \in E} G_E(x+a, u) \omega_u \wedge \overline{\partial}_u G_E(y+a, u) \\
&= \int_{u \in E} G_E(x, u-a) \omega_u \wedge \overline{\partial}_u G_E(y, u-a) \\
&= \int_{u \in E} t_{-a}^* \left( G_E(x, u) \omega_u \wedge \overline{\partial}_u G_E(y, u) \right) \\
&= \int_{u \in E} G_E(x, u) \omega_u \wedge \overline{\partial}_u G_E(y, u) \\
&= R_\omega(x, y),
\end{aligned}$$

où au cours du calcul nous avons utilisé le fait que la forme différentielle $\omega$ est invariante par translation. Cela achève de démontrer la proposition. □

Une conséquence de la proposition 18 est que les fonctions $G_E$ et $R_E$ sont essentiellement des fonctions d'une variable sur la courbe elliptique $E$. Rappelons maintenant la définition du dilogarithme elliptique, due à Bloch [12]. Choisissons un isomorphisme

$$\eta : E \xrightarrow{\cong} \frac{\mathbf{C}}{\mathbf{Z} + \tau \mathbf{Z}}, \tag{1.36}$$

avec $\tau \in \mathbf{C}$ vérifiant $\Im(\tau) > 0$. Composons $\eta$ avec l'application $z \mapsto \exp(2\pi i z)$. Nous obtenons un isomorphisme de groupes

$$E \cong \mathbf{C}^*/q^{\mathbf{Z}}, \tag{1.37}$$

où nous avons posé $q = \exp(2\pi i \tau)$.

**Définition 19.** *Le dilogarithme elliptique $D_{E,\eta}$ associé à $E$ et $\eta$ est la fonction définie par*

$$\begin{aligned}
D_{E,\eta} : E &\to \mathbf{R} \\
[x] &\mapsto \sum_{n=-\infty}^{\infty} D(xq^n),
\end{aligned} \tag{1.38}$$

*où nous avons utilisé l'identification (1.37).*

La série (1.38) définissant $D_{E,\eta}$, vue comme série de fonctions de la variable $x \in \mathbf{C}^*$, converge uniformément sur tout compact. La fonction $D_{E,\eta}$ est de classe $\mathcal{C}^\infty$ sur $E - \{0\}$, et continue en 0. Elle vérifie les relations suivantes, dites *relations de distribution*

$$D_{E,\eta}(nP) = n \sum_{Q \in E[n]} D_{E,\eta}(P+Q) \qquad (P \in E,\ n \in \mathbf{Z},\ n \neq 0), \tag{1.39}$$

où $E[n]$ désigne le groupe des points de $n$-torsion de $E$. En particulier, la fonction $D_{E,\eta}$ est impaire.

*Remarque* 20. La fonction $D_{E,\eta}$ dépend du choix de l'isomorphisme $\eta$. Lorsque la courbe elliptique $E$ est définie sur $\mathbf{R}$, il existe un dilogarithme elliptique $D_E$ bien défini au signe près. Pour voir cela, choisissons une orientation[1] de $E(\mathbf{R})$. Notons $H_1^+(E(\mathbf{C}), \mathbf{Z})$ le sous-groupe de $H_1(E(\mathbf{C}), \mathbf{Z})$ constitué des éléments invariants par la conjugaison complexe agissant sur $E(\mathbf{C})$. C'est un groupe abélien libre de rang 1. L'orientation de $E(\mathbf{R})$ induit un générateur canonique $\gamma_1$ de $H_1^+(E(\mathbf{C}), \mathbf{Z})$. Soit $\gamma_2$ un élément de $H_1(E(\mathbf{C}), \mathbf{Z})$ tel que $(\gamma_1, \gamma_2)$ forme une base directe de $H_1(E(\mathbf{C}), \mathbf{Z})$ pour le produit d'intersection[2]. Nous posons alors

$$\tau = \frac{\int_{\gamma_2} \omega}{\int_{\gamma_1} \omega} \in \mathcal{H} \qquad \text{et} \qquad \eta : P \mapsto \left[ \frac{\int_0^P \omega}{\int_{\gamma_1} \omega} \right], \tag{1.40}$$

où $\omega$ désigne une forme différentielle holomorphe non nulle quelconque sur $E(\mathbf{C})$. Le nombre $q = \exp(2\pi i \tau)$ est un nombre réel non nul. Il vérifie $-1 < q < 1$ et ne dépend que de $E$. Nous vérifions que la fonction

$$D_E = D_{E,\eta} : E(\mathbf{C}) \to \mathbf{R} \tag{1.41}$$

ne dépend que du choix de l'orientation de $E(\mathbf{R})$ ; changer cette orientation revient à changer le signe de la fonction $D_E$. Nous avons encore les relations de distribution

$$D_E(nP) = n \sum_{Q \in E[n]} D_E(P+Q) \qquad (P \in E(\mathbf{C}),\ n \in \mathbf{Z},\ n \neq 0), \tag{1.42}$$

Les fonctions $G_E$, $R_E$ et $D_{E,\eta}$ sont définies respectivement sur $E \times E - \Delta_E$, $E \times E$ et $E$. Il est naturel et intéressant de chercher leurs développements en séries de Fourier. D'une part, cela nous permettra d'établir le lien entre le dilogarithme elliptique $D_{E,\eta}$ et la fonction $R_E$ (voir la proposition 26). D'autre part, cela nous mettra sur la voie de la propriété différentielle satisfaite par la fonction $R_\omega$ (voir la proposition 28).

Il est commode et naturel d'indexer les développements de Fourier des fonctions définies sur $E$ par le groupe d'homologie $H_1(E, \mathbf{Z})$. Rappelons que l'*application d'Abel-Jacobi*

$$\begin{aligned} H_1(E, \mathbf{Z}) &\to \mathrm{Hom}_{\mathbf{C}}(\Omega^{1,0}(E), \mathbf{C}) \\ \lambda &\mapsto \left( \omega \mapsto \int_\lambda \omega \right) \end{aligned} \tag{1.43}$$

induit un isomorphisme d'espaces vectoriels réels

$$H_1(E, \mathbf{R}) \cong \mathrm{Hom}_{\mathbf{C}}(\Omega^{1,0}(E), \mathbf{C}). \tag{1.44}$$

On a alors un isomorphisme

---

[1] En général, si $E$ est donnée par l'équation de Weierstraß $y^2 + a_1 xy + a_3 y = x^3 + a_2 x^2 + a_4 x + a_6$, on peut orienter $E(\mathbf{R})$ dans le sens des $y$ croissants, mais un tel choix dépend du choix de l'équation...

[2] L'orientation canonique de $\mathbf{C}$ induit via $\eta$ une orientation canonique de $E(\mathbf{C})$, qui détermine à son tour le signe du produit d'intersection sur $H_1(E(\mathbf{C}), \mathbf{Z})$.

$$E \xrightarrow{\cong} H_1(E, \mathbf{R}/\mathbf{Z}) \tag{1.45}$$
$$P \mapsto \left[\omega \mapsto \int_0^P \omega\right].$$

Notons $\langle \cdot, \cdot \rangle$ le produit d'intersection sur $H_1(E, \mathbf{Z})$, à valeurs dans $\mathbf{Z}$ (son signe est déterminé par l'orientation canonique de $E$). Il induit une application bilinéaire

$$\langle \cdot, \cdot \rangle : H_1(E, \mathbf{Z}) \otimes_{\mathbf{Z}} H_1(E, \mathbf{R}/\mathbf{Z}) \to \mathbf{R}/\mathbf{Z}. \tag{1.46}$$

*Notation.* Pour tout $\lambda \in H_1(E, \mathbf{Z})$, notons $\chi_\lambda : E \to \mathbf{C}^*$ le caractère de $E$ défini par

$$\chi_\lambda(P) = \exp\bigl(2\pi i \langle \lambda, P \rangle\bigr) \qquad (P \in E), \tag{1.47}$$

où nous avons utilisé l'identification (1.45).

Puisque le produit d'intersection est une forme bilinéaire non dégénérée sur $H_1(E, \mathbf{Z})$, les caractères $\chi_\lambda$ forment une base de Fourier de $E$, lorsque $\lambda$ parcourt $H_1(E, \mathbf{Z})$.

L'isomorphisme $\eta$ en (1.36) induit un isomorphisme

$$\eta : H_1(E, \mathbf{Z}) \xrightarrow{\cong} \mathbf{Z} + \tau \mathbf{Z} \subset \mathbf{C}. \tag{1.48}$$

**Théorème 21** (Bloch, [12]). *Le développement de Fourier du dilogarithme elliptique $D_{E,\eta}$ est donné par*

$$D_{E,\eta}(P) = \frac{i\Im(\tau)^2}{\pi} {\sum_{\lambda \in H_1(E, \mathbf{Z})}}' \frac{\Im(\eta(\lambda))}{|\eta(\lambda)|^4} \chi_\lambda(P) \tag{1.49}$$

$$= -\frac{\Im(\tau)^2}{\pi} \Re\left( {\sum_{\lambda \in H_1(E, \mathbf{Z})}}' \frac{\chi_\lambda(P)}{\eta(\lambda)^2 \overline{\eta(\lambda)}} \right) \qquad (P \in E), \tag{1.50}$$

*où le symbole $\sum'$ indique que la somme porte sur les $\lambda \neq 0$.*

Ce type de série est connu classiquement sous le nom de série d'Eisenstein-Kronecker. Le théorème 21 est essentiellement démontré dans [12]. Expliquons brièvement comment compléter la démonstration. Notons $([1], [\tau])$ la base de $H_1(E, \mathbf{Z})$ déduite de l'isomorphisme (1.48). Fixons un entier $C \geq 1$ et un point de $C$-torsion $P \in E$. La fonction

$$f : \mathbf{Z}^2 \to \mathbf{R} \tag{1.51}$$
$$(m, n) \mapsto \Im(\chi_{m[1]+n[\tau]}(P))$$

est impaire, et $C$-périodique en chacune de ses variables. Utilisons le théorème [12, Thm. 10.2.1, p. 77] avec cette fonction $f$. En utilisant le lemme [12, Lem. 10.2.2, p. 77] et en prenant la partie imaginaire, nous obtenons après calculs l'identité (1.50) au point $P$. L'égalité (1.50) est donc vraie pour tous les points de torsion de $E$ ; par continuité, elle est vraie pour tout point de $E$.

Notons que l'expression (1.50) du développement de Fourier de $D_{E,\eta}$ légitime l'introduction d'une fonction $J_{E,\eta}$ définie en remplaçant $\Re$ par $\Im$ dans (1.50). Cette dernière fonction est liée à la fonction $J_q$ de Bloch [12, (8.1.4), p. 62] ; elle coïncide avec l'opposée de la fonction $J(q; x)$ définie par Zagier [79, p. 616].

Le produit d'intersection sur $H_1(E, \mathbf{Z})$ s'étend par $\mathbf{C}$-bilinéarité en une forme $\mathbf{C}$-bilinéaire alternée sur $H_1(E, \mathbf{C})$. Considérons la décomposition de Hodge $H_1(E, \mathbf{C}) = H^{-1,0} \oplus H^{0,-1}$. Pour un élément $\lambda \in H_1(E, \mathbf{C})$, notons $\lambda^{-1,0}$ et $\lambda^{0,-1}$ ses composantes via cette décomposition.

*Notation.* Pour $\lambda \in H_1(E, \mathbf{Z})$, nous définissons

$$\|\lambda\| = \left(2\pi i \left\langle \lambda^{0,-1}, \lambda^{-1,0} \right\rangle\right)^{\frac{1}{2}}. \tag{1.52}$$

Pour montrer que l'expression (1.52) a bien un sens, nous allons la calculer en termes de $\eta(\lambda)$ et $\tau$. Rappelons que l'isomorphisme (1.48) induit une base de $H_1(E, \mathbf{Z})$ que nous notons $([1], [\tau])$. Nous la considérerons également comme une base de l'espace vectoriel complexe $H_1(E, \mathbf{C})$. Un calcul simple montre que

$$\lambda^{-1,0} = \eta(\lambda) \left( \frac{-\overline{\tau}}{\tau - \overline{\tau}} \cdot [1] + \frac{1}{\tau - \overline{\tau}} \cdot [\tau] \right),$$

$$\lambda^{0,-1} = \overline{\eta(\lambda)} \left( \frac{\tau}{\tau - \overline{\tau}} \cdot [1] + \frac{-1}{\tau - \overline{\tau}} \cdot [\tau] \right).$$

Nous en déduisons

$$2\pi i \left\langle \lambda^{0,-1}, \lambda^{-1,0} \right\rangle = \frac{\pi}{\Im(\tau)} |\eta(\lambda)|^2 \geq 0,$$

ce qui justifie la définition (1.52).

Nous allons maintenant exprimer la fonction de Green $G_E$ au moyen d'une série d'Eisenstein-Kronecker.

**Proposition 22.** *La distribution sur $E$ définie par la série de Fourier*

$$\sum_{\lambda \in H_1(E, \mathbf{Z})}' \frac{\chi_\lambda}{\|\lambda\|^2} \tag{1.53}$$

*est de classe $\mathcal{C}^\infty$ sur $E - \{0\}$. Pour tous points $P, Q \in E$ avec $P \neq Q$, nous avons*

$$G_E(P, Q) = -\frac{1}{2} \sum_{\lambda \in H_1(E, \mathbf{Z})}' \frac{\chi_\lambda(P - Q)}{\|\lambda\|^2}. \tag{1.54}$$

*Démonstration.* Nous allons utiliser la caractérisation de la fonction $G_E$ par les propriétés 3 et 7 (pp. 16 et 17). Des calculs sans difficulté montrent que

$$\partial \chi_\lambda = \frac{\pi}{\Im(\tau)} \overline{\eta(\lambda)} \cdot \chi_\lambda \cdot \eta^* dz, \tag{1.55}$$

$$\overline{\partial} \chi_\lambda = -\frac{\pi}{\Im(\tau)} \eta(\lambda) \cdot \chi_\lambda \cdot \eta^* d\overline{z}, \tag{1.56}$$

$$\|\lambda\|^2 = \frac{\pi}{\Im(\tau)} |\eta(\lambda)|^2, \tag{1.57}$$

$$\mathrm{vol}_E = \frac{i}{2\Im(\tau)} \cdot \eta^*(dz \wedge d\overline{z}). \tag{1.58}$$

Nous en déduisons

$$\partial\overline{\partial}\bigg(\sum_{\lambda}{}' \frac{\chi_\lambda}{\|\lambda\|^2}\bigg) = \sum_{\lambda}{}' \frac{\Im(\tau)}{\pi|\eta(\lambda)|^2} \cdot \frac{-\pi^2|\eta(\lambda)|^2}{\Im(\tau)^2} \cdot \chi_\lambda \cdot \eta^*(dz \wedge d\bar z)$$

$$= \frac{-\pi}{\Im(\tau)}\bigg(\sum_{\lambda}{}' \chi_\lambda\bigg)\frac{2\Im(\tau)}{i} \cdot \mathrm{vol}_E$$

$$= 2\pi i(\delta_0 - \mathrm{vol}_E),$$

où $\delta_0$ désigne le courant d'intégration en 0 sur $E$ et vérifie

$$\bigg(\sum_{\lambda \in H_1(E,\mathbf{Z})} \chi_\lambda\bigg)\mathrm{vol}_E = \delta_0. \tag{1.59}$$

D'après la propriété 7 (p. 17) appliquée à $G_E$ et puisque toute distribution harmonique sur $E$ est constante, il existe une constante $C \in \mathbf{C}$ telle que

$$G_E(0,\cdot) = C - \frac{1}{2}\sum_{\lambda \in H_1(E,\mathbf{Z})}{}' \frac{\chi_\lambda}{\|\lambda\|^2}. \tag{1.60}$$

D'après la propriété 3 (p. 16) appliquée à $G_E$ et puisque $\int_E \chi_\lambda \, \mathrm{vol}_E = 0$ pour $\lambda \neq 0$, il vient $C = 0$. Nous remarquons enfin que

$$G_E(P,Q) = G_E(P-Q,0) = G_E(0,P-Q),$$

ce qui permet de conclure. □

*Remarque* 23. Le recours aux distributions dans la proposition 22 s'explique par le fait que la série intervenant dans le membre de droite de (1.54) n'est pas absolument convergente. Il est possible de contourner ce problème en utilisant des procédés de sommation dûs à Eisenstein et Kronecker [75].

La fonction $R_E$ s'exprime elle aussi au moyen d'une série d'Eisenstein-Kronecker.

**Proposition 24.** *Pour tous points $P, Q \in E$, nous avons*

$$R_E(P,Q) = \frac{\pi i}{2}\sum_{\lambda \in H_1(E,\mathbf{Z})}{}' \frac{\chi_\lambda(P-Q)}{\|\lambda\|^4} \cdot [\lambda], \tag{1.61}$$

*où nous notons $[\lambda]$ l'image de $\lambda \in H_1(E,\mathbf{Z})$ par l'application d'Abel-Jacobi (1.43).*

*Démonstration.* Posons $\omega = \eta^*dz$. Puisque l'espace vectoriel $\Omega^{1,0}(E)$ est de dimension 1, il suffit de vérifier l'égalité obtenue à partir de (1.61) en appliquant $\omega$. Notons que

$$[\lambda](\omega) = \int_\lambda \omega = \int_\lambda \eta^*dz = \int_{\eta_*\lambda} dz = \eta(\lambda),$$

de sorte que nous sommes ramenés à démontrer l'égalité

$$R_\omega(P,Q) = \frac{\pi i}{2}\sum_{\lambda \in H_1(E,\mathbf{Z})}{}' \frac{\chi_\lambda(P-Q)}{\|\lambda\|^4} \cdot \eta(\lambda). \tag{1.62}$$

Cela résulte d'un calcul sans difficulté combinant la définition de $R_\omega$ et la proposition 22.

$$\begin{aligned} R_\omega(P,Q) &= \int_{t\in E} G_E(P,t)\omega_t \wedge \overline{\partial}_t G_E(Q,t) \\ &= \frac{1}{4}\int_{t\in E} {\sum_\lambda}' {\sum_\mu}' \frac{\chi_\lambda(P-t)}{\|\lambda\|^2}(\eta^*dz)_t \wedge \overline{\partial}_t\left(\frac{\chi_\mu(Q-t)}{\|\mu\|^2}\right) \\ &= \frac{1}{4}\int_E {\sum_\lambda}' {\sum_\mu}' \frac{\chi_\lambda(P)\chi_\mu(Q)}{\|\lambda\|^2\|\mu\|^2}\chi_{-\lambda}\cdot \eta^*dz \wedge \overline{\partial}\chi_{-\mu} \\ &= \frac{\pi}{4\Im(\tau)}\int_E {\sum_\lambda}' {\sum_\mu}' \frac{\chi_\lambda(P)\chi_\mu(Q)}{\|\lambda\|^2\|\mu\|^2}\eta(\mu)\cdot \chi_{-\lambda-\mu}\cdot \eta^*(dz\wedge d\overline{z}). \end{aligned}$$

En intervertissant les signes $\sum$ et $\int$, il vient

$$\begin{aligned} R_\omega(P,Q) &= \frac{\pi}{4\Im(\tau)}\cdot \frac{2\Im(\tau)}{i}{\sum_\lambda}' {\sum_\mu}' \frac{\chi_\lambda(P)\chi_\mu(Q)}{\|\lambda\|^2\|\mu\|^2}\eta(\mu)\int_E \chi_{-\lambda-\mu}\,\mathrm{vol}_E \\ &= -\frac{\pi i}{2}{\sum_\lambda}' {\sum_\mu}' \frac{\chi_\lambda(P)\chi_\mu(Q)}{\|\lambda\|^2\|\mu\|^2}\eta(\mu)\delta_{-\lambda-\mu,0} \\ &= -\frac{\pi i}{2}{\sum_\lambda}' \frac{\chi_\lambda(P)\chi_{-\lambda}(Q)}{\|\lambda\|^2\|-\lambda\|^2}\eta(-\lambda) \\ &= \frac{\pi i}{2}{\sum_\lambda}' \frac{\chi_\lambda(P-Q)}{\|\lambda\|^4}\eta(\lambda), \end{aligned}$$

où $\delta$ désigne le symbole de Kronecker. $\square$

*Remarque* 25. La série intervenant dans (1.61) est normalement convergente en tant que série de fonctions sur $E\times E$. En conséquence, la fonction $R_E$ est continue sur $E\times E$.

**Proposition 26.** *Soit $E$ une courbe elliptique sur $\mathbf{C}$ munie d'un isomorphisme $\eta : E \xrightarrow{\cong} \mathbf{C}/(\mathbf{Z}+\tau\mathbf{Z})$. Le dilogarithme elliptique $D_{E,\eta}$ s'exprime en fonction de $R_E$ au moyen de la formule*

$$D_{E,\eta}(P) = 2\cdot \Im(R_\omega(P,0)) \qquad (P\in E), \tag{1.63}$$

*où nous avons posé $\omega = \eta^*dz$. En particulier, supposons la courbe elliptique $E$ définie sur $\mathbf{R}$, et prenons l'isomorphisme $\eta$ défini en (1.40). Alors nous avons*

$$R_\omega(P,Q) = \frac{i}{2}D_E(P-Q) \qquad (P,Q\in E(\mathbf{R})), \tag{1.64}$$

*et $\omega = \eta^*dz \in \Omega^{1,0}(E(\mathbf{C}))$ est l'unique forme différentielle vérifiant $\int_{E^0(\mathbf{R})}\omega = 1$, où $E^0(\mathbf{R})$ est la composante neutre de $E(\mathbf{R})$.*

*Démonstration.* Nous allons utiliser l'identité (1.62) obtenue au cours de la démonstration de la proposition 24. Un calcul sans difficulté nous donne

$$\begin{aligned}
2\cdot \Im(R_\omega(P,0)) &= 2\cdot \Im\left(\frac{\pi i}{2}\underset{\lambda\in H_1(E,\mathbf{Z})}{{\sum}'}\frac{\chi_\lambda(P)}{\|\lambda\|^4}\cdot \eta(\lambda)\right)\\
&= \pi\cdot \frac{\Im(\tau)^2}{\pi^2}\Re\left(\underset{\lambda\in H_1(E,\mathbf{Z})}{{\sum}'}\frac{\chi_\lambda(P)}{|\eta(\lambda)|^4}\cdot \eta(\lambda)\right)\\
&= \frac{\Im(\tau)^2}{\pi}\Re\left(\underset{\lambda\in H_1(E,\mathbf{Z})}{{\sum}'}\frac{\chi_\lambda(P)}{\eta(\lambda)\overline{\eta(\lambda)}^2}\right)\\
&= \frac{\Im(\tau)^2}{\pi}\Re\left(\underset{\lambda\in H_1(E,\mathbf{Z})}{{\sum}'}\frac{\overline{\chi_\lambda(P)}}{\eta(\lambda)^2\overline{\eta(\lambda)}}\right)\\
&= \frac{\Im(\tau)^2}{\pi}\Re\left(\underset{\lambda\in H_1(E,\mathbf{Z})}{{\sum}'}\frac{\chi_{-\lambda}(P)}{\eta(\lambda)^2\overline{\eta(\lambda)}}\right)\\
&= -\frac{\Im(\tau)^2}{\pi}\Re\left(\underset{\lambda\in H_1(E,\mathbf{Z})}{{\sum}'}\frac{\chi_\lambda(P)}{\eta(\lambda)^2\overline{\eta(\lambda)}}\right)\\
&= D_{E,\eta}(P).
\end{aligned}$$

Supposons $E$ définie sur $\mathbf{R}$ et soient $P, Q \in E(\mathbf{R})$. Notons $c : E(\mathbf{C}) \to E(\mathbf{C})$ la conjugaison complexe. Elle correspond à la conjugaison complexe usuelle sur $\mathbf{C}/(\mathbf{Z}+\tau\mathbf{Z})$ via l'isomorphisme $\eta$. Nous avons donc $c^*\overline{\omega} = \omega$. Un calcul simple utilisant la définition de $G_E$ montre que

$$G_E(P, c(t)) = \overline{G_E(P,t)} \qquad (t\in E(\mathbf{C}),\, t\neq P). \tag{1.65}$$

D'après la définition (1.20) de $R_\omega$, nous obtenons

$$\begin{aligned}
\overline{R_\omega(P,Q)} &= \int_{t\in E(\mathbf{C})} \overline{G_E(P,t)\cdot \overline{\omega}\wedge \partial_t G_E(Q,t)}\\
&= -\int_{u\in E(\mathbf{C})} G_E(P,c(u))\cdot c^*\omega \wedge \overline{\partial}_u G_E(Q,c(u))\\
&= -\int_{u\in E(\mathbf{C})} G_E(P,u)\cdot c^*\omega \wedge \overline{\partial}_u G_E(Q,u).
\end{aligned}$$

Le signe $-$ apparaît car le changement de variables $t = c(u)$ renverse l'orientation. Puisque $c^*\overline{\omega} = \omega$, il vient $\overline{R_\omega(P,Q)} = -R_\omega(P,Q)$, c'est-à-dire $R_\omega(P,Q) \in i\mathbf{R}$. L'invariance par translation $R_\omega(P,Q) = R_\omega(P-Q,0)$ et (1.63) permettent alors d'en déduire (1.64). Enfin, nous avons

$$\int_{E^0(\mathbf{R})} \omega = \int_{E^0(\mathbf{R})} \eta^* dz = \int_{\eta_* E^0(\mathbf{R})} dz = 1,$$

puisque $\eta_* E^0(\mathbf{R})$ coïncide avec la classe du segment $[0,1]\subset \mathbf{C}$, orienté dans le sens positif.  $\square$

*Remarque* 27. Nous pouvons déduire de la formule (1.63) la dépendance de la fonction $D_{E,\eta}$ en l'isomorphisme $\eta$.

Nous pouvons maintenant faire l'observation suivante. Elle m'a permis de pressentir le théorème 29.

**Proposition 28.** *Pour toute forme différentielle $\omega \in \Omega^{1,0}(E)$, nous avons la propriété différentielle suivante sur $E - \{0\}$*

$$\partial R_\omega(\cdot, 0) = -\pi i\, G_E(\cdot, 0)\, \omega. \tag{1.66}$$

*Démonstration.* Par linéarité, il suffit de démontrer (1.66) pour $\omega = \eta^* dz$. Nous avons vu au cours de la démonstration de la proposition 24 que

$$R_\omega(P, Q) = \frac{\pi i}{2} \sum_{\lambda \in H_1(E, \mathbf{Z})}{}' \frac{\chi_\lambda(P-Q)}{\|\lambda\|^4} \cdot \eta(\lambda). \tag{1.67}$$

Un calcul de la même veine que les démonstrations précédentes donne

$$\begin{aligned}
\partial R_\omega(\cdot, 0) &= \frac{\pi i}{2} \partial \bigg( \sum_\lambda{}' \frac{\chi_\lambda}{\|\lambda\|^4} \cdot \eta(\lambda) \bigg) \\
&= \frac{\pi i}{2} \sum_\lambda{}' \frac{\eta(\lambda)}{\|\lambda\|^4} \cdot \partial \chi_\lambda \\
&= \frac{\pi i}{2} \sum_\lambda{}' \frac{\eta(\lambda)}{\|\lambda\|^4} \cdot \frac{\pi}{\Im(\tau)} \cdot \overline{\eta(\lambda)} \cdot \chi_\lambda \cdot \omega \\
&= \frac{\pi i}{2} \bigg( \sum_\lambda{}' \frac{\chi_\lambda}{\|\lambda\|^2} \bigg) \cdot \omega \qquad \text{d'après (1.57)} \\
&= -\pi i\, G_E(\cdot, 0)\, \omega.
\end{aligned}$$

$\square$

La signification et la portée de l'identité (1.66) sont les suivantes. La fonction $G_E$ est l'analogue elliptique du logarithme univalué $\log|z|$. Nous avons vu également dans la proposition 24 le lien entre la fonction $R_\omega$ et le dilogarithme elliptique. Nous pouvons alors considérer la propriété différentielle (1.66) comme l'analogue de la relation suivante, satisfaite par la fonction de Bloch-Wigner $D$

$$\frac{\partial}{\partial z} D(z) = \frac{i}{2} \left( \frac{\log|z|}{1-z} + \frac{\log|1-z|}{z} \right) \qquad (z \in \mathbf{C} - \{0, 1\}). \tag{1.68}$$

## 1.3  Propriété différentielle et caractérisation

Nous retournons maintenant au cas général : $X$ est une surface de Riemann compacte connexe non vide, de genre $g \geq 1$. Le but de cette section est d'établir le résultat suivant.

**Théorème 29.** *Fixons une forme différentielle $\omega \in \Omega^{1,0}(X)$ et un point $x \in X$. La fonction $y \mapsto R_\omega(x, y)$ est de classe $\mathcal{C}^\infty$ sur $X - \{x\}$ et satisfait la propriété différentielle*

$$\partial_y R_\omega(x, y) = \pi i \big( G_X(x, y) \cdot \omega_y - \alpha_y \big) \qquad (y \in X - \{x\}), \tag{1.69}$$

*où la forme différentielle $\alpha \in \Omega^{1,0}(X)$ (qui dépend de $\omega$ et $x$) est déterminée de manière unique par la condition*

$$\int_X \alpha \wedge \overline{\beta} = \int_{y \in X} G_X(x, y) \cdot \omega_y \wedge \overline{\beta}_y \qquad (\beta \in \Omega^{1,0}(X)). \tag{1.70}$$

*Démonstration.* L'idée directrice est de dériver la quantité $R_\omega(x, y)$ par rapport à $y$ sous le signe intégrale. La difficulté provient de la singularité de l'intégrand au point $y$. Fixons une forme différentielle $\omega \in \Omega^{1,0}(X)$ et des points $x, y_0 \in X$ avec $x \neq y_0$. Nous allons étudier la fonction $y \mapsto R_\omega(x, y)$ au voisinage de $y_0$. Pour cela, introduisons une coordonnée locale holomorphe $u(y)$ au point $y_0$, de telle sorte que $u(y_0) = 0$. Pour $\epsilon > 0$, posons

$$D_\epsilon = \{y \in X;\ |u(y)| < \epsilon\}, \tag{1.71}$$
$$X_\epsilon = X - \overline{D_\epsilon}. \tag{1.72}$$

Nous raisonnerons toujours en supposant $\epsilon$ suffisamment petit. Définissons des fonctions $R_1$ et $R_2$ sur $D_\epsilon$ par

$$R_1(y) = \int_{z \in X_\epsilon} G_X(x, z) \cdot \omega_z \wedge \overline{\partial}_z G_X(y, z) \qquad (y \in D_\epsilon), \tag{1.73}$$
$$R_2(y) = \int_{z \in D_\epsilon} G_X(x, z) \cdot \omega_z \wedge \overline{\partial}_z G_X(y, z) \qquad (y \in D_\epsilon). \tag{1.74}$$

Nous avons $R_\omega(x, y) = R_1(y) + R_2(y)$ pour $y \in D_\epsilon$.

La seule singularité de l'intégrand définissant $R_1(y)$ se trouve en $z = x$, et cette singularité est indépendante de $y$. La fonction $G_X$ étant de classe $\mathcal{C}^\infty$ sur $X \times X - \Delta_X$, nous en déduisons que la fonction $R_1$ est de classe $\mathcal{C}^\infty$ sur $D_\epsilon$ avec

$$\partial_y R_1(y) = \int_{z \in X_\epsilon} G_X(x, z) \cdot \omega_z \wedge \overline{\partial}_z \partial_y G_X(y, z) \qquad (y \in D_\epsilon). \tag{1.75}$$

Choisissons une base orthonormale $(\omega_j)_{1 \leq j \leq g}$ de $\Omega^{1,0}(X)$ pour le produit scalaire (1.3). D'après la propriété 8 (p. 17) de la fonction $G_X$, nous avons l'expression

$$\overline{\partial}_z \partial_y G_X(y, z) = -\partial_y \overline{\partial}_z G_X(y, z)$$
$$= -\pi \sum_{j=1}^g \omega_{j,y} \wedge \overline{\omega_{j,z}}$$
$$= \pi \sum_{j=1}^g \overline{\omega_{j,z}} \wedge \omega_{j,y}. \tag{1.76}$$

Nous en déduisons

$$\partial_y R_1(y) = \pi \sum_{j=1}^g \left( \int_{X_\epsilon} G_X(x, \cdot) \cdot \omega \wedge \overline{\omega_j} \right) \omega_{j,y} \qquad (y \in D_\epsilon). \tag{1.77}$$

La seule singularité de l'intégrand définissant $R_2(y)$ se trouve en $z = y$. Elle dépend de $y$, mais peut être contrôlée grâce à la propriété 5 (p. 17) de la fonction $G_X$. Cette propriété assure l'existence d'une fonction $\phi : D_\epsilon \times D_\epsilon \to \mathbf{R}$ de classe $\mathcal{C}^\infty$ telle que

$$G_X(y, z) = \log|u(y) - u(z)| + \phi(y, z) \qquad (y, z \in D_\epsilon). \tag{1.78}$$

Notons $\overline{u}(z) = \overline{u(z)}$ pour $z \in D_\epsilon$. Nous avons

$$\overline{\partial}_z \log|u(y) - u(z)| = \frac{1}{2} \frac{d\overline{u}(z)}{\overline{u}(z) - \overline{u}(y)} \qquad (y, z \in D_\epsilon). \tag{1.79}$$

Nous pouvons donc écrire $R_2(y) = R_3(y) + R_4(y)$, avec

$$R_3(y) = \frac{1}{2} \int_{z \in D_\epsilon} G_X(x, z) \cdot \omega_z \wedge \frac{d\overline{u}(z)}{\overline{u}(z) - \overline{u}(y)} \qquad (y \in D_\epsilon) \tag{1.80}$$

$$R_4(y) = \int_{z \in D_\epsilon} G_X(x, z) \cdot \omega_z \wedge \overline{\partial}_z \phi(y, z) \qquad (y \in D_\epsilon). \tag{1.81}$$

Nous déduisons du lemme de Poincaré pour l'opérateur $\partial$ que la fonction $R_3$ est de classe $\mathcal{C}^\infty$ sur $D_\epsilon$, avec

$$\partial_y R_3(y) = 2\pi i \cdot \frac{1}{2} G_X(x, y) \cdot \omega_y = \pi i G_X(x, y) \cdot \omega_y \qquad (y \in D_\epsilon). \tag{1.82}$$

La fonction $R_4$ est de classe $\mathcal{C}^\infty$ sur $D_\epsilon$ et

$$\partial_y R_4(y) = \int_{z \in D_\epsilon} G_X(x, z) \cdot \omega_z \wedge \overline{\partial}_z \partial_y \phi(y, z) \qquad (y \in D_\epsilon). \tag{1.83}$$

Si nous faisons la somme des expressions (1.77), (1.82) et (1.83), nous obtenons

$$\begin{aligned}\partial_y R_\omega(x, y) =& \pi i G_X(x, y) \cdot \omega_y + \pi \sum_{j=1}^{g} \left( \int_{X_\epsilon} G_X(x, \cdot) \cdot \omega \wedge \overline{\omega_j} \right) \omega_{j,y} \\ &+ \int_{z \in D_\epsilon} G_X(x, z) \cdot \omega_z \wedge \overline{\partial}_z \partial_y \phi(y, z) \qquad (y \in D_\epsilon),\end{aligned} \tag{1.84}$$

Évaluons maintenant la forme différentielle (1.84) en $y = y_0$, et faisons tendre $\epsilon$ vers 0. La fonction $\phi$ étant de classe $\mathcal{C}^\infty$ sur $D_\epsilon \times D_\epsilon$ et la mesure de $D_\epsilon$ tendant vers 0 lorsque $\epsilon$ tend vers 0, nous avons

$$\left( \int_{z \in D_\epsilon} G_X(x, z) \cdot \omega_z \wedge \overline{\partial}_z \partial_y \phi(y, z) \right)_{y = y_0} \xrightarrow[\epsilon \to 0]{} 0.$$

Nous en déduisons

$$\partial_y R_\omega(x, y) = \pi i \big( G_X(x, y) \cdot \omega_y - \alpha_y \big) \qquad (y \in X - \{x\}),$$

où $\alpha \in \Omega^{1,0}(X)$ est la forme différentielle définie par

$$\alpha = i \sum_{j=1}^{g} \left( \int_X G_X(x, \cdot) \cdot \omega \wedge \overline{\omega_j} \right) \omega_j. \tag{1.85}$$

Il nous reste à montrer que $\alpha$ vérifie la propriété (1.70). Il suffit de montrer cette dernière propriété pour $\beta = \omega_k$ avec $1 \leq k \leq g$. Nous avons

$$\int_X \alpha \wedge \overline{\omega_k} = i \sum_{j=1}^{g} \left( \int_X G_X(x,\cdot) \cdot \omega \wedge \overline{\omega_j} \right) \int_X \omega_j \wedge \overline{\omega_k}$$
$$= i \left( \int_X G_X(x,\cdot) \cdot \omega \wedge \overline{\omega_k} \right) \cdot (-i)$$
$$= \int_X G_X(x,\cdot) \cdot \omega \wedge \overline{\omega_k}.$$

Cela achève de démontrer le théorème 29. □

*Remarque* 30. Le fait que la forme différentielle $\alpha$ introduite en (1.85) vérifie la propriété (1.70) résulte également du calcul suivant utilisant la formule de Stokes

$$\int_X \alpha \wedge \overline{\beta} = \int_{y \in X} G_X(x,y) \omega_y \wedge \overline{\beta}_y - \frac{1}{\pi i} \int_{y \in X} \partial_y R_\omega(x,y) \wedge \overline{\beta}_y \tag{1.86}$$
$$= \int_{y \in X} G_X(x,y) \omega_y \wedge \overline{\beta}_y - \frac{1}{\pi i} \int_X d\big(R_\omega(x,\cdot) \cdot \overline{\beta}\big)$$
$$= \int_{y \in X} G_X(x,y) \omega_y \wedge \overline{\beta}_y \qquad (\beta \in \Omega^{1,0}(X)).$$

Il est utile d'introduire une notation pour la forme différentielle $\alpha$ apparaissant dans la propriété différentielle (1.69).

**Définition 31.** *Soit $S \subset X$ un ensemble fini et $\eta$ une forme différentielle $\mathcal{C}^\infty$ complexe de type $(1,0)$ sur $X - S$ telle que pour tout $x \in S$, la forme différentielle $\eta$ croisse au plus logarithmiquement au point $x$. Nous appelons* projection holomorphe *de $\eta$, et nous notons $\pi_{\mathrm{hol}}(\eta)$, la forme différentielle holomorphe $\pi_{\mathrm{hol}}(\eta)$ sur $X$ caractérisée par la propriété*

$$\int_X \pi_{\mathrm{hol}}(\eta) \wedge \overline{\beta} = \int_X \eta \wedge \overline{\beta} \qquad (\beta \in \Omega^{1,0}(X)). \tag{1.87}$$

La propriété différentielle (1.69) se réécrit alors de la manière suivante

$$\partial R_\omega(x,\cdot) = \pi i \big( G_X(x,\cdot)\, \omega - \pi_{\mathrm{hol}}(G_X(x,\cdot)\, \omega) \big) \qquad (x \in X), \tag{1.88}$$

l'égalité étant valable sur $X - \{x\}$. Puisque la fonction $R_\omega$ est antisymétrique (proposition 12), nous obtenons facilement l'expression de $\partial R_\omega(\cdot, y)$.

$$\partial R_\omega(\cdot, y) = -\partial R_\omega(y, \cdot)$$
$$= -\pi i \big( G_X(y,\cdot) \cdot \omega - \pi_{\mathrm{hol}}(G_X(y,\cdot) \cdot \omega) \big)$$
$$= -\pi i \big( G_X(\cdot, y) \cdot \omega - \pi_{\mathrm{hol}}(G_X(\cdot, y) \cdot \omega) \big) \qquad (y \in X). \tag{1.89}$$

Il n'est en fait pas très difficile de montrer que la fonction $R_X$ est de classe $\mathcal{C}^\infty$ hors de la diagonale de $X$.

**Proposition 32.** *La fonction $R_X : X \times X \to \mathrm{Hom}_{\mathbf{C}}(\Omega^{1,0}(X), \mathbf{C})$ est de classe $\mathcal{C}^\infty$ sur $X \times X - \Delta_X$.*

*Démonstration.* Il suffit bien sûr de fixer une forme différentielle $\omega \in \Omega^{1,0}(X)$ et de montrer le résultat pour $R_\omega$. Soient $x_0 \neq y_0$ des points de $X$, et considérons la fonction $R_\omega$ au voisinage de $(x_0, y_0)$. Introduisons des coordonnées locales en $x_0$ et $y_0$ :

$$x = x(u) \qquad y = y(v) \qquad (|u|, |v| \leq \epsilon),$$

avec $x(0) = x_0$ et $y(0) = y_0$. Nous supposons $\epsilon$ suffisamment petit pour que

$$\overline{D}(x_0, \epsilon) = \{x(u), \; |u| \leq \epsilon\} \quad \text{et} \quad \overline{D}(y_0, \epsilon) = \{y(v), \; |v| \leq \epsilon\}$$

soient disjoints. Notons alors

$$X_\epsilon = X - \big(\overline{D}(x_0, \epsilon) \cup \overline{D}(y_0, \epsilon)\big).$$

Nous avons

$$R_\omega(x(u), y(v)) = R_1(u,v) + R_2(u,v) + R_3(u,v) \qquad (|u|, |v| \leq \epsilon)$$

avec

$$R_1(u,v) = \int_{z \in \overline{D}(x_0,\epsilon)} G_X(x(u), z) \cdot \omega_z \wedge \overline{\partial}_z G_X(y(v), z)$$

$$R_2(u,v) = \int_{z \in \overline{D}(y_0,\epsilon)} G_X(x(u), z) \cdot \omega_z \wedge \overline{\partial}_z G_X(y(v), z)$$

$$R_3(u,v) = \int_{z \in X_\epsilon} G_X(x(u), z) \cdot \omega_z \wedge \overline{\partial}_z G_X(y(v), z).$$

D'après les théorèmes de régularité des intégrales à paramètres, la fonction $R_3$ est de classe $\mathcal{C}^\infty$ sur le domaine $|u|, |v| < \epsilon$. Considérons maintenant la fonction $R_1$ (le cas de la fonction $R_2$ se traite de la même façon). La variable d'intégration s'écrit $z = x(t)$ avec $|t| \leq \epsilon$. Grâce à la propriété 5 de la proposition 9, nous avons

$$G_X(x(u), x(t)) = \log|u - t| + \phi_1(u, t) \qquad (|u|, |t| \leq \epsilon) \qquad (u \neq t),$$

la fonction $\phi_1$ étant de classe $\mathcal{C}^\infty$ sur le domaine $|u|, |t| \leq \epsilon$. Par conséquent

$$R_1(u,v) = \int_{|t| \leq \epsilon} \log|u - t| \cdot \phi_2(v, t) dt \wedge d\overline{t} + \int_{|t| \leq \epsilon} \phi_1(u, t) \cdot \phi_2(v, t) dt \wedge d\overline{t}, \qquad (1.90)$$

où $\phi_2$ est une fonction de classe $\mathcal{C}^\infty$ sur le domaine $|v|, |t| \leq \epsilon$. D'après les théorèmes de régularité des intégrales à paramètres, le second terme de (1.90) est une fonction de classe $\mathcal{C}^\infty$ sur le domaine $|u|, |v| \leq \epsilon$.

Occupons-nous maintenant du premier terme. Pour tout $z_0 \in \mathbf{C}$ et $r > 0$, nous noterons $B(z_0, r)$ (resp. $\overline{B}(z_0, r)$) la boule ouverte (resp. fermée) de centre $z_0$ et de rayon $r$ dans $\mathbf{C}$. Fixons $|u_0| < \epsilon$ et choisissons $\delta > 0$ tel que $\overline{B}(u_0, 2\delta) \subset B(0, \epsilon)$. Soit $\psi : \mathbf{C} \to \mathbf{R}$ une fonction de classe $\mathcal{C}^\infty$ vérifiant

$$\psi = 1 \text{ sur } \overline{B}(u_0, \delta) \qquad \psi = 0 \text{ en dehors de } B(u_0, 2\delta).$$

En écrivant $1 = \psi + (1 - \psi)$, le premier terme de (1.90) s'écrit devient

$$\int_{|t-u_0|\leq 2\delta} \log|u-t| \cdot \phi_2(v,t)\psi(t)dt \wedge d\bar{t} \tag{1.91}$$

$$+ \int_{\substack{|t|\leq \epsilon \\ |t-u_0|>\delta}} \log|u-t| \cdot \phi_2(v,t)(1-\psi(t))dt \wedge d\bar{t}. \tag{1.92}$$

La fonction $\phi_3 : (v,t) \mapsto \phi_2(v,t)\psi(t)$ s'étend en une fonction de classe $\mathcal{C}^\infty$ sur $\overline{B}(0,\epsilon) \times \mathbf{C}$. L'intégrale (1.91) s'écrit alors

$$\int_{|t-u_0|\leq 2\delta} \log|u-t| \cdot \phi_2(v,t)\psi(t)dt \wedge d\bar{t} = \int_{t\in\mathbf{C}} \log|u-t| \cdot \phi_3(v,t)dt \wedge d\bar{t}$$
$$= \int_{t\in\mathbf{C}} \log|t| \cdot \phi_3(v,u-t)dt \wedge d\bar{t}.$$

Remarquons que $|t| \geq 2\epsilon$ entraîne $\phi_3(v,u-t) = 0$. En posant $t = re^{i\theta}$, la dernière intégrale s'écrit

$$-2i\int_0^{2\epsilon}\int_0^{2\pi} r\log r \cdot \phi_3(v,u-re^{i\theta})dr \wedge d\theta,$$

et définit une fonction $\mathcal{C}^\infty$ sur le domaine $|u|,|v| \leq \epsilon$. Enfin, l'intégrale (1.92) est une fonction de classe $\mathcal{C}^\infty$ sur le domaine $(u,v) \in B(u_0,\delta) \times \overline{B}(0,\epsilon)$. La fonction $R_1$ est donc de classe $\mathcal{C}^\infty$ sur $B(u_0,\delta) \times \overline{B}(0,\epsilon)$. En conclusion, la fonction $R_1$ est de classe $\mathcal{C}^\infty$ sur $B(0,\epsilon) \times \overline{B}(0,\epsilon)$, et la fonction $R_\omega$ est de classe $\mathcal{C}^\infty$ sur $X \times X - \Delta_X$. $\square$

Nous montrons maintenant que la propriété différentielle (1.88) caractérise la fonction $R_\omega$.

**Théorème 33.** *Fixons une forme différentielle $\omega \in \Omega^{1,0}(X)$. La fonction $R_\omega : X \times X \to \mathbf{C}$ est l'unique fonction continue sur $X \times X$, nulle sur $\Delta_X$ et vérifiant, pour tout $x \in X$ fixé, la propriété différentielle*

$$\partial R_\omega(x,\cdot) = \pi i \big(G_X(x,\cdot)\,\omega - \pi_{\mathrm{hol}}\,(G_X(x,\cdot)\,\omega)\big) \qquad \text{sur } X - \{x\},$$

*où $\pi_{\mathrm{hol}}$ désigne la projection holomorphe (définition 31).*

*Démonstration.* Montrons la continuité de $R_\omega$. D'après la proposition 32, il suffit de le faire en tout point de la diagonale $\Delta_X$. Soit donc $x_0 \in X$ et $x = x(u)$ une coordonnée locale au voisinage de $x_0$. Pour $\eta > 0$ suffisamment petit, rappelons la notation

$$\overline{D}(x_0,\eta) = \{x(u),\ |u| \leq \eta\},$$

et posons $X_\eta = X - \overline{D}(x_0,\eta)$. Nous voulons montrer que pour tout $\epsilon > 0$, il existe $\eta > 0$ tel que

$$|u|,|v| \leq \eta \quad \text{entraîne} \quad |R_\omega(x(u),x(v))| \leq \epsilon.$$

Pour $\alpha > 0$ suffisamment petit, nous pouvons écrire

$$R_\omega(x(u),x(v)) = R_1(u,v) + R_2(u,v) \qquad (|u|,|v| \leq \alpha),$$

avec

$$R_1(u,v) = \int_{z \in \overline{D}(x_0,\alpha)} G_X(x(u),z) \cdot \omega_z \wedge \overline{\partial}_z G_X(x(v),z) \tag{1.93}$$

$$R_2(u,v) = \int_{z \in X_\alpha} G_X(x(u),z) \cdot \omega_z \wedge \overline{\partial}_z G_X(x(v),z). \tag{1.94}$$

Pour traiter l'intégrale (1.93), nous écrivons

$$G_X(x(u),x(t)) = \log|u-t| + \phi_1(u,t) \qquad (|u|,|t| \leq \alpha) \qquad (u \neq t), \tag{1.95}$$

la fonction $\phi_1$ étant de classe $\mathcal{C}^\infty$ sur le domaine $|u|, |t| \leq \alpha$. Posant $\phi_2 = \partial \phi_1 / \partial \overline{t}$, nous en déduisons

$$R_1(u,v) = \int_{|t| \leq \alpha} \big(\log|u-t| + \phi_1(u,t)\big)\Big(\frac{-1}{2(\overline{v}-\overline{t})} + \phi_2(v,t)\Big)\phi_3(t) dt \wedge d\overline{t}, \tag{1.96}$$

la fonction $\phi_3$ étant de classe $\mathcal{C}^\infty$ sur le domaine $|t| \leq \alpha$. Nous allons trouver $\alpha > 0$ tel que $|R_1(u,v)| \leq \frac{\epsilon}{2}$ pour tout $|u|, |v| \leq \alpha$. L'intégrale (1.96) s'écrit comme somme de quatre termes. Les fonctions $\phi_1$, $\phi_2$ et $\phi_3$ étant bornées sur leurs domaines de définition respectifs, il suffit de considérer les quantités

$$I_1(u) = \int_{|t| \leq \alpha} \big|\log|u-t| \cdot dt \wedge d\overline{t}\big| \qquad (|u| \leq \alpha)$$

$$I_2(v) = \int_{|t| \leq \alpha} \Big|\frac{dt \wedge d\overline{t}}{v-t}\Big| \qquad (|v| \leq \alpha)$$

$$I_3(u,v) = \int_{|t| \leq \alpha} \Big|\frac{\log|u-t| \cdot dt \wedge d\overline{t}}{v-t}\Big| \qquad (|u|,|v| \leq \alpha),$$

et de montrer qu'elles tendent vers 0 lorsque $\alpha$ tend vers 0, uniformément en $u$ et $v$. Nous allons traiter le cas de $I_3(u,v)$, les cas de $I_1(u)$ et $I_2(v)$ étant plus simples. Par changement de variables, nous avons

$$I_3(u,v) = \int_{t \in \overline{B}(u,\alpha)} \Big|\frac{\log|t| \cdot dt \wedge d\overline{t}}{t+v-u}\Big| \leq \int_{|t| \leq 2\alpha} \Big|\frac{\log|t| \cdot dt \wedge d\overline{t}}{t-x}\Big|,$$

où $x = u - v$ vérifie $|x| \leq 2\alpha$. Il suffit donc de considérer l'intégrale

$$I_3'(x) = \int_{|t| \leq \alpha} \Big|\frac{\log|t| \cdot dt \wedge d\overline{t}}{t-x}\Big| \qquad (|x| \leq \alpha),$$

et de montrer qu'elle tend vers 0 lorsque $\alpha$ tend vers 0, uniformément en $x$. Pour cela, nous écrivons

$$I_3'(x) = \int_{|t| \leq \min(\alpha,2|x|)} \Big|\frac{\log|t| \cdot dt \wedge d\overline{t}}{t-x}\Big| + \int_{2|x| \leq |t| \leq \alpha} \Big|\frac{\log|t| \cdot dt \wedge d\overline{t}}{t-x}\Big|.$$

La fonction $r \mapsto r|\log r|$ étant croissante sur l'intervalle $]0, \frac{1}{e}]$, nous avons

$$\int_{|t|\leq\min(\alpha,2|x|)}\left|\frac{\log|t|\cdot dt\wedge d\bar{t}}{t-x}\right|\leq 2|x\log 2|x||\int_{|t|\leq\min(\alpha,2|x|)}\left|\frac{dt\wedge d\bar{t}}{t(t-x)}\right|$$

$$\leq 2|\log 2|x||\int_{|t|\leq\min(\alpha,2|x|)}\left|\left(\frac{1}{t-x}-\frac{1}{t}\right)\cdot dt\wedge d\bar{t}\right|$$

$$\leq 4|\log 2|x||\int_{|t|\leq 3|x|}\left|\frac{dt\wedge d\bar{t}}{t}\right|.$$

Cette dernière intégrale se calcule grâce aux coordonnées polaires, et donne lieu à la majoration

$$\int_{|t|\leq\min(\alpha,2|x|)}\left|\frac{\log|t|\cdot dt\wedge d\bar{t}}{t-x}\right|\leq 48\pi|x\log 2|x|| \qquad (|x|\leq\alpha). \tag{1.97}$$

Pour traiter la seconde intégrale, nous remarquons que $|t|\geq 2|x|$ entraîne $|t-x|\geq\frac{|t|}{2}$, d'où

$$\int_{2|x|\leq|t|\leq\alpha}\left|\frac{\log|t|\cdot dt\wedge d\bar{t}}{t-x}\right|\leq 2\int_{2|x|\leq|t|\leq\alpha}\left|\frac{\log|t|\cdot dt\wedge d\bar{t}}{t}\right|$$

$$\leq 8\pi\int_0^\alpha|\log r|dr \qquad (|x|\leq\alpha). \tag{1.98}$$

D'après (1.97) et (1.98), la quantité $I_3'(x)$ tend vers 0 lorsque $\alpha$ tend vers 0, uniformément pour $|x|\leq\alpha$. Par suite, il existe $\alpha_0>0$ tel que pour tout choix de $\alpha\in]0,\alpha_0]$ et tout $|u|,|v|\leq\alpha$, nous ayons $|R_1(u,v)|\leq\frac{\epsilon}{2}$ pour tout $|u|,|v|\leq\alpha$. Considérons maintenant l'intégrale $R_2(u,v)$ (1.94). La fonction $R_2$ est continue sur le domaine $|u|,|v|<\alpha$. Or

$$R_2(0,0)=\int_{z\in X_\alpha}G_X(x_0,z)\cdot\omega_z\wedge\overline{\partial}_zG_X(x_0,z)$$

$$=-\frac{1}{2}\int_{z\in X_\alpha}d\bigl(G_X(x_0,z)^2\cdot\omega_z\bigr)$$

$$=\frac{1}{2}\int_{\partial\overline{D}(x_0,\alpha)}G_X(x_0,z)^2\cdot\omega_z.$$

En utilisant (1.95), il vient

$$R_2(0,0)=\frac{1}{2}\int_{|t|=\alpha}\bigl(\log|t|+\phi_1(0,t)\bigr)^2\cdot\phi_3(t)dt$$

$$=\frac{1}{2}\int_{|t|=\alpha}\bigl(\log\alpha+\phi_1(0,t)\bigr)^2\cdot\phi_3(t)dt.$$

Les fonctions $\phi_1$ et $\phi_3$ étant continues en 0, il est facile de voir que cette dernière quantité tend vers 0 lorsque $\alpha$ tend vers 0. Il existe donc $\alpha\in]0,\alpha_0]$ tel que $|R_2(0,0)|\leq\frac{\epsilon}{4}$. Pour ce choix de $\alpha$, la continuité de la fonction $R_2$ implique l'existence d'un $\eta\in]0,\alpha]$ tel que $|R_2(u,v)|\leq\frac{\epsilon}{2}$ pour tout $|u|,|v|\leq\eta$. Nous avons finalement

$$|R_\omega(x(u),x(v))|\leq\epsilon \qquad (|u|,|v|\leq\eta),$$

ce qui achève de montrer la continuité de la fonction $R_\omega$.

Montrons finalement l'unicité de $R_\omega$. Soit $R'_\omega : X \times X \to \mathbf{C}$ une autre fonction vérifiant les propriétés de l'énoncé. Fixons $x_0 \in X$ et considérons la fonction

$$\psi : y \in X \mapsto R'_\omega(x_0, y) - R_\omega(x_0, y).$$

La fonction $\psi$ est continue sur $X$ et de classe $\mathcal{C}^\infty$ sur $X - \{x_0\}$. De plus, elle vérifie $\partial \psi = 0$ sur $X - \{x_0\}$. Par suite, la fonction $\psi$ est antiholomorphe sur $X$, donc constante. Mais $\psi(x_0) = 0$, d'où $\psi = 0$. Ceci étant valable pour tout point $x_0$, nous avons $R'_\omega = R_\omega$. □

*Remarque* 34. Dans la démonstration du théorème 33, nous n'avons en réalité utilisé que la continuité et le caractère $\mathcal{C}^\infty$ par rapport à la seconde variable de la fonction $R_\omega$.

## 1.4  Comportement vis-à-vis des morphismes finis

Le dilogarithme elliptique vérifie les relations de distribution (1.39). Il est naturel de chercher leur analogue dans le cadre de la fonction $R_X$.

Soit $E$ une courbe elliptique sur $\mathbf{C}$. Fixons un entier $n \in \mathbf{Z} - \{0\}$ et notons $\phi : E \to E$ l'isogénie donnée par la multiplication par $n$. La relation (1.39) peut se récrire sous la forme

$$D_{E,\eta}(P) = n D_{E,\eta}(\phi^* P) \qquad (P \in E), \tag{1.99}$$

où $\phi^* P$ est le diviseur sur $E$ défini par

$$\phi^* P = \sum_{\substack{Q \in E \\ \phi(Q) = P}} [Q], \tag{1.100}$$

et où nous avons convenu de définir le dilogarithme elliptique d'un diviseur (ici $\phi^* P$) par linéarité. Lorsque la courbe elliptique $E$ est définie sur $\mathbf{R}$, nous obtenons en particulier

$$D_E(P) = n D_E(\phi^* P) \qquad (P \in E(\mathbf{C})). \tag{1.101}$$

Remarquons également que la multiplication par $n$ intervenant dans (1.99) et (1.101) n'est autre que le morphisme $\phi_*$ induit par $\phi$ sur $H_1(E(\mathbf{C}), \mathbf{R})$. Le lien (1.63) entre les fonctions $D_{E,\eta}$ et $R_\omega$ semble donc indiquer le point de vue suivant. Donnons-nous un morphisme fini

$$\phi : X \to Y \tag{1.102}$$

entre deux surfaces de Riemann compactes connexes. Par image réciproque, $\phi$ induit une application injective

$$\phi^* : \Omega^{1,0}(Y) \to \Omega^{1,0}(X). \tag{1.103}$$

Notons $g_X$ (resp. $g_Y$) le genre de $X$ (resp. $Y$). Nous supposerons dans tout le reste de la section que $g_Y \geq 1$, de sorte que nous avons également $g_X \geq 1$.

Nous étendons la définition de la fonction de Green $G_X$ aux diviseurs sur $X$ en posant

$$G_X(l, m) = \sum_{x \in (l)} \sum_{y \in (m)} \mathrm{ord}_x(l) \, \mathrm{ord}_y(m) G_X(x, y) \tag{1.104}$$

pour tous diviseurs $l, m \in \mathrm{Div}(X)$ tels que $(l) \cap (m) = \emptyset$. Le morphisme fini $\phi$ induit par image réciproque une application

$$\phi^* : \text{Div}(Y) \to \text{Div}(X). \tag{1.105}$$

Dans un premier temps, comparons les fonctions de Green $G_X$ et $G_Y$. Une telle comparaison est sans doute connue, mais je n'ai pas trouvé de référence.

**Proposition 35.** *Il existe une unique fonction $\Phi : X \to \mathbf{R}$ de classe $\mathcal{C}^\infty$ vérifiant les conditions*

$$\frac{1}{\pi i}\partial\overline{\partial}\Phi = \phi^*\operatorname{vol}_Y - (\deg\phi)\operatorname{vol}_X \qquad et \qquad \int_X \Phi\operatorname{vol}_X = 0. \tag{1.106}$$

*Pour tout $x \in X$ et $y \in Y$ tels que $\phi(x) \neq y$, nous avons alors*

$$G_X(x, \phi^*y) = G_Y(\phi(x), y) + \Phi(x) + \frac{\Phi(\phi^*y)}{\deg\phi} + C, \tag{1.107}$$

*où nous avons convenu de définir $\Phi(\phi^*y)$ par linéarité, et avons posé*

$$C = \frac{i}{\pi\deg\phi}\int_X \Phi\partial\overline{\partial}\Phi \leq 0. \tag{1.108}$$

*Supposons de plus que $g_X = g_Y$. Alors $\Phi = 0$, de sorte que pour tout $x \in X$ et $y \in Y$ tels que $\phi(x) \neq y$, nous avons*

$$G_X(x, \phi^*y) = G_Y(\phi(x), y). \tag{1.109}$$

*Démonstration.* D'après [15, Prop. C.5.1, p. 147] et le calcul

$$\int_X (\phi^*\operatorname{vol}_Y - (\deg\phi)\operatorname{vol}_X) = (\deg\phi)\left(\int_Y \operatorname{vol}_Y\right) - \deg\phi = 0,$$

il existe une fonction $\Phi : X \to \mathbf{C}$ de classe $\mathcal{C}^\infty$ vérifiant la première condition de (1.106). D'après [15, Prop. C.5.1, p. 147], la seconde condition de (1.106) détermine $\Phi$ de manière unique. En appliquant la conjugaison complexe aux identités (1.106), il n'est pas difficile de voir que $\overline{\Phi} = \Phi$, donc que $\Phi$ est à valeurs réelles.

Plaçons-nous dans le cas particulier $g_X = g_Y$, et montrons que $\Phi = 0$. Pour des raisons de dimension, l'application (1.103) est un isomorphisme. Soit $(\omega_j)_{1 \leq j \leq g_Y}$ une base orthonormale de $\Omega^{1,0}(Y)$ pour le produit scalaire (1.3). Le calcul

$$i\int_X \phi^*\omega_j \wedge \overline{\phi^*\omega_k} = i(\deg\phi)\int_Y \omega_j \wedge \overline{\omega_k} = (\deg\phi)\delta_{j,k}$$

montre que la famille $(\phi^*\omega_j/\sqrt{\deg\phi})_{1 \leq j \leq g_Y}$ constitue une base orthonormale de $\Omega^{1,0}(X)$ pour le produit scalaire (1.3). Par suite, la forme volume sur $X$ est donnée par

$$\operatorname{vol}_X = \frac{i}{g_X}\sum_{j=1}^{g_Y} \frac{\phi^*\omega_j \wedge \overline{\phi^*\omega_j}}{\deg\phi} = \frac{\phi^*\operatorname{vol}_Y}{\deg\phi}. \tag{1.110}$$

D'après la caractérisation (1.106) de $\Phi$, il suit $\Phi = 0$.

Fixons maintenant un point $y$ de $Y$. Considérons la fonction $\psi_y$ de classe $\mathcal{C}^\infty$ définie par

$$\begin{aligned}\psi_y : X - (\phi^*y) &\to \mathbf{R} \\ x &\mapsto G_X(x, \phi^*y) - G_Y(\phi(x), y).\end{aligned}$$

En utilisant la propriété (1.13) appliquée aux fonctions de Green $G_X$ et $G_Y$, nous voyons que la fonction $\psi$ se prolonge en une fonction de classe $\mathcal{C}^\infty$

$$\widetilde{\psi}_y : X \to \mathbf{R}.$$

D'après les propriétés (1.9) et (1.11) appliquées à $G_X$ et $G_Y$, nous avons sur $X - (\phi^*y)$

$$\begin{aligned}\frac{1}{\pi i}\partial\overline{\partial}\psi_y &= -(\deg\phi^*y)\operatorname{vol}_X - \phi^*(-\operatorname{vol}_Y) \\ &= \phi^*\operatorname{vol}_Y - (\deg\phi)\operatorname{vol}_X.\end{aligned} \quad (1.111)$$

D'après (1.106) et (1.111), il existe une constante $A(y) \in \mathbf{R}$ telle que

$$\widetilde{\psi}_y(x) = \Phi(x) + A(y) \qquad (x \in X). \quad (1.112)$$

Considérons un point $z$ de $Y$ distinct de $y$. Utilisons (1.112) avec le diviseur $x = \phi^*z$. En convenant de définir par linéarité les quantités qui suivent, nous obtenons

$$G_X(\phi^*z, \phi^*y) = (\deg\phi)G_Y(z, y) + \Phi(\phi^*z) + (\deg\phi)A(y).$$

En utilisant la symétrie (1.13) des fonctions de Green, nous en déduisons que la quantité $\Phi(\phi^*z) + (\deg\phi)A(y)$ est symétrique en $y$ et $z$. Il existe par conséquent une constante absolue $C \in \mathbf{R}$ telle que

$$A(y) = \frac{\Phi(\phi^*y)}{\deg\phi} + C \qquad (y \in Y),$$

ce qui montre (1.107).

Déterminons enfin la constante $C$. Considérons $\operatorname{vol}_X$ (resp. $\operatorname{vol}_Y$) comme une forme différentielle sur $X \times Y$ au moyen de la première (resp. seconde) projection. Multiplions chacun des membres de (1.107) par $\operatorname{vol}_X \wedge \operatorname{vol}_Y$ et intégrons sur $X \times Y$. Nous obtenons

$$\begin{aligned}C =& \int_{X \times Y} G_X(x, \phi^*y)\operatorname{vol}_X \wedge \operatorname{vol}_Y - \int_{X \times Y} G_Y(\phi(x), y)\operatorname{vol}_X \wedge \operatorname{vol}_Y \\ &- \int_{X \times Y} \Phi(x)\operatorname{vol}_X \wedge \operatorname{vol}_Y - \int_{X \times Y} \frac{\Phi(\phi^*y)}{\deg\phi}\operatorname{vol}_X \wedge \operatorname{vol}_Y.\end{aligned}$$

Puisque les intégrales des formes différentielles considérées sont absolument convergentes, nous pouvons utiliser le théorème de Fubini. Grâce à la propriété (1.10) des fonctions de Green et à la seconde condition de (1.106), il en découle

$$\begin{aligned}C &= -\frac{1}{\deg\phi}\int_{y\in Y}\Phi(\phi^*y)\operatorname{vol}_Y \quad (1.113)\\ &= -\frac{1}{\deg\phi}\int_X \Phi \cdot \phi^*\operatorname{vol}_Y \\ &= -\frac{1}{\pi i \deg\phi}\int_X \Phi\,\partial\overline{\partial}\Phi - \int_X \Phi \cdot \operatorname{vol}_X \\ &= \frac{i}{\pi\deg\phi}\int_X \Phi\,\partial\overline{\partial}\Phi,\end{aligned}$$

en utilisant (1.106). Le fait que $C \leq 0$ est démontré dans [15, p. 148]. Cela achève de démontrer la proposition 35. $\square$

Cherchons maintenant à comparer les fonctions $R_X$ et $R_Y$. Nous convenons d'étendre la définition de la fonction $R_X$ aux diviseurs sur $X \times X$, de la même manière qu'en (1.104).

**Proposition 36.** *Pour toute forme différentielle $\omega_Y \in \Omega^{1,0}(Y)$ et tout diviseur $m$ de degré 0 sur $Y$, il existe une constante $C(\omega_Y, m) \in \mathbf{C}$ telle que*

$$R_{\phi^*\omega_Y}(x, \phi^*m) = R_{\omega_Y}(\phi(x), m) + C(\omega_Y, m) \qquad (x \in X). \tag{1.114}$$

*Supposons de plus que $g_X = g_Y$. Alors l'égalité (1.114) est valable pour tout diviseur $m$ sur $Y$, et nous avons $C(\omega_Y, m) = 0$.*

*Démonstration.* Considérons une forme différentielle $\omega_Y \in \Omega^{1,0}(Y)$ et un diviseur $m$ de degré 0 sur $Y$. Pour démontrer (1.114), nous allons utiliser la caractérisation différentielle des fonctions $R_{\omega_X}$ et $R_{\omega_Y}$. Posons $\omega_X = \phi^*\omega_Y$. Utilisons également la fonction $\Phi$ définie à la proposition 35. D'après cette proposition et le fait que $m$ est de degré 0, nous avons

$$\begin{aligned}\frac{1}{\pi i}\partial R_{\omega_X}(\cdot, \phi^*m) &= G_X(\cdot, \phi^*m)\,\omega_X - \pi_{\text{hol}}\left(G_X(\cdot, \phi^*m)\,\omega_X\right) \\ &= \left(G_Y(\phi(\cdot), m) + \frac{\Phi(\phi^*m)}{\deg \phi}\right)\omega_X \\ &\quad - \pi_{\text{hol}}\left(\left(G_Y(\phi(\cdot), m) + \frac{\Phi(\phi^*m)}{\deg \phi}\right)\omega_X\right) \\ &= \phi^*G_Y(\cdot, m)\,\omega_X - \pi_{\text{hol}}\left(\phi^*G_Y(\cdot, m)\,\omega_X\right).\end{aligned}$$

D'autre part, nous avons

$$\begin{aligned}\frac{1}{\pi i}\partial R_{\omega_Y}(\phi(\cdot), m) &= \phi^*(\frac{1}{\pi i}\partial R_{\omega_Y}(\cdot, m)) \\ &= \phi^*\left(G_Y(\cdot, m)\,\omega_Y - \pi_{\text{hol}}(G_Y(\cdot, m)\,\omega_Y)\right) \\ &= \phi^*G_Y(\cdot, m)\,\omega_X - \phi^*\pi_{\text{hol}}\left(G_Y(\cdot, m)\,\omega_Y\right).\end{aligned}$$

Considérons la fonction continue

$$\begin{aligned}\psi : X &\to \mathbf{C} \\ x &\mapsto R_{\omega_X}(x, \phi^*m) - R_{\omega_Y}(\phi(x), m).\end{aligned}$$

La fonction $\psi$ est de classe $\mathcal{C}^\infty$ sur $X - (\phi^*m)$. D'après les calculs précédents, il existe une forme différentielle holomorphe $\alpha \in \Omega^{1,0}(X)$ telle que

$$\partial \psi = \alpha|_{X-(\phi^*m)}.$$

Pour toute forme différentielle $\beta \in \Omega^{1,0}(X)$, nous avons

$$\begin{aligned}\int_X \alpha \wedge \overline{\beta} &= \int_{X-(\phi^*m)} \partial\psi \wedge \overline{\beta} \tag{1.115} \\ &= \int_{X-(\phi^*m)} d(\psi\overline{\beta}) \\ &= 0\end{aligned}$$

d'après la formule de Stokes, et puisque $\psi$ est continue sur $X$. Puisque (1.115) est vrai pour toute forme différentielle $\beta$, nous en déduisons $\alpha = 0$. La fonction $\psi$ étant anti-holomorphe sur $X - (\phi^*m)$ et continue sur $X$, elle est nécessairement constante sur $X$, d'où (1.114).

Supposons maintenant que $g_X = g_Y$, et considérons un diviseur quelconque $m$ sur $Y$. En utilisant l'égalité (1.109) fournie par la proposition 35, on se convainc que le raisonnement ci-dessus permet encore de montrer (1.114). Il nous reste donc à montrer $C(\omega_Y, m) = 0$. Multiplions les deux membres de l'égalité (1.114) par $\operatorname{vol}_X$ et intégrons sur $X$. D'après (1.110), nous avons

$$C(\omega_Y, m) = \int_X R_{\omega_X}(\cdot, \phi^*m) \operatorname{vol}_X - \frac{1}{\deg \phi} \int_X \phi^* R_{\omega_Y}(\cdot, m) \, \phi^* \operatorname{vol}_Y$$
$$= \int_X R_{\omega_X}(\cdot, \phi^*m) \operatorname{vol}_X - \int_Y R_{\omega_Y}(\cdot, m) \operatorname{vol}_Y.$$

Pour tout $z \in Y$, posons

$$\Psi(z) = \int_{y \in Y} R_{\omega_Y}(y, z) \operatorname{vol}_Y. \tag{1.116}$$

Il nous suffit de montrer $\Psi = 0$. Par définition de $R_{\omega_Y}$, nous avons

$$\Psi(z) = \int_{y \in Y} \left( \int_{t \in Y} G_Y(y, t) \omega_{Y,t} \wedge \overline{\partial}_t G_Y(z, t) \right) \operatorname{vol}_{Y,y}$$
$$= \int_{t \in Y} \left( \int_{y \in Y} G_Y(y, t) \operatorname{vol}_{Y,y} \right) \omega_{Y,t} \wedge \overline{\partial}_t G_Y(z, t)$$

d'après le théorème de Fubini, qui s'applique car la forme différentielle en question est intégrable sur $Y \times Y$. D'après la propriété (1.11) appliquée à $G_Y$, nous concluons $\Psi = 0$. Nous avons donc bien $C(\omega_Y, m) = 0$. □

La proposition 36 admet la conséquence agréable suivante pour les fonctions $R_X$ et $R_Y$. Par dualité, le morphisme $\phi^*$ de (1.103) induit un morphisme

$$\phi_* : \operatorname{Hom}_{\mathbf{C}}(\Omega^{1,0}(X), \mathbf{C}) \to \operatorname{Hom}_{\mathbf{C}}(\Omega^{1,0}(Y), \mathbf{C}). \tag{1.117}$$

**Corollaire 37.** *Soit $l$ (resp. $m$) un diviseur de degré $0$ sur $X$ (resp. $Y$). Nous avons*

$$\phi_* R_X(l, \phi^*m) = R_Y(\phi_* l, m). \tag{1.118}$$

*Supposons de plus que $g_X = g_Y$. Alors l'égalité (1.118) est valable pour tous les diviseurs $l$ (resp. $m$) sur $X$ (resp. $Y$).*

*Remarque* 38. La formule (1.118) peut être interprétée comme un analogue, pour les fonctions $R_X$ et $R_Y$, de la formule de projection en $K$-théorie algébrique.

D'après le théorème de Hurwitz [34, Cor. 2.4, p. 301] et l'hypothèse $g_Y \geq 1$, la condition $g_X = g_Y$ est réalisée si et seulement si $\phi$ est un isomorphisme ou $g_X = g_Y = 1$. Dans ce dernier cas, le morphisme $\phi$ est une isogénie entre courbes elliptiques. Nous détaillons ces deux cas dans les deux corollaires qui suivent.

**Corollaire 39.** *Soit $\sigma : X \to X$ un automorphisme de $X$. Pour toute forme différentielle $\omega \in \Omega^{1,0}(X)$, nous avons*

$$R_{\sigma^*\omega}(x, y) = R_\omega(\sigma(x), \sigma(y)) \qquad (x, y \in X). \tag{1.119}$$

*Démonstration.* Il suffit d'appliquer la proposition 36 avec $X = Y$, $\phi = \sigma$, $\omega_Y = \omega$ et $m = [\sigma(y)]$. □

**Corollaire 40.** *Soit $\phi : E \to E'$ une isogénie entre deux courbes elliptiques définies sur $\mathbf{C}$. Nous avons alors*

$$R_{E'}(\phi(P), 0) = \phi_* \sum_{Q \in \text{Ker}\,\phi} R_E(P + Q, 0) \qquad (P \in E). \qquad (1.120)$$

*En particulier, lorsque $E = E'$ et $\phi$ est la multiplication par $n \in \mathbf{Z} - \{0\}$, nous obtenons*

$$R_E(nP, 0) = n \sum_{Q \in E[n]} R_E(P + Q, 0) \qquad (P \in E), \qquad (1.121)$$

*où $E[n]$ désigne le groupe des points de $n$-torsion de $E$.*

*Démonstration.* Utilisons le corollaire 37 avec $X = E$, $Y = E'$, $l = [P]$ et $m = [0]$. Il vient

$$R_{E'}(\phi(P), 0) = \phi_* \sum_{Q \in \text{Ker}\,\phi} R_E(P, Q).$$

D'après la proposition 18, nous avons $R_E(P, Q) = R_E(P - Q, 0)$, ce qui montre (1.120). Lorsque $E = E'$ et $\phi$ est la multiplication par $n \in \mathbf{Z} - \{0\}$, les morphismes $\phi^*$ en (1.103) et $\phi_*$ en (1.117) ne sont autres que la multiplication par $n$, ce qui permet de conclure. □

Notons que le corollaire 40 joint à la proposition 26 permettent de redémontrer les relations de distribution pour le dilogarithme elliptique $D_E$.

## 1.5 Questions et perspectives

Nous concluons ce chapitre en proposant plusieurs pistes de recherche pouvant prolonger notre étude de la fonction $R_\omega$.

Les relations de Steinberg (1.31) et les relations issues des morphismes finis (1.118) peuvent être vues comme des équations fonctionnelles pour la fonction $R_X$. En existe-t-il d'autres ?

La proposition 35 permet de décrire le comportement des fonctions de Green vis-à-vis des correspondances. Plus précisément, elle entraîne le résultat suivant (comparer avec [33, (2.6), p. 329]).

**Proposition 41.** *Soit $T$ une correspondance entre deux surfaces de Riemann compactes connexes $X$ et $Y$ de genre $\geq 1$, définie par le diagramme*

$$\begin{array}{c} & Z & \\ {}^p\swarrow & & \searrow^q \\ X \leftarrow\!-\!-\!\stackrel{T}{-}\!-\!-\!\rightarrow & Y, \end{array}$$

*où $Z$ est une surface de Riemann compacte connexe, et $p$ et $q$ sont des morphismes finis. Posons*

$$T_* = q_* \circ p^* : \text{Div}(X) \to \text{Div}(Y) \quad et \quad T^* = p_* \circ q^* : \text{Div}(Y) \to \text{Div}(X).$$

*Pour tous diviseurs $l$ et $m$, respectivement sur $X$ et $Y$, de degré $0$ et vérifiant $(p^*l) \cap (q^*m) = \emptyset$, nous avons l'égalité*

$$G_X(l, T^*m) = G_Y(T_*l, m). \tag{1.122}$$

Peut-on obtenir un énoncé similaire pour la fonction $R_\omega$ ?

Soit $X$ une courbe projective lisse définie sur $\mathbf{Q}$, géométriquement connexe et de genre $g \geq 1$. Notons $X(\mathbf{C})$ la surface de Riemann des points complexes de $X$. Les conjectures de Beĭlinson [5, 52] prédisent que la valeur spéciale $L(H^1(X), 2)$ s'exprime en termes du régulateur de Beĭlinson vu comme morphisme $K_2^{(2)}(X) \to \mathrm{Hom}_{\mathbf{C}}(\Omega^{(1,0)}(X(\mathbf{C})), \mathbf{C})$. D'après la proposition 17, les valeurs prises par ce morphisme sont combinaisons $\mathbf{Q}$-linéaires de valeurs de la fonction $R_{X(\mathbf{C})}$. Peut-on prédire quelle(s) combinaison(s) linéaire(s) prendre dans l'expression de $L(H^1(X), 2)$ ? Il serait intéressant d'écrire une version (conjecturale) précise de ce lien, à la lumière du cas elliptique, qui est traité dans [32] et [76].

La fonction $R_X$ a-t-elle droit au titre de dilogarithme associé à $X$ ? Cette question en apparence naïve peut en réalité se formuler précisément. Pour toute surface de Riemann compacte $X$ et tout point-base $x_0 \in X$, il est possible de définir un polylogarithme associé au couple $(X, x_0)$. Cet objet appartient à la catégorie des variations de structures de Hodge sur $X - \{x_0\}$. Nous renvoyons à [8] pour une définition. Dans quelle mesure la fonction $R_X$ permet-elle de construire une partie du polylogarithme mentionné ci-dessus ?

# Chapitre 2

# Extension à la jacobienne

Soit $X$ une surface de Riemann compacte connexe non vide, de genre $g \geq 1$. Dans le premier chapitre, nous avons introduit une fonction dilogarithme

$$R_X : X \times X \to V := \mathrm{Hom}_{\mathbf{C}}(\Omega^{1,0}(X), \mathbf{C}). \tag{2.1}$$

Notons $J = \mathrm{Jac}\, X$ la jacobienne de $X$ ; c'est une variété abélienne complexe de dimension $g$. Nous disposons d'une application

$$X \times X \to J \tag{2.2}$$
$$(x, y) \mapsto x - y$$

où $x - y \in J$ désigne la classe du diviseur $[x] - [y]$ de degré 0 sur $X$. L'application (2.2) ne dépend pas du choix d'un point-base de $X$. Une question se pose naturellement : est-il possible de prolonger la fonction $R_X$ à $J$ ? Nous allons voir qu'une tel prolongement est effectivement possible, dans le sens suivant.

**Théorème 6.** *Il existe une fonction $R_J : J \to V$ et une fonction continue $\Phi_X : X \to V$ telles que*

$$R_J(x - y) = R_X(x, y) + \Phi_X(x) - \Phi_X(y) \qquad (x, y \in X). \tag{2.3}$$

*Remarque* 42. Dans le cas où $X$ est une courbe elliptique $E$, nous avons un isomorphisme canonique $E \cong J$ (prendre $y = 0$ dans (2.2)). Le théorème résulte alors immédiatement de la proposition 24 : nous pouvons prendre $\Phi_E = 0$ et $R_J$ donnée par (1.61). En fait, nous verrons que la démonstration du théorème 6 fournit précisément ces fonctions $\Phi_E$ et $R_J$ (voir la remarque 55). Dans le cas elliptique, la fonction $R_J$ définie de manière générale à la section 2.2 coïncide donc avec le dilogarithme elliptique :

$$R_J(P - Q) = R_E(P, Q) \qquad (P, Q \in E). \tag{2.4}$$

Expliquons brièvement l'idée ayant donné naissance à la fonction $R_J$. Commençons par considérer une courbe elliptique $E$ définie sur $\mathbf{C}$. Bloch [12] a construit des éléments dans le groupe de $K$-théorie algébrique $K_2(E) \otimes \mathbf{Q}$ à partir des points de torsion de $E$. Plus précisément, il associe à tout point de torsion $P \in E$ un élément $\gamma_P \in K_2(E) \otimes \mathbf{Q}$. Un aspect important de cette construction est que le régulateur de l'élément $\gamma_P$ est essentiellement donné par le dilogarithme elliptique $D_E$ évalué en $P$. L'ensemble des points de torsion étant dense dans $E$, la fonction continue $D_E : E \to \mathbf{R}$ est caractérisée par la collection des régulateurs associés aux éléments $\gamma_P$. Nous nous sommes inspirés de cette remarque pour construire la fonction $R_J$, dans le cas où $J$ est la jacobienne d'une surface de Riemann compacte $X$. Nous commençons par construire des éléments dans le groupe $K_2^{(2)}(J)$ à partir des points de $X$ qui sont de torsion dans $J$. La détermination de l'image de ces éléments par le régulateur de Beĭlinson nous amène alors à une définition naturelle pour $R_J$.

Voici maintenant le plan du chapitre. Dans la première section, nous considérons une courbe $X$ projective lisse définie sur $\mathbf{Q}$, et nous construisons des éléments dans le groupe $K_2^{(2)}(J) = H^2_{\mathcal{M}}(J, \mathbf{Q}(2))$ associé à la jacobienne $J$ de $X$. L'intérêt arithmétique de ces éléments apparaîtra dans le chapitre 3. Dans les sections 2.2 à 2.4, nous supposons que $X$ est une surface de Riemann compacte ; $J$ est donc une variété abélienne complexe. Dans la section 2.2, nous introduisons une nouvelle fonction $R_J$, définie sur $J$ et à valeurs dans un espace vectoriel complexe de dimension finie. Dans la section 2.3, nous établissons un lien entre cette fonction et la fonction $R_X$ du premier chapitre. Ce lien utilise de façon essentielle la description de $J$ en termes de la puissance symétrique $g$-ième de $X$, où $g$ est le genre de $X$. Dans la section 2.4, nous indiquons quelques propriétés supplémentaires de la fonction $R_J$. Dans la section 2.6, nous supposons à nouveau que $X$ est une courbe projective lisse définie sur $\mathbf{Q}$, et nous calculons le régulateur de Beĭlinson des éléments construits dans la première section, au moyen de la fonction $R_J$. Nous terminons le chapitre en mentionnant brièvement deux axes de recherche qui nous semblent importants. Le premier concerne la généralisation éventuelle de la construction précédente aux points de torsion de $J$, voire d'une variété abélienne. Le second soulève la question de l'existence d'un lien entre la fonction $R_J$ et les polylogarithmes abéliens définis par Wildeshaus [77].

## 2.1 Construction d'éléments dans le $K_2$ d'une jacobienne $J$

Soit $X$ une courbe projective, lisse, définie sur $\mathbf{Q}$, et géométriquement connexe[1]. Nous supposons que le genre $g$ de $X$ est $\geq 1$. Soit $J = \mathrm{Pic}^0 X$ la jacobienne de $X$ ; c'est une variété abélienne définie sur $\mathbf{Q}$, de dimension $g$.

Nous allons nous intéresser au groupe de $K$-théorie de Quillen $K_2^{(2)}(J)$, défini comme sous-espace propre de $K_2(J) \otimes \mathbf{Q}$ pour les opérations d'Adams [65, §3]. D'après les conjectures de Beĭlinson, ce groupe est relié à la valeur spéciale $L(J, 2) = L(H^1(X), 2)$, ce qui justifie une étude particulière de ce groupe. Notons $\mathbf{Q}(J)$ le corps des fonctions rationnelles de $J$. La localisation en $K$-théorie algébrique permet d'écrire une suite exacte [52, (7.4), p. 23]

$$0 \longrightarrow K_2^{(2)}(J) \overset{\eta}{\longrightarrow} K_2(\mathbf{Q}(J)) \otimes \mathbf{Q} \overset{\partial}{\longrightarrow} \bigoplus_{V \subset J} \mathbf{Q}(V)^* \otimes \mathbf{Q}, \qquad (2.5)$$

où la somme directe porte sur les hypersurfaces irréductibles $V$ de $J$ (points de codimension 1 du schéma $J$). L'application $\partial = \bigoplus \partial_V$, appelée *symbole modéré*, est définie par

$$\partial_V\{F, G\} = (-1)^{\mathrm{ord}_V(F)\,\mathrm{ord}_V(G)} \left(\frac{F^{\mathrm{ord}_V(G)}}{G^{\mathrm{ord}_V(F)}}\right)\Big|_V \qquad (F, G \in \mathbf{Q}(J)^*). \qquad (2.6)$$

Nous supposons que $X$ possède un point $\mathbf{Q}$-rationnel $x_0 \in X(\mathbf{Q})$. Nous notons $i = i_{x_0} : X \hookrightarrow J$

---

[1] La construction qui suit étant de nature géométrique, il serait plus naturel de considérer une courbe définie sur $\mathbf{C}$. Compte tenu des applications arithmétiques (notamment le calcul du régulateur de Beĭlinson), nous avons cependant préféré nous placer sur le corps des rationnels.

l'immersion fermée associée[2], définie sur les points $\overline{\mathbf{Q}}$-rationnels par

$$i : X(\overline{\mathbf{Q}}) \to J(\overline{\mathbf{Q}}) \qquad (2.7)$$
$$x \mapsto [x - x_0].$$

Nous considérerons souvent (2.7) comme une inclusion. Notons $J_{\text{tors}}$ l'ensemble des points de torsion de $J(\overline{\mathbf{Q}})$, et posons $X_{\text{tors}} = X(\overline{\mathbf{Q}}) \cap J_{\text{tors}}$. L'ensemble $X_{\text{tors}}$ contient $x_0$ et est muni d'une action du groupe de Galois $\operatorname{Gal}(\overline{\mathbf{Q}}/\mathbf{Q})$. Lorsque le genre $g$ de $X$ est $\geq 2$, l'ensemble $X_{\text{tors}}$ est fini, d'après la conjecture de Manin-Mumford, démontrée par Raynaud [54]. Dans le cas où $X$ est une courbe modulaire (complétée) et $x_0$ une pointe rationnelle de $X$, l'ensemble $X_{\text{tors}}$ contient l'ensemble des pointes, d'après le théorème de Manin-Drinfel'd [28, 29]. En nous inspirant d'une construction de Bloch [24, Thm 5.1], nous allons construire des éléments dans $K_2^{(2)}(J)$ à partir de diviseurs $\mathbf{Q}$-rationnels supportés dans $X_{\text{tors}}$ (proposition 44).

Notons $I_X$ le groupe des diviseurs (à coefficients entiers) sur $X_{\text{tors}}$, de degré 0 et invariants sous l'action de $\operatorname{Gal}(\overline{\mathbf{Q}}/\mathbf{Q})$. Notons $\Theta \subset J$ le diviseur thêta, défini comme l'image de la composition

$$X^{g-1} \xrightarrow{i \times \cdots \times i} J^{g-1} \xrightarrow{\Sigma} J. \qquad (2.8)$$

Pour tout $a \in J(\overline{\mathbf{Q}})$, notons $t_a : x \mapsto x + a$ la translation par $a$ dans $J(\overline{\mathbf{Q}})$, et posons $\Theta_a = \Theta + a \subset J(\overline{\mathbf{Q}})$.

**Lemme 43.** *Pour tout diviseur $l \in I_X$, il existe une fonction rationnelle $F \in \mathbf{Q}(J)^* \otimes \mathbf{Q}$ telle que*

$$\operatorname{div} F = \sum_{x \in (l)} \operatorname{ord}_x(l) \, \Theta_x. \qquad (2.9)$$

*Démonstration.* Considérons le diviseur de Weil

$$D = \sum_{x \in (l)} \operatorname{ord}_x(l) \, \Theta_x.$$

Puisque $l$ est $\mathbf{Q}$-rationnel, $D$ définit bien un diviseur de Weil de $J$. Nous allons montrer qu'il existe un entier $N \geq 1$ tel que $N \cdot D$ soit un diviseur principal. Nous avons une identification canonique [49, Thm 6.6]

$$J(\overline{\mathbf{Q}}) \xrightarrow{\cong} \operatorname{Pic}^0(J_{\overline{\mathbf{Q}}}) \qquad (2.10)$$
$$a \mapsto [\Theta_a - \Theta],$$

où $J_{\overline{\mathbf{Q}}}$ désigne l'extension des scalaires de $J$ à $\overline{\mathbf{Q}}$. Puisque le support de $l$ est inclus dans $X_{\text{tors}}$, il existe un entier $N \geq 1$ tel que

$$N \sum_{x \in (l)} \operatorname{ord}_x(l) \cdot x = 0$$

dans $J(\overline{\mathbf{Q}})$. L'application (2.10) étant un homomorphisme de groupes, nous en déduisons que $N \cdot D$ est un diviseur principal, soit $N \cdot D = \operatorname{div} F_0$ avec $F_0 \in \mathbf{Q}(J)^*$. La fonction rationnelle $F = F_0 \otimes \frac{1}{N}$ répond alors à la question. $\square$

---

[2]Plus généralement, il suffit de se donner un diviseur $\mathbf{Q}$-rationnel de degré 1 sur $X$.

**Proposition 44.** *Soient $l, m \in I_X$ des diviseurs et $F, G \in \mathbf{Q}(J)^* \otimes \mathbf{Q}$ des fonctions rationnelles associées respectivement à $l$ et $m$ par le procédé du lemme 43. Soit $M \subset J_{\text{tors}}$ un sous-groupe fini, défini sur $\mathbf{Q}$, contenant les supports des diviseurs $l$ et $m$. Alors*

$$\gamma = \frac{1}{\operatorname{Card} M} \sum_{a \in M} t_a^* \{F, G\} \tag{2.11}$$

*définit un élément de $K_2^{(2)}(J)$ qui ne dépend pas du choix de $F$ et $G$.*

*Remarques 45.*    1. Dans le cas d'une variété abélienne $A$ qui est la puissance d'une courbe elliptique, l'idée de considérer la moyenne relative à un sous-groupe fini de $A$ est présente dans la construction par Beĭlinson du *symbole d'Eisenstein* [6]. On pourra également se référer à [22, Sect. 8], [48, 1.3.1, p. 215], [61, 5, p. 309], [63] et [64].

2. Dokchitser, de Jeu et Zagier [27] ont construit des éléments dans le $K_2$ d'une courbe à partir de points de torsion de cette courbe. Il serait intéressant de comparer les deux constructions.

*Démonstration.* Notons $L$ le corps engendré par les points de $M$. Soit $J_L$ l'extension des scalaires de $J$ à $L$, et $L(J)$ le corps des fonctions rationnelles de $J_L$. Nous avons un diagramme commutatif

$$\begin{array}{ccccccc}
0 & \longrightarrow & K_2^{(2)}(J_L) & \longrightarrow & K_2(L(J)) \otimes \mathbf{Q} & \xrightarrow{\partial_{J_L}} & \bigoplus_{V \subset J_L} K_1(L(V)) \otimes \mathbf{Q} \\
& & \uparrow & & \uparrow & & \uparrow \alpha \\
0 & \longrightarrow & K_2^{(2)}(J) & \longrightarrow & K_2(\mathbf{Q}(J)) \otimes \mathbf{Q} & \xrightarrow{\partial_J} & \bigoplus_{V \subset J} K_1(\mathbf{Q}(V)) \otimes \mathbf{Q}.
\end{array} \tag{2.12}$$

Les translations $t_a$, $a \in M$ sont définies sur $L$. Posons

$$\gamma_L := \frac{1}{\operatorname{Card} M} \sum_{a \in M} t_a^* \{F, G\} \in K_2(L(J)) \otimes \mathbf{Q}.$$

Le sous-groupe $M$ étant défini sur $\mathbf{Q}$, l'élément $\gamma_L$ est invariant sous l'action du groupe de Galois $\operatorname{Gal}(L(J)/\mathbf{Q}(J)) = \operatorname{Gal}(L/\mathbf{Q})$. D'après [37, Théorème 2.2 (a)] appliqué à l'extension $L(J)/\mathbf{Q}(J)$, l'application naturelle

$$K_2(\mathbf{Q}(J)) \otimes \mathbf{Q} \to \bigl(K_2(L(J)) \otimes \mathbf{Q}\bigr)^{\operatorname{Gal}(L/\mathbf{Q})}$$

est un isomorphisme, donc $\gamma_L$ définit un élément $\gamma$ de $K_2(\mathbf{Q}(J)) \otimes \mathbf{Q}$. L'application $\alpha$ du diagramme (2.12) est injective. Pour montrer que $\gamma$ définit un élément de $K_2^{(2)}(J)$, il suffit donc de montrer que $\partial_{J_L}(\gamma_L) = 1$. Par définition de $F$ et $G$, on a $\partial_V \gamma_L = 1$ pour toute hypersurface $V$ de $J_L$ distincte des $\Theta_a$, $a \in M$. D'autre part, puisque $\gamma_L$ est invariant par les translations $t_a$, $a \in M$, nous avons

$$\partial_\Theta \gamma_L = t_a^* (\partial_{\Theta_a} \gamma_L) \qquad (a \in M). \tag{2.13}$$

Pour montrer $\partial_{J_L} \gamma_L = 1$, il suffit donc de montrer $\partial_\Theta \gamma_L = 1$, c'est-à-dire

$$\prod_{a \in M} t_a^* \bigl(\partial_{\Theta_a} \{F, G\}\bigr) = 1. \tag{2.14}$$

Posons $S = (l) \cup (m) \subset M$. Nous avons $\partial_{\Theta_a} \{F, G\} = 1$ pour $a \notin S$. Nous sommes finalement ramenés à montrer

$$H := \prod_{z \in S} t_z^* \big(\partial_{\Theta_z}\{F,G\}\big) = 1. \tag{2.15}$$

Lorsque $g = 1$, l'égalité (2.15) n'est autre que la loi de réciprocité de Weil [68, Chap III, Prop 6] appliquée au symbole $\{F,G\}$. Supposons maintenant $g \geq 2$. Par définition, nous avons

$$\partial_{\Theta_z}\{F,G\} = (-1)^{\operatorname{ord}_z(l)\operatorname{ord}_z(m)} \left(\frac{F^{\operatorname{ord}_z(m)}}{G^{\operatorname{ord}_z(l)}}\right)\bigg|_{\Theta_z} \qquad (z \in S), \tag{2.16}$$

d'où nous déduisons

$$\operatorname{div} \partial_{\Theta_z}\{F,G\} = \Big(\operatorname{ord}_z(m) \sum_{x \in (l)} \operatorname{ord}_x(l) \Theta_x - \operatorname{ord}_z(l) \sum_{y \in (m)} \operatorname{ord}_y(m) \Theta_y\Big) \cap \Theta_z.$$

En appliquant $t_z^*$, puis en prenant le produit sur $z \in S$, nous obtenons

$$\operatorname{div} H = \sum_{z \in S} \sum_{x \in S} \big(\operatorname{ord}_z(m) \operatorname{ord}_x(l) - \operatorname{ord}_z(l) \operatorname{ord}_x(m)\big) \cdot \Theta_{x-z} \cap \Theta$$

$$= \sum_{\substack{z,x \in S \\ z \neq x}} \big(\operatorname{ord}_z(m) \operatorname{ord}_x(l) - \operatorname{ord}_z(l) \operatorname{ord}_x(m)\big) \cdot \Theta_{x-z} \cap \Theta.$$

Remarquons que pour $z \neq x$, l'intersection $\Theta_{x-z} \cap \Theta$ définit bien un diviseur de Weil de $\Theta$. Soit $\omega \in \Omega^1(X)$ une forme différentielle algébrique non nulle, de diviseur

$$\operatorname{div} \omega = \sum_{j=1}^{2g-2} [x_j] \qquad \big(x_j \in X(\overline{\mathbf{Q}})\big), \tag{2.17}$$

et définissons

$$\kappa = \sum_{j=1}^{2g-2} i(x_j) \in J(\mathbf{Q}). \tag{2.18}$$

Le point $\kappa$ ne dépend pas du choix de $\omega$. Notons $W$ l'image de la composition

$$X^{g-2} \xrightarrow{i \times \cdots \times i} J^{g-2} \xrightarrow{\Sigma} J, \tag{2.19}$$

et posons $W' = \kappa - W$. Pour tout $a \in J(\overline{\mathbf{Q}})$, notons $W_a = t_a(W)$ et $W'_a = t_a(W')$. D'après [50, Lemma, p. 76] ou [11, Chap. 11, Prop. 9.1], nous avons

$$\Theta_{x-z} \cap \Theta = W_x + W'_{-z} \qquad \big(x, z \in X(\overline{\mathbf{Q}}),\ x \neq z\big). \tag{2.20}$$

Nous pouvons alors finir le calcul du diviseur de $H$

$$\operatorname{div} H = \sum_{\substack{z,x \in S \\ z \neq x}} \big(\operatorname{ord}_z(m) \operatorname{ord}_x(l) - \operatorname{ord}_z(l) \operatorname{ord}_x(m)\big) \cdot (W_x + W'_{-z})$$

$$= \sum_{z \in S} \sum_{x \in S} \big(\operatorname{ord}_z(m) \operatorname{ord}_x(l) - \operatorname{ord}_z(l) \operatorname{ord}_x(m)\big) \cdot (W_x + W'_{-z})$$

$$= 0,$$

puisque les diviseurs $l$ et $m$ sont de degré 0. En conséquence, la fonction $H$ est constante. Nous allons montrer qu'elle vaut 1 en l'évaluant en un point de $\Theta$ bien choisi. Soit $e \in \Theta(\overline{\mathbf{Q}})$, arbitraire pour l'instant. D'après [11, Thm 11.2.4], nous avons

$$\Theta = \kappa - \Theta. \tag{2.21}$$

Posons $e' = \kappa - e \in \Theta(\overline{\mathbf{Q}})$, et écrivons

$$e' = z_1 + z_2 + \cdots + z_{g-1} \qquad \bigl(z_i \in X(\overline{\mathbf{Q}})\bigr). \tag{2.22}$$

Pour $e$ n'appartenant pas à un fermé de codimension 1 dans $\Theta(\overline{\mathbf{Q}})$, nous avons

$$X(\overline{\mathbf{Q}}) \not\subset \Theta_{x-e} \qquad (x \in S). \tag{2.23}$$

Faisant cette hypothèse sur $e$, le théorème de Riemann [15, Thm III.5.1] donne alors

$$i^* \Theta_{x-e} = [x] + \sum_{i=1}^{g-1} [z_i] \qquad (x \in S). \tag{2.24}$$

Les fonctions $f = (t_e^* F)|_X$ et $g = (t_e^* G)|_X$ sont donc bien définies dans $\overline{\mathbf{Q}}(X)^* \otimes \mathbf{Q}$. Le diviseur de $f$ vaut

$$\operatorname{div} f = i^* \operatorname{div} t_e^* F = i^* \sum_{x \in (l)} \operatorname{ord}_x(l) \, \Theta_{x-e} = \sum_{x \in (l)} \operatorname{ord}_x(l) \cdot [x].$$

Nous voyons donc que $\operatorname{div} f = l$; de même $\operatorname{div} g = m$. Supposons encore que pour tout $z \in S$, le point $e$ n'appartient pas au support de la fonction rationnelle $t_z^*\bigl(\partial_{\Theta_z}\{F, G\}\bigr)$. Nous pouvons alors évaluer $H$ en $e$ par

$$\begin{aligned}
H(e) &= \prod_{z \in S} (-1)^{\operatorname{ord}_z(l) \operatorname{ord}_z(m)} \Bigl(\frac{(t_z^* F)^{\operatorname{ord}_z(m)}}{(t_z^* G)^{\operatorname{ord}_z(l)}}\Bigr)(e) \\
&= \prod_{z \in S} (-1)^{\operatorname{ord}_z(l) \operatorname{ord}_z(m)} \Bigl(\frac{(t_e^* F)^{\operatorname{ord}_z(m)}}{(t_e^* G)^{\operatorname{ord}_z(l)}}\Bigr)(z) \\
&= \prod_{z \in S} (-1)^{\operatorname{ord}_z(l) \operatorname{ord}_z(m)} \Bigl(\frac{f^{\operatorname{ord}_z(m)}}{g^{\operatorname{ord}_z(l)}}\Bigr)(z) \\
&= \prod_{z \in S} \partial_z \{f, g\} \\
&= 1
\end{aligned}$$

d'après la loi de réciprocité de Weil appliquée au symbole $\{f, g\}$. Ainsi $H = 1$ et $\gamma$ définit un élément de $K_2^{(2)}(J)$. Montrons que $\gamma$ ne dépend pas du choix de $F$ (le raisonnement est le même pour $G$). Remplaçons $F$ par $F' = \lambda F$, avec $\lambda \in \mathbf{Q}^* \otimes \mathbf{Q}$. Alors $\gamma$ est changé en

$$\begin{aligned}
\gamma' &= \gamma + \frac{1}{\operatorname{Card} M} \sum_{a \in M} t_a^* \{\lambda, G\} \\
&= \gamma + \frac{1}{\operatorname{Card} M} \Bigl\{\lambda, \prod_{a \in M} t_a^* G\Bigr\}.
\end{aligned}$$

Par définition de $G$, nous avons

$$\operatorname{div} \prod_{a \in M} t_a^* G = \sum_{a \in M} \sum_{y \in (m)} \operatorname{ord}_y(m) \Theta_{y-a} = 0.$$

Donc $\gamma$ est modifié par un élément de $K_2(\mathbf{Q}) \otimes \mathbf{Q} \subset K_2^{(2)}(J)$. Le groupe $K_2(\mathbf{Q})$ étant de torsion [73], nous en déduisons $\gamma' = \gamma$ dans $K_2^{(2)}(J)$. □

Dans la section 2.6, nous calculerons le régulateur de Beĭlinson de $\gamma$. Ce calcul nous amène à la conjecture suivante, concernant la dépendance de $\gamma$ en le sous-groupe $M \subset J_{\mathrm{tors}}$.

**Conjecture 46.** *L'élément $\gamma$ ne dépend pas du choix du sous-groupe $M$.*

La conjecture 46 est entraînée par la conjecture suivante.

**Conjecture 47.** *Soit $A$ une variété abélienne sur un corps de nombres et $a$ un point rationnel de $A$. La translation $t_a : A \to A$ induit l'identité sur $K_2^{(2)}(A)$.*

*Remarque* 48. Dans le cas d'une courbe elliptique, les conjectures 46 et 47 découlent de résultats non triviaux de Goncharov et Levin [32, Prop 3.10]. Indiquons que dans le cas général, la conjecture 47 semble être conséquence de l'injectivité (non démontrée) du régulateur de Beĭlinson.

## 2.2  Définition de la fonction $R_J$

Soit $X$ une surface de Riemann compacte connexe, de genre $g \geq 1$. La jacobienne $J$ de $X$ est un tore complexe de dimension $g$, donné par le conoyau de l'application (injective) d'Abel-Jacobi

$$\begin{aligned} H_1(X, \mathbf{Z}) &\to \operatorname{Hom}_{\mathbf{C}}(\Omega^{1,0}(X), \mathbf{C}) \\ \lambda &\mapsto \left( \omega \mapsto \int_\lambda \omega \right). \end{aligned} \qquad (2.25)$$

Choisissons une base symplectique $B$ de $H_1(X, \mathbf{Z})$. Elle induit [15, III.1.2, p. 183] une matrice de périodes $\Omega \in \mathcal{H}_g$ et un isomorphisme

$$J \cong \frac{\mathbf{C}^g}{\mathbf{Z}^g + \Omega \mathbf{Z}^g}. \qquad (2.26)$$

La *fonction thêta* (relative à $B$) est par définition [15, (III.4.1), p. 192]

$$\begin{aligned} \theta : \mathbf{C}^g &\to \mathbf{C} \\ z &\mapsto \sum_{n \in \mathbf{Z}^g} e^{\pi i \, {}^t n \cdot \Omega \cdot n + 2\pi i \, {}^t n \cdot z}. \end{aligned} \qquad (2.27)$$

En suivant [15, p. 192], posons $Y = \Im(\Omega)^{-1}$ et

$$\|\theta(z)\| = |\theta(z)| \cdot e^{-\pi \, {}^t \Im(z) \cdot Y \cdot \Im(z)} \qquad (z \in \mathbf{C}^g). \qquad (2.28)$$

L'expression (2.28) ne dépend que de la classe de $z$ modulo $\mathbf{Z}^g + \Omega \mathbf{Z}^g$, et définit donc une fonction $\|\theta\| : J \to \mathbf{R}$ à valeurs positives, qui dépend a priori de $B$.

**Définition 49.** *La fonction $R_J$ associée à $J$ est définie par*

$$R_J : J \to \mathrm{Hom}_{\mathbf{C}}(\Omega^{g,g-1}(J), \mathbf{C})$$
$$x \mapsto \left( \alpha \mapsto \int_{u \in J} \log\|\theta(u)\| \cdot \overline{\partial}_u \log\|\theta(u-x)\| \wedge \alpha_u \right). \tag{2.29}$$

Montrons que la fonction $R_J$ est bien définie et ne dépend pas du choix de la base symplectique $B$. Le lieu des zéros de la fonction holomorphe $\theta$ est invariant par translation par les éléments de $\mathbf{Z}^g + \Omega\mathbf{Z}^g$, et définit donc une hypersurface irréductible $\Theta_0 \subset J$. La fonction $\log\|\theta\|$ est de classe $\mathcal{C}^\infty$ sur $J - \Theta_0$ et a une singularité logarithmique le long de $\Theta_0$. Nous en déduisons que l'intégrale (2.29) converge absolument. Soit maintenant $B'$ une autre base symplectique de $H_1(X, \mathbf{Z})$. Notons $\theta'$ la fonction thêta relative à $B'$, et $R'_J$ la fonction définie au moyen de $\theta'$. D'après la formule de transformation pour les fonctions thêta [11, 8.6.1], il existe un point de 2-torsion $a \in J[2]$ et une constante $C \in \mathbf{R}$ tels que

$$\log\|\theta'(z)\| = \log\|\theta(z-a)\| + C \qquad (z \in J). \tag{2.30}$$

Notons $t_a : u \mapsto u + a$ la translation par $a$ dans $J$. Nous avons

$$\begin{aligned}
\langle R'_J(x), \alpha \rangle &= \int_{u \in J} \big(\log\|\theta(u-a)\| + C\big) \cdot \overline{\partial}_u \log\|\theta(u-x-a)\| \wedge \alpha_u \\
&= \int_{u \in J} \log\|\theta(u)\| \cdot \overline{\partial}_u \log\|\theta(u-x)\| \wedge t_a^* \alpha_u \\
&\quad + \int_{u \in J} d\big(C \log\|\theta(u-x-a)\|\alpha_u\big) \\
&= \langle R_J(x), \alpha \rangle
\end{aligned}$$

puisque $\alpha$ est invariante par translation, et d'après la formule de Stokes. Ainsi $R'_J = R_J$, ce qui montre l'indépendance de $R_J$ vis-à-vis du choix de $B$.

*Remarque* 50. La définition ci-dessus vaut en fait pour toute variété abélienne complexe $A$ principalement polarisée. En effet, la polarisation principale de $A$ induit une forme bilinéaire alternée sur $H_1(A, \mathbf{Z})$. Le choix d'une base symplectique de $H_1(A, \mathbf{Z})$ pour cette forme alternée induit alors une fonction $\theta$ sur le revêtement universel de $A$. Il est alors possible de montrer que la définition (2.29) ne dépend pas du choix de la base symplectique, et donne donc lieu à une fonction $R_A$ bien définie sur $A$, à valeurs dans $\mathrm{Hom}_{\mathbf{C}}(\Omega^{d,d-1}(A), \mathbf{C})$, où $d$ est la dimension de $A$.

En vue de comparer les fonctions $R_X$ et $R_J$ (voir la section suivante), nous avons besoin de considérer $R_J$ comme une fonction à valeurs dans l'espace vectoriel complexe $\mathrm{Hom}_{\mathbf{C}}(\Omega^{1,0}(X), \mathbf{C})$. Pour cela, il suffit de définir un isomorphisme entre $\Omega^{1,0}(X)$ et $\Omega^{g,g-1}(J)$. Comme en (2.7), le choix d'un point-base $x_0 \in X$ induit une immersion fermée $i : X \hookrightarrow J$. L'application

$$i^* : \Omega^{1,0}(J) \xrightarrow{\cong} \Omega^{1,0}(X) \tag{2.31}$$

est un isomorphisme qui ne dépend pas du choix de $x_0$, ce qui nous permet d'identifier les espaces $\Omega^{1,0}(X)$ et $\Omega^{1,0}(J)$. Soit $\omega_{\Theta_0} \in \Omega^{1,1}(J)$ la forme différentielle duale de Poincaré de l'hypersurface $\Theta_0$, définie par la propriété

$$\int_J \omega_{\Theta_0} \wedge \beta = \int_{\Theta_0} \beta \qquad (\beta \in \Omega^{g-1,g-1}(J)). \tag{2.32}$$

La forme différentielle $\omega_{\Theta_0}$ ne dépend pas du choix de $B$. La formule de Stokes permet en fait de montrer

$$\omega_{\Theta_0} = i \sum_{j=1}^{g} \omega_j \wedge \overline{\omega_j}. \tag{2.33}$$

Remarquons au passage l'identité $i^*\omega_{\Theta_0} = g \operatorname{vol}_X$, qui fournit une interprétation de la forme volume $\operatorname{vol}_X$.

*Notation.* Nous convenons d'identifier $\Omega^{1,0}(X)$ et $\Omega^{g,g-1}(J)$ au moyen de l'isomorphisme

$$\Omega^{1,0}(X) \xrightarrow{\cong} \Omega^{g,g-1}(J) \tag{2.34}$$
$$\omega \mapsto \frac{1}{(g-1)!} \cdot \omega_{\Theta_0}^{g-1} \wedge \omega.$$

Nous notons encore $R_J : J \to \operatorname{Hom}_{\mathbf{C}}(\Omega^{1,0}(X), \mathbf{C})$ la fonction qui en résulte.

## 2.3 Lien entre les fonctions $R_J$ et $R_X$

Le but de cette section est d'établir un lien entre la fonction $R_J$ et la fonction dilogarithme $R_X$, lorsque $J$ est la jacobienne de $X$. Dans un premier temps, nous relions les fonctions $\log\|\theta\|$ et $G_X$ (la première joue le rôle de fonction de Green pour l'hypersurface $\Theta_0$). Ce lien est classique, voir par exemple [40, Chap. 13, Thm 5.2]. Dans un second temps, nous relions les fonctions $R_J$ et $R_X$, en utilisant la description de $J$ en termes de la puissance symétrique $g$-ième de $X$.

Pour les besoins de la démonstration, nous fixons une base symplectique $B$ de $H_1(X, \mathbf{Z})$ pour le produit d'intersection, ce qui induit une fonction $\theta$ sur le revêtement universel de $J$, une hypersurface $\Theta_0 \subset J$ définie comme le lieu des zéros de $\theta$, et une fonction $\log\|\theta\| : J - \Theta_0 \to \mathbf{R}$ de classe $\mathcal{C}^\infty$, définie grâce à (2.28). Nous choisissons également une base orthonormale $(\omega_1, \ldots, \omega_g)$ de $\Omega^{1,0}(X) \cong \Omega^{1,0}(J)$ pour le produit scalaire hermitien (1.3).

**Lemme 51.** *Sur $J - \Theta_0$, nous avons l'égalité*

$$\partial\overline{\partial}\log\|\theta\| = -\pi \sum_{j=1}^{g} \omega_j \wedge \overline{\omega_j}. \tag{2.35}$$

*Démonstration.* Nous conservons la notation (2.26). Notons $(z_1, \ldots, z_g)$ les coordonnées sur $\mathbf{C}^g$. Soit $\widetilde{\Theta}_0$ l'image réciproque de $\Theta_0$ dans $\mathbf{C}^g$. Puisque la fonction $z \mapsto \log|\theta(z)|$ est harmonique sur $\mathbf{C}^g - \widetilde{\Theta}_0$, et d'après la définition (2.28) de $\|\theta\|$, nous avons sur $J - \Theta_0$

$$\begin{aligned}
\partial\overline{\partial}\log\|\theta(z)\| &= \partial\overline{\partial}\bigl(-\pi\,{}^t\Im(z) \cdot Y \cdot \Im(z)\bigr) \\
&= -\pi \sum_{k,l=1}^{g} Y_{k,l} \cdot \partial\overline{\partial}\Bigl(\frac{z_k - \overline{z_k}}{2i} \cdot \frac{z_l - \overline{z_l}}{2i}\Bigr) \\
&= \frac{\pi}{4} \sum_{k,l=1}^{g} Y_{k,l} \cdot \partial\overline{\partial}(z_k z_l - z_k \overline{z_l} - \overline{z_k} z_l + \overline{z_k}\,\overline{z_l}) \\
&= \frac{\pi}{4} \sum_{k,l=1}^{g} Y_{k,l} \cdot (-dz_k \wedge d\overline{z_l} - dz_l \wedge d\overline{z_k}) \\
&= -\frac{\pi}{2} \sum_{k,l=1}^{g} Y_{k,l} \cdot dz_k \wedge d\overline{z_l},
\end{aligned}$$

la dernière égalité résultant du fait que la matrice $Y$ est symétrique. Soit $A = (a_{k,m})_{1 \leq k,m \leq g} \in \mathcal{M}_g(\mathbf{C})$ la matrice définie par

$$dz_k = \sum_{m=1}^{g} a_{k,m} \omega_m \qquad (1 \leq k \leq g). \tag{2.36}$$

Notons que la matrice $A$ est inversible. Nous avons alors

$$(dz_k, dz_l) = \sum_{m,n=1}^{g} a_{k,m} \overline{a_{l,n}} (\omega_m, \omega_n) = \sum_{m=1}^{g} a_{k,m} \overline{a_{l,m}} = (A {}^t\overline{A})_{k,l} \qquad (1 \leq k,l \leq g).$$

Écrivons maintenant $B = (a_1, \ldots, a_g, b_1, \ldots, b_g)$. Notons $\delta$ le symbole de Kronecker. D'après [15, A.2.5], nous avons également

$$\begin{aligned}(dz_k, dz_l) &= i \int_X dz_k \wedge d\overline{z_l} \\ &= i \sum_{j=1}^{g} \int_{a_j} dz_k \cdot \int_{b_j} d\overline{z_l} - \int_{a_j} d\overline{z_l} \cdot \int_{b_j} dz_k \\ &= i \sum_{j=1}^{g} \delta_{j,k} \cdot \overline{\Omega}_{l,j} - \Omega_{k,j} \cdot \delta_{j,l} \\ &= 2 \cdot \Im(\Omega)_{k,l},\end{aligned}$$

la dernière égalité résultant du fait que $\Omega$ est symétrique. Par identification, nous avons $A {}^t\overline{A} = 2Y^{-1}$, d'où ${}^tA = 2\overline{A}^{-1} Y^{-1}$ et ${}^t A Y \overline{A} = 2 I_g$. Nous avons enfin

$$\begin{aligned}-\frac{\pi}{2} \sum_{k,l=1}^{g} Y_{k,l} \cdot dz_k \wedge d\overline{z_l} &= -\frac{\pi}{2} \sum_{k,l=1}^{g} Y_{k,l} \sum_{m,n=1}^{g} a_{k,m} \overline{a_{l,n}} \cdot \omega_m \wedge \overline{\omega_n} \\ &= -\frac{\pi}{2} \sum_{m,n=1}^{g} ({}^t A Y \overline{A})_{m,n} \cdot \omega_m \wedge \overline{\omega_n} \\ &= -\pi \sum_{m,n=1}^{g} \delta_{m,n} \cdot \omega_m \wedge \overline{\omega_n},\end{aligned}$$

ce qui montre bien $\partial \overline{\partial} \log \|\theta\| = -\pi \sum_{j=1}^{g} \omega_j \wedge \overline{\omega_j}$ sur $J - \Theta_0$. □

La proposition suivante établit le lien entre la fonction $\log \|\theta\|$ et la fonction de Green sur $X$. Nous remarquons que cette proposition ne fait pas référence à un point-base particulier de $X$.

**Proposition 52.** *Soit $e \in \Theta_0$ tel que la fonction*

$$(x,y) \in X \times X \mapsto \|\theta(e + x - y)\| \tag{2.37}$$

*n'est pas identiquement nulle. Alors il existe des $(g-1)$-uplets $(x_1, \ldots, x_{g-1})$ et $(y_1, \ldots, y_{g-1})$ de points de $X$, bien déterminés à permutation près, et une constante $C_e \in \mathbf{R}$, tels que nous ayons*

$$\log \|\theta(e + x - y)\| = G_X(x,y) + \sum_{i=1}^{g-1} G_X(x, x_i) + \sum_{j=1}^{g-1} G_X(y, y_j) + C_e, \tag{2.38}$$

*pour tout couple $(x,y) \in X \times X$ vérifiant $x \neq y$, $x \neq x_i$ et $y \neq y_j$.*

*Démonstration.* Le premier rôle revient à l'ensemble

$$F_e = \{(x,y) \in X \times X \mid \|\theta(e+x-y)\| = 0\}, \tag{2.39}$$

qui est un fermé algébrique strict de $X \times X$. Nous choisissons un point-base $x_0 \in X$ et notons $i : X \hookrightarrow J$ l'immersion fermée associée (prendre $y = x_0$ dans (2.2)). Nous verrons $i$ comme une inclusion. Pour tout $a \in J$, notons $\Theta_a = \Theta_0 + a \subset J$ le translaté de $\Theta_0$ par $a$. Rappelons le théorème de Riemann [15, Thm III.5.1] : il existe $\Delta \in J$ tel que

$$\Theta_\Delta = \Theta_0 + \Delta = \{z_1 + \cdots + z_{g-1} \in J \mid z_1, \ldots, z_{g-1} \in X\}. \tag{2.40}$$

Le point $\Delta$ dépend de la base symplectique $B$ et du point-base $x_0$. Prendre comme point-base $x_0' \in X$ revient à changer $\Delta$ en

$$\Delta' = \Delta - (g-1)(x_0' - x_0). \tag{2.41}$$

L'identité (2.40) entraîne l'existence d'un $(g-1)$-uplet $(y_1, \ldots, y_{g-1})$ de points de $X$ tel que

$$e + \Delta = y_1 + \cdots + y_{g-1}. \tag{2.42}$$

Par parité de $\Theta_0$, nous avons $-e \in \Theta_0$, d'où l'existence d'un $(g-1)$-uplet $(x_1, \ldots, x_{g-1})$ de points de $X$ tel que

$$-e + \Delta = x_1 + \cdots + x_{g-1}. \tag{2.43}$$

D'après [51, Chap. II, Lemma 3.4], chacun de ces $(g-1)$-uplets est bien déterminé à permutation près, et l'identité (2.41) montre qu'ils ne dépendent pas du point-base $x_0$. Toujours d'après [51, Chap. II, Lemma 3.4], nous avons la description suivante de $F_e$ :

$$\|\theta(e+x-y)\| = 0 \quad \Leftrightarrow \quad \begin{cases} x = y \\ \text{ou} \quad x \in \{x_1, \ldots, x_{g-1}\} \\ \text{ou} \quad y \in \{y_1, \ldots, y_{g-1}\}. \end{cases} \tag{2.44}$$

Fixons maintenant $y \in X$ avec $y \notin \{y_1, \ldots, y_{g-1}\}$. Par hypothèse sur $y$, la courbe $X$ n'est pas incluse dans $\Theta_{y-e}$. Toujours d'après le théorème de Riemann [15, Thm III.5.1], l'image réciproque par $i$ de l'hypersurface $\Theta_{y-e}$ est donnée par le diviseur

$$i^*\Theta_{y-e} = [y] + \sum_{i=1}^{g-1}[x_i]. \tag{2.45}$$

En conséquence, la fonction $\phi_y$ définie par

$$\phi_y(x) = \log\|\theta(e+x-y)\| - G_X(x,y) - \sum_{i=1}^{g-1} G_X(x,x_i) \tag{2.46}$$

s'étend en une fonction de classe $\mathcal{C}^\infty$ sur $X$ : cela résulte de la propriété (1.10) de la fonction de Green $G_X$. Montrons que $\phi_y$ est constante sur $X$. Définissons

$$\begin{aligned} t_{e-y} : J &\to J \\ a &\mapsto a + e - y. \end{aligned} \tag{2.47}$$

D'après (1.9) et le lemme 51, nous avons

$$\begin{aligned}
\partial\overline{\partial}\,\phi_y &= \left(\partial\overline{\partial}\,i^*t^*_{e-y}\log\|\theta\|\right) - \pi i g\,\mathrm{vol}_X \\
&= \left(i^*t^*_{e-y}\partial\overline{\partial}\log\|\theta\|\right) - \pi i g\,\mathrm{vol}_X \\
&= -\pi \cdot i^*t^*_{e-y}\Big(\sum_{j=1}^{g}\omega_j\wedge\overline{\omega_j}\Big) - \pi i g\,\mathrm{vol}_X \\
&= -\pi\sum_{j=1}^{g}\omega_j\wedge\overline{\omega_j} - \pi i g\,\mathrm{vol}_X \\
&= 0,
\end{aligned}$$

la dernière égalité résultant de la définition (1.4) de $\mathrm{vol}_X$. Donc $\phi_y$ est constante (réelle) sur $X$. Déterminons cette constante. Par parité de $\theta$, nous avons

$$\|\theta(e + x - y)\| = \|\theta(-e + y - x)\| \qquad (x,y\in X).$$

En utilisant ce que nous avons déjà démontré, il existe, pour tout $x \in X - \{x_1,\ldots,x_{g-1}\}$, une constante $\psi_x \in \mathbf{R}$ telle que

$$\log\|\theta(e + x - y)\| = G_X(y,x) + \sum_{j=1}^{g-1} G_X(y,y_j) + \psi_x \tag{2.48}$$

pour tout $y \in X$ tel que $y \notin \{x,y_1,\ldots,y_{g-1}\}$. La conjonction des équations (2.46) et (2.48) entraîne aisément l'égalité désirée, pour tout $(x,y) \notin F_e$. $\square$

*Remarque* 53. L'ensemble

$$F = \{e \in \Theta_0 \mid (2.37)\text{ est identiquement nulle}\} \tag{2.49}$$

est un fermé algébrique strict de $\Theta_0$ (cela résulte du fait que $J = \Theta_0 + X$), ce qui assure que la proposition précédente est non vide. Il serait intéressant de préciser la dépendance de la constante $C_e$ de la formule (2.38), en le point $e \in \Theta_0 - F$.

Nous sommes maintenant en mesure d'établir le lien entre les fonctions $R_J$ et $R_X$. Rappelons que nous considérons $R_J$ comme une fonction à valeurs dans $\mathrm{Hom}_{\mathbf{C}}(\Omega^{1,0}(X),\mathbf{C})$, grâce à l'isomorphisme (2.34).

*Démonstration du théorème 6.* Définissons la fonction $\Phi_X$ par

$$\begin{aligned}
\Phi_X : X &\to \mathrm{Hom}_{\mathbf{C}}(\Omega^{1,0}(X),\mathbf{C}) \\
x &\mapsto \Big(\omega \mapsto -i\sum_{j=1}^{g}\int_{z\in X} R_{\omega_j}(x,z)\cdot\omega\wedge\overline{\omega_j}\Big),
\end{aligned} \tag{2.50}$$

et notons $\Phi_\omega = \langle\Phi_X,\omega\rangle$ pour $\omega \in \Omega^{1,0}(X)$. On vérifie que $\Phi_X$ ne dépend pas du choix de la base orthonormale $(\omega_k)_{1\leq k\leq g}$. Les fonctions $R_{\omega_j}$ étant continues sur la variété compacte $X\times X$, la fonction $\Phi_X$ est continue sur $X$. Nous allons montrer que la fonction $R_J$ satisfait (2.3). Il suffit de tester cette égalité sur chacune des formes différentielles $\omega_k$, $1 \leq k \leq g$. Notons $\mathcal{S}_g$ le groupe symétrique d'ordre $g$. De (2.33), nous déduisons

$$\omega_{\Theta_0}^{g-1} = i^{g-1} \sum_{\tau \in \mathcal{S}_g} (\omega_{\tau(1)} \wedge \overline{\omega_{\tau(1)}}) \wedge \cdots \wedge (\omega_{\tau(g-1)} \wedge \overline{\omega_{\tau(g-1)}})$$

$$= i^{g-1}(g-1)! \sum_{k=1}^{g} \bigwedge_{\substack{j=1 \\ j \neq k}}^{g} \omega_j \wedge \overline{\omega_j}.$$

L'image $\alpha_k$ de $\omega_k$ par l'isomorphisme (2.34) est donc donnée par

$$\alpha_k = i^{g-1} \cdot (\omega_1 \wedge \overline{\omega_1}) \wedge \cdots \wedge \overset{\overset{k-\text{ième facteur}}{\downarrow}}{\omega_k} \wedge \cdots \wedge (\omega_g \wedge \overline{\omega_g}) \in \Omega^{g,g-1}(J). \tag{2.51}$$

En utilisant la définition (2.29), il s'agit donc de montrer

$$\langle R_J(x-y), \alpha_k \rangle = R_{\omega_k}(x,y) + \Phi_{\omega_k}(x) - \Phi_{\omega_k}(y) \qquad (x, y \in X). \tag{2.52}$$

L'idée-clé est d'utiliser la description de $J$ en termes de la puissance symétrique $g$-ième de $X$. Fixons un point-base $x_0 \in X$. Nous verrons $X$ incluse dans $J$ au moyen de ce point-base. Le morphisme

$$\psi : X^g \to J \tag{2.53}$$
$$(z_1, \ldots, z_g) \mapsto z_1 + \cdots + z_g$$

est fini, de degré $(g!)$. Soit $\Delta \in J$ tel que (2.40) soit vérifié. Nous avons

$$\langle R_J(x-y), \alpha_k \rangle = \int_{u \in J} \log \|\theta(u)\| \cdot \overline{\partial}_u \log \|\theta(u-x+y)\| \wedge \alpha_k$$
$$= \int_{u \in J} \log \|\theta(u-\Delta-y)\| \cdot \overline{\partial}_u \log \|\theta(u-\Delta-x)\| \wedge \alpha_k$$

puisque le changement de variables $u \mapsto u - \Delta - y$ laisse invariante la forme différentielle $\alpha_k$. En utilisant le changement de variables $\psi$, nous obtenons

$$\langle R_J(x-y), \alpha_k \rangle = \frac{1}{g!} \int_{X^g} \log \|\theta(z_1 + \cdots + z_g - \Delta - y)\| \cdot \overline{\partial} \log \|\theta(z_1 + \cdots + z_g - \Delta - x)\| \wedge \psi^* \alpha_k,$$

où l'intégrale porte sur les variables $(z_1, \ldots, z_g) \in X^g$. Nous allons maintenant calculer $\psi^* \alpha_k$. Pour tout $1 \leq l \leq g$, notons $p_l : X^g \to X$ la projection sur le $l$-ème facteur. Puisque

$$\psi^* \omega_j = \sum_{l=1}^{g} p_l^* \omega_j \qquad (1 \leq j \leq g), \tag{2.54}$$

nous en déduisons

$$\psi^* \alpha_k = \sum_{\sigma \in \mathcal{S}_g} \sum_{\tau \in \mathcal{S}_g} \alpha_{\sigma, \tau}, \tag{2.55}$$

où nous avons posé

$$\alpha_{\sigma,\tau} = i^{g-1} \cdot (p_{\sigma(1)}^* \omega_1 \wedge p_{\tau(1)}^* \overline{\omega_1}) \wedge \cdots \wedge \overset{\overset{k-\text{ième facteur}}{\downarrow}}{p_{\sigma(k)}^* \omega_k} \wedge \cdots \wedge (p_{\sigma(g)}^* \omega_g \wedge p_{\tau(g)}^* \overline{\omega_g}). \tag{2.56}$$

Posons

$$I_{\sigma,\tau} = \int_{X^g} \log\|\theta(z_1 + \cdots + z_g - \Delta - y)\| \cdot \overline{\partial}\log\|\theta(z_1 + \cdots + z_g - \Delta - x)\| \wedge \alpha_{\sigma,\tau}. \qquad (2.57)$$

Remarquons que dans l'intégrale (2.57), nous pouvons remplacer le symbole $\overline{\partial}$ par une dérivée partielle par rapport à la variable $z_{\tau(k)}$. Dans la suite, nous poserons pour abréger

$$e = \sum_{\substack{l=1 \\ l \neq \tau(k)}}^{g} z_l - \Delta \in \Theta_0, \qquad (2.58)$$

$$z = z_{\tau(k)} \in X, \qquad (2.59)$$

de sorte que $I_{\sigma,\tau}$ s'écrit

$$I_{\sigma,\tau} = \int_{X^g} \log\|\theta(e + z - y)\| \cdot \overline{\partial}_z \log\|\theta(e + z - x)\| \wedge \alpha_{\sigma,\tau}. \qquad (2.60)$$

Nous reconnaissons ici des quantités similaires à celle qui intervient dans le membre de gauche de (2.38). L'ensemble des $(z_1, \ldots, z_g) \in X^g$ tels que $e$ vérifie les hypothèses de la proposition 52 est un ouvert Zariski non vide de $X^g$. Nous pouvons sans danger remplacer (2.60) par l'intégrale sur cet ouvert, et appliquer la proposition en question à l'intégrand. Dans la suite, nous supposerons donc toujours que $e$ vérifie les hypothèses de la proposition 52. Nous noterons $z'_1(e), \ldots, z'_{g-1}(e)$ les points de $X$, bien déterminés à permutation près, tels que

$$-e + \Delta = \sum_{j=1}^{g-1} z'_j(e). \qquad (2.61)$$

D'après la formule (2.38), nous pouvons écrire

$$\log\|\theta(e + z - y)\| = G_X(z,y) + \sum_{j=1}^{g-1} G_X(z, z'_j(e)) + \sum_{\substack{l=1 \\ l \neq \tau(k)}}^{g} G_X(y, z_l) + C_e \qquad (2.62)$$

$$\text{et } \overline{\partial}_z \log\|\theta(e + z - x)\| = \overline{\partial}_z G_X(z,x) + \overline{\partial}_z\left(\sum_{j=1}^{g-1} G_X(z, z'_j(e))\right). \qquad (2.63)$$

L'intégrale $I_{\sigma,\tau}$ se développe alors comme une somme de 8 termes. Quatre d'entre eux disparaissent grâce à la formule de Stokes :

$$\int_{X^g} \left(\sum_{\substack{l=1 \\ l \neq \tau(k)}}^{g} G_X(y, z_l) + C_e\right) \cdot \overline{\partial}_z \log\|\theta(e + z - x)\| \wedge \alpha_{\sigma,\tau}$$

$$= \int_{X^g} d\left(\left(\sum_{\substack{l=1 \\ l \neq \tau(k)}}^{g} G_X(y, z_l) + C_e\right) \cdot \log\|\theta(e + z - x)\| \cdot \alpha_{\sigma,\tau}\right)$$

$$= 0,$$

puisque la forme différentielle $\alpha_{\sigma,\tau}$ est "saturée" en les variables $z_l$, $l \neq \tau(k)$. D'autre part, pour $1 \leq j_1, j_2 \leq g-1$, la formule de Stokes montre que la quantité

$$\int_{X^g} G_X(z, z'_{j_1}(e)) \cdot \overline{\partial}_z G_X(z, z'_{j_2}(e)) \wedge \alpha_{\sigma,\tau}$$

est antisymétrique en $j_1$ et $j_2$. Donc

$$\int_{X^g} \left( \sum_{j=1}^{g-1} G_X(z, z'_j(e)) \right) \cdot \overline{\partial}_z \left( \sum_{j=1}^{g-1} G_X(z, z'_j(e)) \right) \wedge \alpha_{\sigma,\tau} = 0.$$

Il reste donc

$$\begin{aligned}
I_{\sigma,\tau} &= \int_{X^g} G_X(z,y) \cdot \overline{\partial}_z G_X(z,x) \wedge \alpha_{\sigma,\tau} \\
&+ \int_{X^g} \sum_{j=1}^{g-1} G_X(z, z'_j(e)) \cdot \overline{\partial}_z G_X(z,x) \wedge \alpha_{\sigma,\tau} \\
&+ \int_{X^g} G_X(z,y) \cdot \overline{\partial}_z \left( \sum_{j=1}^{g-1} G_X(z, z'_j(e)) \right) \wedge \alpha_{\sigma,\tau}.
\end{aligned} \tag{2.64}$$

Posons

$$R_{\sigma,\tau}(x,y) = \int_{X^g} G_X(z,y) \cdot \overline{\partial}_z G_X(z,x) \wedge \alpha_{\sigma,\tau} \qquad (x,y \in X) \tag{2.65}$$

$$\Phi_{\sigma,\tau}(x) = \int_{X^g} \sum_{j=1}^{g-1} G_X(z, z'_j(e)) \cdot \overline{\partial}_z G_X(z,x) \wedge \alpha_{\sigma,\tau} \qquad (x \in X). \tag{2.66}$$

D'après la formule de Stokes, le troisième terme de la somme (2.64) est égal à $-\Phi_{\sigma,\tau}(y)$. Nous avons donc

$$\langle R_J(x-y), \alpha_k \rangle = \frac{1}{g!} \sum_{\sigma \in \mathcal{S}_g} \sum_{\tau \in \mathcal{S}_g} \left( R_{\sigma,\tau}(x,y) + \Phi_{\sigma,\tau}(x) - \Phi_{\sigma,\tau}(y) \right). \tag{2.67}$$

Maintenant, occupons-nous de l'intégrale $R_{\sigma,\tau}$. Puisque cette intégrale est à variables séparées, elle se calcule aisément grâce au théorème de Fubini. Nous pouvons l'écrire comme un produit de $g$ intégrales étendues à $X$. Par orthogonalité des $\omega_j$, nous voyons immédiatement que $R_{\sigma,\tau} = 0$ lorsque $\sigma \neq \tau$. Lorsque $\sigma = \tau$, nous pouvons écrire

$$R_{\sigma,\sigma}(x,y) = \left( \int_{z \in X} G_X(z,y) \cdot \overline{\partial}_z G_X(z,x) \wedge \omega_k \right) \cdot i^{g-1} \prod_{\substack{j=1 \\ j \neq k}}^{g} \int_X \omega_j \wedge \overline{\omega_j}.$$

D'après (1.20) et (1.21), et puisque $i \int_X \omega_j \wedge \overline{\omega_j} = 1$, nous obtenons finalement

$$R_{\sigma,\sigma}(x,y) = R_{\omega_k}(x,y) \qquad (x,y \in X).$$

Il nous reste maintenant à calculer $\Phi_{\sigma,\tau}$, ce qui est sensiblement plus technique. Par permutation des variables dans l'intégrale (2.66), nous voyons que $\Phi_{\sigma,\tau} = \Phi_{\sigma_0\sigma,\sigma_0\tau}$ pour tout $\sigma_0 \in \mathcal{S}_g$. Cela implique l'identité

$$\frac{1}{g!}\sum_{\sigma\in\mathcal{S}_g}\sum_{\tau\in\mathcal{S}_g}\Phi_{\sigma,\tau} = \frac{1}{(g-1)!}\sum_{\sigma\in\mathcal{S}_g}\sum_{\substack{\tau\in\mathcal{S}_g\\\tau(k)=g}}\Phi_{\sigma,\tau}. \tag{2.68}$$

Considérons le morphisme

$$\rho : X^g \to \Theta_0 \times X \tag{2.69}$$

$$(z_1,\ldots,z_g) \mapsto (e,z) = \left(\sum_{l=1}^{g-1} z_l - \Delta,\, z_g\right)$$

et le morphisme $\rho'$ défini par le diagramme commutatif

$$\begin{array}{ccc} X^g \xrightarrow{\rho} \Theta_0 \times X & & (e,z) \\ {\scriptstyle \rho'}\searrow \quad \downarrow {\scriptstyle -\mathrm{id}\times\mathrm{id}} & & \downarrow \\ \Theta_0 \times X & & (-e,z). \end{array} \tag{2.70}$$

Ces deux morphismes sont finis, de degré $(g-1)!$. Notons $\mathrm{pr}_1 : \Theta_0 \times X \to \Theta_0$ et $\mathrm{pr}_2 : \Theta_0 \times X \to X$ les deux projections naturelles. Toute forme différentielle $\Omega$ sur $J$ induit par restriction une forme différentielle sur $\Theta_0$, que nous notons encore $\Omega$.

**Lemme 54.** *La forme différentielle $\nu_k$ sur $\Theta_0 \times X$ définie par*

$$\nu_k = \mathrm{pr}_1^*\left(\bigwedge_{\substack{l=1\\l\neq k}}^{g}\omega_l\wedge\overline{\omega_l}\right)\wedge \mathrm{pr}_2^*\omega_k - \sum_{\substack{j=1\\j\neq k}}^{g}\mathrm{pr}_1^*\left(\omega_k\wedge\overline{\omega_j}\wedge\bigwedge_{\substack{l=1\\l\neq j,k}}^{g}\omega_l\wedge\overline{\omega_l}\right)\wedge \mathrm{pr}_2^*\omega_j \tag{2.71}$$

*vérifie*

$$\rho^*\nu_k = \frac{1}{i^{g-1}}\sum_{\sigma\in\mathcal{S}_g}\sum_{\substack{\tau\in\mathcal{S}_g\\\tau(k)=g}}\alpha_{\sigma,\tau}. \tag{2.72}$$

*Démonstration.* Pour toute forme différentielle $\omega \in \Omega^{1,0}(X) = \Omega^{1,0}(J)$, nous avons les formules

$$\rho^*\mathrm{pr}_1^*\omega = \sum_{l=1}^{g-1} p_l^*\omega \tag{2.73}$$

$$\rho^*\mathrm{pr}_2^*\omega = p_g^*\omega. \tag{2.74}$$

Nous en déduisons par exemple

$$\rho^*\mathrm{pr}_1^*\left(\bigwedge_{\substack{l=1\\l\neq k}}^{g}\omega_l\wedge\overline{\omega_l}\right) = \sum_\sigma\sum_\tau\bigwedge_{\substack{l=1\\l\neq k}}^{g}p_{\sigma(l)}^*\omega_l \wedge p_{\tau(l)}^*\overline{\omega_l},$$

où $\sigma$ et $\tau$ parcourent les applications injectives de $\{1,\ldots,\widehat{k},\ldots,g\}$ dans $\{1,\ldots,g-1\}$, c'est-à-dire les permutations de $\mathcal{S}_g$ vérifiant $\sigma(k) = \tau(k) = g$. Nous trouvons ainsi que l'image réciproque par $\rho$ du premier terme de (2.71) vaut

$$\frac{1}{i^{g-1}} \sum_{\substack{\sigma \in \mathcal{S}_g \\ \sigma(k)=g}} \sum_{\substack{\tau \in \mathcal{S}_g \\ \tau(k)=g}} \alpha_{\sigma,\tau}. \tag{2.75}$$

En procédant de la même manière, nous voyons que l'image réciproque par $\rho$ du second terme de (2.71) est

$$\begin{aligned}
&-\sum_{\substack{j=1 \\ j\neq k}}^{g} \sum_{\substack{\sigma \in \mathcal{S}_g \\ \sigma(j)=g}} \sum_{\substack{\tau \in \mathcal{S}_g \\ \tau(k)=g}} p^*_{\sigma(k)}\omega_k \wedge p^*_{\tau(j)}\overline{\omega_j} \wedge \left( \bigwedge_{\substack{l=1 \\ l\neq j,k}}^{g} p^*_{\sigma(l)}\omega_l \wedge p^*_{\tau(l)}\overline{\omega_l} \right) \wedge p^*_g \omega_j \\
&= \sum_{\substack{j=1 \\ j\neq k}}^{g} \sum_{\substack{\sigma \in \mathcal{S}_g \\ \sigma(j)=g}} \sum_{\substack{\tau \in \mathcal{S}_g \\ \tau(k)=g}} p^*_{\sigma(j)}\omega_j \wedge p^*_{\tau(j)}\overline{\omega_j} \wedge \left( \bigwedge_{\substack{l=1 \\ l\neq j,k}}^{g} p^*_{\sigma(l)}\omega_l \wedge p^*_{\tau(l)}\overline{\omega_l} \right) \wedge p^*_{\sigma(k)}\omega_k \\
&= \frac{1}{i^{g-1}} \sum_{\substack{j=1 \\ j\neq k}}^{g} \sum_{\substack{\sigma \in \mathcal{S}_g \\ \sigma(j)=g}} \sum_{\substack{\tau \in \mathcal{S}_g \\ \tau(k)=g}} \alpha_{\sigma,\tau}.
\end{aligned} \tag{2.76}$$

En additionnant (2.75) et (2.76), nous trouvons (2.72). □

*Suite de la démonstration du théorème 6.* D'après l'identité (2.68), nous avons

$$\frac{1}{g!} \sum_{\sigma \in \mathcal{S}_g} \sum_{\tau \in \mathcal{S}_g} \Phi_{\sigma,\tau}(x) = \frac{1}{(g-1)!} \int_{X^g} \sum_{j=1}^{g-1} G_X(z, z'_j(e)) \cdot \overline{\partial}_z G_X(z,x) \wedge \sum_{\sigma \in \mathcal{S}_g} \sum_{\substack{\tau \in \mathcal{S}_g \\ \tau(k)=g}} \alpha_{\sigma,\tau}.$$

Puisque dans cette dernière intégrale $(e,z) = \rho(z_1, \ldots, z_g)$, nous pouvons utiliser le changement de variables $\rho$ et, grâce au lemme 54, écrire

$$\frac{1}{g!} \sum_{\sigma \in \mathcal{S}_g} \sum_{\tau \in \mathcal{S}_g} \Phi_{\sigma,\tau}(x) = i^{g-1} \int_{(e,z) \in \Theta_0 \times X} \sum_{j=1}^{g-1} G_X(z, z'_j(e)) \cdot \overline{\partial}_z G_X(z,x) \wedge \nu_k.$$

Utilisons maintenant le changement de variables $\rho'$. Nous avons

$$(\rho')^* \nu_k = \rho^* \left( -\operatorname{id} \times \operatorname{id} \right)^* \nu_k. \tag{2.77}$$

Mais $\nu_k$ est combinaison linéaire de termes de la forme $\operatorname{pr}_1^* \Omega \wedge \operatorname{pr}_2^* \omega$, avec $\Omega \in \Omega^{g-1,g-1}(J)$ et $\omega \in \Omega^{1,0}(X)$, et nous avons $(-\operatorname{id})^* \Omega = \Omega$, puisque le degré de la forme différentielle $\Omega$ vaut $2g-2$. Par conséquent $(\rho')^* \nu_k = \rho^* \nu_k$. Enfin, en utilisant des notations évidentes, nous avons

$$(z'_1(-e), \ldots, z'_{g-1}(-e)) = (z_1, \ldots, z_{g-1}),$$

ce qui permet d'écrire

$$\frac{1}{g!} \sum_{\sigma \in \mathcal{S}_g} \sum_{\tau \in \mathcal{S}_g} \Phi_{\sigma,\tau}(x) = \frac{i^{g-1}}{(g-1)!} \int_{X^g} \sum_{l=1}^{g-1} G_X(z_g, z_l) \cdot \overline{\partial}_{z_g} G_X(z_g, x) \wedge \sum_{\sigma \in \mathcal{S}_g} \sum_{\substack{\tau \in \mathcal{S}_g \\ \tau(k)=g}} \alpha_{\sigma,\tau}. \tag{2.78}$$

L'intégrale (2.78) peut maintenant se simplifier grâce au théorème de Fubini, pourvu que nous fassions attention à l'ordre des termes. Posons

$$\Psi^{(l)}_{\sigma,\tau}(x) = \int_{X^g} G_X(z_g, z_l) \cdot \overline{\partial}_{z_g} G_X(z_g, x) \wedge \alpha_{\sigma,\tau} \qquad (1 \leq l \leq g-1). \tag{2.79}$$

D'après le théorème de Fubini et la propriété d'orthogonalité des $\omega_i$, nous voyons que $\Psi^{(l)}_{\sigma,\tau} = 0$ à moins que $\sigma^{-1}$ et $\tau^{-1}$ coïncident sur l'ensemble $\{1, \ldots, \widehat{l}, \ldots, g-1\}$. Pour $\tau$ fixé, il y a exactement deux éléments $\sigma$ qui vérifient cette condition, à savoir $\sigma = \tau$ et $\sigma = t \circ \tau$, où $t$ désigne la transposition échangeant $l$ et $g$. Posons $j = \tau^{-1}(l)$. Lorsque $\sigma = \tau$, nous avons

$$\Psi^{(l)}_{\tau,\tau}(x) = \int_{X^g} G_X(z_g, z_l) \overline{\partial}_{z_g} G_X(z_g, x) \wedge (\omega_k)_{z_g} \wedge \bigwedge_{\substack{m=1 \\ m \neq k}}^{g} p^*_{\tau(l)}(\omega_l \wedge \overline{\omega_l})$$

$$= \frac{1}{i^{g-2}} \int_{z \in X} R_{\omega_k}(x,z) \cdot \omega_j \wedge \overline{\omega_j}. \tag{2.80}$$

Lorsque $\sigma = t \circ \tau$, nous avons

$$\Psi^{(l)}_{t \circ \tau, \tau}(x) = \int_{X^g} G_X(z_g, z_l) \overline{\partial}_{z_g} G_X(z_g, x) \wedge (\omega_k)_{z_l} \wedge (\omega_j)_{z_g} \wedge (\overline{\omega_j})_{z_l} \wedge \bigwedge_{\substack{m=1 \\ m \neq j,k}}^{g} p^*_{\tau(m)}(\omega_m \wedge \overline{\omega_m}),$$

d'où nous déduisons

$$\Psi^{(l)}_{t \circ \tau, \tau}(x) = \frac{1}{i^{g-2}} \int_{z \in X} R_{\omega_j}(z,x) \cdot \omega_k \wedge \overline{\omega_j}. \tag{2.81}$$

En sommant (2.80) sur $l$ et $\tau$, nous obtenons (rappelons que $j = \tau^{-1}(l)$)

$$\sum_{l=1}^{g-1} \sum_{\substack{\tau \in \mathcal{S}_g \\ \tau(k)=g}} \Psi^{(l)}_{\tau,\tau}(x) = \frac{1}{i^{g-2}} \sum_{\substack{\tau \in \mathcal{S}_g \\ \tau(k)=g}} \sum_{l=1}^{g-1} \int_{z \in X} R_{\omega_k}(x,z) \cdot \omega_{\tau^{-1}(l)} \wedge \overline{\omega_{\tau^{-1}(l)}}$$

$$= \frac{(g-1)!}{i^{g-2}} \sum_{\substack{j=1 \\ j \neq k}}^{g} \int_{z \in X} R_{\omega_k}(x,z) \cdot \omega_j \wedge \overline{\omega_j}.$$

Par définition de $\mathrm{vol}_X$, nous avons $\sum_{j \neq k} \omega_j \wedge \overline{\omega_j} = -ig\, \mathrm{vol}_X - \omega_k \wedge \overline{\omega_k}$. Au cours de la démonstration de la proposition 36, nous avons montré que

$$\int_{z \in X} R_\omega(x,z) \cdot \mathrm{vol}_X = 0 \qquad (\omega \in \Omega^{1,0}(X),\ x \in X). \tag{2.82}$$

Il en résulte finalement

$$\sum_{l=1}^{g-1} \sum_{\substack{\tau \in \mathcal{S}_g \\ \tau(k)=g}} \Psi^{(l)}_{\tau,\tau}(x) = -\frac{(g-1)!}{i^{g-2}} \int_{z \in X} R_{\omega_k}(x,z) \cdot \omega_k \wedge \overline{\omega_k}. \tag{2.83}$$

En sommant (2.81) sur $l$ et $\tau$, nous obtenons d'autre part

$$\sum_{l=1}^{g-1} \sum_{\substack{\tau \in \mathcal{S}_g \\ \tau(k)=g}} \Psi_{t \circ \tau,\tau}^{(l)}(x) = \frac{1}{i^{g-2}} \sum_{\substack{\tau \in \mathcal{S}_g \\ \tau(k)=g}} \sum_{l=1}^{g-1} \int_{z \in X} R_{\omega_{\tau^{-1}(l)}}(z,x) \cdot \omega_k \wedge \overline{\omega_{\tau^{-1}(l)}}$$

$$= -\frac{(g-1)!}{i^{g-2}} \sum_{\substack{j=1 \\ j \neq k}}^{g} \int_{z \in X} R_{\omega_j}(x,z) \cdot \omega_k \wedge \overline{\omega_j}. \tag{2.84}$$

En reportant (2.83) et (2.84) dans (2.78), il vient

$$\frac{1}{g!} \sum_{\sigma \in \mathcal{S}_g} \sum_{\tau \in \mathcal{S}_g} \Phi_{\sigma,\tau}(x) = \frac{i^{g-1}}{(g-1)!} \sum_{l=1}^{g-1} \sum_{\substack{\tau \in \mathcal{S}_g \\ \tau(k)=g}} \Psi_{\tau,\tau}^{(l)}(x) + \Psi_{t \circ \tau,\tau}^{(l)}(x)$$

$$= -i \sum_{j=1}^{g} \int_{z \in X} R_{\omega_j}(x,z) \cdot \omega_k \wedge \overline{\omega_j}$$

$$= \Phi_{\omega_k}(x).$$

En revenant à (2.67), nous avons donc finalement

$$\langle R_J(x-y), \alpha_k \rangle = R_{\omega_k}(x,y) + \Phi_{\omega_k}(x) - \Phi_{\omega_k}(y),$$

ce qui montre (2.52). □

*Remarque* 55. Dans le cas où $X = E$ est une courbe elliptique, l'identité (2.82) montre que $\Phi_E = 0$. Par suite, la fonction $R_J$ coïncide avec le dilogarithme elliptique :

$$R_J(P-Q) = R_E(P,Q) \qquad (P,Q \in E), \tag{2.85}$$

en utilisant l'identification $E \cong J$. L'identité (2.85) résulte également assez directement de la définition des fonctions $R_J$ et $R_E$.

*Remarque* 56. Les fonctions $R_J$ et $\Phi_X$ ci-dessus ne dépendent d'aucun choix. Néanmoins, il n'est pas clair que l'identité (2.3) détermine les fonctions $R_J$ et $\Phi_X$ de manière unique : par exemple (2.3) reste vérifiée lorsque l'on ajoute une constante à $\Phi_X$. Supposons d'abord que $X = E$ est une courbe elliptique. Donnons-nous des fonctions $R_J$ et $\Phi_E$ vérifiant (2.3). Cette dernière identité implique formellement

$$\Phi_E(x-y) - \Phi_E(0) = R_J(x-y-0) - R_E(x-y,0)$$
$$= R_J(x-y) - R_E(x,y)$$
$$= R_E(x,y) + \Phi_E(x) - \Phi_E(y) - R_E(x,y)$$
$$= \big(\Phi_E(x) - \Phi_E(0)\big) - \big(\Phi_E(y) - \Phi_E(0)\big) \qquad (x,y \in E).$$

Donc $\Phi_E - \Phi_E(0)$ est un homomorphisme de groupes. Puisque $\Phi_E$ est continue, on a $\Phi_E - \Phi_E(0) = 0$, donc $\Phi_E$ est constante. Par suite, la fonction $R_J$ est donnée par (2.85). Nous voyons donc que (2.3) détermine $R_J$ de manière unique, et $\Phi_E$ à une constante près. Supposons maintenant que $X$ est de genre 2. L'application (2.2) étant surjective, l'identité (2.3) entraîne formellement l'imparité de $R_J$. Mais elle ne semble pas déterminer $R_J$ de manière unique. Dans ce cas, ainsi que dans le cas général, il serait intéressant de chercher à caractériser les fonctions $R_J$ et $\Phi_X$ par des propriétés supplémentaires.

*Remarque* 57. Il serait intéressant de définir l'analogue de la fonction $R_J$ pour les fonctions $G_{a,b}$ définies par Goncharov [31, Def 9.1, p. 390]. La fonction $R_J$ correspond au cas de la fonction $G_{1,2}$.

## 2.4 Propriétés de la fonction $R_J$

Nous nous plaçons dans le même cadre que les deux sections précédentes.

**Proposition 58.** *La fonction $R_J$ est impaire, c'est-à-dire*

$$\langle R_J(-x), \alpha \rangle = -\langle R_J(x), \alpha \rangle \qquad (x \in J, \ \alpha \in \Omega^{g,g-1}(J)). \tag{2.86}$$

*Démonstration.* Cela résulte du calcul suivant :

$$\begin{aligned}
\langle R_J(-x), \alpha \rangle &= \int_{u \in J} \log\|\theta(u)\| \cdot \overline{\partial}_u \log\|\theta(u+x)\| \wedge \alpha_u \\
&= \int_{u \in J} \log\|\theta(-u)\| \cdot \overline{\partial}_u \log\|\theta(-u+x)\| \wedge (-1)^* \alpha_u \\
&= -\int_{u \in J} \log\|\theta(u)\| \cdot \overline{\partial}_u \log\|\theta(u-x)\| \wedge \alpha_u \\
&= -\langle R_J(x), \alpha \rangle,
\end{aligned}$$

par parité de la fonction $\theta$ et d'après $(-1)^* \alpha = -\alpha$. Noter qu'il est également possible de démontrer (2.86) en utilisant la formule de Stokes et en faisant un calcul analogue à (1.21). □

Soit maintenant $\sigma : X \to X$ un automorphisme (holomorphe) de $X$. Il induit, par fonctorialité d'Albanese, un automorphisme $\sigma_J$ de $J$, défini par

$$\sigma_J : J \to J \tag{2.87}$$
$$\left[\sum_i n_i \cdot x_i\right] \mapsto \left[\sum_i n_i \cdot \sigma(x_i)\right].$$

**Proposition 59.** *Nous avons la formule*

$$\langle R_J(\sigma_J(x)), \alpha \rangle = \langle R_J(x), \sigma_J^* \alpha \rangle \qquad (x \in J, \ \alpha \in \Omega^{g,g-1}(J)). \tag{2.88}$$

*Démonstration.* Nous reprenons les notations de la section 2.2. Soit $\Theta_0 \subset J$ l'hypersurface induite par le lieu des zéros de la fonction $\theta$. Fixant un point-base $x_0 \in X$ et utilisant (2.40), nous voyons que

$$\sigma_J(\Theta_0) + \sigma_J(\Delta) = \left\{ \sigma(z_1) + \cdots + \sigma(z_{g-1}) - (g-1)\sigma(x_0) \mid z_1, \ldots, z_{g-1} \in X \right\}.$$

Puisque $\sigma$ est une bijection de $X$, nous voyons que $\sigma_J(\Theta_0)$ coïncide avec un translaté de $\Theta_0$, et il en va de même de $\sigma_J^{-1}(\Theta_0)$. Il existe donc $a \in J$ tel que

$$H := \sigma_J^{-1}(\Theta_0) = \Theta_{-a} = \Theta_0 - a. \tag{2.89}$$

Notons $t_a : x \mapsto x + a$ la translation par $a$ dans $J$ et considérons la fonction

$$\phi = \log\|\theta \circ \sigma_J\| - \log\|\theta \circ t_a\| \tag{2.90}$$

définie et de classe $\mathcal{C}^\infty$ sur $J - H$. En considérant une coordonnée locale holomorphe $\lambda$ pour l'hypersurface $H$, nous voyons que $\log\|\theta \circ \sigma_J\|$ et $\log\|\theta \circ t_a\|$ s'écrivent toutes deux localement comme la somme de $\log|\lambda|$ et d'une fonction $\mathcal{C}^\infty$. En conséquence $\phi$ se prolonge en une fonction de classe $\mathcal{C}^\infty$ sur $J$. En appliquant le lemme 51, nous obtenons sur $J - H$

$$\partial\overline{\partial}\phi = -\pi \sum_{j=1}^g \sigma_J^*(\omega_j \wedge \overline{\omega_j}) - t_a^*(\omega_j \wedge \overline{\omega_j}),$$

où nous rappelons que $(\omega_j)_{1 \leq j \leq g}$ est une base de l'espace $\Omega^{1,0}(X)$ identifié à $\Omega^{1,0}(J)$. Par invariance par translation, nous avons $t_a^*(\omega_j \wedge \overline{\omega_j}) = \omega_j \wedge \overline{\omega_j}$. D'autre part, le calcul

$$\sigma(x) - x_0 = \big(\sigma(x) - \sigma(x_0)\big) + \big(\sigma(x_0) - x_0\big)$$
$$= \sigma_J(x - x_0) + \big(\sigma(x_0) - x_0\big) \qquad (x \in X)$$

montre que $i \circ \sigma = t \circ \sigma_J \circ i$, où $t$ est une translation de $J$. Par conséquent $\sigma^* i^* = i^* \sigma_J^* t^* = i^* \sigma_J^*$ sur l'espace des formes différentielles holomorphes, ce qui permet d'écrire $\sigma_J^* \omega_j = \sigma^* \omega_j$. Or, la famille $(\sigma^* \omega_j)_{1 \leq j \leq g}$ est encore une base orthonormale de $\Omega^{1,0}(X)$. Par indépendance de la forme volume (1.4) vis-à-vis du choix de la base orthonormale, il vient

$$\sum_{j=1}^g \sigma^*(\omega_j \wedge \overline{\omega_j}) = \sum_{j=1}^g \omega_j \wedge \overline{\omega_j},$$

et par suite $\partial\overline{\partial}\phi = 0$ sur $J$. Donc $\phi = C$ est constante (réelle) sur $J$.

Soit maintenant $x \in J$ et $\alpha \in \Omega^{g,g-1}(J)$. Il vient

$$\langle R_J(\sigma_J(x)), \alpha \rangle = \int_{u \in J} \log\|\theta(u)\| \cdot \overline{\partial}_u \log\|\theta(u - \sigma_J(x))\| \wedge \alpha_u$$
$$= \int_{u \in J} \log\|\theta \circ \sigma_J(u)\| \cdot \overline{\partial}_u \log\|\theta \circ \sigma_J(u - x)\| \wedge (\sigma_J^* \alpha)_u$$
$$= \int_{u \in J} \big(C + \log\|\theta \circ t_a(u)\|\big) \cdot \overline{\partial}_u \log\|\theta \circ t_a(u - x)\| \wedge (\sigma_J^* \alpha)_u.$$

La constante $C$ disparaît grâce à la formule de Stokes. Le changement de variables $t_a$ dans l'intégrale permet alors d'obtenir la formule désirée. □

*Remarque* 60. Une question naturelle est de savoir si (2.88) reste vrai pour un automorphisme $\sigma_J$ de $J$ ne provenant pas nécessairement de $X$. En fait, la démonstration de (2.88) n'utilise que le fait qu'il existe $a \in J$ tel que la fonction $\phi$ définie par (2.90) soit constante. On peut montrer que cette hypothèse équivaut au fait que $\sigma_J(\Theta_0)$ soit un translaté de $\Theta_0$. Dans ce cas, $\sigma_J$ laisse invariante la forme $\sum_{j=1}^g \omega_j \wedge \overline{\omega_j}$. Cette hypothèse me semble assez restrictive.

## 2.5 Rappels sur la cohomologie de Deligne

Soit $A$ une variété abélienne définie sur $\mathbf{Q}$, de dimension $d \geq 1$. Considérons le *régulateur de Beĭlinson*

$$r_A : K_2^{(2)}(A) \to H_{\mathcal{D}}^2(A_\mathbf{R}, \mathbf{R}(2)), \tag{2.91}$$

tel qu'il est défini par exemple dans [65, p. 30] et [52, (5.4) et (7.1)]. Notre but est de calculer explicitement (en un sens que nous allons bientôt préciser) l'application régulateur (2.91). Nous noterons $A_{\mathbf{R}}$ (resp. $A_{\mathbf{C}}$) l'extension des scalaires de $A$ à $\mathbf{R}$ (resp. $\mathbf{C}$). Nous noterons $\mathbf{Q}(A)$ et $\mathbf{C}(A)$ les corps des fonctions rationnelles respectifs de $A$ et $A_{\mathbf{C}}$. Nous avons une inclusion naturelle $\mathbf{Q}(A) \subset \mathbf{C}(A)$. Le groupe de $K$-théorie de Quillen $K_2^{(2)}(A)$ est décrit par la suite exacte [52, (7.4), p. 23]

$$0 \to K_2^{(2)}(A) \xrightarrow{\eta^*} K_2(\mathbf{Q}(A)) \otimes \mathbf{Q} \xrightarrow{\partial} \bigoplus_{V \subset A} \mathbf{Q}(V)^* \otimes \mathbf{Q}, \qquad (2.92)$$

où la somme directe porte sur les hypersurfaces irréductibles $V$ de $A$, $\mathbf{Q}(V)$ désigne le corps des fonctions de $V$, et $\partial = (\partial_V)_{V \subset A}$ est défini par

$$\partial_V : K_2(\mathbf{Q}(A)) \to \mathbf{Q}(V)^* \qquad (2.93)$$
$$\{F, G\} \mapsto (-1)^{\operatorname{ord}_V(F) \operatorname{ord}_V(G)} \left( \frac{F^{\operatorname{ord}_V(G)}}{G^{\operatorname{ord}_V(F)}} \right) \Big|_V.$$

Nous allons donner une formule explicite pour $r_A$, dans le sens suivant : nous allons définir une application explicite $\widehat{r}_A$ sur $K_2(\mathbf{Q}(A)) \otimes \mathbf{Q}$, et montrer que la composition de cette application avec l'application $\eta^*$ de (2.92) coïncide avec le régulateur de Beĭlinson $r_A$.

Commençons par décrire explicitement l'espace d'arrivée $H^2_{\mathcal{D}}(A_{\mathbf{R}}, \mathbf{R}(2))$ de l'application régulateur, en terme de la cohomologie singulière de $A(\mathbf{C})$. Le lemme suivant se déduit de [65, (*), p. 9].

**Lemme 61.** *Les groupes de cohomologie de Deligne $H^2_{\mathcal{D}}(A_{\mathbf{C}}, \mathbf{R}(2))$ et $H^2_{\mathcal{D}}(A_{\mathbf{R}}, \mathbf{R}(2))$ sont décrits par le diagramme commutatif*

$$\begin{array}{ccc} H^2_{\mathcal{D}}(A_{\mathbf{C}}, \mathbf{R}(2)) & \xrightarrow{\cong} & H^1(A(\mathbf{C}), \mathbf{R}(1)) \\ \cup | & & \cup | \\ H^2_{\mathcal{D}}(A_{\mathbf{R}}, \mathbf{R}(2)) & \xrightarrow{\cong} & H^1(A(\mathbf{C}), \mathbf{R}(1))^- \end{array} \qquad (2.94)$$

*où les lignes sont des isomorphismes. Dans ce diagramme, nous avons posé $\mathbf{R}(1) = 2\pi i \mathbf{R} \subset \mathbf{C}$, et $(\cdot)^-$ désigne le sous-espace des éléments anti-invariants par le morphisme $c^*$ induit par la conjugaison complexe $c : A(\mathbf{C}) \to A(\mathbf{C})$.*

*Démonstration.* La suite exacte en bas de [65, p. 8] s'écrit

$$0 \to F^2 H^1_{\mathrm{DR}}(A_{\mathbf{C}}) \to H^1(A(\mathbf{C}), \mathbf{R}(1)) \to H^2_{\mathcal{D}}(A_{\mathbf{C}}, \mathbf{R}(2)) \to 0, \qquad (2.95)$$

où $H^1_{\mathrm{DR}}(A_{\mathbf{C}})$ désigne la cohomologie de de Rham algébrique de $A_{\mathbf{C}}$. Pour toute variété projective lisse $X$ définie sur $\mathbf{C}$, la cohomologie de de Rham algébrique de $X$ et sa filtration sont données par

$$H^i_{\mathrm{DR}}(X) \cong H^i(X(\mathbf{C}), \mathbf{C}) \qquad (0 \leq i \leq 2 \dim_{\mathbf{C}} X), \qquad (2.96)$$
$$F^{p_0} H^i_{\mathrm{DR}}(X) \cong \bigoplus_{\substack{p+q=i \\ p \geq p_0}} H^{p,q}(X(\mathbf{C}), \mathbf{C}) \qquad (p_0 \in \mathbf{Z}). \qquad (2.97)$$

Il est donc clair que $F^2 H^1_{\mathrm{DR}}(A_{\mathbf{C}}) = 0$, ce qui montre le premier isomorphisme du diagramme (2.94). D'après [65, p. 8], le groupe de cohomologie de Deligne $H^2_{\mathcal{D}}(A_{\mathbf{R}}, \mathbf{R}(2))$ se définit comme le sous-espace de $H^2_{\mathcal{D}}(A_{\mathbf{C}}, \mathbf{R}(2))$ fixé par la "conjugaison de de Rham". Cette conjugaison correspond à l'involution $[\omega] \mapsto [c^*\overline{\omega}] = -[c^*\omega]$ sur $H^1(A(\mathbf{C}), \mathbf{R}(1))$. L'espace des invariants sous cette involution est précisément $H^1(A(\mathbf{C}), \mathbf{R}(1))^-$, d'où le second isomorphisme du diagramme (2.94). □

Introduisons maintenant une notation. L'espace des formes différentielles complexes invariantes par translation de type $(d, d-1)$ sur $A(\mathbf{C})$ sera noté $\Omega^{d,d-1}(A(\mathbf{C}))$. C'est un espace vectoriel complexe de dimension $d$ muni d'une involution anti-linéaire

$$\sigma : \Omega^{d,d-1}(A(\mathbf{C})) \to \Omega^{d,d-1}(A(\mathbf{C})) \tag{2.98}$$
$$\alpha \mapsto c^*\overline{\alpha},$$

où $c : A(\mathbf{C}) \to A(\mathbf{C})$ désigne la conjugaison complexe. Nous noterons $\Omega^{d,d-1}_{\mathbf{R}}(A(\mathbf{C}))$ le sous-espace de $\Omega^{d,d-1}(A(\mathbf{C}))$ fixé par $\sigma$ (notation non standard). C'est un espace vectoriel réel de dimension $d$, et l'inclusion induit un isomorphisme canonique

$$\Omega^{d,d-1}_{\mathbf{R}}(A(\mathbf{C})) \otimes_{\mathbf{R}} \mathbf{C} \cong \Omega^{d,d-1}(A(\mathbf{C})). \tag{2.99}$$

**Lemme 62.** *Nous avons un isomorphisme*

$$H^1(A(\mathbf{C}), \mathbf{R}(1)) \to \mathrm{Hom}_{\mathbf{C}}(\Omega^{d,d-1}(A(\mathbf{C})), \mathbf{C}) \tag{2.100}$$
$$[\omega] \mapsto \left(\alpha \mapsto \int_{A(\mathbf{C})} \omega \wedge \alpha\right).$$

*Via cet isomorphisme, le sous-espace $H^1(A(\mathbf{C}), \mathbf{R}(1))^-$ s'identifie à*

$$\mathrm{Hom}_{\mathbf{R}}(\Omega^{d,d-1}_{\mathbf{R}}(A(\mathbf{C})), \mathbf{R}(d)) \subset \mathrm{Hom}_{\mathbf{C}}(\Omega^{d,d-1}(A(\mathbf{C})), \mathbf{C}), \tag{2.101}$$

*où la dernière inclusion est induite par l'isomorphisme (2.99).*

*Démonstration.* Nous vérifions sans peine que l'application (2.100) est bien définie ; montrons que c'est un isomorphisme. Les espaces de départ et d'arrivée ayant même dimension réelle $2d$, il est suffisant de montrer l'injectivité. Soit $[\omega] \in H^1(A(\mathbf{C}), \mathbf{R}(1))$ un élément vérifiant

$$\int_{A(\mathbf{C})} \omega \wedge \alpha = 0 \qquad (\alpha \in \Omega^{d,d-1}(A(\mathbf{C}))). \tag{2.102}$$

La décomposition de Hodge

$$H^{2d-1}(A(\mathbf{C}), \mathbf{C}) = \Omega^{d,d-1}(A(\mathbf{C})) \oplus \overline{\Omega^{d,d-1}(A(\mathbf{C}))} \tag{2.103}$$

et le fait que $\overline{\omega} = -\omega$ montrent que (2.102) reste valable pour tout $\alpha \in H^{2d-1}(A(\mathbf{C}), \mathbf{C})$. Grâce au théorème de dualité de Poincaré à coefficients complexes, appliqué à la variété compacte $A(\mathbf{C})$, nous en déduisons $\omega = 0$. L'application (2.100) est donc un isomorphisme.

Soit à présent $\omega$ une forme différentielle fermée sur $A(\mathbf{C})$ avec $[\omega] \in H^1(A(\mathbf{C}), \mathbf{R}(1))^-$. Soit également $\alpha \in \Omega^{d,d-1}_{\mathbf{R}}(A(\mathbf{C}))$. Nous avons

$$\overline{\int_{A(\mathbf{C})} \omega \wedge \alpha} = \int_{A(\mathbf{C})} \overline{\omega} \wedge \overline{\alpha} = -\int_{A(\mathbf{C})} \omega \wedge c^*\alpha$$

$$= -(-1)^d \int_{A(\mathbf{C})} c^*\omega \wedge \alpha = (-1)^d \int_{A(\mathbf{C})} \omega \wedge \alpha, \quad (2.104)$$

ce qui montre que (2.100) induit une application injective

$$H^1(A(\mathbf{C}), \mathbf{R}(1))^- \to \mathrm{Hom}_{\mathbf{R}}(\Omega^{d,d-1}_{\mathbf{R}}(A(\mathbf{C})), \mathbf{R}(d)). \quad (2.105)$$

Puisque les deux espaces ci-dessus ont même dimension réelle $d$, cette application est un isomorphisme. □

Nous déduisons des lemmes 61 et 62 les identifications

$$\begin{array}{ccc}
H^2_{\mathcal{D}}(A_{\mathbf{C}}, \mathbf{R}(2)) & \xrightarrow{\cong} & \mathrm{Hom}_{\mathbf{C}}(\Omega^{d,d-1}(A(\mathbf{C})), \mathbf{C}) \\
\bigcup & & \bigcup \\
H^2_{\mathcal{D}}(A_{\mathbf{R}}, \mathbf{R}(2)) & \xrightarrow{\cong} & \mathrm{Hom}_{\mathbf{R}}(\Omega^{d,d-1}_{\mathbf{R}}(A(\mathbf{C})), \mathbf{R}(d))
\end{array} \quad (2.106)$$

Nous allons les utiliser pour prolonger le régulateur de Beĭlinson $r_A$.

**Définition 63.** *Pour toutes fonctions rationnelles $F, G \in \mathbf{C}(A)^*$, notons $\eta(F,G)$ la forme différentielle définie par*

$$\eta(F,G) = \log|F| \cdot (\partial - \overline{\partial}) \log|G| - \log|G| \cdot (\partial - \overline{\partial}) \log|F|. \quad (2.107)$$

*Remarque* 64. La définition de $\eta(F,G)$ est directement issue de la formule pour le cup-produit en cohomologie de Deligne [52, (7.3.2), p. 23]. La forme différentielle $\eta(F,G)$ est un représentant privilégié du cup-produit des classes des fonctions $F$ et $G$.

La forme différentielle $\eta(F,G)$ est définie sur l'ouvert $U = A(\mathbf{C}) - S$, où $S$ est la réunion des hypersurfaces intervenant dans les diviseurs de $F$ et $G$. Remarquons que

$$\overline{\eta(F,G)} = -\eta(F,G). \quad (2.108)$$

D'autre part, si $c : U \to U$ désigne la conjugaison complexe, nous avons $c^*\overline{F}, c^*\overline{G} \in \mathbf{C}(A)^*$ et

$$c^*\eta(F,G) = -\eta(c^*\overline{F}, c^*\overline{G}). \quad (2.109)$$

Enfin, nous avons

$$d\eta(F,G) = \frac{1}{2}\big(\mathrm{dlog}\, F \wedge \mathrm{dlog}\, G - \overline{\mathrm{dlog}\, F} \wedge \overline{\mathrm{dlog}\, G}\big) = \pi_1(\mathrm{dlog}\, F \wedge \mathrm{dlog}\, G), \quad (2.110)$$

en posant $\pi_1(\omega) = \frac{1}{2}(\omega - \overline{\omega})$ pour toute forme forme différentielle $\omega$. Lorsque $d = 1$, la forme différentielle $\eta(F,G)$ est fermée ; lorsque $d > 1$, ce n'est pas le cas en général.

**Définition 65.** *Soit $\widehat{r}_{A_{\mathbf{C}}}$ l'application définie par*

$$\widehat{r}_{A_{\mathbf{C}}} : K_2(\mathbf{C}(A)) \to H^2_{\mathcal{D}}(A_{\mathbf{C}}, \mathbf{R}(2))$$

$$\{F, G\} \mapsto \left(\alpha \mapsto \int_{A(\mathbf{C})} \eta(F, G) \wedge \alpha\right), \quad (2.111)$$

*où nous utilisons l'identification $H^2_{\mathcal{D}}(A_{\mathbf{C}}, \mathbf{R}(2)) \cong \mathrm{Hom}_{\mathbf{C}}(\Omega^{d,d-1}(A(\mathbf{C})), \mathbf{C})$ de (2.106).*

Montrons que l'application $\widehat{r}_{A_\mathbf{C}}$ est bien définie. Puisque les fonctions $\log|F|$ et $\log|G|$ sont à croissance logarithmique le long des supports des diviseurs de $F$ et $G$, l'intégrale (2.111) définissant $\widehat{r}_A$ converge absolument. La forme différentielle $\eta(F,G)$, et donc l'intégrale (2.111), sont clairement bilinéaires en $F, G \in \mathbf{C}(A)^*$. Pour justifier la définition 65, il reste à montrer que l'intégrale (2.111) annule les relations de Steinberg. Soit $F \in \mathbf{C}(A)$ une fonction rationnelle, avec $F \neq 0, 1$. Un petit calcul montre que

$$\eta(F, 1-F) = i\, d(D \circ F), \tag{2.112}$$

où $D$ est la fonction de Bloch-Wigner. Il en découle

$$\int_{A(\mathbf{C})} \eta(F, 1-F) \wedge \alpha = i \int_{A(\mathbf{C})} d\big((D \circ F)\, \alpha\big) = 0 \qquad (\alpha \in \Omega^{d,d-1}(A(\mathbf{C}))), \tag{2.113}$$

d'après la formule de Stokes, ce qui achève de montrer que $\widehat{r}_{A_\mathbf{C}}$ est bien définie.

Notons la propriété suivante d'invariance par translation : pour toutes fonctions rationnelles $F, G \in \mathbf{C}(A)^*$ et tout $a \in A(\mathbf{C})$, nous avons

$$\widehat{r}_{A_\mathbf{C}}(\{t_a^*F, t_a^*G\}) = \widehat{r}_{A_\mathbf{C}}(\{F, G\}). \tag{2.114}$$

**Proposition 66.** *L'application $\widehat{r}_{A_\mathbf{C}}$ définie en (2.111) induit une application*

$$\widehat{r}_A : K_2(\mathbf{Q}(A)) \to H^2_\mathcal{D}(A_\mathbf{R}, \mathbf{R}(2)). \tag{2.115}$$

*Démonstration.* Rappelons que l'inclusion naturelle $\mathbf{Q}(A) \subset \mathbf{C}(A)$ induit un homomorphisme $K_2(\mathbf{Q}(A)) \to K_2(\mathbf{C}(A))$. Considérons des fonctions rationnelles $F, G \in \mathbf{Q}(A)^*$. Puisque $c^*\overline{F} = F$ et $c^*\overline{G} = G$, nous avons d'après (2.108) et (2.109) l'identité $c^*\overline{\eta(F,G)} = \eta(F,G)$. Un calcul similaire à (2.104) montre que

$$\int_{A(\mathbf{C})} \eta(F, G) \wedge \alpha \in \mathbf{R}(d) \qquad (\alpha \in \Omega^{d,d-1}_\mathbf{R}(A(\mathbf{C}))). \tag{2.116}$$

L'élément $\widehat{r}_A(\{F,G\})$ appartient donc bien à $H^2_\mathcal{D}(A_\mathbf{R}, \mathbf{R}(2))$. □

Si nous tensorisons l'application $\widehat{r}_A$ en (2.115) par $\mathbf{Q}$, nous obtenons l'application suivante, que nous notons encore $\widehat{r}_A$ :

$$\widehat{r}_A : K_2(\mathbf{Q}(A)) \otimes \mathbf{Q} \to H^2_\mathcal{D}(A_\mathbf{R}, \mathbf{R}(2)). \tag{2.117}$$

**Proposition 67.** *Le régulateur de Beĭlinson $r_A$ de (2.91) s'obtient comme la composée*

$$K_2^{(2)}(A) \xrightarrow{\eta^*} K_2(\mathbf{Q}(A)) \otimes \mathbf{Q} \xrightarrow{\widehat{r}_A} H^2_\mathcal{D}(A_\mathbf{R}, \mathbf{R}(2)). \tag{2.118}$$

*Démonstration.* Soit $\gamma \in K_2^{(2)}(A)$. En utilisant l'identification (2.94), écrivons $r_A(\gamma) = [\omega]$, où $\omega$ est une 1-forme différentielle fermée sur $A(\mathbf{C})$ (nous pouvons supposer que $\omega$ est invariante par translation). Nous avons donc

$$r_A(\gamma) = \left(\alpha \mapsto \int_{A(\mathbf{C})} \omega \wedge \alpha\right) \tag{2.119}$$

en utilisant l'identification (2.106).

D'autre part, écrivons $\eta^*\gamma = \frac{1}{N}\sum_i \{F_i, G_i\}$ avec $F_i, G_i \in \mathbf{Q}(A)^*$ et $N \geq 1$. Soit $U$ l'ouvert complémentaire dans $A$ de la réunion des supports des fonctions $F_i$ et $G_i$. Notons $j: U \hookrightarrow A$ l'inclusion. Nous avons un diagramme commutatif

$$\begin{array}{ccc}
H^2_{\mathcal{D}}(A_{\mathbf{R}}, \mathbf{R}(2)) & \xrightarrow{j^*} & H^2_{\mathcal{D}}(U_{\mathbf{R}}, \mathbf{R}(2)) \\
{\scriptstyle r_A}\uparrow & & {\scriptstyle r_U}\uparrow \\
K_2^{(2)}(A) & \xrightarrow{j^*} K_2^{(2)}(U) \xrightarrow{\eta_U^*} & K_2(\mathbf{Q}(A)) \otimes \mathbf{Q},
\end{array} \qquad (2.120)$$

$$\underbrace{\phantom{K_2^{(2)}(A) \xrightarrow{j^*} K_2^{(2)}(U) \xrightarrow{\eta_U^*} K_2(\mathbf{Q}(A))}}_{\eta^*}$$

où $r_U$ désigne le régulateur de Beĭlinson associé à $U$, tel qu'il est défini dans [52, (7.1), p. 21] par exemple. La commutativité du carré du diagramme (2.120) exprime la compatibilité du régulateur aux images réciproques [30, Section 2]. Grâce au cup-produit en cohomologie motivique, nous pouvons former l'élément $\gamma_U = \frac{1}{N}\sum_i \{F_i, G_i\} \in K_2^{(2)}(U)$. Il vérifie $\eta_U^*\gamma_U = \eta^*\gamma$. Or, l'application $\eta_U^*$ est injective, pour la même raison qu'en (2.5). Nous avons donc $\gamma_U = j^*\gamma$. Par commutativité du diagramme, il vient alors $j^*r_A(\gamma) = r_U(\gamma_U)$.

Nous avons maintenant besoin d'une description du groupe de cohomologie de Deligne $H^2_{\mathcal{D}}(U_{\mathbf{R}}, \mathbf{R}(2))$ similaire à (2.94). Notons $\mathcal{C}^\infty(U(\mathbf{C}), \mathbf{R}(1))^-$ l'espace des fonctions $\phi$ de classe $\mathcal{C}^\infty$ sur $A(\mathbf{C})$, à valeurs dans $\mathbf{R}(1)$, et satisfaisant $\phi \circ c = -\phi$, où $c: U(\mathbf{C}) \to U(\mathbf{C})$ désigne la conjugaison complexe. Notons encore $\mathcal{A}^1(U(\mathbf{C}))^-$ l'espace des 1-formes différentielles réelles $\nu$ sur $U(\mathbf{C})$ vérifiant $c^*\nu = -\nu$. Notons enfin $\Omega^2_{\log}(U(\mathbf{C}))$ l'espace des 2-formes différentielles complexes fermées sur $U(\mathbf{C})$ qui sont à croissance logarithmique le long de $A(\mathbf{C}) - U(\mathbf{C})$ [21]. L'identification suivante se trouve dans [52, (7.3.1)].

$$H^2_{\mathcal{D}}(U_{\mathbf{R}}, \mathbf{R}(2)) \xrightarrow{\cong} \frac{\left\{\nu \in \mathcal{A}^1(U(\mathbf{C}))^- \otimes_{\mathbf{R}} \mathbf{R}(1) \mid d\nu \in \pi_1\bigl(\Omega^2_{\log}(U(\mathbf{C}))\bigr)\right\}}{d\bigl(\mathcal{C}^\infty(U(\mathbf{C}), \mathbf{R}(1))^-\bigr)}, \qquad (2.121)$$

Puisque le régulateur de Beĭlinson est compatible aux cup-produits [52, p. 23], et d'après la formule du cup-produit en cohomologie de Deligne [52, (7.3.2), p. 23], nous avons

$$r_U(\gamma_U) = \frac{1}{N}\sum_i [\eta(F_i, G_i)]. \qquad (2.122)$$

Nous avons montré $j^*[\omega] = r_U(\gamma_U)$. L'application $j^* : H^2_{\mathcal{D}}(A_{\mathbf{R}}, \mathbf{R}(2)) \to H^2_{\mathcal{D}}(U_{\mathbf{R}}, \mathbf{R}(2))$ du haut du diagramme (2.120) est donnée, via les identifications (2.94) et (2.121), par la restriction des formes différentielles. Il existe donc une fonction $\phi \in \mathcal{C}^\infty(U(\mathbf{C}))^- \otimes_{\mathbf{R}} \mathbf{R}(1)$ telle que nous ayons

$$\omega|_{U(\mathbf{C})} = \frac{1}{N}\sum_i \eta(F_i, G_i) + d\phi. \qquad (2.123)$$

En particulier, la forme différentielle $\sum_i \eta(F_i, G_i)$ est fermée sur $U(\mathbf{C})$. La fonction $\phi$ s'écrit partout localement comme une primitive de la forme différentielle fermée $\omega - \frac{1}{N}\sum_i \eta(F_i, G_i)$, qui est à croissance logarithmique le long de $A(\mathbf{C}) - U(\mathbf{C})$. Il en résulte que $\phi$ croît au plus logarithmiquement le long de $A(\mathbf{C}) - U(\mathbf{C})$. Un calcul utilisant la formule de Stokes donne, pour toute forme différentielle $\alpha \in \Omega^{d,d-1}_{\mathbf{R}}(A(\mathbf{C}))$,

$$\int_{A(\mathbf{C})} \omega \wedge \alpha = \frac{1}{N} \sum_i \int_{U(\mathbf{C})} \eta(F_i, G_i) \wedge \alpha + \int_{U(\mathbf{C})} d\phi \wedge \alpha$$
$$= \frac{1}{N} \sum_i \widehat{r}_A(\{F_i, G_i\})(\alpha) + \int_{U(\mathbf{C})} d(\phi\,\alpha)$$
$$= \widehat{r}_A(\eta^*\gamma)(\alpha).$$

En conséquence $r_A(\gamma) = \widehat{r}_A(\eta^*\gamma)$. □

*Question.* Soit $U$ un ouvert Zariski non vide de $A$ et $j_U : U \hookrightarrow A$ l'inclusion. Considérons le morphisme

$$j_U^* : H^2_{\mathcal{D}}(A_{\mathbf{R}}, \mathbf{R}(2)) \to H^2_{\mathcal{D}}(U_{\mathbf{R}}, \mathbf{R}(2)).$$

Est-il possible de construire un projecteur

$$p_U : H^2_{\mathcal{D}}(U_{\mathbf{R}}, \mathbf{R}(2)) \to H^2_{\mathcal{D}}(A_{\mathbf{R}}, \mathbf{R}(2)), \tag{2.124}$$

c'est-à-dire un morphisme "naturel" vérifiant $p_U \circ j_U^* = \mathrm{id}$ ? L'application $\widehat{r}_A$ définie en (2.115) pourrait alors, de manière plus conceptuelle, être définie par la formule

$$\widehat{r}_A = \varinjlim_{U \subset A} p_U \circ r_U, \tag{2.125}$$

où la limite inductive porterait sur les ouverts Zariski non vides $U$ de $A$.

## 2.6 Calcul du régulateur de Beĭlinson

Nous retournons maintenant au cadre présenté au début de la section 2.1 (p. 44), et calculons le régulateur de Beĭlinson de l'élément $\gamma$ construit dans la proposition 44. Grâce à l'identification (2.106), le régulateur de Beĭlinson s'écrit

$$r_J : K_2^{(2)}(J) \to \mathrm{Hom}_{\mathbf{C}}(\Omega^{g,g-1}(J(\mathbf{C})), \mathbf{C}). \tag{2.126}$$

**Théorème 68.** *Soient $l, m \in I_X$ des diviseurs et $M \subset J_{\mathrm{tors}}$ un sous-groupe fini, défini sur $\mathbf{Q}$, contenant les supports de $l$ et $m$. Notons $\gamma \in K_2^{(2)}(J)$ l'élément associé à $l$, $m$ et $M$, obtenu dans la proposition 44. Nous avons*

$$r_J(\gamma) = 2 \sum_{x \in (l)} \sum_{y \in (m)} \mathrm{ord}_x(l)\,\mathrm{ord}_y(m) R_J(x - y), \tag{2.127}$$

*où nous notons $x - y \in J_{\mathrm{tors}}$ la classe du diviseur $[x] - [y]$.*

*Démonstration.* Grâce aux propositions 66 et 67, ainsi qu'à l'invariance par translation (2.114), le régulateur de Beĭlinson de $\gamma$ est donné par

$$r_J(\gamma) = \left(\alpha \in \Omega^{g,g-1}(J(\mathbf{C})) \mapsto \int_{J(\mathbf{C})} \eta(F, G) \wedge \alpha\right). \tag{2.128}$$

Nous allons maintenant exprimer $\eta(F, G)$ en termes de le fonction $\log\|\theta\|$. Reprenons les notations de la section 2.3. Le lieu des zéros de $\|\theta\|$ est noté $\Theta_0$, et il existe $\Delta \in J(\mathbf{C})$ tel que $\Theta_0 + \Delta$ soit le diviseur thêta donné par (2.8). Soit $S$ la réunion des hypersurfaces $\Theta_x$, $x \in (l)$. Considérons la fonction $\phi_F$, de classe $\mathcal{C}^\infty$ sur $J(\mathbf{C}) - S$, définie par

$$\phi_F(u) = \log|F(u)| - \sum_{x \in (l)} \mathrm{ord}_x(l) \cdot \log\|\theta(u - \Delta - x)\| \qquad (u \in J(\mathbf{C}) - S). \tag{2.129}$$

Puisque la fonction $\log\|\theta\|$ est à singularité logarithmique le long de $\Theta_0$, la fonction $\phi_F$ s'étend en une fonction de classe $\mathcal{C}^\infty$ sur $J(\mathbf{C})$. De plus, le lemme 51 entraîne l'identité suivante sur $J(\mathbf{C}) - S$

$$\partial\overline{\partial}\phi_F = \pi \sum_{x \in (l)} \mathrm{ord}_x(l) \sum_{j=1}^g \omega_j \wedge \overline{\omega_j} = 0, \tag{2.130}$$

puisque le diviseur $l$ est de degré 0. Donc $\phi_F$ est constante sur $J(\mathbf{C})$ : il existe une constante $C_F$ telle que

$$\log|F(u)| = C_F + \sum_{x \in (l)} \mathrm{ord}_x(l) \cdot \log\|\theta(u - \Delta - x)\| \qquad (u \notin S). \tag{2.131}$$

De même, il existe une constante $C_G$ telle que

$$\log|G(u)| = C_G + \sum_{y \in (m)} \mathrm{ord}_y(m) \cdot \log\|\theta(u - \Delta - y)\| \qquad (u \notin \Theta_y,\ y \in (m)). \tag{2.132}$$

Soit $\alpha \in \Omega^{g,g-1}(J(\mathbf{C}))$. Nous avons

$$\begin{aligned}\eta(F,G) \wedge \alpha &= \bigl(-\log|F| \cdot \overline{\partial}\log|G| + \log|G| \cdot \overline{\partial}\log|F|\bigr) \wedge \alpha \\ &= -\log|F| \cdot d\bigl(\log|G| \cdot \alpha\bigr) + \log|G| \cdot d\bigl(\log|F| \cdot \alpha\bigr) \\ &= -2\log|F| \cdot d\bigl(\log|G| \cdot \alpha\bigr) + d\bigl(\log|F| \cdot \log|G| \cdot \alpha\bigr),\end{aligned}$$

d'où nous déduisons

$$\int_{J(\mathbf{C})} \eta(F,G) \wedge \alpha = -2 \int_{J(\mathbf{C})} \log|F| \cdot \overline{\partial}\log|G| \wedge \alpha.$$

Remplaçons $\log|F|$ et $\log|G|$ par leurs expressions (2.131) et (2.132). D'après la formule de Stokes $\int_{J(\mathbf{C})} \overline{\partial}\log|G| \wedge \alpha = 0$, et il vient

$$\int_{J(\mathbf{C})} \eta(F,G) \wedge \alpha = -2 \sum_{x \in (l)} \sum_{y \in (m)} \mathrm{ord}_x(l)\,\mathrm{ord}_y(m)\cdot$$
$$\int_{u \in J(\mathbf{C})} \log\|\theta(u - \Delta - x)\| \cdot \overline{\partial}_u \log\|\theta(u - \Delta - y)\| \wedge \alpha_u.$$

En effectuant le changement de variables $u' = u - \Delta - x$ dans l'intégrale, il vient facilement

$$\begin{aligned}\int_{J(\mathbf{C})} \eta(F,G) \wedge \alpha &= -2 \sum_{x \in (l)} \sum_{y \in (m)} \mathrm{ord}_x(l)\,\mathrm{ord}_y(m) R_J(y - x) \\ &= 2 \sum_{x \in (l)} \sum_{y \in (m)} \mathrm{ord}_x(l)\,\mathrm{ord}_y(m) R_J(x - y)\end{aligned}$$

d'après la proposition 58. Le résultat découle alors de (2.128). $\qquad\square$

*Conclusion.* La fonction $R_J$ calcule donc le régulateur de Beĭlinson de certains éléments du groupe de cohomologie motivique $K_2^{(2)}(J) = H_{\mathcal{M}}^2(J, \mathbf{Q}(2))$, ce qui justifie son introduction. Les conjectures de Beĭlinson prédisent que la valeur spéciale $L(J, 2)$ s'exprime à l'aide de ce régulateur. Le théorème 113 donne un exemple explicite de lien entre $L(J, 2)$ et $R_J$, pour $J = J_1(13)$.

*Remarque* 69. Nous pouvons maintenant donner un argument heuristique permettant de mieux comprendre la fonction $R_J$. Cette dernière est définie grâce à une intégrale qui n'est autre qu'une moyenne sur le groupe compact $J(\mathbf{C})$. Cette moyenne est un analogue continu de la moyenne discrète (2.11) qui a permis de construire des éléments dans la cohomologie motivique de $J$. Il est utile de penser à ces opérations de moyenne comme à des projections, respectivement vers $H_{\mathcal{D}}^2(J_{\mathbf{R}}, \mathbf{R}(2))$ et $K_2^{(2)}(J)$.

## 2.7  Questions et perspectives

Plaçons-nous dans le cadre de la section 2.1. La formulation des conjectures de Beĭlinson pour la valeur spéciale $L(J, 2)$ nécessite l'introduction d'un *sous-espace entier* $K_2^{(2)}(J)_{\mathbf{Z}} \subset K_2^{(2)}(J)$, défini à l'aide d'un modèle entier convenable de $J$, cf. [52]. Admettant l'existence d'un tel sous-espace, à quelle(s) condition(s) les éléments construits dans la proposition 44 sont-ils entiers ?

Il serait intéressant de généraliser la construction (2.11) au cas d'une variété abélienne quelconque. Soit $A$ une variété abélienne de dimension $d \geq 1$, définie sur $\mathbf{Q}$, et $N \geq 1$ un entier. Notons $A[N]$ le groupe des points de $N$-torsion de $A(\overline{\mathbf{Q}})$. Considérons le morphisme (abstrait)

$$p_N = \frac{1}{N^{2d}} \sum_{a \in A[N]} t_a^*, \tag{2.133}$$

En faisant agir $p_N$ sur les groupes $K_2^{(2)}(A)$ et $K_2(\mathbf{Q}(A)) \otimes \mathbf{Q}$, nous obtenons un diagramme commutatif

$$\begin{array}{ccccc} 0 & \longrightarrow & K_2^{(2)}(A) & \xrightarrow{\eta^*} & K_2(\mathbf{Q}(A)) \otimes \mathbf{Q} \\ & & {\scriptstyle p_N'}\downarrow & & \downarrow{\scriptstyle p_N''} \\ 0 & \longrightarrow & K_2^{(2)}(A) & \xrightarrow{\eta^*} & K_2(\mathbf{Q}(A)) \otimes \mathbf{Q}. \end{array} \tag{2.134}$$

La conjecture 47 entraîne que $p_N'$ est l'identité. En revanche, l'application $p_N''$ n'est pas injective, comme le montrent des considérations de symbole modéré. Il semble donc intéressant de trouver des combinaisons linéaires de symboles de Milnor

$$\gamma = \sum_i \{F_i, G_i\} \qquad \bigl(F_i, G_i \in \mathbf{Q}(A)^*\bigr), \tag{2.135}$$

telles que $p_N''(\gamma)$ appartienne à l'image de $\eta^*$. Nous avons vu dans la section 2.1 que de tels symboles existent dans le cas d'une jacobienne[3]. En général, nous n'avons pas réussi à construire de tels symboles, qui seraient l'analogue des symboles d'Eisenstein dans le cas elliptique. Enfin, il serait très intéressant d'utiliser le morphisme $p_N$ pour étudier d'autres groupes de cohomologie motivique associés à $A$. Cela nécessiterait d'écrire des suites exactes de localisation en $K$-théorie algébrique et de démontrer des lois de réciprocité pour les homomorphismes de bord associés.

---

[3]Il faut également supposer que la courbe possède suffisamment de points de torsion.

Il est tentant de relier la fonction $R_J$ au polylogarithme abélien associé à $J$. L'existence et la formulation d'un tel lien ne sont cependant pas claires. Beĭlinson et Levin [7] ont défini les faisceaux logarithme et polylogarithme dans le cas d'une courbe elliptique. Wildeshaus [77] a étendu ces définitions au cas d'une variété abélienne ; en particulier, le polylogarithme abélien est défini comme extension de faisceaux. D'autre part, Levin [43] a défini des "courants polylogarithmiques" sur les variétés abéliennes complexes polarisées. Le lien entre les approches de Wildeshaus et Levin est étudié dans la thèse de D. Blottière [14]. Enfin, Kings [38] a construit des éléments dans la cohomologie motivique des schémas abéliens, dont les réalisations étale et de Hodge sont les éléments définis par Wildeshaus. Il est à noter que le polylogarithme abélien vit dans un groupe d'extensions $\mathrm{Ext}^{2g-1}$, avec $g = \dim J$ [38, 3.1]. Le théorème 68 indique quant à lui que la fonction $R_J$ calcule la réalisation de certains éléments de $H^2_{\mathcal{M}}(J, \mathbf{Q}(2))$. Or, le groupe de cohomologie de Deligne $H^2_{\mathcal{D}}(J_{\mathbf{C}}, \mathbf{R}(2))$ s'interprète comme un groupe d'extensions $\mathrm{Ext}^1$ dans la catégorie $\mathcal{MH}_{\mathbf{R}}$ des $\mathbf{R}$-structures de Hodge mixtes [52, (3.4.1) et (7.1)], groupe de nature a priori différente de $\mathrm{Ext}^{2g-1}$ pour $g > 1$ [4]. Pour établir un lien éventuel entre la fonction $R_J$ et le polylogarithme abélien, une possibilité serait de mettre en lumière un théorème de dualité entre les groupes d'extensions mentionnés ci-dessus. Un autre point de départ serait de chercher les propriétés différentielles de la fonction $R_J$.

---

[4] Notons également que pour $g > 1$, nous avons $\mathrm{Ext}^{2g-1} = 0$ dans $\mathcal{MH}_{\mathbf{R}}$ [52, (2.4)].

# Chapitre 3

# Calculs explicites dans le cas modulaire

## 3.1 Formulaire pour les séries d'Eisenstein

Dans cette section, nous suivons de près l'exposition remarquable de Siegel [69, pp. 1–73]. Le lecteur pourra y trouver les démonstrations que nous avons omises.

Soit $N$ un entier $\geq 1$. Le groupe $\mathrm{SL}_2(\mathbf{Z})$ et les sous-groupes

$$\Gamma_0(N) = \left\{ \begin{pmatrix} a & b \\ c & d \end{pmatrix} \in \mathrm{SL}_2(\mathbf{Z}) \,|\, c \equiv 0 \pmod{N} \right\} \text{ et} \tag{3.1}$$

$$\Gamma_1(N) = \left\{ \begin{pmatrix} a & b \\ c & d \end{pmatrix} \in \mathrm{SL}_2(\mathbf{Z}) \,|\, c \equiv 0 \pmod{N}, a \equiv d \equiv 1 \pmod{N} \right\}, \tag{3.2}$$

opèrent sur le demi-plan de Poincaré $\mathcal{H} = \{z \in \mathbf{C} \,|\, \Im(z) > 0\}$ par homographies :

$$\begin{pmatrix} a & b \\ c & d \end{pmatrix} \cdot z = \frac{az+b}{cz+d} \qquad \left( \begin{pmatrix} a & b \\ c & d \end{pmatrix} \in \mathrm{SL}_2(\mathbf{Z}),\, z \in \mathcal{H} \right). \tag{3.3}$$

Nous allons définir, au moyen de *séries d'Eisenstein*, des fonctions analytiques-réelles sur $\mathcal{H}$, invariantes sous l'action de $\Gamma_1(N)$ ou d'un sous-groupe conjugué dans $\mathrm{SL}_2(\mathbf{Z})$. Pour tout $z \in \mathcal{H}$, nous noterons $y = \Im(z)$. Pour $(u, v) \in (\frac{\mathbf{Z}}{N\mathbf{Z}})^2$, définissons

$$E_{(u,v)}(z,s) = \sideset{}{'}\sum_{\substack{m \equiv u \,(N) \\ n \equiv v \,(N)}} \frac{y^s}{|mz+n|^{2s}} \qquad (z \in \mathcal{H},\, \Re(s) > 1), \tag{3.4}$$

où le symbole $\sum'$ indique que la somme est étendue aux $(m, n) \in \mathbf{Z}^2$, $(m, n) \neq (0, 0)$. Nous allons aussi définir certaines combinaisons linéaires des séries $E_{(u,v)}$. Pour $(a, b) \in (\frac{\mathbf{Z}}{N\mathbf{Z}})^2$, définissons

$$\zeta_{a,b}(z,s) = \sideset{}{'}\sum_{(m,n) \in \mathbf{Z}^2} \frac{e^{\frac{2\pi i}{N}(ma+nb)} \cdot y^s}{|mz+n|^{2s}} \qquad (z \in \mathcal{H},\, \Re(s) > 1). \tag{3.5}$$

Les fonctions $\zeta_{a,b}(z, s)$ sont un cas très particulier de *fonctions zêta d'Epstein*. Remarquons l'identité

$$\zeta_{a,b}(z,s) = \sum_{(u,v) \in (\frac{\mathbf{Z}}{N\mathbf{Z}})^2} e^{\frac{2\pi i}{N}(au+bv)} E_{(u,v)}(z,s), \tag{3.6}$$

ainsi que la formule de transformée de Fourier inverse

$$E_{(u,v)}(z,s) = \frac{1}{N^2} \sum_{(a,b) \in (\frac{\mathbf{Z}}{N\mathbf{Z}})^2} e^{-\frac{2\pi i}{N}(au+bv)} \zeta_{a,b}(z,s). \tag{3.7}$$

Pour $(a,b) \in (\frac{\mathbf{Z}}{N\mathbf{Z}})^2$ et $z \in \mathcal{H}$ fixé, la fonction $s \mapsto \zeta_{a,b}(z,s)$ admet un prolongement méromorphe au plan complexe [69, Thm 3, p. 69]. Lorsque $(a,b) = (0,0)$, le prolongement a un unique pôle en $s = 1$, ce pôle est simple et son résidu (indépendant de $z$!) est égal à $\pi$. Lorsque $(a,b) \neq (0,0)$, le prolongement est holomorphe sur $\mathbf{C}$. Nous définissons pour $z \in \mathcal{H}$

$$\zeta_{a,b}^*(z) = \begin{cases} \lim_{s \to 1} \bigl( \zeta_{a,b}(z,s) - \frac{\pi}{s-1} \bigr) & \text{si } (a,b) = (0,0), \\ \zeta_{a,b}(z,1) & \text{si } (a,b) \neq (0,0). \end{cases} \tag{3.8}$$

Nous déduisons des affirmations précédentes le prolongement méromorphe des fonctions $s \mapsto E_{(u,v)}(z,s)$ au plan complexe, avec un unique pôle en $s = 1$, simple et de résidu égal à $\frac{\pi}{N^2}$. En accord avec la notation (3.8), nous définissons

$$E_x^*(z) = \lim_{s \to 1} \Bigl( E_x(z,s) - \frac{\pi}{N^2(s-1)} \Bigr) \qquad \bigl( x \in (\mathbf{Z}/N\mathbf{Z})^2,\ z \in \mathcal{H} \bigr). \tag{3.9}$$

Passons maintenant aux deux *formules-limite de Kronecker*. Ces formules donnent une expression de $\zeta_{0,0}^*(z)$ (resp. $\zeta_{a,b}^*(z)$ avec $(a,b) \neq (0,0)$). La première formule-limite de Kronecker [69, Thm 1, p. 17] s'écrit

$$\zeta_{0,0}^*(z) = 2\pi \bigl( \gamma - \log 2 - \log \sqrt{y} - 2 \log|\eta(z)| \bigr) \qquad (z \in \mathcal{H}), \tag{3.10}$$

où $\gamma$ désigne la constante d'Euler et

$$\eta(z) = e^{\frac{\pi i z}{12}} \prod_{n=1}^{\infty} (1 - e^{2\pi i n z}) \qquad (z \in \mathcal{H}). \tag{3.11}$$

La deuxième formule-limite de Kronecker [69, Thm 2, p. 40] s'écrit de la façon suivante. Soit $(a,b) \in (\frac{\mathbf{Z}}{N\mathbf{Z}})^2$ avec $(a,b) \neq (0,0)$. Choisissons un couple de représentants $(\tilde{a}, \tilde{b})$ de $(a,b)$ dans $\mathbf{Z}^2$. Alors

$$\zeta_{a,b}^*(z) = \frac{2\pi^2 \tilde{b}^2}{N^2} y - 2\pi \log \Bigl| \vartheta \Bigl( \frac{\tilde{a} - \tilde{b} z}{N}, z \Bigr) \Bigr|, \tag{3.12}$$

où nous avons posé, pour $w \in \mathbf{C}$ et $z \in \mathcal{H}$ :

$$\vartheta(w,z) = e^{\frac{\pi i z}{6}} (e^{\pi i w} - e^{-\pi i w}) \prod_{n=1}^{\infty} (1 - e^{2\pi i (w+nz)})(1 - e^{-2\pi i (w-nz)}). \tag{3.13}$$

Les séries $E_x$ vérifient la propriété de modularité suivante

$$E_x(gz, s) = E_{xg}(z, s) \qquad \bigl( x \in (\mathbf{Z}/N\mathbf{Z})^2,\ g \in \mathrm{SL}_2(\mathbf{Z}),\ z \in \mathcal{H},\ \Re(s) > 1 \bigr), \tag{3.14}$$

où $xg \in (\frac{\mathbf{Z}}{N\mathbf{Z}})^2$ désigne le produit du vecteur ligne $x$ par la matrice $g$. Nous en déduisons

$$E_x^*(gz) = E_{xg}^*(z) \qquad \bigl( x \in (\mathbf{Z}/N\mathbf{Z})^2,\ g \in \mathrm{SL}_2(\mathbf{Z}),\ z \in \mathcal{H} \bigr). \tag{3.15}$$

En particulier, pour tout $x \in (\frac{\mathbf{Z}}{N\mathbf{Z}})^2$, la fonction $z \mapsto E_x^*(z)$ est invariante par la transformation $z \mapsto z + N$, et admet donc un développement de Fourier (non holomorphe) en la variable $e^{\frac{2\pi i z}{N}}$. Il en va donc de même des fonctions $\zeta^*_{a,b}$, dont le développement de Fourier se déduit des formules-limite de Kronecker. Posons $q = e^{2\pi i z}$ et $q^{\frac{1}{N}} = e^{\frac{2\pi i z}{N}}$ pour $z \in \mathcal{H}$. Pour tout entier $r \geq 1$, notons $\sigma(r)$ la somme des diviseurs positifs de $r$. Nous avons alors

$$\zeta^*_{0,0}(z) = \frac{\pi^2 y}{3} - \pi \log y + 2\pi \Big( \gamma - \log 2 + \sum_{r=1}^{\infty} \frac{\sigma(r)}{r} (q^r + \overline{q}^r) \Big). \tag{3.16}$$

Écrivons ensuite le développement de Fourier de $\zeta^*_{a,0}$, avec $a \in \frac{\mathbf{Z}}{N\mathbf{Z}}$, $a \neq 0$. Notons $\zeta_N = e^{\frac{2\pi i}{N}}$. Nous avons alors

$$\zeta^*_{a,0}(z) = \frac{\pi^2 y}{3} - 2\pi \log|1 - \zeta_N^a| + \pi \sum_{r=1}^{\infty} \Big( \sum_{k \mid r} \frac{\zeta_N^{ka} + \zeta_N^{-ka}}{k} \Big)(q^r + \overline{q}^r). \tag{3.17}$$

Écrivons enfin le développement de Fourier de $\zeta^*_{a,b}$, avec $a \in \frac{\mathbf{Z}}{N\mathbf{Z}}$ et $b \in \frac{\mathbf{Z}}{N\mathbf{Z}}$, $b \neq 0$. Notons $B_2(X) = X^2 - X + \frac{1}{6}$ le deuxième polynôme de Bernoulli, et définissons une fonction 1-périodique $\overline{B_2}$ sur $\mathbf{R}$ par

$$\overline{B_2}(x) = B_2(x - \lfloor x \rfloor) \qquad (x \in \mathbf{R}), \tag{3.18}$$

où $\lfloor x \rfloor$ désigne le plus grand entier $\leq x$. Alors la quantité $\overline{B_2}(\frac{\tilde{b}}{N})$ ne dépend pas du représentant $\tilde{b}$ de $b$ dans $\mathbf{Z}$, et nous avons

$$\zeta^*_{a,b}(z) = 2\pi^2 \overline{B_2}\Big(\frac{\tilde{b}}{N}\Big) y + \pi \sum_{r=1}^{\infty} \alpha_r \cdot q^{\frac{r}{N}} + \overline{\alpha_r} \cdot \overline{q}^{\frac{r}{N}}, \tag{3.19}$$

où les coefficients $\alpha_r$ sont donnés par la formule

$$\alpha_r = \sum_{\substack{k \mid r \\ \frac{r}{k} \equiv b \ (N)}} \frac{\zeta_N^{-ka}}{k} + \sum_{\substack{k \mid r \\ \frac{r}{k} \equiv -b \ (N)}} \frac{\zeta_N^{ka}}{k} \qquad (r \geq 1). \tag{3.20}$$

Notons $E_N \subset (\frac{\mathbf{Z}}{N\mathbf{Z}})^2$ l'ensemble des éléments d'ordre $N$ du groupe additif $(\frac{\mathbf{Z}}{N\mathbf{Z}})^2$. Nous avons une bijection

$$\Gamma_1(N) \backslash \mathrm{SL}_2(\mathbf{Z}) \xrightarrow{\cong} E_N \tag{3.21}$$
$$\left[ \begin{pmatrix} a & b \\ c & d \end{pmatrix} \right] \mapsto (c, d) \mod N.$$

Pour tout $x \in E_N$, notons $g_x \in \mathrm{SL}_2(\mathbf{Z})$ un représentant de l'image réciproque de $x$ par la bijection (3.21). Il est alors facile de voir que la fonction $E^*_x$, définie sur $\mathcal{H}$, est invariante sous l'action du groupe $g_x^{-1} \Gamma_1(N) g_x$. Lorsque $x = (u, v) \notin E_N$, notons $d = (u, v, N) > 1$ et $x' = (\frac{u}{d}, \frac{v}{d}) \in (\mathbf{Z}/\frac{N}{d}\mathbf{Z})^2$. Alors $x' \in E_{\frac{N}{d}}$ et $E^*_x = d^{-2} E^*_{x'}$. En particulier, la fonction $E^*_x$ sur $\mathcal{H}$ est invariante sous l'action d'un sous-groupe conjugué de $\Gamma_1(\frac{N}{d})$ dans $\mathrm{SL}_2(\mathbf{Z})$; elle est *a fortiori* invariante sous l'action d'un sous-groupe conjugué de $\Gamma_1(N)$.

Le comportement des séries d'Eisenstein $E^*_x$ vis-à-vis de la conjugaison complexe est le suivant. Notons $c : z \mapsto -\overline{z}$ la conjugaison complexe sur $\mathcal{H}$. Pour tout $x = (u, v) \in (\frac{\mathbf{Z}}{N\mathbf{Z}})^2$, notons $x^c = (-u, v)$. Nous avons alors

$$E^*_x\big(c(z)\big) = E^*_{x^c}(z) \qquad \big(x \in (\mathbf{Z}/N\mathbf{Z})^2, \ z \in \mathcal{H}\big). \tag{3.22}$$

**Définition 70.** *Pour toute fonction* $f : \frac{\mathbf{Z}}{N\mathbf{Z}} \to \mathbf{C}$, *nous définissons les séries d'Eisenstein* $E_f$ *et* $E_f^*$ *par*

$$E_f(z,s) = \sum_{v \in \frac{\mathbf{Z}}{N\mathbf{Z}}} f(v)\, E_{(0,v)}(z,s) \tag{3.23}$$

$$E_f^*(z) = \sum_{v \in \frac{\mathbf{Z}}{N\mathbf{Z}}} f(v)\, E_{(0,v)}^*(z). \tag{3.24}$$

Remarquons que les notations $E_f$ et $E_f^*$ sont en accord avec le procédé (3.9). L'application $f \mapsto E_f^*$ est **C**-linéaire. Remarquons que (3.22) entraîne

$$E_f^*\big(c(z)\big) = E_f^*(z) \qquad (z \in \mathcal{H}). \tag{3.25}$$

La transformée de Fourier $\widehat{f} : \frac{\mathbf{Z}}{N\mathbf{Z}} \to \mathbf{C}$ de $f$ est définie par

$$\widehat{f}(b) = \sum_{v \in \frac{\mathbf{Z}}{N\mathbf{Z}}} f(v) \cdot e^{-\frac{2\pi i b v}{N}} \qquad (b \in \frac{\mathbf{Z}}{N\mathbf{Z}}). \tag{3.26}$$

Rappelons que $f$ et $\widehat{f}$ induisent des fonctions $N$-périodiques de **Z** dans **C**. Le développement de Fourier de $E_f^*$ se calcule aisément à l'aide des formules (3.16), (3.17) et (3.19).

**Proposition 71.** *Soit* $f : \frac{\mathbf{Z}}{N\mathbf{Z}} \to \mathbf{C}$ *une fonction de somme nulle. La série d'Eisenstein* $E_f^*$ *admet le développement de Fourier*

$$E_f^*(z) = \Big(\sideset{}{'}\sum_{n \in \mathbf{Z}} \frac{f(n)}{n^2}\Big) y + \frac{\pi}{N^2} \sum_{r=1}^{\infty} \frac{1}{r} \Big( \sum_{k \mid r} k\big(\widehat{f}(k) + \widehat{f}(-k)\big) \Big)(q^r + \overline{q}^r). \tag{3.27}$$

*Démonstration.* Puisque $f$ est de somme nulle, nous avons $\widehat{f}(0) = 0$. Soit $z \in \mathcal{H}$. Nous avons

$$\begin{aligned}
E_f^*(z) &= \sum_{v \in \frac{\mathbf{Z}}{N\mathbf{Z}}} f(v)\, E_{(0,v)}^*(z) \\
&= \sum_{v \in \frac{\mathbf{Z}}{N\mathbf{Z}}} f(v)\, \frac{1}{N^2} \sum_{(a,b) \in (\frac{\mathbf{Z}}{N\mathbf{Z}})^2} e^{-\frac{2\pi i b v}{N}} \zeta_{a,b}^*(z) \\
&= \frac{1}{N^2} \sum_{b \in \frac{\mathbf{Z}}{N\mathbf{Z}}} \widehat{f}(b) \Big( \sum_{a \in \frac{\mathbf{Z}}{N\mathbf{Z}}} \zeta_{a,b}^*(z) \Big).
\end{aligned}$$

Nous déduisons de (3.19) que pour $b \in \frac{\mathbf{Z}}{N\mathbf{Z}}$, $b \neq 0$, nous avons

$$\sum_{a \in \frac{\mathbf{Z}}{N\mathbf{Z}}} \zeta_{a,b}^*(z) = 2\pi^2 N \overline{B}_2(\frac{\tilde{b}}{N}) \cdot y + \pi \sum_{r=1}^{\infty} \frac{1}{r} \Big( \sum_{\substack{k \mid r \\ k \equiv b\ (N)}} k + \sum_{\substack{k \mid r \\ k \equiv -b\ (N)}} k \Big)(q^r + \overline{q}^r).$$

Il en résulte

$$E_f^*(z) = \frac{2\pi^2}{N} \Big( \sum_{b \in \frac{\mathbf{Z}}{N\mathbf{Z}}} \widehat{f}(b) \overline{B}_2(\frac{\tilde{b}}{N}) \Big) \cdot y + \frac{\pi}{N^2} \sum_{r=1}^{\infty} \frac{1}{r} \Big( \sum_{k \mid r} k\big(\widehat{f}(k) + \widehat{f}(-k)\big) \Big)(q^r + \overline{q}^r).$$

La fonction $\overline{B_2}$ est donnée par la série de Fourier suivante, qui converge normalement sur $\mathbf{R}$ [18, (1.56), p. 14]

$$\overline{B_2}(x) = \frac{1}{2\pi^2} \sum_{\substack{n \in \mathbf{Z} \\ n \neq 0}} \frac{e^{2\pi i n x}}{n^2} \qquad (x \in \mathbf{R}).$$

Nous en déduisons aisément

$$\sum_{b \in \frac{\mathbf{Z}}{N\mathbf{Z}}} \widehat{f}(b) \overline{B_2}(\frac{\tilde{b}}{N}) = \frac{N}{\pi^2} \sideset{}{'}\sum_{n \in \mathbf{Z}} \frac{f(n)}{n^2},$$

ce qui achève de montrer (3.27). □

Il est amusant de constater que le développement de Fourier de $E_f^*$ fait intervenir naturellement la transformée de Fourier de $f$. Typiquement, nous utiliserons la série d'Eisenstein $E_f^*$ dans le cas où $f$ est un caractère de Dirichlet modulo $N$, ou la transformée de Fourier d'un tel caractère.

## 3.2  Calcul d'une intégrale par la méthode de Rankin-Selberg

Soit $N$ un entier $\geq 1$. Notons $S_2(\Gamma_1(N))$ l'espace des formes paraboliques de poids 2 pour le groupe $\Gamma_1(N)$. Il s'identifie canoniquement à l'espace des formes différentielles holomorphes sur $X_1(N)(\mathbf{C})$, au moyen de l'application $f \mapsto \omega_f := 2\pi i f(z) dz$. Cette section a pour but de calculer l'intégrale

$$\int_{X_1(N)(\mathbf{C})} E_\chi^* \cdot \omega_f \wedge \overline{\partial} E_{\widehat{\chi'}}^*, \qquad (3.28)$$

où $f$ est une forme primitive de poids 2 pour $\Gamma_1(N)$, et $\chi$ (resp. $\chi'$) est un caractère de Dirichlet pair modulo $N$ (resp. modulo un diviseur $M$ de $N$). D'après la section précédente, les séries d'Eisenstein $E_\chi^*$ et $E_{\widehat{\chi'}}^*$ sont des fonctions sur $\mathcal{H}$ respectivement invariantes sous l'action de $\Gamma_1(N)$ et $\Gamma_1(M)$, et donnent donc lieu à des fonctions sur $Y_1(N)(\mathbf{C})$. Nous obtenons une expression exacte pour (3.28) en termes de valeurs spéciales de la fonction $L$ de $f$ (éventuellement tordue). Signalons que Scholl [66, Thm 4.6.3] a également obtenu une formule exacte pour une intégrale analogue, en utilisant le langage adélique. Je n'ai pas encore pu comparer les deux résultats.

La *fonction L* (ou *série L*) d'une forme parabolique $f(z) = \sum_{n \geq 1} a_n e^{2\pi i n z}$, $z \in \mathcal{H}$ est définie comme la série de Dirichlet

$$L(f, s) = \sum_{n=1}^{\infty} \frac{a_n}{n^s} \qquad \left(\Re(s) > \frac{3}{2}\right). \qquad (3.29)$$

Elle admet un prolongement holomorphe au plan complexe. La fonction $L$ de $f$ tordue par un caractère de Dirichlet $\chi$ arbitraire

$$L(f, \chi, s) = \sum_{n=1}^{\infty} \frac{a_n \chi(n)}{n^s} \qquad \left(\Re(s) > \frac{3}{2}\right) \qquad (3.30)$$

possède la même propriété.

**Théorème 72.** *Sous les hypothèses précédentes, notons $\psi$ le caractère de $f$ et $\chi'_N$ le caractère modulo $N$ induit par $\chi'$. Nous avons*

$$\int_{X_1(N)(\mathbf{C})} E^*_\chi \cdot \omega_f \wedge \overline{\partial} E^*_{\widehat{\chi'}} = \begin{cases} -i\pi \frac{\varphi(N)}{M} \cdot L(f,2) L(f, \chi', 1) & \text{si } \psi = \chi \overline{\chi'_N}, \\ 0 & \text{sinon.} \end{cases} \qquad (3.31)$$

Il est intéressant de reformuler le théorème précédent en utilisant l'algèbre de Hecke et la forme modulaire universelle. Cette démarche m'a conduit à une formule plus agréable et naturelle. Rappelons que l'*algèbre de Hecke* $\mathbf{T} \subset \mathrm{End}_{\mathbf{C}} S_2(\Gamma_1(N))$ est le sous-anneau engendré par tous les opérateurs de Hecke $T_n$, $n \geq 1$ et les opérateurs diamants $\langle d \rangle$, $d \in (\frac{\mathbf{Z}}{N\mathbf{Z}})^*$. Nous avons un isomorphisme canonique [47, Lemma 9]

$$\begin{aligned} \mathbf{T} \otimes \mathbf{C} &\xrightarrow{\cong} \mathrm{Hom}_{\mathbf{C}}(S_2(\Gamma_1(N)), \mathbf{C}) \\ T &\mapsto \bigl(f \mapsto a_1(Tf)\bigr), \end{aligned} \qquad (3.32)$$

où $a_1(\cdot)$ désigne le premier coefficient de Fourier d'une forme modulaire. La fonction $L$ (éventuellement tordue) de l'algèbre de Hecke est définie par

$$L(\mathbf{T}, s) = \sum_{n=1}^\infty T_n \otimes \frac{1}{n^s} \qquad L(\mathbf{T}, \chi, s) = \sum_{n=1}^\infty T_n \otimes \frac{\chi(n)}{n^s} \qquad \left(\Re(s) > \frac{3}{2}\right). \qquad (3.33)$$

Elle est à valeurs dans $\mathbf{T} \otimes \mathbf{C}$ et admet un prolongement holomorphe au plan complexe. Elle s'interprète comme la fonction $L$ de la *forme modulaire universelle* [47]

$$\Omega = 2\pi i \sum_{n=1}^\infty T_n \cdot e^{2\pi i n z} dz \in S_2(\Gamma_1(N)) \otimes \mathbf{T}. \qquad (3.34)$$

Pour tout caractère de Dirichlet $\psi$ modulo $N$, notons $\mathbf{T}^\psi$ la *composante $\psi$-isotypique* de $\mathbf{T} \otimes \mathbf{C}$

$$\mathbf{T}^\psi = \{T \in \mathbf{T} \otimes \mathbf{C} \mid \forall d \in (\mathbf{Z}/N\mathbf{Z})^*, \ T \circ \langle d \rangle = \psi(d) \cdot T\}. \qquad (3.35)$$

Nous avons une décomposition canonique de $\mathbf{T} \otimes \mathbf{C}$ en produit de sous-algèbres

$$\mathbf{T} \otimes \mathbf{C} \cong \prod_\psi \mathbf{T}^\psi, \qquad (3.36)$$

le produit étant étendu aux caractères de Dirichlet $\psi$ pairs modulo $N$. Les projections de $L(\mathbf{T}, s)$ et $L(\mathbf{T}, \chi, s)$ sur $\mathbf{T}^\psi$ seront notées respectivement $L(\mathbf{T}^\psi, s)$ et $L(\mathbf{T}^\psi, \chi, s)$.

Nous démontrerons le résultat suivant, qui comprend le théorème 72 comme cas particulier (appliquer une forme primitive $f \in S_2(\Gamma_1(N))$).

**Théorème 73.** *Soient $\chi$ un caractère de Dirichlet pair modulo $N$ et $\chi'$ un caractère de Dirichlet pair modulo un diviseur $M$ de $N$. Posons $\psi = \chi \overline{\chi'_N}$. Nous avons alors*

$$\int_{X_1(N)(\mathbf{C})} E^*_\chi \cdot \Omega \wedge \overline{\partial} E^*_{\widehat{\chi'}} = -i\pi \frac{\varphi(N)}{M} \cdot L(\mathbf{T}^\psi, 2) L(\mathbf{T}^\psi, \chi', 1). \qquad (3.37)$$

Faisons maintenant la somme de (3.37) sur les caractères $\chi$ pairs modulo $N$. La série d'Eisenstein $E^*_\chi$ étant nulle lorsque $\chi$ est impair, nous avons

$$\sum_{\chi \text{ pair}} E^*_\chi = \sum_\chi E^*_\chi = \varphi(N) \cdot E^*_{(0,1)},$$

où $E^*_{(0,1)}$ est la série d'Eisenstein associée à $(0,1) \in E_N$. Sommer le membre de gauche de (3.37) sur $\chi$ équivaut à sommer le membre de droite sur $\psi$, ce qui conduit à la forme suivante de notre résultat.

**Théorème 74.** *Soit $M$ un diviseur de $N$ et $\chi$ un caractère de Dirichlet pair modulo $M$. Nous avons l'identité*

$$\int_{X_1(N)(\mathbf{C})} E^*_{(0,1)} \cdot \Omega \wedge \overline{\partial} E^*_{\widehat{\chi}} = -\frac{\pi i}{M} \cdot L(\mathbf{T}, 2) L(\mathbf{T}, \chi, 1), \tag{3.38}$$

*où $E^*_{(0,1)}$ est la série d'Eisenstein associée à $(0,1) \in E_N$.*

*Remarque* 75. L'identité (3.38) est plus intrinsèque que (3.37). Curieusement, il semble que nous soyons obligés de passer par le théorème 73 pour démontrer le théorème 74. Une démonstration plus directe serait la bienvenue.

*Remarque* 76. Il n'y a pas de raison a priori de se limiter à la torsion par un caractère dans (3.38) : il est possible de tordre la fonction $L$ par n'importe quelle application $\alpha : \frac{\mathbf{Z}}{M\mathbf{Z}} \to \mathbf{C}$ ; d'ailleurs la formule (3.38) est linéaire en l'application $\chi$. Il est raisonnable d'espérer que (3.38) reste valable pour toute application $\alpha$ paire, mais nous n'avons pas de démonstration. Il serait intéressant d'envisager une version (et une démonstration) du théorème de Beĭlinson dans le cas où la fonction $L$ de $f$ est tordue par un caractère additif.

*Démonstration du théorème 73.* La démonstration consiste en deux grandes étapes. La première étape, de nature globale, utilise la méthode de Rankin-Selberg et exprime l'intégrale en termes d'une convolution de séries de Dirichlet (formule (3.45)). Pour une introduction à la méthode de Rankin-Selberg, voir [80, 3. B]. La seconde étape, de nature locale, exprime la série de Dirichlet précédente comme un produit eulérien (lemme 77). Il est à noter que jusqu'au bout du calcul, nous tiendrons compte des facteurs locaux aux mauvaises places, c'est-à-dire aux nombres premiers divisant $N$.

Notons $I$ le membre de gauche de (3.37). Montrons que $I$ appartient à $\mathbf{T}^\psi$. Soit $d \in (\frac{\mathbf{Z}}{N\mathbf{Z}})^*$. Si nous effectuons le changement de variables $z \to \langle d \rangle z$ dans l'intégrale $I$, nous obtenons grâce à (3.74) et (3.75)

$$I = \overline{\chi}(d) \, \chi'_N(d) \int_{X_1(N)(\mathbf{C})} E^*_\chi \cdot \langle d \rangle^* \Omega \wedge \overline{\partial} E^*_{\widehat{\chi'}},$$

où $\langle d \rangle^* \Omega$ désigne l'image réciproque de la forme différentielle $\Omega$ par l'automorphisme $\langle d \rangle$. L'isomorphisme (3.32) étant compatible à l'action des opérateurs diamants, nous avons $\langle d \rangle^* \Omega = \Omega \otimes \langle d \rangle$, où à droite $\langle d \rangle$ agit par multiplication sur le facteur $\mathbf{T}$ du produit tensoriel. Il en résulte

$$I \circ \langle d \rangle = \chi(d) \, \overline{\chi'_N}(d) \, I = \psi(d) \, I \qquad (d \in (\mathbf{Z}/N\mathbf{Z})^*),$$

c'est-à-dire exactement $I \in \mathbf{T}^\psi$. Pour obtenir la composante de caractère $\psi$ de l'intégrale $I$, il suffit de remplacer $\Omega$ par sa composante de caractère $\psi$

$$\Omega^\psi \in S_2(\Gamma_1(N)) \otimes_{\mathbf{C}} \mathbf{T}^\psi, \tag{3.39}$$

ce que nous ferons désormais. Remarquons que $\langle d \rangle^* \Omega^\psi = \psi(d) \cdot \Omega^\psi$, $d \in (\frac{\mathbf{Z}}{N\mathbf{Z}})^*$, c'est-à-dire $\Omega^\psi \in S_2(\Gamma_1(N), \psi) \otimes_{\mathbf{C}} \mathbf{T}^\psi$, où $S_2(\Gamma_1(N), \psi)$ désigne le sous-espace des formes de caractère $\psi$.

Pour calculer $I$, nous pouvons considérer l'intégrale étendue au domaine $\Gamma_1(N) \backslash \mathcal{H}$. Pour utiliser la méthode de Rankin-Selberg, nous allons remplacer $E^*_\chi$ par une somme indexée par $\Gamma_\infty \backslash \Gamma_0(N)$, où nous rappelons que $\Gamma_\infty$ est le sous-groupe de $SL_2(\mathbf{Z})$ engendré par les matrices $\begin{pmatrix} -1 & 0 \\ 0 & -1 \end{pmatrix}$ et $T = \begin{pmatrix} 1 & 1 \\ 0 & 1 \end{pmatrix}$. La fonction $E_\chi^*$ est définie en appliquant le procédé (3.9) à la fonction $E_\chi$. Remarquons que le résidu en $s=1$ de la fonction $s \mapsto E_\chi(z,s)$ est indépendant de $z \in \mathcal{H}$. Or

$$\int_{Y_1(N)(\mathbf{C})} \Omega^\psi \wedge \overline{\partial} E_{\widehat{\chi}'}^* = \int_{Y_1(N)(\mathbf{C})} d\bigl(-E_{\widehat{\chi}'}^* \cdot \Omega^\psi\bigr) = 0,$$

d'après la formule de Stokes et le fait que la fonction $E_{\widehat{\chi}'}^*$ est à croissance modérée aux pointes. Dans l'intégrale $I$, nous pouvons donc remplacer $E_\chi^*$ par $E_\chi(\cdot, s)$, puis faire $s=1$ :

$$I = \bigl(I(s)\bigr)_{s=1} = \left(\int_{Y_1(N)(\mathbf{C})} E_\chi(\cdot, s) \cdot \Omega^\psi \wedge \overline{\partial} E_{\widehat{\chi}'}^*\right)_{s=1}.$$

Le caractère analytique de la fonction $I(s)$ pour $s \neq 1$, et la possibilité d'intervertir le signe $\int$ et l'opération $(\cdot)_{s=1}$, résultent du fait qu'en chaque pointe de $X_1(N)(\mathbf{C})$, la fonction $z \mapsto E_\chi(z,s)$ possède un développement de Fourier convergeant uniformément sur tout compact par rapport à $s$. Maintenant, nous avons

$$E_\chi(z,s) = \sum_{v \in (\frac{\mathbf{Z}}{N\mathbf{Z}})^*} \chi(v)\, E_{(0,v)}(z,s) = \sum_{\substack{m \equiv 0 \ (N) \\ (n,N)=1}} \frac{\chi(n)\, y^s}{|mz+n|^{2s}}.$$

Nous pouvons récrire cette somme en distinguant les valeurs du p.g.c.d. $d=(m,n)$, qui est $\geq 1$ et premier à $N$. Nous obtenons

$$\begin{aligned}
E_\chi(z,s) &= \sum_{\substack{d \geq 1 \\ (d,N)=1}} \sum_{\substack{m \equiv 0 \ (N) \\ (n,N)=1 \\ (m,n)=d}} \frac{\chi(n)\, y^s}{|mz+n|^{2s}} \\
&= \sum_{\substack{d \geq 1 \\ (d,N)=1}} \sum_{\substack{(\mu,\nu)=1 \\ d\mu \equiv 0 \ (N) \\ (d\nu, N)=1}} \frac{\chi(d)}{d^{2s}} \cdot \frac{\chi(\nu)\, y^s}{|\mu z + \nu|^{2s}} \\
&= L(\chi, 2s) \sum_{\substack{(\mu,\nu)=1 \\ \mu \equiv 0 \ (N) \\ (\nu,N)=1}} \frac{\chi(\nu)\, y^s}{|\mu z + \nu|^{2s}}.
\end{aligned}$$

Or nous avons une bijection

$$\Gamma_\infty \backslash \Gamma_0(N) \xrightarrow{\cong} \bigl\{(\mu,\nu)=1 \mid \mu \equiv 0 \ (N)\bigr\}/\pm 1$$
$$\left[\begin{pmatrix} a & b \\ c & d \end{pmatrix}\right] \mapsto [(c,d)].$$

Puisque $-1$ agit sans point fixe sur l'ensemble des couples $(\mu,\nu)$ ci-dessus, il vient

$$E_\chi(z,s) = 2L(\chi, 2s) \sum_{\gamma \in \Gamma_\infty \backslash \Gamma_0(N)} \chi(\gamma)\, \Im(\gamma z)^s,$$

où nous avons posé $\chi(\gamma) = \chi(d)$ pour $\gamma = \begin{pmatrix} a & b \\ c & d \end{pmatrix} \in \Gamma_0(N)$. Pour $z \in \mathcal{H}$, notons $z = x + iy$ et posons

$$\Omega^\psi \wedge \overline{\partial} E^*_{\widehat{\chi'}} = F(z) \cdot \frac{dx \wedge dy}{y^2}, \tag{3.40}$$

avec $F : \mathcal{H} \to \mathbf{T}^\psi$ de classe $\mathcal{C}^\infty$. La forme différentielle $\Omega^\psi \wedge \overline{\partial} E^*_{\widehat{\chi'}}$ est de caractère $\chi = \psi\overline{\chi'_N}$, tandis que $dx \wedge dy/y^2$ est invariante sous l'action de $\mathrm{SL}_2(\mathbf{R})$. Nous avons donc

$$F(\gamma z) = \chi(\gamma) F(z) \qquad (\gamma \in \Gamma_0(N), z \in \mathcal{H}). \tag{3.41}$$

Nous en déduisons

$$I(s) = 2L(\chi, 2s) \int_{\Gamma_1(N)\backslash\mathcal{H}} \sum_{\gamma \in \Gamma_\infty\backslash\Gamma_0(N)} \chi(\gamma)\, \Im(\gamma z)^s \cdot F(z) \cdot \frac{dx \wedge dy}{y^2}$$

$$= 2L(\chi, 2s) \int_{\Gamma_1(N)\backslash\mathcal{H}} \sum_{\gamma \in \Gamma_\infty\backslash\Gamma_0(N)} \Im(\gamma z)^s \cdot F(\gamma z) \cdot \frac{dx \wedge dy}{y^2}.$$

L'espace $S_2(\Gamma_1(N))$ étant trivial pour $N = 1$ ou $2$, nous pouvons supposer $N \geq 3$; par suite le morphisme $\Gamma_1(N)\backslash\mathcal{H} \to \Gamma_0(N)\backslash\mathcal{H}$ est fini, de degré $\varphi(N)/2$. Par conséquent

$$I(s) = \varphi(N) L(\chi, 2s) \int_{\Gamma_0(N)\backslash\mathcal{H}} \sum_{\gamma \in \Gamma_\infty\backslash\Gamma_0(N)} \Im(\gamma z)^s \cdot F(\gamma z) \cdot \frac{dx \wedge dy}{y^2}$$

$$= \varphi(N) L(\chi, 2s) \int_{\Gamma_\infty\backslash\mathcal{H}} \Im(z)^s \cdot F(z) \cdot \frac{dx \wedge dy}{y^2}.$$

La dernière égalité est le point-clé de la méthode de Rankin-Selberg.

Développons maintenant $F$ en série de Fourier

$$F(x + iy) = \sum_{m \in \mathbf{Z}} \widehat{F}_m(y)\, e^{2\pi i m x} \qquad (x + iy \in \mathcal{H}). \tag{3.42}$$

Un calcul simple utilisant la définition de $\Omega^\psi$ et $E^*_{\widehat{\chi'}}$, ainsi que le développement de Fourier (3.27), donne

$$\widehat{F}_0(y) = -\frac{16\pi^3 i}{M} y^2 \sum_{n=1}^\infty c_n\, e^{-4\pi n y} \qquad (y > 0), \tag{3.43}$$

$$\text{avec } c_n = T_n^\psi \cdot \sum_{k \mid n} k\chi'(k) \in \mathbf{T}^\psi \qquad (n \geq 1), \tag{3.44}$$

où nous notons $T_n^\psi$ l'image de $T_n$ dans $\mathbf{T}^\psi$ pour tout $n \geq 1$. Noter que dans l'unique cas $M = 1$, la formule (3.27) ne s'applique pas à $f = \widehat{\chi'}$, mais (3.16) permet quand même de mener le calcul, aboutissant au même résultat.

$$I(s) = \varphi(N) L(\chi, 2s) \int_0^1 \int_0^\infty y^s \sum_{m \in \mathbf{Z}} \widehat{F}_m(y) \, e^{2\pi i m x} \cdot \frac{dx \wedge dy}{y^2}$$

$$= \varphi(N) L(\chi, 2s) \int_0^\infty y^s \, \widehat{F}_0(y) \, \frac{dy}{y^2}$$

$$= -\frac{16\pi^3 i \varphi(N)}{M} L(\chi, 2s) \sum_{n=1}^\infty c_n \int_0^\infty y^s e^{-4\pi n y} \, dy.$$

Puisque $\int_0^\infty y^s e^{-4\pi n y} \, dy = \frac{\Gamma(s+1)}{(4\pi n)^{s+1}}$, nous obtenons

$$I(s) = -\frac{16\pi^3 i \varphi(N)}{M} \cdot \frac{\Gamma(s+1)}{(4\pi)^{s+1}} L(\chi, 2s) \sum_{n=1}^\infty \frac{c_n}{n^{s+1}}, \tag{3.45}$$

les coefficients $c_n$ étant donnés par (3.44). Nous voici arrivés au terme de la première étape du calcul.

**Lemme 77** (Une convolution de séries de Dirichlet). *Soient $\psi$ un caractère de Dirichlet modulo $N$ et $\chi_1$, $\chi_2$ des caractères de Dirichlet arbitraires. Posons*

$$\sigma_{\chi_1, \chi_2}(n) = \sum_{d \mid n} d \, \chi_1(d) \, \chi_2\Big(\frac{n}{d}\Big) \qquad (n \geq 1). \tag{3.46}$$

*Nous avons alors pour $s \in \mathbf{C}$, $\Re(s) > \frac{5}{2}$ :*

$$\sum_{n=1}^\infty T_n^\psi \cdot \frac{\sigma_{\chi_1, \chi_2}(n)}{n^s} = \frac{L(\mathbf{T}^\psi, \chi_2, s) \cdot L(\mathbf{T}^\psi, \chi_1, s-1)}{L(\psi \chi_1 \chi_2, 2s - 2)}, \tag{3.47}$$

*où nous avons posé $L(\psi \chi_1 \chi_2, s) = \sum_{n=1}^\infty \frac{\psi(n) \chi_1(n) \chi_2(n)}{n^s}$.*

*Démonstration.* Nous avons les estimations $T_n = O(n^{\frac{1}{2}+\epsilon})$ et $\sigma_{\chi_1, \chi_2}(n) = O(n^{1+\epsilon})$ pour tout $\epsilon > 0$, ce qui montre la convergence absolue de la série du membre de gauche de (3.47) pour $\Re(s) > \frac{5}{2}$. La fonction arithmétique $\sigma_{\chi_1, \chi_2}$ est convolution de deux fonctions multiplicatives. Elle est donc faiblement multiplicative *i. e.* vérifie

$$\sigma_{\chi_1, \chi_2}(mn) = \sigma_{\chi_1, \chi_2}(m) \, \sigma_{\chi_1, \chi_2}(n) \qquad \big((m,n) = 1\big).$$

Il en va de même de la fonction $n \mapsto T_n^\psi \cdot \frac{1}{n^s}$. Il suit que le membre de gauche de (3.47) admet l'expression en produit eulérien

$$\prod_{p \text{ premier}} \bigg( \sum_{a=0}^\infty T_{p^a}^\psi \cdot \frac{\sigma_{\chi_1, \chi_2}(p^a)}{p^{as}} \bigg). \tag{3.48}$$

D'autre part, nous avons formellement

$$L_p(\mathbf{T}^\psi, X) := \sum_{a=0}^\infty T_{p^a}^\psi \cdot X^a = \frac{1}{1 - T_p^\psi \cdot X + p \psi(p) \cdot X^2} \in \mathbf{T}^\psi[[X]], \tag{3.49}$$

où 1 désigne l'élément unité de $\mathbf{T}^\psi$. Nous pouvons calculer $\sigma_{\chi_1, \chi_2}(p^a)$ grâce à la multiplicativité de $\chi_1$ et $\chi_2$. Nous trouvons

$$\sigma_{\chi_1,\chi_2}(p^a) = \begin{cases} \frac{\chi_2(p)^{a+1} - (p\chi_1(p))^{a+1}}{\chi_2(p) - p\chi_1(p)} & \text{si } \chi_1(p) \neq 0 \text{ ou } \chi_2(p) \neq 0 \, ; \\ 1 & \text{si } \chi_1(p) = \chi_2(p) = 0 \text{ et } a = 0 \, ; \\ 0 & \text{si } \chi_1(p) = \chi_2(p) = 0 \text{ et } a \geq 1. \end{cases}$$

Il en résulte que, pour $p$ premier tel que $\chi_1(p) \neq 0$ ou $\chi_2(p) \neq 0$, le facteur local en $p$ du produit eulérien (3.48) est donné par

$$\frac{\chi_2(p)}{\chi_2(p) - p\chi_1(p)} \cdot \frac{1}{1 - T_p^\psi \cdot \chi_2(p) p^{-s} + \psi(p) \chi_2(p)^2 p^{1-2s}}$$
$$- \frac{p\chi_1(p)}{\chi_2(p) - p\chi_1(p)} \cdot \frac{1}{1 - T_p^\psi \cdot \chi_1(p) p^{1-s} + \psi(p) \chi_1(p)^2 p^{3-2s}},$$

soit après simplifications

$$(1 - \psi(p)\chi_1(p)\chi_2(p) \cdot p^{2-2s}) \cdot L_p(\mathbf{T}^\psi, \chi_2(p) p^{-s}) \cdot L_p(\mathbf{T}^\psi, \chi_1(p) p^{1-s}), \tag{3.50}$$

et ce dernier résultat est encore valable lorsque $\chi_1(p) = \chi_2(p) = 0$. Pour tout caractère de Dirichlet $\mu$, nous avons

$$\prod_{p \text{ premier}} L_p(\mathbf{T}^\psi, \mu(p) p^{-s}) = L(\mathbf{T}^\psi, \mu, s) \qquad \left(\Re(s) > \frac{3}{2}\right). \tag{3.51}$$

En prenant le produit sur tous les nombres premiers à partir de l'expression (3.50), nous obtenons le résultat souhaité. □

*Suite et fin de la démonstration du théorème 73.* Reprenons l'égalité (3.45). Utilisons le lemme 77 avec $\chi_1 = \chi'$ (de niveau $M$) et $\chi_2 = 1$ (de niveau 1). Il vient

$$I(s) = -\frac{16\pi^3 i \varphi(N)}{M} \cdot \frac{\Gamma(s+1)}{(4\pi)^{s+1}} L(\chi, 2s) \frac{L(\mathbf{T}^\psi, s+1) \cdot L(\mathbf{T}^\psi, \chi', s)}{L(\psi\chi', 2s)}$$
$$= -\frac{16\pi^3 i \varphi(N)}{M} \cdot \frac{\Gamma(s+1)}{(4\pi)^{s+1}} L(\mathbf{T}^\psi, s+1) \cdot L(\mathbf{T}^\psi, \chi', s),$$

puisque $L(\psi\chi', s) = L(\psi\chi'_N, s) = L(\chi, s)$. La fonction

$$s \mapsto N^{s/2} (2\pi)^{-s} \Gamma(s) L(\mathbf{T}, s)$$

admettant un prolongement *holomorphe* au plan complexe, il en va de même de la fonction $s \mapsto I(s)$. En évaluant en $s = 1$, il vient finalement

$$I = I(1) = -\frac{\pi i \varphi(N)}{M} \cdot L(\mathbf{T}^\psi, 2) L(\mathbf{T}^\psi, \chi', 1).$$

□

## 3.3 Rappels sur les unités modulaires

Dans cette section, nous faisons le lien entre les séries d'Eisenstein définies à la section 3.1 et certaines unités modulaires. Ces résultats sont bien connus [39], mais il nous est utile de les rappeler pour fixer les notations.

L'ensemble quotient

$$Y_1(N)(\mathbf{C}) = \Gamma_1(N)\backslash\mathcal{H} \tag{3.52}$$

peut être muni d'une structure de surface de Riemann (non compacte), appelée *courbe modulaire ouverte*. L'ensemble quotient

$$X_1(N)(\mathbf{C}) = \Gamma_1(N)\backslash(\mathcal{H} \sqcup \mathbf{P}^1(\mathbf{Q})), \tag{3.53}$$

obtenu à partir de $Y_1(N)(\mathbf{C})$ en rajoutant l'ensemble

$$P_N = \Gamma_1(N)\backslash\mathbf{P}^1(\mathbf{Q}), \tag{3.54}$$

peut être muni d'une structure de surface de Riemann compacte, appelée *courbe modulaire complétée*. L'ensemble $P_N$ est appelé *ensemble des pointes* de $X_1(N)(\mathbf{C})$. La *pointe infinie* est par définition la classe de $\infty \in \mathbf{P}^1(\mathbf{Q})$; elle est encore notée $\infty$. Notons $\Gamma_\infty$ le sous-groupe de $\operatorname{SL}_2(\mathbf{Z})$ engendré par les matrices $T = \begin{pmatrix} 1 & 1 \\ 0 & 1 \end{pmatrix}$ et $\begin{pmatrix} -1 & 0 \\ 0 & -1 \end{pmatrix}$. Le groupe $\Gamma_\infty$ opère par multiplication à droite sur $\operatorname{SL}_2(\mathbf{Z})$ et sur $E_N$, et nous avons

$$\begin{aligned} \operatorname{SL}_2(\mathbf{Z})/\Gamma_\infty &\xrightarrow{\cong} \mathbf{P}^1(\mathbf{Q}) \\ \left[\begin{pmatrix} a & b \\ c & d \end{pmatrix}\right] &\mapsto \frac{a}{c}. \end{aligned} \tag{3.55}$$

Nous en déduisons la suite de bijections

$$P_N \cong \Gamma_1(N)\backslash\mathbf{P}^1(\mathbf{Q}) \cong \Gamma_1(N)\backslash\operatorname{SL}_2(\mathbf{Z})/\Gamma_\infty \cong E_N/\Gamma_\infty. \tag{3.56}$$

En particulier, l'ensemble $P_N$ est fini. Pour tout $(u,v) \in E_N$, nous notons $[u,v]$ la classe de $(u,v)$ dans $P_N = E_N/\Gamma_\infty$. Nous avons donc $\infty = [0,1]$.

**Lemme 78.** *Pour toute application $f : \frac{\mathbf{Z}}{N\mathbf{Z}} \to \mathbf{C}$ de somme nulle, la série d'Eisenstein $E_f^*$ induit une fonction $Y_1(N)(\mathbf{C}) \to \mathbf{C}$ de classe $\mathcal{C}^\infty$, qui vérifie $\partial\bar{\partial}E_f^* = 0$.*

*Démonstration.* D'après (3.15), la fonction $E_f^* : \mathcal{H} \to \mathbf{C}$ est invariante sous l'action de $\Gamma_1(N)$, donc induit une fonction sur $Y_1(N)(\mathbf{C})$. Montrons que cette fonction est de classe $\mathcal{C}^\infty$. L'identité (3.27) nous permet d'écrire $E_f^* = E_1 + E_2$, où $E_1$ (resp. $E_2$) est une fonction holomorphe (resp. antiholomorphe) sur $\mathcal{H}$. Notons $\pi : \mathcal{H} \to Y_1(N)(\mathbf{C})$ la projection naturelle. Soit $z_0 \in \mathcal{H}$. D'après [15, Ex. i), p. 75], nous pouvons trouver une coordonnée locale holomorphe $u$ (resp. $v$) au point $z_0 \in \mathcal{H}$ (resp. $\pi(z_0) \in Y_1(N)(\mathbf{C})$), de telle sorte que la fonction $\pi$ soit donnée au voisinage de $z_0$ par $v = \pi(u) = u^n$, où $n$ est un entier $\geq 1$ (l'indice de ramification de $\pi$ en $z_0$). Dans ces coordonnées, nous avons donc

$$E_f^*(v) = E_1(u) + E_2(u). \tag{3.57}$$

Soit $\zeta_n = e^{\frac{2\pi i}{n}}$. D'après l'équation (3.57), nous avons $E_1(u\zeta_n) + E_2(u\zeta_n) = E_1(u) + E_2(u)$. Par conséquent, la fonction

$$u \mapsto E_1(u\zeta_n) - E_1(u) = E_2(u\zeta_n) - E_2(u)$$

est holomorphe et antiholomorphe, donc constante au voisinage de 0. Cette constante vaut $E_1(0) - E_1(0) = 0$, d'où $E_1(u\zeta_n) = E_1(u)$ et $E_2(u\zeta_n) = E_2(u)$. Donc $E_1$ (resp. $E_2$) induit une fonction holomorphe (resp. antiholomorphe) de $v$. D'après (3.57), la fonction $E_f^*$ est de classe $\mathcal{C}^\infty$ sur un voisinage de $\pi(z_0)$ dans $Y_1(N)(\mathbf{C})$ ; de plus, elle vérifie $\partial\bar{\partial}E_f^* = 0$ sur ce voisinage. Le résultat étant vrai localement, il est vrai sur $Y_1(N)(\mathbf{C})$. □

Nous notons $\mathbf{C}(X_1(N))$ le corps des fonctions méromorphes de la surface de Riemann compacte $X_1(N)(\mathbf{C})$. Les éléments de $\mathbf{C}(X_1(N))$ sont appelés *fonctions modulaires* pour $\Gamma_1(N)$. Par définition, le groupe $\mathcal{O}^*(Y_1(N)(\mathbf{C}))$ des *unités modulaires* est le sous-groupe de $\mathbf{C}(X_1(N))^*$ formé des fonctions $u$ dont le diviseur est à support dans $P_N$, c'est-à-dire vérifiant $(u) \subset P_N$. Notons $\mathrm{Div}(P_N)$ (resp. $\mathrm{Div}^0(P_N)$) le groupe des diviseurs (resp. diviseurs de degré 0) sur $P_N$. Le diviseur d'une unité modulaire $u \in \mathcal{O}^*(Y_1(N)(\mathbf{C}))$ vérifie $\mathrm{div}\, u \in \mathrm{Div}^0(P_N)$. Le théorème de Manin-Drinfel'd [28, 29] énonce que l'application naturelle

$$\mathrm{div} \otimes \mathbf{Q} : \mathcal{O}^*(Y_1(N)(\mathbf{C})) \otimes \mathbf{Q} \to \mathrm{Div}^0(P_N) \otimes \mathbf{Q} \tag{3.58}$$

est surjective. Nous en déduisons la suite exacte

$$0 \to \mathbf{C}^* \otimes \mathbf{C} \to \mathcal{O}^*(Y_1(N)(\mathbf{C})) \otimes \mathbf{C} \to \mathrm{Div}^0(P_N) \otimes \mathbf{C} \to 0 \tag{3.59}$$

où la flèche de gauche est déduite de l'inclusion naturelle $\mathbf{C} \subset \mathbf{C}(X_1(N))$. Toute unité modulaire $u \in \mathcal{O}^*(Y_1(N)(\mathbf{C}))$ induit une fonction holomorphe $u : \mathcal{H} \to \mathbf{C}$ invariante par la transformation $z \mapsto z + 1$ et admettant donc un développement de Fourier de la forme

$$u(z) = \sum_{n=n_0}^{\infty} a_n q^n \qquad (z \in \mathcal{H},\ q = e^{2\pi i z}) \tag{3.60}$$

avec $n_0 \in \mathbf{Z}$ et $a_{n_0} \neq 0$. Nous définissons alors

$$\widehat{u}(\infty) = a_{n_0}. \tag{3.61}$$

L'application $u \mapsto \widehat{u}(\infty)$ est un homomorphisme de groupes de $\mathcal{O}^*(Y_1(N)(\mathbf{C}))$ vers $\mathbf{C}^*$ qui, après tensorisation par $\mathbf{C}$, scinde la suite exacte (3.59). Toute unité modulaire $u \in \mathcal{O}^*(Y_1(N)(\mathbf{C}))$ induit une fonction $\log|u| : \mathcal{H} \to \mathbf{R}$, invariante sous l'action de $\Gamma_1(N)$. L'application $u \mapsto \log|u|$ est un homomorphisme de groupes et s'étend par linéarité à $\mathcal{O}^*(Y_1(N)(\mathbf{C})) \otimes \mathbf{C}$ : pour tout $u \in \mathcal{O}^*(Y_1(N)(\mathbf{C})) \otimes \mathbf{C}$, nous disposons donc d'une fonction $\log|u| : \mathcal{H} \to \mathbf{C}$, invariante sous l'action de $\Gamma_1(N)$. Enfin, nous étendons par linéarité les définitions de $\mathrm{div}\, u$ et $\mathrm{ord}_P(u)$, $P \in P_N$, au cas où $u \in \mathcal{O}^*(Y_1(N)(\mathbf{C})) \otimes \mathbf{C}$.

**Proposition 79.** *Soit $f : \frac{\mathbf{Z}}{N\mathbf{Z}} \to \mathbf{C}$ une fonction de somme nulle. Il existe une unique unité modulaire $u_f \in \mathcal{O}^*(Y_1(N)(\mathbf{C})) \otimes \mathbf{C}$ vérifiant*

$$\log|u_f| = \frac{1}{\pi} \cdot E_f^* \qquad et \qquad \widehat{u_f}(\infty) = 1 \in \mathbf{C}^* \otimes \mathbf{C}. \tag{3.62}$$

*Le diviseur de $u_f$ peut se décrire de la manière suivante : pour toute pointe $P = [u, v] \in P_N$, avec $(u, v) \in E_N$, nous avons*

$$\mathrm{ord}_P(u_f) = -\frac{1}{N \cdot (u, N)} \sum_{(a,b) \in (\frac{\mathbf{Z}}{N\mathbf{Z}})^2} \widehat{f}(au + bv) \cdot \overline{B_2}(\frac{\tilde{b}}{N}). \tag{3.63}$$

*De plus, l'application $f \mapsto u_f$ ainsi définie est $\mathbf{C}$-linéaire.*

*Démonstration.* Soit $P = [u,v] \in P_N$ une pointe et $g = \begin{pmatrix} \alpha & \beta \\ \gamma & \delta \end{pmatrix} \in \mathrm{SL}_2(\mathbf{Z})$ un représentant de $P$ via la bijection (3.56). Nous avons donc $\gamma \equiv u \pmod{N}$ et $\delta \equiv v \pmod{N}$. Soit $z \in \mathcal{H}$. D'après (3.15), nous avons

$$
\begin{aligned}
E_f^*(gz) &= \sum_{w \in \frac{\mathbf{Z}}{N\mathbf{Z}}} f(w)\, E_{(0,w)}^*(gz) \\
&= \sum_{w \in \frac{\mathbf{Z}}{N\mathbf{Z}}} f(w)\, E_{(uw,vw)}^*(z) \\
&= \sum_{w \in \frac{\mathbf{Z}}{N\mathbf{Z}}} f(w) \frac{1}{N^2} \sum_{(a,b) \in (\frac{\mathbf{Z}}{N\mathbf{Z}})^2} e^{-\frac{2\pi i}{N}(auw+bvw)} \cdot \zeta_{a,b}^*(z) \\
&= \frac{1}{N^2} \sum_{(a,b) \in (\frac{\mathbf{Z}}{N\mathbf{Z}})^2} \widehat{f}(au+bv) \cdot \zeta_{a,b}^*(z).
\end{aligned}
$$

Puisque $\widehat{f}(0) = 0$, nous pouvons omettre le terme $(a,b) = (0,0)$ dans la somme précédente. D'après les développements de Fourier (3.17) et (3.19), nous voyons que $E_f^*(gz)$ admet un développement de Fourier de la forme

$$E_f^*(gz) = K_{f,P} \cdot y + \alpha_0 + \sum_{r=1}^{\infty} \alpha_r q^{\frac{r}{N}} + \beta_r \overline{q}^{\frac{r}{N}}, \tag{3.64}$$

avec $K_{f,P}, \alpha_r, \beta_r \in \mathbf{C}$. La constante $K_{f,P}$ est donnée par

$$K_{f,P} = \frac{2\pi^2}{N^2} \sum_{(a,b) \in (\frac{\mathbf{Z}}{N\mathbf{Z}})^2} \widehat{f}(au+bv) \cdot \overline{B_2}(\frac{\tilde{b}}{N}) \tag{3.65}$$

(cette expression ne dépend que de la classe de $(u,v)$ dans $P_N$). Rappelons qu'un paramètre local en la pointe $P \in P_N$ est donné par

$$q_P = q^{\frac{(u,N)}{N}} = e^{\frac{2\pi i (u,N)}{N} z}. \tag{3.66}$$

Soit $G_{X_1(N)(\mathbf{C})}$ la fonction de Green définie dans le premier chapitre. Notons $\pi : \mathcal{H} \to Y_1(N)(\mathbf{C})$ la surjection naturelle. D'après (1.10) et (3.66), nous avons l'estimation

$$
\begin{aligned}
G_{X_1(N)(\mathbf{C})}(P, \pi(gz)) &= \log|q_P| + O_{y \to \infty}(1) \\
&= -\frac{2\pi(u,N)}{N} \cdot y + O_{y \to \infty}(1) \qquad (z \in \mathcal{H}).
\end{aligned} \tag{3.67}
$$

Définissons une fonction $\phi$ sur $Y_1(N)(\mathbf{C})$ par

$$\phi = E_f^* + \frac{N}{2\pi} \sum_{P \in P_N} \frac{K_{f,P}}{(u,N)} \cdot G_{X_1(N)(\mathbf{C})}(P, \cdot). \tag{3.68}$$

D'après le lemme 78, la fonction $\phi$ est de classe $\mathcal{C}^{\infty}$. D'après (3.64) et (3.67), la fonction $\phi$ s'étend en une fonction de classe $\mathcal{C}^{\infty}$ sur $X_1(N)(\mathbf{C})$. Nous avons sur $Y_1(N)(\mathbf{C})$ (et donc sur $X_1(N)(\mathbf{C})$)

$$\partial\overline{\partial}\phi = \partial\overline{\partial}E_f^* + \frac{N}{2\pi}\sum_{P\in P_N}\frac{K_{f,P}}{(u,N)}\cdot \pi i\,\mathrm{vol}_{X_1(N)(\mathbf{C})}$$

$$= \frac{N}{2\pi}\sum_{P\in P_N}\frac{K_{f,P}}{(u,N)}\cdot \pi i\,\mathrm{vol}_{X_1(N)(\mathbf{C})}.$$

D'après la formule de Stokes $\int_{X_1(N)(\mathbf{C})} \partial\overline{\partial}\phi = \int_{X_1(N)(\mathbf{C})} d(\overline{\partial}\phi) = 0$. D'autre part $\int \mathrm{vol}_{X_1(N)(\mathbf{C})} = 1$, d'où nous déduisons

$$\sum_{P\in P_N}\frac{K_{f,P}}{(u,N)} = 0. \tag{3.69}$$

Il en résulte $\partial\overline{\partial}\phi = 0$, c'est-à-dire que $\phi$ est constante sur $X_1(N)(\mathbf{C})$. D'après (3.69) et la suite exacte scindée (3.59), il existe une unique unité modulaire $u_f \in \mathcal{O}^*(Y_1(N)(\mathbf{C}))\otimes\mathbf{C}$ telle que

$$\mathrm{div}\,u_f = -\frac{N}{2\pi^2}\sum_{P\in P_N}\frac{K_{f,P}}{(u,N)}\cdot[P] \qquad \text{et} \qquad \widehat{u_f}(\infty) = 1 \in \mathbf{C}^*\otimes\mathbf{C}. \tag{3.70}$$

D'après (1.14), il existe une constante $C \in \mathbf{C}$ telle que

$$\log|u_f| = C - \frac{N}{2\pi^2}\sum_{P\in P_N}\frac{K_{f,P}}{(u,N)}\cdot G_{X_1(N)(\mathbf{C})}(P,\cdot). \tag{3.71}$$

Nous déduisons de (3.68) et (3.71) l'existence d'une constante $C' \in \mathbf{C}$ telle que

$$\log|u_f| = C' + \frac{E_f^*}{\pi}.$$

Pour déterminer la constante $C'$, nous considérons les développements de Fourier des deux membres de l'égalité précédente. D'après la définition (3.61) de $\widehat{u_f}(\infty)$, le terme constant du développement de Fourier de $\log|u_f|$ vaut $\log|\widehat{u_f}(\infty)|$, c'est-à-dire $0$. Or, le terme constant du développement de Fourier (3.27) de $E_f^*$ est nul ; nous avons donc $C' = 0$. L'identité (3.63) résulte de la définition de $u_f$. D'après cette même identité, l'application $f \mapsto \mathrm{div}\,u_f$ est $\mathbf{C}$-linéaire. Or, $u_f$ n'est autre que l'image de $\mathrm{div}\,u_f$ par l'application linéaire

$$\mathrm{Div}^0(P_N)\otimes\mathbf{C} \to \mathcal{O}^*(Y_1(N)(\mathbf{C}))\otimes\mathbf{C}$$

scindant la suite exacte (3.59). Donc $f \mapsto u_f$ est $\mathbf{C}$-linéaire. $\square$

Nous appliquons maintenant la proposition 79 dans le cas où $f$ est un caractère de Dirichlet, ou sa transformée de Fourier. Pour tout $d \in (\frac{\mathbf{Z}}{N\mathbf{Z}})^*$, notons $\langle d\rangle \in \mathrm{SL}_2(\mathbf{Z})$ un représentant de l'image réciproque de $(0,d)$ par la bijection (3.21) ; un tel élément est appelé *opérateur diamant*. Ces derniers induisent une action du groupe $(\frac{\mathbf{Z}}{N\mathbf{Z}})^*/\pm 1$ sur $P_N$, par la règle

$$\langle d\rangle\cdot[u,v] = [du,dv] \qquad \bigl(d\in(\mathbf{Z}/N\mathbf{Z})^*, (u,v)\in E_N\bigr). \tag{3.72}$$

Pour tout $v \in (\frac{\mathbf{Z}}{N\mathbf{Z}})^*/\pm 1$, nous notons $P_v = [0,v] = \langle v\rangle\cdot\infty$. Nous obtenons ainsi une bijection naturelle entre $(\frac{\mathbf{Z}}{N\mathbf{Z}})^*/\pm 1$ et l'orbite de la pointe infinie sous l'action des opérateurs diamants.

Par définition, l'*involution d'Atkin-Lehner* $W_N$ de $X_1(N)(\mathbf{C})$ est induite par l'involution $z \mapsto -\frac{1}{Nz}$ de $\mathcal{H}$.

Soit $\chi$ un *caractère de Dirichlet* modulo $N$, c'est-à-dire un homomorphisme de groupes $\chi : (\frac{\mathbf{Z}}{N\mathbf{Z}})^* \to \mathbf{C}^*$. Rappelons que $\chi$ s'étend (par 0) en une application $\frac{\mathbf{Z}}{N\mathbf{Z}} \to \mathbf{C}$, que nous notons encore $\chi$. La *somme de Gauß* $\tau(\chi)$ de $\chi$ est définie par

$$\tau(\chi) = \sum_{v \in \frac{\mathbf{Z}}{N\mathbf{Z}}} \chi(v) \cdot e^{\frac{2\pi i v}{N}}. \tag{3.73}$$

Nous disposons de séries d'Eisenstein $E_\chi^*$ et $E_{\widehat{\chi}}^*$. Remarquons les relations

$$E_\chi^*(\langle d \rangle \cdot z) = \overline{\chi}(d)\, E_\chi^*(z)\,; \tag{3.74}$$
$$E_{\widehat{\chi}}^*(\langle d \rangle \cdot z) = \chi(d)\, E_{\widehat{\chi}}^*(z) \qquad \big(d \in (\mathbf{Z}/N\mathbf{Z})^*, z \in \mathcal{H}\big). \tag{3.75}$$

Lorsque $\chi$ (resp. $\widehat{\chi}$) vérifie les hypothèses de la proposition 79, nous en déduisons une unité modulaire $u_\chi$ (resp. $u_{\widehat{\chi}}$). Dans la proposition suivante, nous obtenons explicitement le diviseur de $u_\chi$ et $u_{\widehat{\chi}}$.

**Proposition 80.** *Soit $\chi$ un caractère de Dirichlet modulo $N$, pair et non trivial* i. e. *vérifiant $\chi(-1) = 1$ et $\chi \neq 1$. Alors $u_\chi$ est définie et nous avons*

$$\operatorname{div} u_\chi = -\frac{L(\chi, 2)}{\pi^2} \sum_{v \in (\mathbf{Z}/N\mathbf{Z})^*/\pm 1} \overline{\chi}(v) \cdot P_v. \tag{3.76}$$

*Soit $\chi$ un caractère de Dirichlet pair modulo $N > 1$. Alors $u_{\widehat{\chi}}$ est définie. Notons $N_\chi$ le conducteur de $\chi$ et, pour tout entier $d \geq 1$ vérifiant $N_\chi \,|\, d \,|\, N$, notons $\chi_d$ le caractère modulo $d$ induit par $\chi$. Pour toute pointe $P = [u, v] \in P_N$, avec $(u, v) \in E_N$, posons $d = (u, N)\,;$ nous avons alors*

$$\operatorname{ord}_P(u_{\widehat{\chi}}) = \begin{cases} 0 & \textit{si } N_\chi \nmid d\,; \\ -\frac{\varphi(N)/N}{\varphi(d)/d} \chi_d(v_d) \sum_{\beta \in (\mathbf{Z}/d\mathbf{Z})^*} \overline{B_2}(\frac{\tilde{\beta}}{d}) \chi_d(\beta) & \textit{si } N_\chi \,|\, d, \end{cases} \tag{3.77}$$

*où $v_d \in (\mathbf{Z}/d\mathbf{Z})^*$ désigne l'image de $v$ modulo $d$. Lorsque de plus $\chi$ est primitif* i. e. *$N_\chi = N$, nous avons $\operatorname{div} u_{\widehat{\chi}} = \tau(\chi) \cdot \operatorname{div} u_{\overline{\chi}}$.*

*Démonstration.* Supposons d'abord $\chi$ pair et non trivial. La fonction $\chi : \frac{\mathbf{Z}}{N\mathbf{Z}} \to \mathbf{C}$ est paire, et de somme nulle puisque le caractère est non trivial. Calculons $\operatorname{ord}_P(u_\chi)$ pour $P = [u, v] \in P_N$. Nous avons

$$\operatorname{ord}_P(u_\chi) = -\frac{1}{N \cdot (u, N)} \sum_{(a,b) \in (\frac{\mathbf{Z}}{N\mathbf{Z}})^2} \widehat{\chi}(au + bv) \cdot \overline{B_2}(\frac{\tilde{b}}{N})$$

$$= -\frac{1}{N \cdot (u, N)} \sum_{(a,b) \in (\frac{\mathbf{Z}}{N\mathbf{Z}})^2} \sum_{w \in (\frac{\mathbf{Z}}{N\mathbf{Z}})^*} \chi(w) e^{-\frac{2\pi i(au+bv)w}{N}} \cdot \overline{B_2}(\frac{\tilde{b}}{N}).$$

Or nous avons $\sum_{a \in \frac{\mathbf{Z}}{N\mathbf{Z}}} e^{-\frac{2\pi i a u w}{N}} = 0$ si $u \neq 0$. Il en résulte $\operatorname{ord}_P(u_\chi) = 0$ si $u \neq 0$. Supposons maintenant $u = 0$, c'est-à-dire $P = [0, v]$ avec $v \in (\frac{\mathbf{Z}}{N\mathbf{Z}})^*$. Nous obtenons

$$\operatorname{ord}_P(u_\chi) = -\frac{1}{N} \sum_{b \in \frac{\mathbf{Z}}{N\mathbf{Z}}} \sum_{w \in (\frac{\mathbf{Z}}{N\mathbf{Z}})^*} \chi(w) e^{-\frac{2\pi i b v w}{N}} \cdot \overline{B_2}(\frac{\tilde{b}}{N})$$

$$= -\frac{1}{N} \sum_{b \in \frac{\mathbf{Z}}{N\mathbf{Z}}} \widehat{\chi}(bv) \cdot \overline{B_2}(\frac{\tilde{b}}{N})$$

$$= -\overline{\chi}(v) \frac{L(\chi, 2)}{\pi^2},$$

comme dans la démonstration de la proposition 71. Il en résulte (3.76). Passons à la deuxième partie de la proposition. Soit $\chi$ un caractère de Dirichlet pair modulo $N > 1$, de conducteur $N_\chi$. Puisque $\widehat{\widehat{\chi}} = N \cdot \chi$ et d'après l'hypothèse $N > 1$, la fonction $\widehat{\chi} : \frac{\mathbf{Z}}{N\mathbf{Z}} \to \mathbf{C}$ est paire et de somme nulle. Il vient également

$$\operatorname{ord}_P(u_{\widehat{\chi}}) = -\frac{1}{(u,N)} \sum_{(a,b) \in (\frac{\mathbf{Z}}{N\mathbf{Z}})^2} \chi(au+bv) \cdot \overline{B_2}(\frac{\tilde{b}}{N})$$

$$= -\frac{1}{(u,N)} \sum_{w \in (\frac{\mathbf{Z}}{N\mathbf{Z}})^*} \chi(w) \sum_{\substack{(a,b) \in (\frac{\mathbf{Z}}{N\mathbf{Z}})^2 \\ au+bv=w}} \overline{B_2}(\frac{\tilde{b}}{N})$$

$$= -\frac{1}{(u,N)} \sum_{w \in (\frac{\mathbf{Z}}{N\mathbf{Z}})^*} \chi(w) \sum_{\substack{(a,b) \in (\frac{\mathbf{Z}}{N\mathbf{Z}})^2 \\ aw^{-1}u+bw^{-1}v=1}} \overline{B_2}(\frac{\tilde{b}}{N}).$$

Nous allons maintenant calculer la somme

$$S(u,v) = \sum_{\substack{(a,b) \in (\frac{\mathbf{Z}}{N\mathbf{Z}})^2 \\ au+bv=1}} \overline{B_2}(\frac{\tilde{b}}{N}) \qquad ((u,v) \in E_N) \tag{3.78}$$

En posant $d = (u,N)$, nous avons $S(u,v) = S(d,v)$. Soient $a_0, b_0 \in \frac{\mathbf{Z}}{N\mathbf{Z}}$ tels que $a_0 d + b_0 v = 1$, ce qui est possible. La solution générale de l'équation $ad + bv = 1$, $a, b \in \frac{\mathbf{Z}}{N\mathbf{Z}}$ s'écrit

$$\begin{cases} a = a_0 - kv\,; \\ b = b_0 + kd, \quad k \in \frac{\mathbf{Z}}{N\mathbf{Z}}. \end{cases}$$

Par conséquent

$$S(d,v) = \sum_{k=0}^{N-1} \overline{B_2}\Big(\frac{\tilde{b}_0 + kd}{N}\Big) = d \sum_{k=0}^{N/d-1} \overline{B_2}\Big(\frac{\tilde{b}_0/d + k}{N/d}\Big).$$

La fonction $\overline{B_2}$ satisfait la relation de distribution [18, (**), p. 22]

$$\overline{B_2}(x) = r \sum_{k=0}^{r-1} \overline{B_2}\Big(\frac{x+k}{r}\Big) \qquad (x \in \mathbf{R}, r \geq 1). \tag{3.79}$$

Nous en déduisons

$$S(d,v) = \frac{d^2}{N}\overline{B_2}(\frac{\tilde{b_0}}{d}),$$

où $\tilde{b_0} \in \mathbf{Z}$ est un entier quelconque vérifiant $\tilde{b_0}v \equiv 1 \pmod{d}$. Alors

$$\begin{aligned}\operatorname{ord}_P(u_{\widehat{\chi}}) &= -\frac{1}{d}\sum_{w \in (\frac{\mathbf{Z}}{N\mathbf{Z}})^*} \chi(w)\, S(d, w^{-1}v) \\ &= -\frac{d}{N}\sum_{w \in (\frac{\mathbf{Z}}{N\mathbf{Z}})^*} \chi(w)\, \overline{B_2}(\frac{\tilde{w}\tilde{b_0}}{d}),\end{aligned}$$

où $\tilde{w}$ désigne un représentant quelconque de $w$ dans $\mathbf{Z}$. Par conséquent

$$\operatorname{ord}_P(u_{\widehat{\chi}}) = -\frac{d}{N}\sum_{\beta \in (\mathbf{Z}/d\mathbf{Z})^*}\Big(\sum_{\substack{w \in (\frac{\mathbf{Z}}{N\mathbf{Z}})^* \\ w \equiv \beta \,(d)}} \chi(w)\Big)\overline{B_2}(\frac{\tilde{\beta}\tilde{b_0}}{d}), \qquad (3.80)$$

où $\tilde{\beta}$ désigne un représentant quelconque de $\beta$ dans $\mathbf{Z}$. Nous distinguons maintenant deux cas. Si $N_\chi \nmid d$, c'est-à-dire $\chi$ ne se factorise pas par $(\mathbf{Z}/d\mathbf{Z})^*$, la somme intérieure de (3.80) est nulle, donc $\operatorname{ord}_P(u_{\widehat{\chi}}) = 0$. Si $N_\chi \,|\, d$, notons $\chi_d$ le caractère modulo $d$ induit par $\chi$, alors la somme intérieure de (3.80) vaut $\frac{\varphi(N)}{\varphi(d)}\chi_d(\beta)$, d'où

$$\begin{aligned}\operatorname{ord}_P(u_{\widehat{\chi}}) &= -\frac{\varphi(N)/N}{\varphi(d)/d}\sum_{\beta \in (\mathbf{Z}/d\mathbf{Z})^*} \chi_d(\beta)\,\overline{B_2}(\frac{\tilde{\beta}\tilde{b_0}}{d}) \\ &= -\frac{\varphi(N)/N}{\varphi(d)/d}\sum_{\beta' \in (\mathbf{Z}/d\mathbf{Z})^*} \chi_d(\beta' v_d)\,\overline{B_2}(\frac{\tilde{\beta'}}{d}),\end{aligned}$$

grâce au changement de variables $\beta = \beta' v_d$, ce qui montre (3.77). Lorsque $\chi$ est primitif (et pair), l'égalité $\widehat{\chi} = \tau(\chi) \cdot \overline{\chi}$ est classique. L'application $f \mapsto \operatorname{div} u_f$ étant $\mathbf{C}$-linéaire, nous en déduisons $\operatorname{div} u_{\widehat{\chi}} = \tau(\chi) \cdot \operatorname{div} u_{\overline{\chi}}$. $\square$

## 3.4 Version explicite du théorème de Beĭlinson

Nous utilisons les calculs de la section 3.2 pour obtenir une version explicite du théorème de Beĭlinson sur les courbes modulaires. Rappelons que l'application régulateur (1.27)

$$r_N = r_{X_1(N)(\mathbf{C})} : K_2(\mathbf{C}(X_1(N))) \otimes \mathbf{C} \to \operatorname{Hom}_\mathbf{C}(S_2(\Gamma_1(N)), \mathbf{C}) \qquad (3.81)$$

est définie par

$$\langle r_N(\{u,v\}), f\rangle = \int_{X_1(N)(\mathbf{C})} \log|u| \cdot \omega_f \wedge \overline{\partial}\log|v|, \qquad (3.82)$$

pour toutes fonctions rationnelles $u, v \in \mathbf{C}(X_1(N))^* \otimes \mathbf{C}$ et toute forme parabolique $f \in S_2(\Gamma_1(N))$, avec la notation $\omega_f = 2\pi i f(z)dz$.

Le théorème suivant est une version explicite du théorème de Beĭlinson sur les courbes modulaires [5, 62].

**Théorème 81.** *Soit $f$ une forme parabolique primitive de poids $2$ pour $\Gamma_1(N)$, de caractère $\psi$. Soient $\chi$ un caractère de Dirichlet pair modulo $N$ et $\chi'$ un caractère de Dirichlet pair modulo un diviseur $M$ de $N$. Supposons $\chi$ non trivial et $M > 1$. Notons $\alpha : Y_1(N)(\mathbf{C}) \to Y_1(M)(\mathbf{C})$ le morphisme de dégénérescence induit par l'identité sur $\mathcal{H}$, et $\chi'_N$ le caractère modulo $N$ induit par $\chi'$. Nous avons alors*

$$\langle r_N(\{u_\chi, \alpha^* u_{\widehat{\chi'}}\}), f \rangle = \begin{cases} \frac{\varphi(N)}{M\pi i} \cdot L(f,2) L(f, \chi', 1) & si\ \psi = \chi \overline{\chi'_N}, \\ 0 & sinon. \end{cases} \quad (3.83)$$

*Remarques.* 1. Les hypothèses $\chi$ non trivial et $M > 1$ assurent respectivement que les unités modulaires $u_\chi$ et $u_{\widehat{\chi'}}$ sont bien définies (proposition 80).

2. Le facteur de proportionnalité liant régulateur et produit de valeurs spéciales est totalement explicite.

*Démonstration.* Nous avons

$$\langle r_N(\{u_\chi, \alpha^* u_{\widehat{\chi'}}\}), f \rangle = \int_{X_1(N)(\mathbf{C})} \log|u_\chi| \cdot \omega_f \wedge \overline{\partial} \log|\alpha^* u_{\widehat{\chi'}}|$$

$$= \frac{1}{\pi^2} \int_{X_1(N)(\mathbf{C})} E^*_\chi \cdot \omega_f \wedge \overline{\partial} E^*_{\widehat{\chi'}}.$$

Le résultat suit alors du théorème 72. □

**Corollaire** (Théorème 4). *Soit $f$ une forme parabolique primitive de poids $2$ pour $\Gamma_1(N)$, de caractère $\psi$. Pour tout caractère de Dirichlet $\chi$ modulo $N$, pair, distinct de $\overline{\psi}$ et primitif, nous avons*

$$L(f,2) L(f,\chi,1) = \frac{N\pi i}{\varphi(N)} \tau(\chi) \langle r_N(\{u_{\psi\chi}, u_{\overline{\chi}}\}), f \rangle. \quad (3.84)$$

*Démonstration.* On utilise le théorème précédent avec $M = N$. Le caractère $\chi$ étant primitif, nous avons $\widehat{\chi} = \tau(\chi)\overline{\chi}$, d'où $u_{\widehat{\chi}} = u_{\overline{\chi}} \otimes \tau(\chi)$. Par hypothèse, le caractère $\psi\chi$ est non trivial. Il suit

$$L(f,2)L(f,\chi,1) = \frac{N\pi i}{\varphi(N)} \langle r_N(\{u_{\psi\chi}, u_{\widehat{\chi}}\}), f \rangle$$

$$= \frac{N\pi i}{\varphi(N)} \tau(\chi) \langle r_N(\{u_{\psi\chi}, u_{\overline{\chi}}\}), f \rangle.$$

□

*Remarque* 82. Sans l'hypothèse $\chi$ primitif, nous n'avons pas trouvé de formule satisfaisante pour $\langle r_N(\{u_{\psi\chi}, u_{\overline{\chi}}\}), f \rangle$. La formule recherchée fait-elle intervenir $L(f, \widehat{\overline{\chi}}, 1)$ ?

Les courbes modulaires $Y_1(N)(\mathbf{C})$ et $X_1(N)(\mathbf{C})$ admettent des modèles sur $\mathbf{Q}$ [25, II. 8.], notés $Y_1(N)$ et $X_1(N)$. Nous adoptons la convention [25, Variant 9.3.6], grâce à laquelle la pointe infinie de $X_1(N)$ est définie sur $\mathbf{Q}$ (un tel modèle est parfois noté $X_\mu(N)$). La courbe $X_1(N)$ est projective, lisse et géométriquement irréductible [25, Thm 9.3.7]. La courbe $Y_1(N)$ peut être vue comme un ouvert affine de $X_1(N)$. Notons $\mathbf{Q}(X_1(N))$ le corps des fonctions rationnelles de $X_1(N)$. Par extension des scalaires, nous avons une inclusion naturelle $\mathbf{Q}(X_1(N)) \subset \mathbf{C}(X_1(N))$.

**Lemme 83.** *La composition*

$$K_2(\mathbf{Q}(X_1(N))) \to K_2(\mathbf{C}(X_1(N))) \xrightarrow{r_N} \mathbf{T} \otimes \mathbf{C} \qquad (3.85)$$

*est à valeurs dans* $\mathbf{T} \otimes \mathbf{R}(1) = \mathbf{T} \otimes 2\pi i \mathbf{R}$.

*Démonstration.* Soient $f, g \in \mathbf{Q}(X_1(N))^*$ des fonctions rationnelles. L'image du symbole $\{f, g\} \in K_2(\mathbf{Q}(X_1(N)))$ par la composition (3.85) est donnée par

$$r_N(\{f,g\}) = \int_{X_1(N)(\mathbf{C})} \log|f| \cdot \Omega \wedge \overline{\partial} \log|g| \in \mathbf{T} \otimes \mathbf{C},$$

où $\Omega$ est la forme modulaire universelle (3.34). Notons $c : X_1(N)(\mathbf{C}) \to X_1(N)(\mathbf{C})$ la conjugaison complexe, induite par l'involution $z \mapsto -\overline{z}$ sur $\mathcal{H}$. Notons $\overline{\cdot}$ la conjugaison complexe sur les coefficients de $\mathbf{T} \otimes \mathbf{C}$. Nous avons

$$\begin{aligned}\overline{r_N(\{f,g\})} &= \int_{X_1(N)(\mathbf{C})} \log|f| \cdot \overline{\Omega} \wedge \partial \log|g| \\ &= -\int_{X_1(N)(\mathbf{C})} \log|c^* f| \cdot c^* \overline{\Omega} \wedge \overline{\partial} \log|c^* g| \\ &= -\int_{X_1(N)(\mathbf{C})} \log|f| \cdot c^* \overline{\Omega} \wedge \overline{\partial} \log|g|,\end{aligned}$$

puisque $c^* f = \overline{f}$ et $c^* g = \overline{g}$. D'autre part, l'expression (3.34) de $\Omega$ entraîne

$$c^* \overline{\Omega} = -2\pi i \sum_{n=1}^{\infty} T_n \cdot \overline{e^{-2\pi i n \overline{z}} d(-\overline{z})} = \Omega.$$

Par conséquent $\overline{r_N(\{f,g\})} = -r_N(\{f,g\})$, c'est-à-dire $r_N(\{f,g\}) \in \mathbf{T} \otimes \mathbf{R}(1)$. □

Nous nous intéressons au groupe de $K$-théorie de Quillen $K_2(X_1(N))$ associé à la courbe $X_1(N)$. Il est possible de donner une description de ce groupe (au moins après tensorisation par $\mathbf{Q}$) en termes de symboles de Milnor. La localisation en $K$-théorie algébrique permet d'écrire une suite exacte

$$0 \to K_2(X_1(N)) \otimes \mathbf{Q} \xrightarrow{\eta^*} K_2(\mathbf{Q}(X_1(N))) \otimes \mathbf{Q} \xrightarrow{\partial} \bigoplus_{P \in X_1(N)(\overline{\mathbf{Q}})} \overline{\mathbf{Q}}^* \otimes \mathbf{Q}, \qquad (3.86)$$

où l'application $\eta^*$ est la restriction au point générique de $X_1(N)$, et l'application $\partial$ est définie comme suit. Nous avons $\partial = (\partial_P \otimes \mathrm{id})_{P \in X_1(N)(\overline{\mathbf{Q}})}$ et, pour tout $P \in X_1(N)(\overline{\mathbf{Q}})$, l'application $\partial_P$, appelée *symbole modéré en* $P$, est définie par

$$\begin{aligned}\partial_P : K_2(\mathbf{Q}(X_1(N))) &\to \overline{\mathbf{Q}}^* \\ \{f, g\} &\mapsto (-1)^{\mathrm{ord}_P(f) \mathrm{ord}_P(g)} \left(\frac{f^{\mathrm{ord}_P(g)}}{g^{\mathrm{ord}_P(f)}}\right)(P).\end{aligned} \qquad (3.87)$$

En composant l'application $\eta^*$ de (3.86) et l'application (3.85), nous obtenons le régulateur de Beĭlinson, noté encore $r_N$

$$r_N : K_2(X_1(N)) \otimes \mathbf{Q} \to \mathbf{T} \otimes \mathbf{R}(1). \qquad (3.88)$$

Soit $X_1(N)_{\mathbf{Z}}$ un modèle propre et régulier de $X_1(N)$ sur $\mathbf{Z}$ (un tel modèle existe d'après la résolution des singularités [1, 44, 3]). Définissons un sous-groupe $K_2(X_1(N))_{\mathbf{Z}}$ de $K_2(X_1(N))$ par

$$K_2(X_1(N))_{\mathbf{Z}} = \mathrm{Im}\big(K_2(X_1(N)_{\mathbf{Z}}) \to K_2(X_1(N))\big). \tag{3.89}$$

D'après [65, Remark p. 13], ce sous-groupe ne dépend pas du choix du modèle (propre et régulier) $X_1(N)_{\mathbf{Z}}$. L'inclusion $K_2(X_1(N))_{\mathbf{Z}} \subset K_2(X_1(N))$ identifie alors $K_2(X_1(N))_{\mathbf{Z}} \otimes \mathbf{Q}$ à un sous-espace vectoriel de $K_2(X_1(N)) \otimes \mathbf{Q}$. Le régulateur de Beĭlinson s'écrit finalement

$$r_N : K_2(X_1(N))_{\mathbf{Z}} \otimes \mathbf{Q} \to \mathbf{T} \otimes \mathbf{R}(1). \tag{3.90}$$

Notons $\{\mathcal{O}^*(Y_1(N)), \mathcal{O}^*(Y_1(N))\}$ l'image de l'application bilinéaire alternée

$$\mathcal{O}^*(Y_1(N))^{\otimes 2} \to K_2(\mathbf{Q}(X_1(N))) \otimes \mathbf{Q} \tag{3.91}$$
$$u \otimes v \mapsto \{u, v\} \otimes 1,$$

et définissons via la suite exacte (3.86)

$$K_N = \{\mathcal{O}^*(Y_1(N)), \mathcal{O}^*(Y_1(N))\} \cap \mathrm{Ker}\, \partial \subset K_2(X_1(N)) \otimes \mathbf{Q}. \tag{3.92}$$

De manière informelle, $K_N$ est formé des éléments de $K_2(X_1(N)) \otimes \mathbf{Q}$ que l'on peut construire à partir des unités modulaires de la courbe $X_1(N)$. Notons $V_N = \mathbf{T} \otimes \mathbf{R}(1)$ l'espace d'arrivée de l'application régulateur (3.88). Schappacher et Scholl ont soulevé le problème suivant [62, 1.1.3].

*Problème.* Le groupe $r_N(K_N)$ engendre-t-il l'espace vectoriel réel $V_N$ ?

*Remarques* 84.    1. Ce problème admet un analogue naturel pour tout sous-groupe de congruence $\Gamma \subset \mathrm{SL}_2(\mathbf{Z})$ tel que la courbe modulaire associée à $\Gamma$ soit définie sur $\mathbf{Q}$.

2. Schappacher et Scholl ont démontré [62, 1.1.2 (iii)] que $K_N \subset K_2(X_1(N))_{\mathbf{Z}} \otimes \mathbf{Q}$.

Schappacher et Scholl ont remarqué que le problème ci-dessus admet une réponse négative dans le cas du sous-groupe de congruence $\Gamma_0(p)$, où $p$ est un nombre premier tel que le genre de $X_0(p)$ est non nul, c'est-à-dire $p = 11$ ou $p \geq 17$, cf. [62, 1.1.3 (i)]. Pour un sous-groupe de congruence donné, la méthode de Beĭlinson permet (en théorie) de calculer l'image par l'application régulateur des symboles associés aux unités modulaires ; le problème ci-dessus relève alors d'un simple calcul de déterminant. Ainsi, pour les sous-groupes de congruence $\Gamma_0(20)$ et $\Gamma_0(27)$, la réponse est positive [62, 1.1.3 (ii)]. Nous démontrerons le résultat général suivant, annoncé dans l'introduction.

**Théorème 5.** *Pour tout nombre premier $p$, l'espace vectoriel réel $V_p$ est engendré par $r_p(K_p)$.*

La réponse au problème proposé ci-dessus est donc différente pour les sous-groupes de congruence $\Gamma = \Gamma_0(N)$ et $\Gamma = \Gamma_1(N)$. La méthode de Beĭlinson semble mieux adaptée à la courbe modulaire $X_1(N)$. Nous avons également constaté ce phénomène lors de l'application des résultats de Beĭlinson à la conjecture de Zagier explicite pour les courbes elliptiques. Par exemple, nous n'avons pas d'équivalent du théorème 8, qui concerne la courbe elliptique $X_1(11)$, pour la courbe elliptique $X_0(11)$. Il semble intéressant de rapprocher ces observations des travaux de Stevens sur les paramétrisations modulaires. Dans [72], Stevens conjecture que pour toute courbe elliptique $E$ définie sur $\mathbf{Q}$, de conducteur $N$, il existe une paramétrisation $X_1(N) \to E$ telle que la constante de Manin associée soit égale à 1. De plus, il montre que dans toute classe d'isogénie de courbes elliptiques définies sur $\mathbf{Q}$ (de conducteur $N$), il existe une unique courbe paramétrée par $X_1(N)$ de manière optimale. Cette courbe, appelée *courbe optimale*, est conjecturalement la courbe dont la hauteur de Faltings est minimale. Nos résultats semblent donc eux aussi indiquer qu'il est plus naturel de paramétrer une courbe elliptique $E$ par la courbe modulaire $X_1(N)$. L'exemple des courbes modulaires $X_1(11)$ et $X_1(13)$, que nous étudions en détail dans les sections 3.7 et 3.8, est particulièrement instructif.

Le groupe $\mathrm{Aut}(\mathbf{C}/\mathbf{Q})$ agit sur les ensembles $Y_1(N)(\mathbf{C})$ et $X_1(N)(\mathbf{C})$, et donc sur l'ensemble des pointes $P_N$. D'après [62, 3.0.2], nous pouvons décrire explicitement cette action. Pour tout $\sigma \in \mathrm{Aut}(\mathbf{C}/\mathbf{Q})$, définissons $\epsilon(\sigma) \in (\frac{\mathbf{Z}}{N\mathbf{Z}})^*$ par l'égalité $\sigma(\zeta_N) = \zeta_N^{\epsilon(\sigma)}$, avec $\zeta_N = e^{\frac{2\pi i}{N}}$. Alors $\mathrm{Aut}(\mathbf{C}/\mathbf{Q})$ agit sur $P_N$ par la règle

$$[u,v]^\sigma = [\epsilon(\sigma)^{-1}u, v] \qquad \big((u,v) \in E_N,\ \sigma \in \mathrm{Aut}(\mathbf{C}/\mathbf{Q})\big), \tag{3.93}$$

où nous avons utilisé l'identification (3.56). En particulier, la conjugaison complexe $c$ envoie $[u,v]$ sur $[-u,v]$. Remarquons qu'avec cette convention, la pointe infinie (et plus généralement les pointes $P_v$, $v \in (\mathbf{Z}/N\mathbf{Z})^*/\pm 1$) est bien définie sur $\mathbf{Q}$. Le groupe $\mathcal{O}^*(Y_1(N))$ des unités de $Y_1(N)$ est de façon naturelle un sous-groupe de $\mathcal{O}^*(Y_1(N)(\mathbf{C}))$, le groupe des unités modulaires pour $\Gamma_1(N)$.

*Notation.* Soient $u \in \mathcal{O}^*(Y_1(N)(\mathbf{C})) \otimes \mathbf{C}$ et $K, L$ deux sous-corps de $\mathbf{C}$. Nous dirons que $u$ est *définie sur $K$ et à coefficients dans $L$* lorsque

$$u \in \mathcal{O}^*(Y_1(N)_K) \otimes L \subset \mathcal{O}^*(Y_1(N)(\mathbf{C})) \otimes \mathbf{C}.$$

où $Y_1(N)_K$ désigne l'extension des scalaires à $K$.

**Lemme 85.** *Pour toute fonction $f : \frac{\mathbf{Z}}{N\mathbf{Z}} \to \mathbf{C}$ de somme nulle, l'unité modulaire $u_f$ est définie sur $\mathbf{Q}$ et à coefficients dans $\mathbf{Q}(\widehat{f})$, le corps engendré par les valeurs de $\widehat{f}$.*

*Démonstration.* Posons $L = \mathbf{Q}(\widehat{f})$. Notons $\mathrm{Div}^0_{\mathbf{Q}} P_N$ le sous-groupe de $\mathrm{Div}^0 P_N$ formé des diviseurs qui sont définis sur $\mathbf{Q}$. Nous avons un diagramme commutatif

$$\begin{array}{ccccccccc}
0 & \longrightarrow & \mathbf{Q}^* \otimes L & \longrightarrow & \mathcal{O}^*(Y_1(N)) \otimes L & \longrightarrow & \mathrm{Div}^0_{\mathbf{Q}} P_N \otimes L & \longrightarrow & 0 \\
& & \downarrow & & \downarrow & & \downarrow & & \\
0 & \longrightarrow & \mathbf{C}^* \otimes L & \longrightarrow & \mathcal{O}^*(Y_1(N)(\mathbf{C})) \otimes L & \longrightarrow & \mathrm{Div}^0 P_N \otimes L & \longrightarrow & 0,
\end{array}$$

où les flèches verticales sont injectives et les lignes sont exactes (l'exactitude à droite de la ligne du haut résulte du théorème Hilbert 90). L'application $u \in \mathcal{O}^*(Y_1(N)(\mathbf{C})) \mapsto \widehat{u}(\infty) \in \mathbf{C}^*$ scinde la suite exacte du bas du diagramme. Lorsque $u \in \mathcal{O}^*(Y_1(N))$, nous avons $\widehat{u}(\infty) \in \mathbf{Q}^*$, puisque le développement de Fourier (3.60) de $u$ est à coefficients rationnels. L'application $u \in \mathcal{O}^*(Y_1(N)) \mapsto \widehat{u}(\infty) \in \mathbf{Q}^*$ scinde la suite exacte du haut du diagramme, de manière compatible à la scission de celle du bas.

D'après (3.63), nous avons $\mathrm{ord}_P(u_f) \in L$ pour toute pointe $P \in P_N$. Soit $\sigma \in \mathrm{Aut}(\mathbf{C}/\mathbf{Q})$. Le changement de variables $a = \epsilon(\sigma)a'$ dans la formule (3.63) montre que

$$\mathrm{ord}_{P^\sigma}(u_f) = \mathrm{ord}_P(u_f) \qquad (P \in P_N,\ \sigma \in \mathrm{Aut}(\mathbf{C}/\mathbf{Q})).$$

En conséquence $D = \mathrm{div}\, u_f \in \mathrm{Div}^0_{\mathbf{Q}} P_N \otimes L$, et $u_f$ n'est autre que l'image de $D$ par l'une des deux compositions du diagramme commutatif

$$\begin{array}{ccc}
\mathrm{Div}^0_{\mathbf{Q}} P_N \otimes L & \longrightarrow & \mathcal{O}^*(Y_1(N)) \otimes L \\
\downarrow & & \downarrow \\
\mathrm{Div}^0 P_N \otimes L & \longrightarrow & \mathcal{O}^*(Y_1(N)(\mathbf{C})) \otimes L.
\end{array}$$

Par suite, nous avons $u_f \in \mathcal{O}^*(Y_1(N)) \otimes L$. □

Nous allons maintenant chercher à former des éléments dans $K_2(X_1(N)) \otimes \mathbf{Q}$ à partir des symboles de Milnor $\{u_f, u_g\}$, où $u_f, u_g$ sont les unités modulaires définies précédemment (proposition 79). Nous avons la proposition suivante.

**Proposition 86.** *Soient $f, g : \frac{\mathbf{Z}}{N\mathbf{Z}} \to \mathbf{C}$ deux fonctions de somme nulle et à support dans $(\frac{\mathbf{Z}}{N\mathbf{Z}})^*$. Notons $L$ le sous-corps de $\mathbf{C}$ engendré par les valeurs de $\widehat{f}$ et $\widehat{g}$. Alors l'élément $\{u_f, u_g\}$ appartient à $K_2(X_1(N)) \otimes L$.*

*En particulier, pour toutes unités modulaires $u, v \in \mathcal{O}^*(Y_1(N)) \otimes \mathbf{Q}$ à support dans les pointes $P_w$, $w \in (\frac{\mathbf{Z}}{N\mathbf{Z}})^*$ et vérifiant $\widehat{u}(\infty) = \widehat{v}(\infty) = 1$, l'élément $\{u, v\}$ appartient à $K_2(X_1(N)) \otimes \mathbf{Q}$.*

*Démonstration.* Considérons la suite exacte (3.86) tensorisée par $L$. Nous avons

$$u_f, u_g \in \mathcal{O}^*(Y_1(N)) \otimes L \subset \mathbf{Q}(X_1(N))^* \otimes L$$

et le symbole $\{u_f, u_g\}$ définit a priori un élément de $K_2(\mathbf{Q}(X_1(N))) \otimes L$. Il s'agit donc de montrer que son image par $\partial$ est triviale. Par hypothèse, $f$ et $g$ sont combinaisons linéaires à coefficients dans $L$ de caractères de Dirichlet non triviaux modulo $N$. Puisque $(f, g) \mapsto \{u_f, u_g\}$ est bilinéaire, il suffit de montrer le résultat pour $f = \chi$ et $g = \chi'$, avec $\chi$ et $\chi'$ caractères non triviaux modulo $N$. Puisque $u_\chi = 1$ pour un caractère $\chi$ impair, nous pouvons supposer que $\chi$ et $\chi'$ sont pairs. Soit $d \in (\frac{\mathbf{Z}}{N\mathbf{Z}})^*$. La proposition 3.76 montre que

$$\mathrm{div}\langle d \rangle^* u_\chi = \langle d \rangle^* \mathrm{div}\, u_\chi = \overline{\chi}(d) \cdot \mathrm{div}\, u_\chi = \mathrm{div}\bigl(u_\chi \otimes \overline{\chi}(d)\bigr).$$

D'après la suite exacte

$$0 \longrightarrow \mathbf{Q}^* \otimes L \longrightarrow \mathcal{O}^*(Y_1(N)) \otimes L \longrightarrow \mathrm{Div}^0_{\mathbf{Q}} P_N \otimes L \longrightarrow 0,$$

il existe une unique constante $C_d \in \mathbf{Q}^* \otimes L$ telle que

$$\langle d \rangle^* u_\chi = C_d \cdot \bigl(u_\chi \otimes \overline{\chi}(d)\bigr),$$

où nous utilisons la notation multiplicative pour le groupe $\mathcal{O}^*(Y_1(N)) \otimes L$. L'application $d \mapsto C_d$ est alors un homomorphisme de groupes, qui part d'un groupe fini et arrive dans un $\mathbf{Q}$-espace vectoriel. Il en résulte $C_d = 1$ et $\langle d \rangle^* u_\chi = u_\chi \otimes \overline{\chi}(d)$ pour tout $d \in (\frac{\mathbf{Z}}{N\mathbf{Z}})^*$. De même pour $u_{\chi'}$. Puisque $u_\chi$ et $u_{\chi'}$ sont à support dans les pointes $P_w$, $w \in (\frac{\mathbf{Z}}{N\mathbf{Z}})^*$, il suffit de montrer $\partial_{P_w}\{u_\chi, u_{\chi'}\} = 1$ pour tout $w$. Or

$$\partial_{P_w}\{u_\chi, u_{\chi'}\} = \partial_\infty\{\langle w \rangle^* u_\chi, \langle w \rangle^* u_{\chi'}\} = \bigl(\partial_\infty\{u_\chi, u_{\chi'}\}\bigr) \otimes \overline{\chi\chi'}(w).$$

Il suffit donc de montrer que $\partial_\infty\{u_\chi, u_{\chi'}\} = 1$. Pour des unités modulaires $u, v \in \mathcal{O}^*(Y_1(N)(\mathbf{C}))$, nous avons

$$\partial_\infty\{u, v\} = (-1)^{\mathrm{ord}_\infty(u)\,\mathrm{ord}_\infty(v)} \frac{\widehat{u}(\infty)^{\mathrm{ord}_\infty(v)}}{\widehat{v}(\infty)^{\mathrm{ord}_\infty(u)}}.$$

D'après la définition (3.62) de $u_f$, nous avons

$$\widehat{u_\chi}(\infty) = \widehat{u_{\chi'}}(\infty) = 1 \otimes 0 \in \mathbf{Q}^* \otimes L.$$

Il en résulte $\partial_\infty \{u_\chi, u_{\chi'}\} = 1$, comme annoncé. Pour finir, soient $u, v \in \mathcal{O}^*(Y_1(N)) \otimes \mathbf{Q}$ des unités modulaires vérifiant les hypothèses de la seconde partie du théorème. En utilisant (3.76) et le fait que $L(\chi, 2) \neq 0$, il n'est pas difficile de voir que les diviseurs de $u$ et $v$ sont combinaisons linéaires à coefficients complexes des diviseurs div $u_\chi$, où $\chi$ parcourt les caractères pairs non triviaux modulo $N$. Il existe donc des constantes $C, C' \in \mathbf{C}^* \otimes \mathbf{C}$ et $c_\chi, c'_\chi \in \mathbf{C}$ telles que

$$u = C \cdot \prod_\chi u_\chi \otimes c_\chi \qquad \text{et} \qquad v = C' \cdot \prod_\chi u_\chi \otimes c'_\chi.$$

Puisque $\widehat{u}(\infty) = \widehat{v}(\infty) = \widehat{u_\chi}(\infty) = 1$, il vient $C = C' = 1$. La première partie de la démonstration entraîne alors

$$\partial\{u, v\} = 1 \qquad \text{dans} \qquad \bigoplus_{P \in X_1(N)(\mathbf{C})} \mathbf{C}^* \otimes \mathbf{C}.$$

L'application naturelle $\mathbf{C}^* \otimes \mathbf{Q} \to \mathbf{C}^* \otimes \mathbf{C}$ étant injective, le symbole $\{u, v\}$ appartient au noyau de $\partial$ dans la suite exacte (3.86), et définit donc un élément de $K_2(X_1(N)) \otimes \mathbf{Q}$. □

*Remarque* 87. Lorsque $f$ ou $g$ n'est pas à support dans $(\frac{\mathbf{Z}}{N\mathbf{Z}})^*$, tous les symboles modérés de $\{u_f, u_g\}$ ne semblent pas nécessairement triviaux. Soit $P \in P_N$ et $\gamma \in \mathrm{SL}_2(\mathbf{Z})$ un représentant de $P$ via la bijection (3.56). Notons $\alpha_P(f)$ (resp. $\alpha_P(g)$) le terme constant du développement de Fourier de $E_f^*(\gamma z)$ (resp. $E_g^*(\gamma z)$). Nous avons alors

$$\log|\partial_P\{u_f, u_g\}| = \frac{1}{\pi}\bigl(\mathrm{ord}_P(u_g)\alpha_P(f) - \mathrm{ord}_P(u_f)\alpha_P(g)\bigr) \in \mathbf{R} \cdot L \subset \mathbf{C}.$$

Grâce au formulaire sur les séries d'Eisenstein et à la formule (3.63), il est possible d'obtenir une formule complètement explicite pour $\log|\partial_P\{u_f, u_g\}|$. Cette dernière quantité ne semble pas toujours nulle. Une manière de remédier à cela serait d'ajouter à $\{u_f, u_g\}$ des "symboles constants" (voir la remarque à la fin de cette section).

Nous avons enfin besoin des *symboles de Manin*, dont nous rappelons ici brièvement la définition [45]. Pour tous points $\alpha, \beta \in \mathbf{P}^1(\mathbf{Q})$, notons $\{\alpha, \beta\}$ la géodésique du demi-plan de Poincaré reliant $\alpha$ à $\beta$. Pour tout $x \in E_N$ (c'est-à-dire $x = (u, v)$ avec $u, v \in \frac{\mathbf{Z}}{N\mathbf{Z}}$ engendrant $\frac{\mathbf{Z}}{N\mathbf{Z}}$ comme groupe additif), choisissons une matrice $g_x = \begin{pmatrix} a & b \\ c & d \end{pmatrix} \in \mathrm{SL}_2(\mathbf{Z})$ vérifiant $(c, d) \in x$. On définit un cycle relatif $\xi(x) \in H_1(X_1(N)(\mathbf{C}), P_N, \mathbf{C})$, où $P_N$ est l'ensemble des pointes de $X_1(N)(\mathbf{C})$, par

$$\xi(x) = \{g_x 0, g_x \infty\} \qquad (x \in E_N). \tag{3.94}$$

On vérifie que $\xi(x)$ ne dépend pas du choix de la matrice $g_x$.

*Démonstration du théorème 5.* Soit $p$ un nombre premier. Pour tout $\alpha \in (\mathbf{Z}/p\mathbf{Z})^*$, notons $u_\alpha \in \mathcal{O}^*(Y_1(p)) \otimes \mathbf{Q}$ l'unique unité modulaire vérifiant

$$\mathrm{div}\, u_\alpha = P_\alpha - P_1 \qquad \widehat{u_\alpha}(\infty) = 1, \tag{3.95}$$

ce qui est possible d'après le théorème de Manin-Drinfel'd [28, 29]. D'après la proposition 86, nous avons

$$\{u_\alpha, u_\beta\} \in K_2(X_1(p)) \otimes \mathbf{Q} \qquad (\alpha, \beta \in (\mathbf{Z}/p\mathbf{Z})^*).$$

Nous allons démontrer que les $r_p(\{u_\alpha, u_\beta\})$ engendrent l'espace d'arrivée $V_p$. Notons $V \subset V_p = \mathbf{T} \otimes \mathbf{R}(1)$ l'espace engendré par les $r_p(\{u_\alpha, u_\beta\})$. Soient $\chi, \chi'$ deux caractères pairs non triviaux modulo $p$. Par linéarité et d'après (3.76), nous avons la formule suivante dans $\mathbf{T} \otimes \mathbf{C}$

$$r_p(\{u_\chi, u_{\chi'}\}) = \frac{L(\chi, 2)L(\chi', 2)}{\pi^4} \sum_{\alpha, \beta \in (\mathbf{Z}/p\mathbf{Z})^*/\pm 1} \overline{\chi}(\alpha)\overline{\chi'}(\beta) r_p(\{u_\alpha, u_\beta\}) \in V \otimes \mathbf{C}.$$

Il suffit donc de montrer que les éléments $r_p(\{u_\chi, u_{\chi'}\})$ engendrent $\mathbf{T} \otimes \mathbf{C}$ comme espace vectoriel complexe. Remarquons que tous les caractères non triviaux modulo $p$ sont primitifs, et que les formes primitives de poids 2 pour $\Gamma_1(p)$ constituent une base de $S_2(\Gamma_1(p))$. D'après le théorème 81 et le corollaire 3.4, nous avons

$$r_p(\{u_\chi, u_{\chi'}\}) = \frac{p-1}{p\pi i \cdot \tau(\overline{\chi'})} \cdot L(\mathbf{T}^{\chi\chi'}, 2) L(\mathbf{T}^{\chi\chi'}, \overline{\chi'}, 1) \in \mathbf{T}^{\chi\chi'}.$$

Fixons maintenant un caractère pair $\psi$ modulo $p$. Remarquons que $L(\mathbf{T}^\psi, 2)$ est inversible dans $\mathbf{T}^\psi$ : pour chaque forme primitive $f$, on a $L(f, 2) \neq 0$. D'après ce qui précède, il suffit de montrer que les $L(\mathbf{T}^\psi, \chi, 1)$, avec $\chi$ caractère pair modulo $p$, et $\chi \neq 1, \overline{\psi}$, engendrent $\mathbf{T}^\psi$ comme espace vectoriel complexe. Nous noterons $H_1^+(X_1(p)(\mathbf{C}), \psi)$ la composante $\psi$-isotypique du groupe d'homologie $H_1^+(X_1(p)(\mathbf{C}), \mathbf{C})$, où $H_1^+$ indique le sous-espace des invariants pour la conjugaison complexe agissant sur $X_1(p)(\mathbf{C})$. Nous disposons d'isomorphismes

$$\mathbf{T}^\psi \cong \mathrm{Hom}_\mathbf{C}(S_2(\Gamma_1(p), \psi), \mathbf{C}) \cong H_1^+(X_1(p)(\mathbf{C}), \psi), \tag{3.96}$$

le premier isomorphisme étant induit par (3.32) et le second donné par l'intégration. Le groupe d'homologie (3.96) admet une présentation en termes des symboles de Manin [45]. Nous allons exprimer l'image de $L(\mathbf{T}^\psi, \chi, 1)$ par l'isomorphisme (3.96) comme combinaison linéaire explicite de tels symboles. C'est un calcul classique [45, Thm 3.9 et 4.2.b)]. Nous trouvons que $L(\mathbf{T}, \chi, 1)$ est donné par l'intégration le long du cyle (relatif)

$$\theta_\chi = -\frac{\tau(\chi)}{p} \sum_{v \in (\mathbf{Z}/p\mathbf{Z})^*} \overline{\chi}(v) \{\frac{v}{p}, \infty\} \tag{3.97}$$

Rappelons que l'involution d'Atkin-Lehner $W_p$ de $X_1(p)(\mathbf{C})$ est induite par $z \mapsto -\frac{1}{pz}$ sur $\mathcal{H}$. Nous avons $W_p\{\frac{v}{p}, \infty\} = \xi(1, v)$, d'où

$$\theta_\chi = -\frac{\tau(\chi)}{p} W_p \xi(1, \overline{\chi}) \quad \text{avec} \quad \xi(1, \overline{\chi}) := \sum_{v \in (\mathbf{Z}/p\mathbf{Z})^*} \overline{\chi}(v)\, \xi(1, v). \tag{3.98}$$

Pour tout cycle $c$, nous noterons $c^\psi$ sa projection sur la composante $\psi$-isotypique. Il suit de (3.98) que l'élément $L(\mathbf{T}^\psi, \chi, 1)$ est donné par l'intégration le long du cycle (a priori relatif)

$$\theta_\chi^\psi = -\frac{\tau(\chi)}{p} W_p \bigl(\xi(1, \overline{\chi})^{\overline{\psi}}\bigr). \tag{3.99}$$

Nous allons montrer que les cycles $\theta_\chi^\psi$, $\chi \neq 1, \overline{\psi}$ sont fermés et engendrent $H_1^+(X_1(p)(\mathbf{C}), \psi)$ comme espace vectoriel complexe. Notons $1_p$ le caractère trivial modulo $p$. D'après (3.99), il suffit de montrer que les cycles $\xi(1, \chi)^\psi$, $\chi \neq 1_p, \overline{\psi}$ engendrent $H_1^+(X_1(p)(\mathbf{C}), \psi)$. L'ensemble des pointes de $X_1(p)(\mathbf{C})$ est donné par

$$P_v = [0, v] \qquad Q_v = [v, 0] \qquad (v \in (\mathbf{Z}/p\mathbf{Z})^*/\pm 1). \tag{3.100}$$

Considérons la projection naturelle

$$\pi_\psi : H_1(X_1(p)(\mathbf{C}), \text{ptes}, \mathbf{C}) \to H_1(X_1(p)(\mathbf{C}), \text{ptes}, \psi)$$

$$\xi(x) \mapsto \xi(x)^\psi = \frac{1}{p-1} \sum_{\lambda \in (\mathbf{Z}/p\mathbf{Z})^*} \overline{\psi}(\lambda)\xi(\lambda x).$$

Manin [45, 1.6] a démontré que les $\xi(x)$, $x \in E_p$ engendrent l'espace de départ de $\pi_\psi$. En conséquence, l'espace d'arrivée de $\pi_\psi$ est engendré par les cycles

$$\xi(x)^\psi \qquad \big(x = (0,1),\ x = (1,v),\ v \in \mathbf{Z}/p\mathbf{Z}\big).$$

Nous avons $\xi(1,0)^\psi = -\xi(0,1)^\psi$, et le bord des cycles $\xi(x)^\psi$ est donné par

$$\partial \xi(0,1)^\psi = \frac{1}{p-1} \sum_{\lambda \in (\mathbf{Z}/p\mathbf{Z})^*} \overline{\psi}(\lambda)\big(P_\lambda - Q_\lambda\big)$$

$$\partial \xi(1,v)^\psi = \frac{1-\psi(v)}{p-1} \sum_{\lambda \in (\mathbf{Z}/p\mathbf{Z})^*} \overline{\psi}(\lambda)\, Q_\lambda \qquad (v \in (\mathbf{Z}/p\mathbf{Z})^*). \tag{3.101}$$

L'espace $H_1(X_1(p)(\mathbf{C}), \psi)$ est donc contenu dans le sous-espace engendré par les cycles $\xi(1,v)^\psi$, $v \in (\mathbf{Z}/p\mathbf{Z})^*$. Par suite, l'espace $H_1^+(X_1(p)(\mathbf{C}), \psi)$ est contenu dans le sous-espace engendré par les $\xi(1,v)^\psi + \xi(1,-v)^\psi$, $v \in (\mathbf{Z}/p\mathbf{Z})^*$. Or, nous avons la formule

$$\xi(1,\chi)^\psi = \sum_{v \in (\mathbf{Z}/p\mathbf{Z})^*} \chi(v)\, \xi(1,v)^\psi.$$

Par transformée de Fourier inverse, l'espace $H_1^+(X_1(p)(\mathbf{C}), \psi)$ est contenu dans le sous-espace engendré par les $\xi(1,\chi)^\psi$, où $\chi$ parcourt les caractères pairs modulo $p$. D'autre part, nous disposons de la première relation de Manin $\xi(u,v) = -\xi(v,-u)$, ce qui permet d'écrire

$$\xi(1,1_p)^\psi = \frac{1}{p-1} \sum_{v \in (\mathbf{Z}/p\mathbf{Z})^*} \sum_{\lambda \in (\mathbf{Z}/p\mathbf{Z})^*} \overline{\psi}(\lambda)\xi(\lambda, \lambda v)$$

$$= -\frac{1}{p-1} \sum_{v \in (\mathbf{Z}/p\mathbf{Z})^*} \sum_{\lambda \in (\mathbf{Z}/p\mathbf{Z})^*} \overline{\psi}(\lambda)\xi(\lambda v, -\lambda)$$

$$= -\frac{1}{p-1} \sum_{w \in (\mathbf{Z}/p\mathbf{Z})^*} \sum_{\mu \in (\mathbf{Z}/p\mathbf{Z})^*} \overline{\psi}(\mu w)\xi(\mu, \mu w) = -\xi(1,\overline{\psi})^\psi.$$

Nous devons maintenant distinguer deux cas. Si $\psi = 1_p$, alors $\xi(1,1_p)^\psi = 0$; par conséquent $H_1^+(X_1(p)(\mathbf{C}), 1_p)$ est contenu dans l'espace engendré par les $\xi(1,\chi)^\psi$, $\chi \neq 1_p$. Les cycles $\xi(1,\chi)^\psi$ étant de bord nul d'après (3.101), nous obtenons le résultat. Supposons maintenant $\psi \neq 1_p$. Nous savons que $H_1^+(X_1(p)(\mathbf{C}), \psi)$ est contenu dans l'espace engendré par les $\xi(1,\chi)^\psi$, $\chi \neq \overline{\psi}$. Mais le calcul du bord de ces cycles donne

$$\partial \xi(1,\chi)^\psi = \begin{cases} \sum_{\lambda \in (\mathbf{Z}/p\mathbf{Z})^*} \overline{\psi}(\lambda)\, Q_\lambda & \text{si } \chi = 1_p \\ 0 & \text{si } \chi \neq 1_p, \overline{\psi}. \end{cases}$$

Par conséquent $H_1^+(X_1(p)(\mathbf{C}), \psi)$ est l'espace engendré par les cycles $\xi(1,\chi)^\psi$, avec $\chi \neq 1_p, \overline{\psi}$. $\square$

Nous présentons maintenant un argument heuristique concernant l'extension éventuelle du théorème 5 à tout niveau $N$. Le rang $\rho_1(N)$ du $\mathbf{Z}$-module libre $\mathcal{O}^*(Y_1(N))/\mathbf{Q}^*$ est donné par

$$\rho_1(N) = \left\lfloor \frac{N}{2} \right\rfloor \qquad (N \geq 1).$$

Appelons *symboles constants* les éléments de la forme $\{\lambda, u\}$ avec $\lambda \in \mathbf{Q}^*$ et $u \in \mathcal{O}^*(Y_1(N))$. En utilisant le caractère alterné de (3.91) et en étudiant le symbole modéré des symboles constants, il est possible de montrer que $K_N$ est un $\mathbf{Z}$-module libre, avec

$$\operatorname{rg} K_N \leq \frac{\rho_1(N)(\rho_1(N)-1)}{2} \underset{N \to \infty}{\sim} \frac{N^2}{8}. \tag{3.102}$$

D'autre part [25, Ex 9.1.6] montre que l'espace vectoriel réel $V_N$ est de dimension $g_1(N) \leq 1 + \frac{N^2}{24}$. Le rapport asymptotique égal à 3 entre les deux quantités précédentes constitue un indice favorable à la validité du théorème 5 pour tout entier $N$. Le même argument de dimension, appliqué dans le cas du sous-groupe de congruence $\Gamma_0(N)$, permet à l'inverse, pour beaucoup de valeurs de $N$, de répondre négativement à la question de Schappacher et Scholl.

*Question.* La dimension de l'espace d'arrivée du régulateur, jointe à l'injectivité conjecturale de ce dernier, imposent des relations sur les symboles $\{u,v\}$, $u,v \in \mathcal{O}^*(Y_1(N))$. Est-il possible d'expliciter ces relations ?

Nous reformulons enfin le théorème 4 au niveau de la jacobienne $J_1(N)$ de $X_1(N)$, en suivant la démarche proposée dans la section 0.6 de l'introduction. Rappelons que nous pouvons considérer $R_{J_1(N)} = R_{J_1(N)(\mathbf{C})}$ comme une fonction à valeurs dans $\mathbf{T} \otimes \mathbf{C}$.

**Théorème 7.** *Soit $\psi$ un caractère de Dirichlet pair modulo $N$. Pour tout caractère de Dirichlet $\chi$ modulo $N$, pair, primitif et distinct de $\overline{\psi}$, nous avons*

$$L(\mathbf{T}^\psi, 2)\, L(\mathbf{T}^\psi, \chi, 1) = C_{\psi,\chi} \sum_{\lambda,\mu \in (\frac{\mathbf{Z}}{N\mathbf{Z}})^*/\pm 1} \overline{\psi\chi}(\lambda)\, \chi(\mu)\, R_{J_1(N)}(P_\lambda - P_\mu) \tag{3.103}$$

*où la constante $C_{\psi,\chi}$ est donnée par*

$$C_{\psi,\chi} = \frac{N\pi i}{\varphi(N)} \cdot \frac{L(\psi\chi, 2)\, \tau(\chi)\, L(\overline{\chi}, 2)}{\pi^4}. \tag{3.104}$$

*Remarque.* En particulier, le produit $L(\mathbf{T}^\psi, 2)L(\mathbf{T}^\psi, \chi, 1)$ est combinaison linéaire (explicite) de valeurs de la fonction $R_{J_1(N)}$ en des points $\mathbf{Q}$-rationnels du sous-groupe cuspidal de $J_1(N)$.

*Démonstration.* Il s'agit d'utiliser le théorème 4. Puisque nous voulons un résultat portant sur toute l'algèbre de Hecke $\mathbf{T}^\psi$, nous devons reprendre les ingrédients de la démonstration de ce théorème. D'après le théorème 73, nous avons

$$L(\mathbf{T}^\psi, 2)L(\mathbf{T}^\psi, \chi, 1) = -\frac{N}{\pi i \varphi(N)} \int_{X_1(N)(\mathbf{C})} E^*_{\psi\chi} \cdot \Omega \wedge \overline{\partial} E^*_{\overline{\chi}}.$$

Puisque $\chi$ est primitif, nous avons $E^*_{\overline{\chi}} = \tau(\chi)\, E^*_{\overline{\chi}}$. D'après la proposition 79, il vient

$$L(\mathbf{T}^\psi, 2)L(\mathbf{T}^\psi, \chi, 1) = \frac{N\pi i}{\varphi(N)} \cdot \tau(\chi) \int_{X_1(N)(\mathbf{C})} \log|u_{\psi\chi}| \cdot \Omega \wedge \overline{\partial} \log|u_{\overline{\chi}}|.$$

Grâce à (1.28) et (3.76), nous en déduisons

$$L(\mathbf{T}^\psi, 2)L(\mathbf{T}^\psi, \chi, 1) = C_{\psi,\chi} \sum_{\lambda,\mu \in (\frac{\mathbf{Z}}{N\mathbf{Z}})^*/\pm 1} \overline{\psi\chi}(\lambda)\chi(\mu) R_{X_1(N)}(P_\lambda, P_\mu).$$

Par hypothèse $\psi\chi$ et $\overline{\chi}$, vus comme diviseurs sur $\frac{\mathbf{Z}}{N\mathbf{Z}}$, sont de degré 0. En partant de l'identité donnée par le théorème 5, il suit formellement

$$\begin{aligned}
&\sum_{\lambda,\mu \in (\frac{\mathbf{Z}}{N\mathbf{Z}})^*/\pm 1} \overline{\psi\chi}(\lambda)\chi(\mu) R_{X_1(N)}(P_\lambda, P_\mu) \\
&= \sum_{\lambda,\mu \in (\frac{\mathbf{Z}}{N\mathbf{Z}})^*/\pm 1} \overline{\psi\chi}(\lambda)\chi(\mu) \big(R_{J_1(N)}(P_\lambda - P_\mu) - \Phi_{X_1(N)}(P_\lambda) + \Phi_{X_1(N)}(P_\mu)\big) \\
&= \sum_{\lambda,\mu \in (\frac{\mathbf{Z}}{N\mathbf{Z}})^*/\pm 1} \overline{\psi\chi}(\lambda)\chi(\mu) R_{J_1(N)}(P_\lambda - P_\mu)
\end{aligned}$$

ce qui achève de montrer (3.103). □

## 3.5 Une formule en termes de cycles

Le but de cette section est de montrer le théorème 3 de l'introduction. Nous commençons par rappeler des résultats bien connus sur la forme différentielle $\eta(f,g)$ associée à deux fonctions méromorphes $f$ et $g$ d'une surface de Riemann compacte $X$ [56].

Soit donc $X$ une surface de Riemann compacte, connexe, non vide, et $f, g \in \mathbf{C}(X)^*$ deux fonctions méromorphes non nulles. Posons [56]

$$\begin{aligned}
\eta(f,g) &= \log|f| d\arg g - \log|g| d\arg f \\
&= -i\log|f|(\partial - \overline{\partial})\log|g| + i\log|g|(\partial - \overline{\partial})\log|f|.
\end{aligned} \tag{3.105}$$

Rappelons les propriétés de cette forme différentielle. C'est une 1-forme réelle, définie et de classe $\mathcal{C}^\infty$ sur le complémentaire des supports de $f$ et $g$ dans $X$. Elle est fermée, et se comporte de manière bilinéaire et alternée en $f, g$. De plus, elle est "exacte sur les relations de Steinberg" : nous avons

$$\eta(f, 1-f) = d(D \circ f) \qquad (f \in \mathbf{C}(X), f \neq 0, 1), \tag{3.106}$$

où $D$ est la fonction de Bloch-Wigner. Étudions maintenant le comportement de la forme différentielle $\eta(f,g)$ au voisinage d'un point $P \in X$ appartenant à la réunion $S$ des supports des fonctions $f$ et $g$. Soit $u$ une coordonnée locale en $P$. Pour toute fonction $h \in \mathbf{C}(X)^*$, nous pouvons écrire

$$h(u) \sim \widehat{h}(P) u^{\mathrm{ord}_P(h)} \qquad (u \to 0, u \neq 0). \tag{3.107}$$

Il suit

$$\log|h(u)| = \mathrm{ord}_P(h) \log|u| + \log|\widehat{h}(P)| + o_{u \to 0}(1)$$

Il vient de même

$$d\arg h(u) = \mathrm{ord}_P(h) d\arg u + o_{u \to 0}(1) du + o_{u \to 0}(1) d\overline{u}.$$

Nous en tirons

$$\begin{aligned}\eta(f,g) &= \bigl(\mathrm{ord}_P(g)\log|\widehat{f}(P)| - \mathrm{ord}_P(f)\log|\widehat{g}(P)|\bigr)d\arg u \\ &\quad + o_{u\to 0}(\log|u|)\,du + o_{u\to 0}(\log|u|)\,d\overline{u} \\ &= \log|\partial_P\{f,g\}|\,d\arg u + o_{u\to 0}(\log|u|)\,du + o_{u\to 0}(\log|u|)\,d\overline{u}, \end{aligned} \quad (3.108)$$

la dernière égalité résultant de la définition du symbole modéré en $P$. Voici les conséquences de cette estimation. Posons $Y = X - S$. La quantité $\int_\gamma \eta(f,g)$, où $\gamma$ est un chemin non nécessairement fermé de $Y$, est bien définie et ne dépend que de la classe d'homotopie de $\gamma$. D'après l'estimation (3.108), l'intégrale $\int_\gamma \eta(f,g)$ converge absolument lorsque $\gamma$ est un chemin de classe $C^1$ dans $X$, à dérivée non nulle aux points de $S$. D'après cette même égalité, l'intégrale de $\eta(f,g)$ sur un lacet $\gamma_P$ orienté positivement autour d'un point $P \in S$, est donnée par

$$\int_{\gamma_P} \eta(f,g) = 2\pi \log|\partial_P\{f,g\}| \quad (3.109)$$

(comparer avec [56, Lemma p. 12]). Enfin, grâce à la formule de Stokes, le régulateur associé aux fonctions rationnelles $f,g$ et à la forme différentielle $\omega \in \Omega^{1,0}(X)$ s'exprime facilement en termes de la forme différentielle $\eta(f,g)$ :

$$\langle r_X(\{f,g\}), \omega\rangle = -\frac{i}{2}\int_X \omega \wedge \eta(f,g). \quad (3.110)$$

Remarquons que par linéarité, toutes les définitions et résultats ci-dessus s'étendent au cas où $f,g \in \mathbf{C}(X)^* \otimes \mathbf{C}$. Nous nous intéresserons par la suite au cas où $X$ est une courbe modulaire et $f,g$ sont des unités modulaires. Nous disposons par exemple des unités modulaires de la proposition 79. Néanmoins, nous avons besoin d'étendre légèrement la définition de $\eta$ aux séries d'Eisenstein ne provenant pas d'unités modulaires. Cela justifie la définition suivante. Nous identifierons les fonctions $\frac{\mathbf{Z}}{N\mathbf{Z}} \to \mathbf{C}$ aux diviseurs sur $\frac{\mathbf{Z}}{N\mathbf{Z}}$.

**Définition 88.** *Soient $l,m$ deux diviseurs sur $\frac{\mathbf{Z}}{N\mathbf{Z}}$, la forme différentielle $\eta(l,m)$ sur $\mathcal{H}$ est définie par*

$$\eta(l,m) = E_l^* \cdot (\partial - \overline{\partial})E_m^* - E_m^* \cdot (\partial - \overline{\partial})E_l^*. \quad (3.111)$$

Les séries d'Eisenstein $E_l^*$ et $E_m^*$ étant modulaires pour le groupe $\Gamma_1(N)$, la forme différentielle $\eta(l,m)$ est invariante sous l'action de ce groupe, définissant ainsi une forme différentielle de $Y_1(N)(\mathbf{C})$. Dans le cas où $l$ et $m$ sont des diviseurs de degré 0, la proposition 79 et la définition (3.105) montrent que

$$\eta(l,m) = \pi^2 i \cdot \eta(u_l, u_m). \quad (3.112)$$

Dans ce cas, la forme différentielle $\eta(l,m)$ est donc fermée. Ce n'est pas le cas en général, comme le montre le lemme suivant.

**Lemme 89.** *La différentielle de la forme $\eta(l,m)$ est donnée par*

$$d\eta(l,m) = \frac{\pi i}{N^2} E_{D(l,m)}^* \cdot \frac{dx \wedge dy}{y^2} \quad (3.113)$$

*où $D(l,m)$ est le diviseur (de degré 0) défini par*

$$D(l,m) = (\deg m)l - (\deg l)m. \quad (3.114)$$

*Démonstration.* Nous avons

$$d\eta(l,m) = (\partial + \overline{\partial})E_l^* \wedge (\partial - \overline{\partial})E_m^* - (\partial + \overline{\partial})E_m^* \wedge (\partial - \overline{\partial})E_l^*$$
$$+ E_l^* \cdot (-2\partial\overline{\partial}E_m^*) + E_m^* \cdot (2\partial\overline{\partial}E_l^*)$$
$$= -2E_l^* \cdot \partial\overline{\partial}E_m^* + 2E_m^* \cdot \partial\overline{\partial}E_l^*.$$

Par définition de $E_l^*$, nous avons

$$E_l^* = \frac{1}{N^2} \sum_{v \in \frac{\mathbf{Z}}{N\mathbf{Z}}} \operatorname{ord}_v(l) \sum_{a,b \in \frac{\mathbf{Z}}{N\mathbf{Z}}} e^{-\frac{2\pi i b v}{N}} \zeta_{a,b}^*.$$

La considération des développements de Fourier (3.16), (3.17) et (3.19) amène à la formule suivante pour $\partial\overline{\partial}E_l^*$

$$\partial\overline{\partial}E_l^* = \frac{1}{N^2} \sum_{v \in \frac{\mathbf{Z}}{N\mathbf{Z}}} \operatorname{ord}_v(l) \cdot \partial\overline{\partial}\zeta_{0,0}^*$$
$$= -\frac{\pi}{N^2}(\deg l) \cdot \partial\overline{\partial} \log y$$
$$= -\frac{\pi i}{2N^2}(\deg l) \cdot \frac{dx \wedge dy}{y^2},$$

la dernière égalité provenant de (1.32). Le résultat suit. □

Remarquons que l'intégrale (3.28) peut s'exprimer à l'aide de la forme différentielle $\eta(l,m)$. C'est un cas particulier du lemme suivant.

**Lemme 90.** *Pour toute forme parabolique $f$ de poids 2 pour $\Gamma_1(N)$, et tous diviseurs $l, m$ sur $\frac{\mathbf{Z}}{N\mathbf{Z}}$, nous avons*

$$\int_{X_1(N)(\mathbf{C})} E_l^* \cdot \omega_f \wedge \overline{\partial}E_m^* = -\frac{1}{2} \int_{X_1(N)(\mathbf{C})} \omega_f \wedge \eta(l,m). \tag{3.115}$$

*Démonstration.* Par définition de $\eta(l,m)$, nous avons

$$\int_{X_1(N)(\mathbf{C})} \omega_f \wedge \eta(l,m) = \int_{X_1(N)(\mathbf{C})} -E_l^* \cdot \omega_f \wedge \overline{\partial}E_m^* + E_m^* \cdot \omega_f \wedge \overline{\partial}E_l^*.$$

Or, une intégration par parties et la formule de Stokes donnent

$$\int_{X_1(N)(\mathbf{C})} E_m^* \cdot \omega_f \wedge \overline{\partial}E_l^* = -\int_{X_1(N)(\mathbf{C})} E_l^* \cdot \omega_f \wedge \overline{\partial}E_m^*,$$

puisque la forme différentielle $E_l^* E_m^* \cdot \omega_f$ est à croissance modérée aux pointes. □

La forme différentielle $\eta(l,m)$ n'étant pas nécessairement fermée, nous avons besoin de préciser la notation suivante.

*Notation.* Pour tous points $\alpha, \beta \in \mathcal{H} \cup \mathbf{P}^1(\mathbf{Q})$, nous convenons que l'intégrale $\int_\alpha^\beta \eta(l,m)$ est calculée le long d'une géodésique de $\mathcal{H}$ reliant $\alpha$ à $\beta$.

Lorsque $\alpha$ ou $\beta$ appartient à $\mathbf{P}^1(\mathbf{Q})$, la convergence absolue de l'intégrale résulte de (3.15), (3.16), (3.17) et (3.19).

*Notation.* Pour toute forme parabolique $f \in S_2(\Gamma_1(N))$, nous poserons

$$\xi_f(x) = -\frac{1}{2\pi} \int_{\xi(x)} \omega_f = -i \int_{g_x 0}^{g_x \infty} f(z) dz \qquad (x \in E_N). \tag{3.116}$$

Pour $x = (u, v) \in E_N$, posons $x^c = (-u, v)$ et $\xi_f^{\pm}(x) = \frac{1}{2}\bigl(\xi_f(x) \pm \xi_f(x^c)\bigr)$. Posons $\rho = e^{\frac{\pi i}{3}}$. Puisque $\eta(l, m)$ est invariante sous l'action du groupe $\Gamma_1(N)$ et que ce groupe préserve les géodésiques de $\mathcal{H}$, l'intégrale $\int_{g_x\rho}^{g_x\rho^2} \eta(l, m)$ ne dépend que de $x$. Pour tout $g \in \mathrm{SL}_2(\mathbf{Z})$ et toute forme différentielle $\eta$ sur $\mathcal{H}$, posons $\eta \mid g = g^*\eta$. Notons $\sigma$ et $\tau$ les matrices suivantes de $\mathrm{SL}_2(\mathbf{Z})$

$$\sigma = \begin{pmatrix} 0 & -1 \\ 1 & 0 \end{pmatrix} \quad \text{et} \quad \tau = \begin{pmatrix} 0 & -1 \\ 1 & -1 \end{pmatrix}.$$

**Définition 91.** *Pour tous diviseurs $l$ et $m$ sur $\frac{\mathbf{Z}}{N\mathbf{Z}}$, définissons un cycle relatif $c(l, m)$ sur $X_1(N)(\mathbf{C})$ par*

$$c(l, m) = \frac{1}{4} \sum_{x \in E_N} \Bigl( \int_{g_x\rho}^{g_x\rho^2} \eta(l, m) \Bigr) \xi(x). \tag{3.117}$$

Le cycle $c(l, m)$ définit un élément de $H_1(X_1(N)(\mathbf{C}), P_N, \mathbf{C})$, où $P_N$ est l'ensemble des pointes de $X_1(N)(\mathbf{C})$. Lorsque $l$ et $m$ sont de degré 0, le bord du cycle $c(l, m)$ possède une expression agréable en termes du symbole modéré associé aux unités modulaires $u_l$ et $u_m$, comme le montre la proposition suivante.

**Proposition 92.** *Pour tous diviseurs $l$ et $m$ de degré 0 sur $\frac{\mathbf{Z}}{N\mathbf{Z}}$, le bord du cycle $c(l, m)$ est donné par*

$$\partial c(l, m) = -2\pi^3 i \sum_{P \in P_N} \log|\partial_P\{u_l, u_m\}| \cdot [P], \tag{3.118}$$

*où $\partial_P$ désigne le symbole modéré en $P$.*

*Démonstration.* Nous avons

$$\partial c(l, m) = \frac{1}{4} \sum_{x \in E_N} \Bigl( \int_{\rho}^{\rho^2} \eta(l, m) \mid g_x \Bigr) \bigl([g_x \infty] - [g_x 0]\bigr).$$

Mais $[g_x 0] = [g_{x\sigma} \infty]$, d'où

$$\partial c(l, m) = \frac{1}{4} \sum_{x \in E_N} \Bigl( \int_{\rho}^{\rho^2} \eta(l, m) \mid g_x - \int_{\rho}^{\rho^2} \eta(l, m) \mid g_{x\sigma^{-1}} \Bigr) [g_x \infty].$$

Par changement de variable,

$$\int_{\rho}^{\rho^2} \eta(l, m) \mid g_{x\sigma^{-1}} = -\int_{\rho}^{\rho^2} \eta(l, m) \mid g_x,$$

il vient donc

$$\partial c(l, m) = \frac{1}{2} \sum_{x \in E_N} \Bigl( \int_{\rho}^{\rho^2} \eta(l, m) \mid g_x \Bigr) [g_x \infty]$$

$$= \sum_{x \in E_N/\pm 1} \Bigl( \int_{\rho}^{\rho^2} \eta(l, m) \mid g_x \Bigr) [g_x \infty].$$

Pour toute pointe $P \in P_N \cong E_N/\pm\Gamma_\infty$, notons $\widetilde{E}_P$ l'ensemble des préimages de $P$ par la surjection naturelle $E_N/\pm 1 \to E_N/\pm\Gamma_\infty$. Par définition, la largeur $l_P$ de la pointe $P$ est le cardinal de $\widetilde{E}_P$. La matrice $T$ agit sur $\widetilde{E}_P$ par multiplication à droite ; cette action est transitive, cyclique d'ordre $l_P$. Nous pouvons écrire

$$\partial c(l,m) = \sum_{P \in P_N} \Big( \sum_{x \in \widetilde{E}_P} \int_\rho^{\rho^2} \eta(l,m) \mid g_x \Big)[P].$$

Fixons $P \in P_N$. Il s'agit de calculer l'intégrale de $\eta(l,m)$ sur le cycle (a priori relatif) $\gamma_P$ de $Y_1(N)(\mathbf{C})$ défini par

$$\gamma_P = \sum_{x \in \widetilde{E}_P} \{g_x\rho, g_x\rho^2\}.$$

Fixons $x_0 \in \widetilde{E}_P$. D'après $\rho = T\rho^2$, nous avons

$$\begin{aligned}
\gamma_P &= \sum_{k=0}^{l_P-1} \{g_{x_0 T^k} T\rho^2, g_{x_0 T^k}\rho^2\} \\
&= \sum_{k=0}^{l_P-1} \{g_{x_0} T^{k+1}\rho^2, g_{x_0} T^k\rho^2\} \\
&= \{g_{x_0} T^{l_P}\rho^2, g_{x_0}\rho^2\}.
\end{aligned} \tag{3.119}$$

Puisque $g_{x_0} T^{l_P} g_{x_0}^{-1} \in \pm\Gamma_1(N)$, nous voyons déjà que $\gamma_P$ est un cycle fermé de $Y_1(N)(\mathbf{C})$, autrement dit $\gamma_P \in H_1(Y_1(N)(\mathbf{C}), \mathbf{Z})$. Pour la même raison, nous pouvons remplacer $\rho^2$ par n'importe quel point de $\mathcal{H}$ dans (3.119). Par conséquent, $\gamma_P$ n'est autre qu'un lacet orienté négativement autour de $P \in X_1(N)(\mathbf{C})$. D'après (3.112) et (3.109), nous avons alors

$$\int_{\gamma_P} \eta(l,m) = \pi^2 i \int_{\gamma_P} \eta(u_l, u_m) = -2\pi^3 i \log|\partial_P\{u_l, u_m\}|,$$

d'où nous déduisons (3.118). $\square$

Le théorème qui suit est inspiré du théorème C de l'appendice (p. 152), dû à Merel, qui exprime le produit scalaire de Petersson de deux formes paraboliques en fonction de leurs périodes de Manin. La démonstration suit largement celle du théorème C.

Notons $\mathcal{F}$ le domaine fondamental standard du demi-plan de Poincaré :

$$\mathcal{F} = \{z \in \mathcal{H}\,;\, |\Re(z)| \leq \frac{1}{2} \text{ et } |z| \geq 1\}. \tag{3.120}$$

**Théorème 93.** *Soient $f$ une forme parabolique de poids 2 pour $\Gamma_1(N)$ et $l, m$ deux diviseurs sur $\frac{\mathbf{Z}}{N\mathbf{Z}}$. Posons*

$$F_x(z) = \int_\infty^z \omega_f \mid g_x \qquad (x \in E_N, z \in \mathcal{H}). \tag{3.121}$$

*Nous avons alors*

$$\int_{X_1(N)(\mathbf{C})} \omega_f \wedge \eta(l,m) = \int_{c(l,m)} \omega_f - \int_{\mathcal{F}} \sum_{x \in E_N/\pm 1} F_x \cdot d\big(\eta(l,m) \mid g_x\big). \tag{3.122}$$

*Lorsque $l$ et $m$ sont de degré $0$, nous avons en particulier*

$$\int_{X_1(N)(\mathbf{C})} \omega_f \wedge \eta(l,m) = \int_{c(l,m)} \omega_f = -\frac{\pi}{2} \sum_{x \in E_N} \Big(\int_{g_x \rho}^{g_x \rho^2} \eta(l,m)\Big) \xi_f(x). \qquad (3.123)$$

*Remarque* 94. Nous n'avons pas réussi à nous débarasser du second terme de (3.122), ce qui nous obligera à faire des hypothèses sur les caractères de Dirichlet par la suite. Notons cependant que ce terme désagréable ne dépend que de $f$ et du diviseur $D(l,m)$. Nous ne savons pas si la formule (3.123) est valable pour tous diviseurs $l$ et $m$.

*Démonstration.* Nous avons

$$\int_{X_1(N)(\mathbf{C})} \omega_f \wedge \eta(l,m) = \sum_{x \in E_N/\pm 1} \int_{g_x \mathcal{F}} \omega_f \wedge \eta(l,m)$$
$$= \sum_{x \in E_N/\pm 1} \int_{\mathcal{F}} \big(\omega_f \wedge \eta(l,m)\big) \mid g_x.$$

Or, nous avons

$$d\big(F_x \cdot \eta(l,m) \mid g_x\big) = dF_x \wedge \eta(l,m) \mid g_x + F_x \cdot d\big(\eta(l,m) \mid g_x\big)$$
$$= \big(\omega_f \wedge \eta(l,m)\big) \mid g_x + F_x \cdot d\big(\eta(l,m) \mid g_x\big).$$

D'après la formule de Stokes, nous avons donc

$$\int_{\mathcal{F}} \big(\omega_f \wedge \eta(l,m)\big) \mid g_x = \int_{\partial \mathcal{F}} F_x \cdot \eta(l,m) \mid g_x - \int_{\mathcal{F}} F_x \cdot d\big(\eta(l,m) \mid g_x\big).$$

Le bord du domaine $\mathcal{F}$ est le triangle géodésique de sommets $\rho^2$, $\rho$ et $\infty$, orienté dans cet ordre. Puisque la forme $f$ est parabolique, la fonction $F_x$ est à décroissance exponentielle en $\infty$; par ailleurs, la forme différentielle $\eta(l,m) \mid g_x$ est à croissance modérée en $\infty$. Cela nous permet d'écrire

$$\int_{\partial \mathcal{F}} F_x \cdot \eta(l,m) \mid g_x = \Big(\int_{\rho^2}^{\rho} + \int_{\rho}^{\infty} + \int_{\infty}^{\rho^2}\Big) F_x \cdot \eta(l,m) \mid g_x.$$

Nous écrirons $\sum_x$ pour signifier la somme sur les $x \in E_N/\pm 1$. Nous allons démontrer que

$$\sum_x \Big(\int_{\rho}^{\infty} + \int_{\infty}^{\rho^2}\Big) F_x \cdot \eta(l,m) \mid g_x = 0.$$

La matrice $T = \begin{pmatrix} 1 & 1 \\ 0 & 1 \end{pmatrix} \in \mathrm{SL}_2(\mathbf{Z})$ fixe $\infty$ et envoie $\rho^2$ sur $\rho$. Par changement de variable, nous obtenons

$$\int_{\rho}^{\infty} F_x \cdot \eta(l,m) \mid g_x = \int_{\rho^2}^{\infty} F_x \mid T \cdot \eta(l,m) \mid g_x T.$$

Par ailleurs, nous avons

$$F_x(Tz) = \int_{\infty}^{Tz} \omega_f \mid g_x = \int_{\infty}^{z} \omega_f \mid g_x T = F_{xT}(z) \qquad (x \in E_N/\pm 1, z \in \mathcal{H}).$$

En sommant sur les $x$, il vient donc

$$\sum_x \int_\rho^\infty F_x \cdot \eta(l,m) \mid g_x = \sum_x \int_{\rho^2}^\infty F_{xT} \cdot \eta(l,m) \mid g_{xT}$$
$$= \sum_x \int_{\rho^2}^\infty F_x \cdot \eta(l,m) \mid g_x,$$

ce qui montre la simplification annoncée. Utilisons maintenant la matrice $\sigma = \begin{pmatrix} 0 & -1 \\ 1 & 0 \end{pmatrix} \in$ SL$_2(\mathbf{Z})$, qui échange $\rho$ et $\rho^2$. Par changement de variables, nous obtenons

$$\int_{\rho^2}^\rho F_x \cdot \eta(l,m) \mid g_x = \int_\rho^{\rho^2} F_x \mid \sigma \cdot \eta(l,m) \mid g_{x\sigma}.$$

Par ailleurs, puisque $\sigma$ échange $0$ et $\infty$, nous avons

$$F_x(\sigma z) = \int_\infty^{\sigma z} \omega_f \mid g_x = \int_0^z \omega_f \mid g_{x\sigma} = \int_0^\infty \omega_f \mid g_{x\sigma} + F_{x\sigma}(z) = F_{x\sigma}(z) + 2\pi \xi_f(x),$$

la quantité $\xi_f(x)$ étant indépendante de $z$. En remplaçant dans l'intégrale précédente, il vient

$$\int_{\rho^2}^\rho F_x \cdot \eta(l,m) \mid g_x = \int_\rho^{\rho^2} F_{x\sigma} \cdot \eta(l,m) \mid g_{x\sigma} + 2\pi \xi_f(x) \int_\rho^{\rho^2} \eta(l,m) \mid g_{x\sigma}$$

En sommant sur les $x$, nous obtenons

$$\sum_x \int_{\rho^2}^\rho F_x \cdot \eta(l,m) \mid g_x = \frac{1}{2} \sum_x \int_{\rho^2}^\rho F_x \cdot \eta(l,m) \mid g_x + \int_{\rho^2}^\rho F_{x\sigma} \cdot \eta(l,m) \mid g_{x\sigma}$$
$$= \frac{1}{2} \sum_x 2\pi \xi_f(x) \int_\rho^{\rho^2} \eta(l,m) \mid g_{x\sigma}$$
$$= -\pi \sum_x \xi_f(x) \int_\rho^{\rho^2} \eta(l,m) \mid g_x.$$

Nous obtenons donc au final

$$\int_{X_1(N)(\mathbf{C})} \omega_f \wedge \eta(l,m) = \int_{c(l,m)} \omega_f - \int_\mathcal{F} \sum_{x \in E_N/\pm 1} F_x \cdot d\bigl(\eta(l,m) \mid g_x\bigr).$$

Lorsque $l$ et $m$ sont de degré 0, le lemme 89 montre que la forme différentielle $\eta(l,m)$ est fermée, ce qui entraîne (3.123). □

Notons $\mathbf{1}$ le diviseur sur $\frac{\mathbf{Z}}{N\mathbf{Z}}$ valant 1 en $\bar 1 \in \frac{\mathbf{Z}}{N\mathbf{Z}}$ et 0 ailleurs.

**Corollaire** (Théorème 3). *Soit $f$ une forme parabolique primitive de poids 2 pour $\Gamma_1(N)$, de caractère $\psi$. Pour tout caractère de Dirichlet $\chi$ modulo $N$, pair et distinct de $\overline{\psi}$, nous avons*

$$L(f,2)L(f,\chi,1) = \frac{N\,i}{4} \sum_{x \in E_N} \Bigl(\int_{g_x\rho}^{g_x\rho^2} \eta(1,\widehat{\chi})\Bigr) \xi_f^+(x). \tag{3.124}$$

*Démonstration du théorème 3.* D'après le théorème 72 et le lemme 90, nous avons

$$L(f,2)L(f,\chi,1) = -\frac{N}{\pi i \varphi(N)} \int_{X_1(N)(\mathbf{C})} E^*_{\psi\chi} \cdot \omega_f \wedge \overline{\partial} E^*_{\widehat{\chi}}$$
$$= \frac{N}{2\pi i \varphi(N)} \int_{X_1(N)(\mathbf{C})} \omega_f \wedge \eta(\psi\chi, \widehat{\chi}).$$

Appliquons maintenant le théorème 93. Par hypothèse sur $\chi$, le diviseur $\psi\chi$ est de degré 0. Il en va de même pour le diviseur $\widehat{\chi}$, puisque $N > 1$ par hypothèse sur $f$. Nous en déduisons

$$L(f,2)L(f,\chi,1) = \frac{N}{2\pi i \varphi(N)} \int_{c(\psi\chi,\widehat{\chi})} \omega_f. \tag{3.125}$$

Soit maintenant $\epsilon$ un caractère de Dirichlet pair modulo $N$, avec $\epsilon \neq \psi\chi$. Nous allons utiliser des arguments de caractère pour montrer

$$\int_{c(\epsilon,\widehat{\chi})} \omega_f = 0.$$

Fixons $d \in (\frac{\mathbf{Z}}{N\mathbf{Z}})^*$. De (3.74) et (3.75), nous tirons

$$\langle d \rangle^* \eta(\epsilon, \widehat{\chi}) = \overline{\epsilon}\chi(d)\, \eta(\epsilon, \widehat{\chi}).$$

D'autre part, puisque $f$ est de caractère $\psi$, nous avons

$$\xi_f(dx) = \psi(d)\, \xi_f(x) \qquad (x \in E_N).$$

Il suit

$$\int_{c(\epsilon,\widehat{\chi})} \omega_f = -\frac{\pi}{2} \sum_{x \in E_N} \Big( \int_{g_x \rho}^{g_x \rho^2} \eta(\epsilon, \widehat{\chi}) \Big) \xi_f(x)$$
$$= -\frac{\pi}{2} \sum_{x \in E_N} \Big( \int_{g_{dx} \rho}^{g_{dx} \rho^2} \eta(\epsilon, \widehat{\chi}) \Big) \xi_f(dx).$$

Puisque $\langle d \rangle$ préserve les géodésiques de $\mathcal{H}$, nous déduisons

$$\int_{c(\epsilon,\widehat{\chi})} \omega_f = -\frac{\pi}{2} \sum_{x \in E_N} \Big( \int_{g_x \rho}^{g_x \rho^2} \langle d \rangle^* \eta(\epsilon, \widehat{\chi}) \Big) \xi_f(dx)$$
$$= -\frac{\pi}{2} \sum_{x \in E_N} \Big( \int_{g_x \rho}^{g_x \rho^2} \overline{\epsilon}\chi(d)\, \eta(\epsilon, \widehat{\chi}) \Big) \psi(d)\, \xi_f(x)$$
$$= \overline{\epsilon}\chi\psi(d) \cdot \int_{c(\epsilon,\widehat{\chi})} \omega_f.$$

D'où $\int_{c(\epsilon,\widehat{\chi})} \omega_f = 0$. Prenons maintenant la somme sur tous les caractères pairs $\epsilon$ modulo $N$. Nous obtenons

$$L(f,2)L(f,\chi,1) = \frac{N}{2\pi i \varphi(N)} \sum_\epsilon \int_{c(\epsilon,\widehat{\chi})} \omega_f$$
$$= \frac{N}{2\pi i \varphi(N)} \int_{c(\sum_\epsilon \epsilon,\widehat{\chi})} \omega_f.$$

Remarquons que $E_l^*$ est nulle lorsque $l : \frac{\mathbf{Z}}{N\mathbf{Z}} \to \mathbf{C}$ est une fonction impaire. Il en va donc de même pour $\eta(l,m)$ et $c(l,m)$ lorsque $l$ est impaire. Dans la somme ci-dessus, nous pouvons donc faire parcourir à $\epsilon$ tous les caractères de Dirichlet modulo $N$. La somme de tous les caractères de Dirichlet modulo $N$ donnant $\varphi(N) \cdot 1$, nous obtenons

$$L(f,2)L(f,\chi,1) = \frac{N}{2\pi i} \int_{c(1,\widehat{\chi})} \omega_f.$$

En remplaçant $c(1,\widehat{\chi})$ par sa définition, il vient

$$L(f,2)L(f,\chi,1) = \frac{Ni}{4} \sum_{x \in E_N} \Big(\int_{g_x\rho}^{g_x\rho^2} \eta(1,\widehat{\chi})\Big) \xi_f(x). \tag{3.126}$$

En vue de faire apparaître $\xi_f^+(x)$ dans l'équation précédente, utilisons le changement de variable $x \mapsto x^c$. Pour tous diviseurs $l$, $m$, nous avons

$$\int_{g_{x^c}\rho}^{g_{x^c}\rho^2} \eta(l,m) = \int_{c(g_x\rho^2)}^{c(g_x\rho)} \eta(l,m)$$
$$= -\int_{g_x\rho}^{g_x\rho^2} c^*\eta(l,m).$$

D'après (3.25), les séries d'Eisenstein $E_l^*$ vérifient $c^*E_l^* = E_l^*$, il suit

$$c^*\eta(l,m) = c^*E_l^* \cdot (\overline{\partial} - \partial)c^*E_m^* - c^*E_m^* \cdot (\overline{\partial} - \partial)c^*E_m^* = -\eta(l,m).$$

En prenant la demi-somme à partir de (3.126), il vient donc

$$L(f,2)L(f,\chi,1) = \frac{Ni}{4} \sum_{x \in E_N} \Big(\int_{g_x\rho}^{g_x\rho^2} \eta(1,\widehat{\chi})\Big) \xi_f^+(x)$$

c'est-à-dire (3.124). □

*Remarque* 95. Au cours de la démonstration du théorème 3, nous avons obtenu l'expression

$$L(f,2)L(f,\chi,1) = \frac{N}{2\pi i \varphi(N)} \int_{c(\psi\chi,\widehat{\chi})} \omega_f.$$

Dans la pratique, cette expression est peut-être plus avantageuse que (3.124), puisque la définition du cycle $c(\psi\chi,\widehat{\chi})$ fait intervenir la forme différentielle $\eta(\psi\chi,\widehat{\chi})$, qui est fermée. Faisons l'hypothèse supplémentaire $\chi$ *primitif*. Le diviseur $\widehat{\chi}$ est donc à support dans $(\frac{\mathbf{Z}}{N\mathbf{Z}})^*$. Les propositions 86 et 92 entraînent alors que le cycle $c(\psi\chi,\widehat{\chi})$ est *fermé*.

*Remarque* 96. Il serait agréable de s'affranchir de l'hypothèse sur $\chi$ dans le théorème 3. La formule (3.124) nous semble suffisamment naturelle pour être vraie pour tout caractère $\chi$, voire pour toute application paire de $\mathbf{Z}/N\mathbf{Z}$ dans $\mathbf{C}$, mais nous n'avons pas de démonstration. Nous remarquons également que dans le cas $\chi = \overline{\psi}$, l'équation fonctionnelle donne une expression de $L(f,\overline{\psi},1)$ en termes de $L(f,1_N,1)$, où $1_N$ désigne le caractère trivial modulo $N$. Cela conduit à une expression pour $L(f,2)L(f,\overline{\psi},1)$, au moins dans le cas $\psi \neq 1_N$.

*Question.* Soit $f$ une forme primitive de poids 2 pour $\Gamma_1(N)$. Existe-t-il une formule pour $L(f,2)L(f,1)$ analogue à (3.124) ?

## 3.6 Application aux courbes elliptiques

Nous détaillons maintenant les conséquences de la formule de la section précédente pour les courbes elliptiques. Nous montrons en particulier le résultat suivant.

**Théorème 97.** *Soit $N \geq 1$ un entier. Pour toute courbe elliptique $E$ définie sur $\mathbf{Q}$, de conducteur $N$, la valeur spéciale $L(E,2)$ est combinaison linéaire à coefficients rationnels des quantités*

$$\frac{1}{\pi i} \int_{g_x \rho}^{g_x \rho^2} \eta(1,a), \tag{3.127}$$

*où $x$ et $a$ parcourent respectivement $E_N$ et $\frac{\mathbf{Z}}{N\mathbf{Z}}$.*

Il est tout à fait frappant de remarquer que $L(E,2)$ est combinaison linéaire d'un nombre fini de quantités purement analytiques et dépendant exclusivement du conducteur de $E$. On peut rapprocher ce résultat de la conjecture de Zagier pour $L(E,2)$, démontrée par Goncharov et Levin [32], qui exprime $L(E,2)$ comme combinaison linéaire à coefficients rationnels de valeurs de la fonction dilogarithme elliptique en des points algébriques de $E$. Il serait intéressant de pousser plus avant cette analogie.

*Démonstration.* Soit $E$ une courbe elliptique définie sur $\mathbf{Q}$, de conducteur $N$. D'après les travaux de Breuil, Conrad, Diamond, Taylor et Wiles [78, 74, 17], il existe une forme primitive $f \in S_2(\Gamma_0(N))$ telle que $L(E,s) = L(f,s)$. D'après le théorème 3, nous avons pour tout caractère $\chi$ modulo $N$, pair et non trivial

$$L(f,2)L(f,\chi,1) = \frac{Ni}{4} \sum_{x \in E_N} \Big(\int_{g_x \rho}^{g_x \rho^2} \eta(1,\widehat{\chi})\Big) \xi_f^+(x). \tag{3.128}$$

L'idée naturelle consiste à diviser (3.128) par la période réelle de $E$. Notons $\Omega_E^+$ cette période réelle, définie comme la valeur absolue de l'intégrale d'une forme différentielle de Néron de $E$ le long d'un générateur de $H_1^+(E(\mathbf{C}),\mathbf{Z})$. D'après le théorème de Manin-Drinfel'd [28, 29], nous avons

$$\xi_f^+(x) \in \mathbf{Q} \cdot \frac{\Omega_E^+}{2\pi} \qquad (x \in E_N). \tag{3.129}$$

Soit $\alpha$ une application de $\frac{\mathbf{Z}}{N\mathbf{Z}}$ dans $\mathbf{C}$. Notons $L(f,\alpha,s)$ la série $L$ de $f$ tordue par $\alpha$. Un calcul classique montre que

$$L(f,\alpha,1) = -\frac{1}{N} \sum_{a \in \frac{\mathbf{Z}}{N\mathbf{Z}}} \widehat{\alpha}(a) \int_{a/N}^{\infty} \omega_f. \tag{3.130}$$

Cela montre en particulier que $L(f,\chi,1) \in \mathbf{Q}(\widehat{\chi}) \cdot \Omega_E^+$ pour $\chi$ caractère pair modulo $N$, où $\mathbf{Q}(\widehat{\chi})$ est le corps engendré par les valeurs de $\widehat{\chi}$. Nous souhaitons réduire ce corps des coefficients à $\mathbf{Q}$. Notons $V$ le sous-$\mathbf{C}$-espace vectoriel de $\mathbf{C}[\frac{\mathbf{Z}}{N\mathbf{Z}}]$ engendré par les caractères de Dirichlet modulo $N$, pairs et non triviaux. Soit $\widehat{V}$ l'image de $V$ par la transformée de Fourier ; c'est aussi le sous-espace engendré par les $\widehat{\chi}$, avec $\chi$ pair et non trivial modulo $N$. Pour toute application $m$ de $\frac{\mathbf{Z}}{N\mathbf{Z}}$ dans $\mathbf{C}$, posons

$$\int_{m/N}^{\infty} \omega_f = \sum_{a \in \frac{\mathbf{Z}}{N\mathbf{Z}}} m(a) \int_{a/N}^{\infty} \omega_f. \tag{3.131}$$

Par linéarité, l'identité (3.128) entraîne visiblement

$$L(f,2) \int_{m/N}^{\infty} \omega_f = -\frac{N^2 i}{4} \sum_{x \in E_N} \Big( \int_{g_x\rho}^{g_x\rho^2} \eta(1,m) \Big) \xi_f^+(x) \qquad (m \in \widehat{V}). \tag{3.132}$$

Les deux lemmes suivants montrent qu'il existe toujours $m \in \widehat{V} \cap \mathbf{Q}[\frac{\mathbf{Z}}{N\mathbf{Z}}]$ tel que $\int_{m/N}^{\infty} \omega_f \neq 0$, ce qui entraîne le théorème.

**Lemme 98.** *L'espace $\widehat{V}_{\mathbf{Q}} = \widehat{V} \cap \mathbf{Q}[\frac{\mathbf{Z}}{N\mathbf{Z}}]$ est une $\mathbf{Q}$-structure de $\widehat{V}$.*

**Lemme 99.** *Il existe $m \in \widehat{V}$ tel que $\int_{m/N}^{\infty} \omega_f \neq 0$.*

*Démonstration du lemme 98.* Il s'agit de montrer que $\widehat{V}_{\mathbf{Q}}$ engendre $\widehat{V}$ comme espace vectoriel complexe. Notons $W$ le sous-espace de $\mathbf{C}[\frac{\mathbf{Z}}{N\mathbf{Z}}]$ formé des applications paires et à support dans $(\frac{\mathbf{Z}}{N\mathbf{Z}})^*$. L'espace $V$ n'est autre que l'hyperplan de $W$ formé des applications de somme nulle. L'espace $W$ admet le système de générateurs suivants

$$f_a = [a] + [-a] \qquad (a \in (\frac{\mathbf{Z}}{N\mathbf{Z}})^*).$$

Nous avons

$$\widehat{f_a}(n) = e^{-\frac{2\pi i a n}{N}} + e^{\frac{2\pi i a n}{N}} \qquad (a \in (\frac{\mathbf{Z}}{N\mathbf{Z}})^*,\ n \in \frac{\mathbf{Z}}{N\mathbf{Z}}).$$

Notons $\mathbf{Q}(\mu_N)$ le corps cyclotomique engendré par $\zeta_N = e^{\frac{2\pi i}{N}}$. Nous identifions $(\frac{\mathbf{Z}}{N\mathbf{Z}})^*$ au groupe de Galois $G = \mathrm{Gal}(\mathbf{Q}(\mu_N)/\mathbf{Q})$, au moyen de $a \mapsto \sigma_a$, où $\sigma_a \in G$ est défini par $\sigma_a(\zeta_N) = \zeta_N^a$. Soit $\mathbf{Q}(\mu_N)^0$ le sous-espace de $\mathbf{Q}(\mu_N)$ formé des éléments de trace nulle. Soit $(e_1, \ldots e_r)$ une base de $\mathbf{Q}(\mu_N)^0$ sur $\mathbf{Q}$. Nous avons alors

$$\sum_{\sigma \in G} \sigma(e_i) \widehat{f_\sigma}(n) = \mathrm{tr}\big(e_i(\zeta_N^n + \zeta_N^{-n})\big) \qquad (1 \leq i \leq r,\ n \in \frac{\mathbf{Z}}{N\mathbf{Z}}).$$

Les éléments

$$g_i = \sum_{\sigma \in G} \sigma(e_i) f_\sigma \qquad (1 \leq i \leq r)$$

appartiennent à $V$ puisque $e_i \in \mathbf{Q}(\mu_N)^0$, et le calcul précédent entraîne

$$\widehat{g_i} \in \widehat{V}_{\mathbf{Q}} \qquad (1 \leq i \leq r)$$

Montrons que ces éléments engendrent $\widehat{V}$ comme espace vectoriel complexe. Posant

$$e_0 = 1 \quad \text{et} \quad g_0 = \sum_{\sigma \in G} f_\sigma,$$

il est classique de montrer que les $g_i$, $0 \le i \le r$ engendrent $W$. Par suite, les $\widehat{g_i}$, $0 \le i \le r$ engendrent $\widehat{W}$. Puisque $\widehat{V}$ est un hyperplan de $\widehat{W}$, un argument de dimension montre que les $\widehat{g_i}$, $1 \le i \le r$ engendrent $\widehat{V}$. □

*Démonstration du lemme 99.* Il suffit de montrer qu'il existe un caractère de Dirichlet $\chi_N$ modulo $N$, pair et non trivial, tel que $L(f, \chi_N, 1) \ne 0$. Supposons le contraire. D'après le corollaire 2 de l'appendice (p. 145), il existe un caractère $\chi$ pair et primitif, de conducteur divisant $N$, tel que $L(f \otimes \chi, 1) \ne 0$. Notons $\chi_N$ le caractère modulo $N$ induit par $\chi$. La non-nullité de $L(f \otimes \chi, 1)$ équivaut à celle de $L(f, \chi, 1)$, elle même équivalente à celle de $L(f, \chi_N, 1)$. Donc $\chi_N$ est égal au caractère trivial $1_N$. De plus, pour tout caractère $\chi$ pair et primitif de conducteur divisant $N$, nous avons $\Lambda(f \otimes \chi, 1) = 0$ si $\chi$ est non trivial. Dans le théorème A de l'appendice (p. 144), la somme donnant $\xi_f^+(u,v)$ est donc réduite au terme $\chi = 1$, ce qui admet la conséquente suivante

$$\xi_f^+(\lambda u, v) = \xi_f^+(u, \lambda v) = \xi_f^+(u,v) \qquad (\lambda \in (\frac{\mathbf{Z}}{N\mathbf{Z}})^*). \tag{3.133}$$

Remarquons également que $\xi_f^+(u,v)$ est la partie réelle de $\xi_f(u,v)$. Nous allons utiliser l'opérateur de Hecke $T_2$. D'après [47, Thm 2], nous avons

$$T_2\xi(u,v) = \xi(2u,v) + \xi(u,2v) + \xi(2u, u+v) + \xi(u+v, 2v) \qquad \big((u,v) \in E_N\big), \tag{3.134}$$

où par convention $\xi(u',v') = 0$ lorsque $(u',v') \notin E_N$. Nous allons maintenant distinguer les cas suivant la valuation de $N$ au nombre premier 2.

**Premier cas : $N$ impair.** Utilisant (3.134) avec $(u,v) = (0,1)$, nous trouvons

$$T_2\xi(0,1) = \xi(0,1) + \xi(0,2) + \xi(0,1) + \xi(1,2).$$

En appliquant cette égalité à $f$ et en prenant la partie réelle, il vient

$$a_2(f)\xi_f^+(0,1) = 3\xi_f^+(0,1) + \xi_f^+(1,2) = 3\xi_f^+(0,1) + \xi_f^+(1,1)$$

d'après (3.133). Les relations de Manin sur les symboles modulaires donnent $\xi(1,1) = \xi(1,0) + \xi(0,1) = 0$. Par suite

$$(a_2(f) - 3)\xi_f^+(0,1) = 0.$$

Par les bornes de Hasse et Weil $a_2(f) - 3 \ne 0$, d'où une contradiction avec l'hypothèse $L(f,1) \ne 0$.

**Deuxième cas : $N = 2N'$ avec $N'$ impair.** Utilisant (3.134) avec $(u,v) = (2t, 1-2t) \in E_N$, $t \in \mathbf{Z}/N'\mathbf{Z}$, il vient

$$T_2\xi(2t, 1-2t) = \xi(4t, 1-2t) + \xi(2t, 2-4t) + \xi(4t, 1) + \xi(1, 2-4t)$$
$$= \xi(4t, 1-2t) + \xi(4t, 1) + \xi(1, 2-4t)$$

car $(2t, 2-4t) \notin E_N$. En appliquant à $f$ et en prenant la partie réelle, nous déduisons

$$a_2(f)\xi_f^+(2t, 1-2t) = \xi_f^+(4t, 1-2t) + \xi_f^+(4t, 1) + \xi_f^+(1, 2-4t)$$
$$= \xi_f^+(2t, 1-2t) + \xi_f^+(2t, 1) + \xi_f^+(1, 2-4t)$$

car $(4t, N) = (2t, N)$. Avec $t = 0$, nous obtenons

$$\xi_f^+(1, 2) = (a_2(f) - 2)\xi_f^+(0, 1) \neq 0.$$

Avec $t = 1$, nous avons également

$$a_2(f)\xi_f^+(2, -1) = \xi_f^+(2, -1) + \xi_f^+(2, 1) + \xi_f^+(1, -2) = \xi_f^+(2, -1),$$

d'après la première relation de Manin. Comme $\xi_f^+(2, -1) = -\xi_f^+(1, 2) \neq 0$, il vient $a_2(f) = 1$. Pour tout $t$, nous obtenons donc

$$\xi_f^+(2t, 1) + \xi_f^+(1, 2 - 4t) = 0$$

Nous avons alors successivement

$$\xi_f^+(1, 2t) = -\xi_f^+(2t, 1) = \xi_f^+(1, 2 - 4t) = \xi_f^+(1, 4t - 2)$$
$$= \xi_f^+(1, 2 - 4(1 - t)) = -\xi_f^+(2(1 - t), 1) = \xi_f^+(1, 2 - 2t)$$
$$= \xi_f^+(1, 2t - 2) \qquad (t \in \mathbf{Z}/N'\mathbf{Z}).$$

D'où $\xi_f^+(1, 2t) = \xi_f^+(1, 0)$ pour tout $t$. Notons $w(f) = \pm 1$ la valeur propre de $f$ pour l'involution d'Atkin-Lehner $W_N$. Pour tout $a \in \frac{\mathbf{Z}}{N\mathbf{Z}}$, nous avons

$$\xi_f(1, a) = -\frac{w(f)}{2\pi}\int_{W_N\xi(1,a)} \omega_f = -\frac{w(f)}{2\pi}\int_{a/N}^{\infty} \omega_f. \tag{3.135}$$

Notons $1_{N'}$ le caractère trivial modulo $N'$. D'après (3.130) et (3.135), nous avons

$$L(f, 1_{N'}, 1) = -\frac{1}{N'}\sum_{a \in \mathbf{Z}/N'\mathbf{Z}} \widehat{1_{N'}}(a) \int_{a/N'}^{\infty} \omega_f$$
$$= \frac{2\pi w(f)}{N'}\sum_{a \in \mathbf{Z}/N'\mathbf{Z}} \widehat{1_{N'}}(a)\xi_f(1, 2a)$$
$$= \frac{2\pi w(f)}{N'}\sum_{a \in \mathbf{Z}/N'\mathbf{Z}} \widehat{1_{N'}}(a)\xi_f^+(1, 2a)$$
$$= -\frac{w(f)}{N'}\Big(\sum_{a \in \mathbf{Z}/N'\mathbf{Z}} \widehat{1_{N'}}(a)\Big)\xi_f^+(1, 0).$$

Puisque $N' > 1$, la somme intérieure est nulle, ce qui contredit l'hypothèse $L(f, 1) \neq 0$.

**Troisième cas :** $N = 4N'$ **avec** $N'$ **entier.** La courbe elliptique $E$ a donc réduction additive en 2, d'où $a_2(f) = 0$. En utilisant (3.134) avec $(u, v) = (0, 1)$, nous trouvons

$$\xi_f^+(1, 2) = 2\xi_f(1, 0). \tag{3.136}$$

D'autre part, avec $(u, v) = (1, 2)$ et en utilisant les relations de Manin, nous obtenons

$$\begin{aligned} 0 &= \xi_f^+(1,4) + \xi_f^+(2,3) + \xi_f^+(3,4) \\ &= \xi_f^+(1,4) + \xi_f^+(2,1) + \xi_f^+(1,3) + \xi_f^+(3,1) + \xi_f^+(1,4) \\ &= 2\xi_f^+(1,4) - \xi_f^+(1,2), \end{aligned}$$

d'où

$$\xi_f^+(1,4) = \xi_f(1,0). \qquad (3.137)$$

D'après (3.130) et (3.135), nous avons

$$\xi_f(1,0) = \frac{w(f)}{2\pi} L(f, 1). \qquad (3.138)$$

Pour tout $k$ diviseur de $N$, notons $\epsilon_k$ le diviseur suivant sur $\frac{\mathbf{Z}}{N\mathbf{Z}}$

$$\epsilon_k = \sum_{\lambda \in (\mathbf{Z}/\frac{N}{k}\mathbf{Z})^*} [k\lambda].$$

Nous avons $\widehat{\widehat{\epsilon_k}} = N\epsilon_k$, d'où

$$L(f, \widehat{\epsilon_k}, 1) = -\int_{\epsilon_k/N}^{\infty} \omega_f = 2\pi w(f)\xi_f(1, \epsilon_k) = 2\pi w(f)\xi_f^+(1, \epsilon_k).$$

Grâce à (3.133), nous en déduisons

$$L(f, \widehat{\epsilon_k}, 1) = 2\pi w(f)\varphi(N/k)\xi_f^+(1, k). \qquad (3.139)$$

L'application $\widehat{\epsilon_k}$ est $\frac{N}{k}$-périodique et invariante sous l'action de $(\mathbf{Z}/\frac{N}{k}\mathbf{Z})^*$. Pour la déterminer, il suffit donc de calculer $\widehat{\epsilon_k}(d)$ pour $d \mid \frac{N}{k}$, et nous avons

$$\begin{aligned} \widehat{\epsilon_k}(d) &= \sum_{\lambda \in (\mathbf{Z}/\frac{N}{k}\mathbf{Z})^*} e^{-\frac{2\pi ik\lambda d}{N}} = \operatorname{tr}_{\mathbf{Q}(\mu_{N/k})/\mathbf{Q}}(e^{\frac{2\pi ikd}{N}}) \\ &= \frac{\varphi(N/k)}{\varphi(N/(kd))} \operatorname{tr}_{\mathbf{Q}(\mu_{N/(kd)})/\mathbf{Q}}(\zeta_{N/(kd)}) = \frac{\varphi(N/k)}{\varphi(N/(kd))} \mu(N/(kd)), \end{aligned}$$

où $\mu$ est la fonction de Möbius [35, Chap 2, §2], puisque pour tout entier $n \geq 1$, le polynôme cyclotomique d'ordre $n$ débute par les termes $\Phi_n = X^{\varphi(n)} - \mu(n)X^{\varphi(n)-1} + \cdots$. Nous avons alors

$$\begin{aligned} L(f, \widehat{\epsilon_k}, s) &= \sum_{n=1}^{\infty} \frac{a_n(f)\widehat{\epsilon_k}(n)}{n^s} = \sum_{d \mid \frac{N}{k}} \sum_{\substack{n=1 \\ (n, \frac{N}{k})=d}}^{\infty} \frac{a_n(f)\widehat{\epsilon_k}(n)}{n^s} \\ &= \sum_{d \mid \frac{N}{k}} \widehat{\epsilon_k}(d) \sum_{\substack{n=1 \\ (n, \frac{N}{k})=d}}^{\infty} \frac{a_n(f)}{n^s} = \sum_{d \mid \frac{N}{k}} \widehat{\epsilon_k}(d) \sum_{n=1}^{\infty} \frac{a_{dn}(f) \cdot 1_{\frac{N}{kd}}(n)}{(dn)^s} \end{aligned}$$

où nous avons noté $1_{\frac{N}{kd}}$ le caractère trivial modulo $\frac{N}{kd}$. La fonction $1_{\frac{N}{kd}}$ est complètement multiplicative. En décomposant $dn$ en facteurs premiers, on est ramenés à calculer la somme suivante, pour $p$ premier divisant $N$ et $a \geq 0$

$$\sum_{b=0}^{\infty} \frac{a_{p^{a+b}}(f) \cdot 1_{\frac{N}{kd}}(p)^b}{p^{(a+b)s}} = \Big(\frac{a_p(f)}{p^s}\Big)^a \sum_{b=0}^{\infty} \Big(\frac{a_p(f) \cdot 1_{\frac{N}{kd}}(p)}{p^s}\Big)^b$$
$$= \Big(\frac{a_p(f)}{p^s}\Big)^a \frac{1}{1 - a_p(f) \cdot 1_{\frac{N}{kd}}(p) \cdot p^{-s}}.$$

Nous en déduisons

$$L(f, \widehat{\epsilon_k}, s) = \sum_{d \mid \frac{N}{k}} \widehat{\epsilon_k}(d) \cdot \Big( \prod_{p \mid N} \big(\frac{a_p(f)}{p^s}\big)^{v_p(d)} \Big) \cdot L(f, 1_{\frac{N}{kd}}, s)$$
$$= \sum_{d \mid \frac{N}{k}} \widehat{\epsilon_k}(d) \cdot \frac{a_d(f)}{d^s} \cdot L(f, 1_{\frac{N}{kd}}, s).$$

Il suit

$$L(f, \widehat{\epsilon_k}, 1) = \varphi(N/k) \sum_{d \mid \frac{N}{k}} \frac{\mu(\frac{N}{kd})}{\varphi(\frac{N}{kd})} \cdot \frac{a_d(f)}{d} \cdot L(f, 1_{\frac{N}{kd}}, 1).$$

La quantité $L(f, 1_{\frac{N}{kd}}, 1)$ s'obtient à partir de $L(f, 1)$ en supprimant les facteurs d'Euler aux nombres premiers $p$ divisant $\frac{N}{kd}$. Posons

$$R(n) = \prod_{p \mid n} L_p(f, p^{-1})^{-1} \neq 0 \qquad (n \geq 1)$$

où $L_p(f, p^{-s})$ est le facteur local en $p$ de $L(f, s)$ (voir l'appendice p. 143), de telle sorte que

$$L(f, 1_n, 1) = R(n) L(f, 1) \qquad (n \geq 1).$$

Lorsque $p$ divise $N$, nous avons $R(p) = 1 - \frac{a_p(f)}{p}$. Il vient finalement

$$L(f, \widehat{\epsilon_k}, 1) = \varphi(N/k) F(N/k) L(f, 1) \tag{3.140}$$

avec la notation

$$F(n) = \sum_{d \mid n} \frac{a_d(f)}{d} \cdot \frac{\mu(n/d) R(n/d)}{\varphi(n/d)} \qquad (n \geq 1).$$

Appliquons tous ces calculs à $k = 2$ et $k = 4$. Des égalités (3.136) à (3.140), nous déduisons (par hypothèse $L(f, 1) \neq 0$)

$$F(N/2) = 2 \quad \text{et} \quad F(N/4) = 1.$$

La fonction $F$ est faiblement multiplicative puisqu'elle est convolution de telles fonctions. Nous voyons sans difficulté que $F(2^a) = 0$ pour $a \geq 2$. Puisque $F(2N') \neq 0$, nous en déduisons que $N'$ est impair. La propriété de $F$ entraîne alors $F(2N') = F(2)F(N')$, d'où $F(2) = 2$, or la définition de $F$ donne $F(2) = \mu(2)R(2)/\varphi(2) = -1$, d'où la contradiction.

Le théorème 97 se trouve donc démontré. □

*Remarques.* 1. Le caractère inexplicite des appartenances $2\pi\xi_f^+(x) \in \mathbf{Q} \cdot \Omega_E^+$ et $L(f,\chi,1) \in \mathbf{Q}(\widehat{\chi}) \cdot \Omega_E^+$ rend difficile, dans un cadre général, l'explicitation de la combinaison linéaire donnant $L(E,2)$ en termes des quantités (3.127). Pour une courbe elliptique donnée, en revanche, des calculs utilisant les symboles modulaires permettent d'expliciter cette combinaison linéaire.

2. Que peut-on dire des dénominateurs des coefficients de la combinaison linéaire intervenant dans le théorème 97 ?

3. Il est peut-être possible d'améliorer la nature des nombres (3.127) intervenant dans l'énoncé du théorème 97. Par exemple, le lemme 89 entraîne que pour $a \neq 1$, la forme différentielle $\eta(1,a)$ n'est pas fermée. Le théorème reste-t-il vrai en considérant, par exemple, les quantités

$$\frac{1}{\pi i} \int_\gamma \eta(l,m), \tag{3.141}$$

où $\gamma$ parcourt les chemins fermés de $Y_1(N)(\mathbf{C})$, et $l, m$ parcourent les diviseurs de degré 0 sur $\frac{\mathbf{Z}}{N\mathbf{Z}}$ ?

4. Une réponse positive à la question précédente aurait la conséquence intéressante suivante. Pour des diviseurs $l$, $m$ de degré 0 sur $\frac{\mathbf{Z}}{N\mathbf{Z}}$, la forme différentielle $\eta(l,m)$ est fermée sur $Y_1(N)(\mathbf{C})$, donc admet une primitive $F_{l,m}$ sur $\mathcal{H}$. Soit $\gamma \in H_1(Y_1(N)(\mathbf{C}), \mathbf{Z})$ et $\widetilde{\gamma} \in \Gamma_1(N)$ tel que $\gamma$ soit la classe d'un chemin reliant $z_0$ à $\widetilde{\gamma}z_0$ dans $\mathcal{H}$, avec $z_0 \in \mathcal{H}$. Nous avons alors

$$\int_\gamma \eta(l,m) = \int_{z_0}^{\widetilde{\gamma}z_0} dF_{l,m} = F_{l,m}(\widetilde{\gamma}z_0) - F_{l,m}(z_0).$$

En particulier, la valeur spéciale $L(E,2)$ serait combinaison linéaire à coefficients rationnels de valeurs des fonctions $F_{l,m}$, en des points de $\mathcal{H}$ que l'on peut choisir dans la même orbite sous l'action de $\Gamma_1(N)$.

Nous spécialisons encore les résultats précédents aux courbes elliptiques de niveau premier, et montrons les théorèmes 1 et 2 de l'introduction. Soit donc $E$ une courbe elliptique sur $\mathbf{Q}$, de conducteur $p$ premier. Rappelons que nous notons $w(E)$ l'opposé du signe de l'équation fonctionnelle de $L(E,s)$. Pour $v \in \mathbf{Z}$, posons $g_v = \begin{pmatrix} 0 & -1 \\ 1 & v \end{pmatrix} \in \mathrm{SL}_2(\mathbf{Z})$. Pour $\chi$ caractère de Dirichlet modulo $p$, posons

$$\eta_\chi = \eta(\chi,\overline{\chi}) = \sum_{a \in (\frac{\mathbf{Z}}{p\mathbf{Z}})^*} \sum_{b \in (\frac{\mathbf{Z}}{p\mathbf{Z}})^*} \chi(a)\,\overline{\chi}(b)\,\eta(a,b).$$

**Théorème 1.** *Pour tout caractère de Dirichlet $\chi$ modulo $p$, pair et non trivial, nous avons la formule*

$$L(E,2)L(E,\chi,1) = \frac{p\,w(E)\,\tau(\chi)}{8\pi i(p-1)} \sum_{\chi'} c_{\chi,\chi'} L(E,\chi',1) \tag{3.142}$$

*où la somme est étendue aux caractères $\chi'$ modulo $p$, pairs et non triviaux, et les coefficients $c_{\chi,\chi'}$ sont donnés par*

$$c_{\chi,\chi'} = \tau(\overline{\chi'}) \sum_{v=1}^{p-1} \chi'(v) \int_{g_v\rho}^{g_v\rho^2} \eta_\chi. \tag{3.143}$$

*Démonstration.* Notons $f \in S_2(\Gamma_0(p))$ la forme primitive associée à $E$. Soit $\chi$ un caractère de Dirichlet modulo $p$, pair et non trivial. Au cours de la démonstration du théorème 3, nous avons obtenu l'expression suivante

$$L(f,2)L(f,\chi,1) = \frac{p}{2\pi i(p-1)} \int_{c(\chi,\widehat{\chi})} \omega_f.$$

Puisque $\chi$ est primitif, on a $\widehat{\chi} = \tau(\chi)\overline{\chi}$. Par définition de $c(\chi, \overline{\chi})$ il vient alors

$$L(f,2)L(f,\chi,1) = \frac{p\, i\, \tau(\chi)}{4(p-1)} \sum_{x \in E_p} \Big(\int_{g_x\rho}^{g_x\rho^2} \eta(\chi,\overline{\chi})\Big)\xi_f(x).$$

Les quantités $\int_{g_x\rho}^{g_x\rho^2} \eta_\chi$ et $\xi_f^+(x)$ ($x \in E_p$) sont invariantes par $x \mapsto \alpha x$, $\alpha \in (\frac{\mathbf{Z}}{p\mathbf{Z}})^*$. Il suit

$$L(f,2)L(f,\chi,1) = \frac{p\, i\, \tau(\chi)}{4} \sum_{x \in \mathbf{P}^1(\frac{\mathbf{Z}}{p\mathbf{Z}})} \Big(\int_{g_x\rho}^{g_x\rho^2} \eta_\chi\Big)\xi_f^+(x).$$

Soit $x = \infty = \frac{1}{0}$. Nous pouvons prendre $g_x = \begin{pmatrix} 0 & -1 \\ 1 & 0 \end{pmatrix}$. Remarquons que $\eta_\chi$ est fermée (lemme 89) et invariante par $T$. Il suit

$$\int_{g_x\rho}^{g_x\rho^2} \eta_\chi = \int_{\sigma\rho}^{\sigma\rho^2} \eta_\chi = \int_{\rho^2}^{T\rho^2} \eta_\chi = \int_{\infty}^{T\infty} \eta_\chi = 0.$$

De même, pour $x = 0$ on a $\int_{g_x\rho}^{g_x\rho^2} \eta_\chi = 0$. Il reste

$$L(f,2)L(f,\chi,1) = \frac{p\, i\, \tau(\chi)}{4} \sum_{x \in (\frac{\mathbf{Z}}{p\mathbf{Z}})^*} \Big(\int_{g_x\rho}^{g_x\rho^2} \eta_\chi\Big)\xi_f(x). \tag{3.144}$$

La formule suivante, obtenue au cours (p. 155) de la démonstration du théorème D de l'appendice, nous sera utile

$$\xi_f(x) = \frac{w(E)}{p-1} \sum_{\chi \neq 1_p} \frac{\tau(\overline{\chi})}{p} \chi(x) \Lambda(f \otimes \chi, 1) \qquad (x \in (\frac{\mathbf{Z}}{p\mathbf{Z}})^*),$$

la somme portant sur les caractères non triviaux modulo $p$. Pour $\chi$ modulo $p$ et non trivial, la forme primitive $f \otimes \chi$ est de niveau $p^2$, avec $a_p(f \otimes \chi) = 0$ [4]. On a donc

$$\Lambda(f \otimes \chi, 1) = \frac{p}{2\pi} L(f \otimes \chi, 1) = \frac{p}{2\pi} L(f, \chi, 1).$$

On en déduit

$$\xi_f(x) = \frac{w(E)}{2\pi(p-1)} \sum_{\chi \neq 1_p} \tau(\overline{\chi})\chi(x) L(f,\chi,1) \qquad (x \in (\frac{\mathbf{Z}}{p\mathbf{Z}})^*). \tag{3.145}$$

En reportant (3.145) dans (3.144), il vient

$$L(f,2)L(f,\chi,1) = -\frac{p\,w(E)\,\tau(\chi)}{8\pi i(p-1)} \sum_{x \in (\frac{\mathbf{Z}}{p\mathbf{Z}})^*} \Big(\int_{g_x\rho}^{g_x\rho^2} \eta_\chi\Big) \sum_{\chi' \neq 1_p} \tau(\overline{\chi'})\chi'(x)L(f,\chi',1).$$

où la seconde somme porte sur les caractères $\chi'$ non triviaux modulo $p$. Posons

$$c_{\chi,\chi'} = -\tau(\overline{\chi'}) \sum_{x \in (\frac{\mathbf{Z}}{p\mathbf{Z}})^*} \chi'(x)\Big(\int_{g_x\rho}^{g_x\rho^2} \eta_\chi\Big) \qquad (\chi' \neq 1_p),$$

de sorte que

$$L(f,2)L(f,\chi,1) = \frac{p\,w(E)\,\tau(\chi)}{8\pi i(p-1)} \sum_{\chi' \neq 1_p} c_{\chi,\chi'} L(f,\chi',1).$$

Il reste à montrer que $c_{\chi,\chi'}$ est nul pour $\chi'$ impair, et donné par (3.143) lorsque $\chi'$ est pair (non trivial). En effectuant le changement de variable $x \mapsto -x$ dans la somme définissant $c_{\chi,\chi'}$, il vient

$$\begin{aligned}
c_{\chi,\chi'} &= -\tau(\overline{\chi'}) \sum_{x \in (\frac{\mathbf{Z}}{p\mathbf{Z}})^*} \chi'(-x)\Big(\int_{g_{-x}\rho}^{g_{-x}\rho^2} \eta_\chi\Big) \\
&= -\chi'(-1)\tau(\overline{\chi'}) \sum_{x \in (\frac{\mathbf{Z}}{p\mathbf{Z}})^*} \chi'(x)\Big(\int_{c(g_x\rho^2)}^{c(g_x\rho)} \eta_\chi\Big) \\
&= -\chi'(-1)\tau(\overline{\chi'}) \sum_{x \in (\frac{\mathbf{Z}}{p\mathbf{Z}})^*} \chi'(x)\Big(\int_{g_x\rho^2}^{g_x\rho} c^*\eta_\chi\Big) \\
&= \chi'(-1)c_{\chi,\chi'}
\end{aligned}$$

puisque $c^*\eta_\chi = -\eta_\chi$ (cf. démonstration du théorème 3). Donc $c_{\chi,\chi'} = 0$ pour $\chi'$ impair. Supposons enfin $\chi'$ pair non trivial modulo $p$, et fixons $1 \leq v \leq p-1$. Pour $x = [v] \in (\frac{\mathbf{Z}}{p\mathbf{Z}})^*$ nous pouvons prendre $g_x = \begin{pmatrix} 1 & 0 \\ v & 1 \end{pmatrix} = g_{-v}\sigma$, d'où

$$c_{\chi,\chi'} = -\tau(\overline{\chi'})\sum_{v=1}^{p-1} \chi'(v) \int_{g_{-v}\sigma\rho}^{g_{-v}\sigma\rho^2} \eta_\chi = \tau(\overline{\chi'})\sum_{v=1}^{p-1} \chi'(v) \int_{g_{-v}\rho}^{g_{-v}\rho^2} \eta_\chi.$$

Puisque $g_{-v} \in \Gamma_1(p) \cdot g_{p-v}$, nous avons $\int_{g_{-v}\rho}^{g_{-v}\rho^2} \eta_\chi = \int_{g_{p-v}\rho}^{g_{p-v}\rho^2} \eta_\chi$. Le changement de variables $v \mapsto p-v$ donne alors

$$c_{\chi,\chi'} = \tau(\overline{\chi'})\sum_{v=1}^{p-1} \chi'(v) \int_{g_v\rho}^{g_v\rho^2} \eta_\chi$$

c'est-à-dire (3.143). $\square$

Considérons maintenant la série de Dirichlet

$$L(E \otimes E, s) = \sum_{n=1}^{\infty} \frac{a_n^2}{n^s} \qquad (\Re(s) > 2).$$

Elle admet un prolongement méromorphe au plan complexe [19] dont le résidu en $s = 2$ est non nul (cf. appendice, p.154).

**Théorème 2.** *Nous avons la formule*

$$L(E,2) = \frac{p^3\,w(E)}{8(p+1)(p-1)^3\,\pi^2} \frac{\sum_{\chi,\chi'} \lambda_{\chi,\chi'} L(E,\chi,1) L(E,\chi',1)}{\operatorname{Res}_{s=2} L(E \otimes E, s)} \tag{3.146}$$

*où la somme est étendue aux caractères $\chi$ (resp. $\chi'$) modulo $p$, pairs et non triviaux (resp. impairs), et les coefficients $\lambda_{\chi,\chi'}$ sont donnés par*

$$\lambda_{\chi,\chi'} = \sum_{\chi''} \frac{\tau(\chi'')}{\tau(\chi'\chi'')} c_{\chi'',\chi} \tag{3.147}$$

*la dernière somme portant sur les caractères $\chi''$ pairs non triviaux modulo $p$, les nombres $c_{\chi'',\chi}$ étant donnés par (3.143).*

*Démonstration.* Notons encore $f \in S_2(\Gamma_0(p))$ la forme primitive associée à $E$. Le théorème D de l'appendice (p. 154) exprime $\operatorname{Res}_{s=2} L(E \otimes E, s)$ en termes des valeurs $\Lambda(f \otimes \chi, 1)$. Joint à l'identité $\Lambda(f \otimes \chi, 1) = \frac{p}{2\pi} L(f, \chi, 1)$ pour $\chi$ caractère non trivial modulo $p$, il entraîne

$$\operatorname{Res}_{s=2} L(E \otimes E, s) = \frac{p^2\,i}{(p+1)(p-1)^2\,\pi} \sum_{\chi,\chi'} \frac{L(E,\chi,1)L(E,\chi',1)}{\tau(\chi\chi')} \tag{3.148}$$

où la somme est étendue aux caractères $\chi$ (resp. $\chi'$) modulo $p$, pairs et non triviaux (resp. impairs). L'idée consiste maintenant à multiplier (3.142) par $L(E,\chi',1)$ pour tout caractère $\chi'$ impair modulo $p$, puis à prendre une combinaison bilinéaire appropriée sur $\chi$ et $\chi'$, de manière à faire apparaître $\operatorname{Res}_{s=2} L(E \otimes E, s)$. Fixons $\chi$ (resp. $\chi'$) modulo $p$, pair et non trivial (resp. impair). Nous avons

$$L(E,2)L(E,\chi,1)L(E,\chi',1) = \frac{p\,w(E)\,\tau(\chi)}{8\pi i(p-1)} \sum_{\chi''} c_{\chi,\chi''} L(E,\chi'',1) L(E,\chi',1)$$

où la somme est étendue aux caractères $\chi''$ pairs non triviaux modulo $p$. En prenant la combinaison bilinéaire sur $\chi$ et $\chi'$ indiquée par (3.148), il vient

$$L(E,2)\operatorname{Res}_{s=2} L(E \otimes E, s) = \frac{p^3\,w(E)}{8(p+1)(p-1)^3\,\pi^2} \sum_{\chi'',\chi'} \Big(\sum_{\chi} \frac{\tau(\chi)}{\tau(\chi\chi')} c_{\chi,\chi''}\Big) L(E,\chi'',1) L(E,\chi',1),$$

d'où la formule (3.146). □

*Question.* Peut-on simplifier l'expression de $\lambda_{\chi,\chi'}$ ?

*Remarque* 100. En remplaçant $\operatorname{Res}_{s=2} L(E \otimes E, s)$ par son expression (3.148), il suit la formule

$$L(E,2) = \frac{p^2\,i\,w(E)}{8(p-1)\,\pi^2} \frac{\sum_{\chi,\chi'} \lambda_{\chi,\chi'} L(E,\chi,1) L(E,\chi',1)}{\sum_{\chi,\chi'} \tau(\overline{\chi\chi'}) L(E,\chi,1) L(E,\chi',1)} \tag{3.149}$$

En particulier, la valeur spéciale $L(E,2)$ se déduit mécaniquement des quantités $L(E,\chi,1)$, où $\chi$ parcourt les caractères de Dirichlet non triviaux modulo $p$, ce qui justifie la remarque de l'introduction.

## 3.7 Exemple en genre 1

Nous démontrons dans cette section le résultat suivant, annoncé dans l'introduction.

**Théorème 8.** *Soit $E$ la courbe elliptique donnée par l'équation $y^2 + y = x^3 - x^2$, et $P = (0,0)$, point d'ordre 5 de $E(\mathbf{Q})$. Orientons $E(\mathbf{R})$ dans le sens des $y$ croissants. Soit $\chi$ un caractère de Dirichlet modulo 11, pair et non trivial. Posons $\zeta = \chi(3) \in \mu_5$. Nous avons la formule*

$$L(E,2) = \frac{2^2 \cdot 5 \cdot \pi}{11^2} \cdot \frac{1 + 3(\zeta + \overline{\zeta})}{\zeta - \overline{\zeta}} \sum_{a \in \mathbf{Z}/5\mathbf{Z}} \zeta^a D_E(aP). \tag{3.150}$$

Nous déduisons du théorème 8 la preuve d'une identité conjecturée par Bloch et Grayson [13].

**Corollaire 101.** *En conservant les hypothèses du théorème 8, nous avons*

$$L(E,2) = \frac{10}{11} \cdot \pi \cdot D_E(P) \quad et \quad D_E(2P) = \frac{3}{2} D_E(P). \tag{3.151}$$

La seconde des identités (3.151), appelée *relation exotique*, a été démontrée récemment par Bertin [10].

*Démonstration du théorème 8.* La courbe elliptique $E$ ci-dessus n'est autre que la courbe modulaire $X_1(11)$. Notons $f$ l'unique forme primitive de poids 2 pour $\Gamma_1(11)$ et $\omega_f = 2\pi i f(z)dz$ la forme différentielle associée. Nous choisissons l'isomorphisme entre $X_1(11)$ et $E$ de telle sorte que la pointe infinie s'envoie vers 0 et $\omega_f$ s'identifie à la forme différentielle $dx/(2y+1)$. Nous avons alors

$$\Omega_f^+ := \int_{E(\mathbf{R})} \omega_f > 0$$

d'après le choix d'orientation de $E(\mathbf{R})$. O. Lecacheux m'a indiqué comment trouver les coordonnées $(x,y)$ des pointes $P_v$ pour $v \in (\mathbf{Z}/11\mathbf{Z})^*/\pm 1$. Ces pointes sont rationnelles et ont pour coordonnées

$$P_1 = \infty \quad P_2 = (1,0) \quad P_3 = (0,-1) \quad P_4 = (0,0) \quad P_5 = (1,-1). \tag{3.152}$$

Le groupe de Mordell-Weil $E(\mathbf{Q})$ est cyclique d'ordre 5, engendré par $P = P_4$ [20]. Parallèlement, le groupe $(\mathbf{Z}/11\mathbf{Z})^*/\pm 1$ est cyclique d'ordre 5, engendré par la classe de 4, et nous pouvons vérifier la formule

$$P_{4^a} = a \cdot P \qquad (a \in \mathbf{Z}/5\mathbf{Z}), \tag{3.153}$$

Cette formule résulte du fait que chaque opérateur diamant est nécessairement une translation rationnelle de la courbe elliptique $E$.

Utilisons le théorème 7 avec $\psi = 1$ et le caractère primitif $\chi$. Puisque $J_1(11)$ s'identifie canoniquement à $X_1(11)$, il vient

$$\begin{aligned} L(f,2)L(f,\chi,1) &= C_{1,\chi} \sum_{\lambda, \mu \in (\frac{\mathbf{Z}}{11\mathbf{Z}})^*/\pm 1} \overline{\chi}(\lambda)\chi(\mu) \left\langle R_{J_1(11)}(P_\lambda - P_\mu), \omega_f \right\rangle \\ &= C_{1,\chi} \sum_{a,b \in \frac{\mathbf{Z}}{5\mathbf{Z}}} \overline{\chi}(4)^a \chi(4)^b \left\langle R_{J_1(11)}\big((a-b)P\big), \omega_f \right\rangle. \end{aligned}$$

Puisque $\overline{\chi}(4) = \chi(3) = \zeta$, il suit

$$L(f,2)L(f,\chi,1) = 5\, C_{1,\chi} \sum_{a \in \frac{\mathbf{Z}}{5\mathbf{Z}}} \zeta^a \, \langle R_{J_1(11)}(aP), \omega_f \rangle.$$

On a $\omega_f = \Omega_f^+ \cdot \eta^* dz$ en utilisant les notations de la remarque 20. D'autre part $P \in E(\mathbf{R})$. D'après (2.4) et (1.64), nous en déduisons

$$\langle R_{J_1(11)}(aP), \omega_f \rangle = R_{\omega_f}(aP, 0) = \Omega_f^+ R_{\eta^* dz}(aP, 0) = \frac{i\Omega_f^+}{2} D_E(aP).$$

Il vient donc

$$L(f,2)L(f,\chi,1) = \frac{5i\Omega_f^+}{2} C_{1,\chi} \sum_{a \in \frac{\mathbf{Z}}{5\mathbf{Z}}} \zeta^a D_E(aP). \tag{3.154}$$

Des calculs utilisant les symboles modulaires pour $\Gamma_0(11)$ montrent que

$$L(f,\chi,1) = \pm \frac{\alpha}{5\tau(\overline{\chi})} \cdot \Omega_f^+ \quad \text{avec} \quad \alpha = \zeta - \overline{\zeta} + 2(\zeta^2 - \overline{\zeta}^2).$$

Le signe dans cette dernière égalité est déterminé en évaluant numériquement $L(f,\chi,1)$, par exemple grâce au package ComputeL [26] écrit pour le logiciel Pari (nous ne connaissons pas d'autre méthode). Nous trouvons le signe moins. D'autre part, le calcul explicite de $C_{1,\chi}$ donne

$$\begin{aligned}C_{1,\chi} &= \frac{11\pi i\,\tau(\chi)}{10} \frac{L(\chi,2)L(\overline{\chi},2)}{\pi^4} \\ &= \frac{11\pi i\,\tau(\chi)}{10} \left(\frac{2}{11}\right)^4 (7 - \zeta - \overline{\zeta}).\end{aligned}$$

En reportant ces expressions dans (3.154), il vient

$$\begin{aligned}L(f,2) &= -\frac{5\tau(\overline{\chi})}{\alpha} \cdot \frac{5i}{2} \cdot \frac{11\pi i \tau(\chi)}{10} \left(\frac{2}{11}\right)^4 (7 - \zeta - \overline{\zeta}) \sum_{a \in \frac{\mathbf{Z}}{5\mathbf{Z}}} \zeta^a D_E(aP) \\ &= \frac{2^2 \cdot 5 \cdot \pi}{11^2} \cdot \frac{7 - \zeta - \overline{\zeta}}{\alpha} \sum_{a \in \frac{\mathbf{Z}}{5\mathbf{Z}}} \zeta^a D_E(aP)\end{aligned}$$

puisque $\tau(\chi)\tau(\overline{\chi}) = 11$. Posons $\eta = \zeta + \overline{\zeta} \in \mathbf{Q}(\frac{1+\sqrt{5}}{2})$. Nous avons $\eta^2 = 1 - \eta$ d'où

$$\frac{7 - \zeta - \overline{\zeta}}{\alpha} = \frac{7 - \eta}{(\zeta - \overline{\zeta})(1 + 2\eta)} = \frac{1 + 3\eta}{\zeta - \overline{\zeta}}.$$

Puisque $L(E,2) = L(f,2)$ on obtient (3.150). $\square$

*Démonstration du corollaire 101.* Posons $\eta = \zeta + \overline{\zeta}$ comme dans la démonstration du théorème 8. D'après (1.42), nous avons $D_E(4P) = -D_E(P)$, $D_E(3P) = -D_E(2P)$ et $D_E(0) = 0$. De (3.150) nous déduisons

$$\begin{aligned}L(E,2) &= \frac{2^2 \cdot 5 \cdot \pi}{11^2} \cdot \frac{7 - \eta}{\zeta - \overline{\zeta}} \big((\zeta - \overline{\zeta})D_E(P) + (\zeta^2 - \overline{\zeta}^2)D_E(2P)\big) \\ &= \frac{2^2 \cdot 5 \cdot \pi}{11^2} \cdot (7 - \eta)\big(D_E(P) + \eta D_E(2P)\big). \tag{3.155}\end{aligned}$$

Cette identité est valable pour tout $\zeta$ racine primitive 5-ième de l'unité. Soit $\sigma \in \mathrm{Gal}(\mathbf{Q}(\mu_5)/\mathbf{Q})$ l'automorphisme défini par $\sigma(\zeta) = \zeta^2$ (noter que la définition de $\sigma$ est indépendante de $\zeta$). Puisque (3.155) reste vraie si l'on remplace $\eta$ par $\sigma(\eta) = -\eta - 1$, il vient également

$$L(E,2) = \frac{2^2 \cdot 5 \cdot \pi}{11^2} \cdot (8+\eta)\bigl(D_E(P) - (1+\eta)D_E(2P)\bigr). \tag{3.156}$$

En éliminant $\eta$ des équations (3.155) et (3.156), nous obtenons

$$L(E,2) = \frac{2 \cdot 5 \cdot \pi}{11^2} \cdot \bigl(-D_E(P) + 8D_E(2P)\bigr). \tag{3.157}$$

D'autre part, en écrivant l'égalité des membres de droite de (3.155) et (3.156), il vient

$$(7-\eta)\bigl(D_E(P) + \eta D_E(2P)\bigr) = (8+\eta)\bigl(D_E(P) - (1+\eta)D_E(2P)\bigr)$$

ce qui montre, après simplifications, la seconde des identités (3.151). En reportant cette dernière dans (3.157), nous obtenons la première identité. □

*Remarque* 102. Un aspect intéressant est que l'on voit assez clairement au cours du calcul d'où viennent les facteurs premiers 5 et 11 de l'identité (3.151). Le nombre premier 11 est le niveau de $f$ (donc aussi le conducteur de $E$), tandis que le nombre premier 5 est l'ordre du groupe $(\mathbf{Z}/11\mathbf{Z})^*/\pm 1$ (qui est isomorphe au groupe de Mordell-Weil de $E$).

*Remarque* 103. La courbe $X_1(N)$ est de genre 1 pour deux autres valeurs de $N$, à savoir $N = 14$ et $N = 15$. Les méthodes utilisées ici devraient s'adapter sans difficulté à ces cas-là.

## 3.8    Exemple en genre 2

Nous considérons ici la courbe $X_1(13)$, qui est de genre 2. Nous proposons deux formules pour $L(f_\epsilon, 2)$, où $f_\epsilon$ est une forme primitive de poids 2 pour $\Gamma_1(13)$. La première formule utilise les résultats de la section 3.5 et fait intervenir l'intégrale de la forme différentielle $\eta(x,y)$ associée aux fonctions $x$ et $y$ relatives à un modèle convenable de $X_1(13)$. Nous avons constaté que l'utilisation de l'équation fonctionnelle ($L'(f_\epsilon, 0)$ au lieu de $L(f_\epsilon, 2)$) amène à une formule plus naturelle. La seconde formule utilise la théorie développée au chapitre 2, et exprime $L(f_\epsilon, 2)$ en termes de la fonction $R_{J_1(13)(\mathbf{C})}$, évaluée aux points rationnels de $J_1(13)$.

Un modèle plan de $X_1(13)$ a été déterminé par Lecacheux [42, pp. 55-57]. Nous noterons $P_v = [0, v]$ pour $v \in (\mathbf{Z}/13\mathbf{Z})^*/\pm 1$, les pointes de $X_1(13)(\mathbf{C})$ situées au-dessus de la pointe infinie de $X_0(13)(\mathbf{C})$. Soit $W_{13}$ l'involution d'Atkin-Lehner de $X_1(13)$. Les pointes $P_v$ sont rationnelles sur $\mathbf{Q}$, et nous notons $a_1\,a_2\,\ldots\,a_6$ le diviseur $\sum_{i=1}^{6} a_i P_i$ (attention, ces notations diffèrent de [42] par l'application de $W_{13}$). Dans l'article cité précédemment, l'équation suivante est donnée pour $X_1(13)$

$$H^2 h(h-1) + H(-h^3 + h^2 + 2h - 1) - h^2 + h = 0. \tag{3.158}$$

Posant $x = h|W_{13}$ et $y = H|W_{13}$, les fonctions rationnelles $x$ et $y$ sont définies sur $\mathbf{Q}$ et à support dans les pointes $P_v$. En appliquant $W_{13}$ à (3.158), nous trouvons que $X_1(13)$ est birationnelle à la courbe plane de même équation, soit en coordonnées projectives

$$y^2 x(x-z) + y(-x^3 + x^2 z + 2xz^2 - z^3) - xz^2(x-z) = 0. \tag{3.159}$$

Noter que la courbe définie par (3.159) est singulière. Le point de coordonnées projectives $(0 : 1 : 0)$ est un point double ordinaire et la courbe $X_1(13)$ s'identifie à la désingularisée de (3.159) en ce point. Dans la carte $y = 1$, les directions tangentes au point $(0 : 1 : 0)$ sont données par $x = z$ et $x = 0$. Pour la seconde direction, on a $x \sim -z^2$ au voisinage de $(x, z) = (0, 0)$. Les diviseurs des fonctions $x$ et $y$ sont donnés par

$$\begin{array}{rrrrrrr} \operatorname{div} x = & 0 & 1 & 1 & -1 & 0 & -1 \\ \operatorname{div} y = & 1 & -1 & 1 & 1 & -1 & -1. \end{array} \quad (3.160)$$

Nous en déduisons les coordonnées projectives des pointes $P_v$ :

$$\begin{array}{lll} P_1 = (1 : 0 : 1) & P_2 = (0 : 1 : 0)_{x=0} & P_3 = (0 : 0 : 1) \\ P_4 = (1 : 0 : 0) & P_5 = (0 : 1 : 0)_{x=z} & P_6 = (1 : 1 : 0). \end{array} \quad (3.161)$$

**Proposition 104.** *L'élément $\{x, y\} \otimes 1$ appartient à $K_2(X_1(13)) \otimes \mathbf{Q}$.*

Nous donnons deux démonstrations de ce résultat.

*Démonstration directe.* Nous avons a priori $\{x, y\} \in K_2(\mathbf{Q}(X_1(13)))$. D'après la suite exacte (3.86), il suffit de montrer la trivialité des symboles modérés de $\{x, y\}$ aux pointes $P_v$, $1 \leq v \leq 6$. Cela se fait directement grâce à (3.160) et (3.161), et nous ne détaillons pas les calculs. □

*Démonstration utilisant les résultats de ce chapitre.* Cette méthode est plus intéressante car susceptible de généralisation. Les fonctions $x$ et $y$ sont des unités modulaires pour $\Gamma_1(13)$, à support dans les pointes $P_v$, $v \in (\mathbf{Z}/13\mathbf{Z})^* / \pm 1$. La proposition résultera donc de la proposition 86, à condition de montrer les égalités $\widehat{x}(\infty) = \widehat{y}(\infty) = 1$, où nous utilisons la notation (3.61). Il s'agit donc d'étudier le comportement des unités modulaires $x$ et $y$ au voisinage de la pointe infinie $P_1$. Pour l'unité modulaire $x$, l'égalité proposée est immédiate puisque $x(P_1) = 1$. L'unité modulaire $y$ possède en revanche un zéro simple en la pointe $P_1$, d'après (3.160), ce qui nous oblige à calculer le développement de Fourier de $y$ en cette pointe. Notons $q = e^{2\pi i z}$, $z \in \mathcal{H}$ la variable modulaire. En utilisant [42, pp. 56-57], on exprime $y$ en termes de la fonction $\wp$ de Weierstraß, puis on calcule le développement de Fourier en $q$ grâce à [71, I.6.2a], ce qui donne

$$y(q) = -q + O(q^2).$$

Nous avons donc $\widehat{y}(P_1) = -1$. D'après la proposition 86, nous avons $\{x, y^2\} \in K_2(X_1(13)) \otimes \mathbf{Q}$. Il en va donc de même de $\{x, y\} = \{x, y^2\} \otimes \frac{1}{2}$. □

*Remarque* 105. D'après les résultats de Schappacher et Scholl [62, 7.3.1], nous avons même $\{x, y\} \in K_2(X_1(13))_\mathbf{Z} \otimes \mathbf{Q}$.

Le groupe $(\mathbf{Z}/13\mathbf{Z})^* / \pm 1$ est cyclique d'ordre 6, engendré par la classe de 2. La donnée d'un caractère de Dirichlet $\chi$ pair modulo 13 équivaut à celle d'une racine sixième de l'unité $\chi(2) \in \mu_6$. Nous poserons $\zeta_6 = e^{\frac{2\pi i}{6}}$ et noterons $\epsilon$ le caractère de Dirichlet pair modulo 13 défini par $\epsilon(2) = \zeta_6$. Les caractères de Dirichlet pairs modulo 13 sont 1, $\epsilon$, $\epsilon^2$, $\epsilon^3$, $\overline{\epsilon}^2$ et $\overline{\epsilon}$ ; ils sont d'ordres respectifs 1, 6, 3, 2, 3 et 6. L'espace $S_2(\Gamma_1(13))$ est de dimension 2 sur $\mathbf{C}$, une base étant donnée par les formes primitives $f_\epsilon$ et $f_{\overline{\epsilon}}$, de caractères respectifs $\epsilon$ et $\overline{\epsilon}$. Nous allons maintenant calculer le régulateur (3.81) de l'élément $\{x, y\}$. Posons

$$r_\epsilon = \langle r_{13}(\{x, y\}), \omega_{f_\epsilon} \rangle \quad \text{et} \quad r_{\overline{\epsilon}} = \langle r_{13}(\{x, y\}), \omega_{f_{\overline{\epsilon}}} \rangle. \quad (3.162)$$

**Théorème 106.** *Nous avons*

$$r_\epsilon = \frac{13\sqrt{13}}{4\pi i}(1 + \zeta_6) \cdot L(f_\epsilon, 2)L(f_\epsilon, \epsilon^2, 1) \quad \text{et} \quad r_{\overline{\epsilon}} = -\overline{r_\epsilon}. \quad (3.163)$$

Pour démontrer ce théorème, nous utiliserons le lemme suivant.

**Lemme 107.** *Dans $\mathcal{O}^*(Y_1(13)) \otimes \mathcal{O}^*(Y_1(13)) \otimes \mathbf{C}$, nous avons l'identité*

$$x \otimes y = -\frac{13^2\sqrt{13}}{2^4 \cdot 3}\big((1+\zeta_6)(u_{\widehat{\epsilon}^2} \otimes u_{\epsilon^3}) + (2-\zeta_6)(u_{\widehat{\overline{\epsilon}}^2} \otimes u_{\epsilon^3})\big). \tag{3.164}$$

*Démonstration.* Le principe est simple : décomposer les diviseurs (3.160) de $x$ et $y$ comme combinaisons linéaires de caractères de Dirichlet, puis appliquer la proposition 79. Nous avons donc besoin de calculer $L(\chi, 2)/\pi^2$ pour $\chi$ caractère modulo 13, pair et non trivial. Posons [18, (1.81)]

$$B_{2,\chi} = 13 \sum_{a \in \frac{\mathbf{Z}}{13\mathbf{Z}}} \chi(a)\, \overline{B_2}(\frac{\widetilde{a}}{13}) \qquad (\chi \neq 1).$$

L'équation fonctionnelle [18, (3.87)] pour $L(\chi, s)$ et l'expression [18, (1.80)] donnent alors

$$\frac{L(\chi, 2)}{\pi^2} = -\frac{2\tau(\chi)}{13^2} L(\overline{\chi}, -1) = \frac{\tau(\chi)}{13^2} B_{2,\overline{\chi}} = \frac{B_{2,\overline{\chi}}}{13\tau(\overline{\chi})} \qquad (\chi \neq 1), \tag{3.165}$$

compte tenu de $\tau(\chi)\tau(\overline{\chi}) = 13$. Grâce à la proposition 79, nous avons

$$\operatorname{div} u_{\epsilon^3} = -\frac{L(\epsilon^3, 2)}{\pi^2} \cdot (\begin{array}{cccccc} 1 & -1 & 1 & 1 & -1 & -1 \end{array}) = -\frac{4\sqrt{13}}{13^2} \operatorname{div} y.$$

En considérant $u_{\epsilon^3}$ et $y$ comme des éléments de $\mathcal{O}^*(Y_1(13)) \otimes \mathbf{C}$, nous avons $\widehat{u_{\epsilon^3}}(\infty) = \widehat{y}(\infty) = 1$ (voir la démonstration de la proposition 104), d'où

$$y \otimes 1 = u_{\epsilon^3} \otimes -\frac{13\sqrt{13}}{4}. \tag{3.166}$$

Le diviseur de $x$ étant de degré 0 et invariant par l'opérateur diamant $\langle 5 \rangle$, il est combinaison linéaire de $\operatorname{div} u_{\epsilon^2}$ et $\operatorname{div} u_{\overline{\epsilon}^2}$. Nous trouvons explicitement

$$\operatorname{div} x = \frac{1-2\zeta_6}{3}\Big(\frac{\operatorname{div} u_{\epsilon^2}}{L(\epsilon^2, 2)/\pi^2} - \frac{\operatorname{div} u_{\overline{\epsilon}^2}}{L(\overline{\epsilon}^2, 2)/\pi^2}\Big)$$

$$= \frac{13}{12}\big((2-\zeta_6)\tau(\overline{\epsilon}^2)\operatorname{div} u_{\epsilon^2} + (1+\zeta_6)\tau(\epsilon^2)\operatorname{div} u_{\overline{\epsilon}^2}\big).$$

D'après la proposition 80, nous avons $u_{\widehat{\overline{\epsilon}}^2} = u_{\epsilon^2} \otimes \tau(\overline{\epsilon}^2)$ et $u_{\widehat{\epsilon}^2} = u_{\overline{\epsilon}^2} \otimes \tau(\epsilon^2)$, d'où

$$\operatorname{div} x = \frac{13}{12}\big((2-\zeta_6)\operatorname{div} u_{\widehat{\overline{\epsilon}}^2} + (1+\zeta_6)\operatorname{div} u_{\widehat{\epsilon}^2}\big).$$

De $\widehat{x}(\infty) = 1$, nous déduisons comme précédemment

$$x \otimes 1 = \frac{13}{12}\big(u_{\widehat{\epsilon}^2} \otimes (1+\zeta_6) + u_{\widehat{\overline{\epsilon}}^2} \otimes (2-\zeta_6)\big). \tag{3.167}$$

Le lemme résulte alors de (3.166) et (3.167). □

*Démonstration du théorème 106.* D'après le lemme précédent, nous avons

$$r_\epsilon = -\frac{13^2\sqrt{13}}{2^4 \cdot 3}\langle r_{13}\big((1+\zeta_6)\{u_{\widehat{\epsilon}^2}, u_{\epsilon^3}\} + (2-\zeta_6)\{u_{\widehat{\overline{\epsilon}}^2}, u_{\epsilon^3}\}\big), \omega_{f_\epsilon}\rangle$$

$$= \frac{13^2\sqrt{13}}{2^4 \cdot 3}\langle r_{13}\big((1+\zeta_6)\{u_{\epsilon^3}, u_{\widehat{\epsilon}^2}\} + (2-\zeta_6)\{u_{\epsilon^3}, u_{\widehat{\overline{\epsilon}}^2}\}\big), \omega_{f_\epsilon}\rangle.$$

Utilisons le théorème 81 avec $N = 13$, $f = f_\epsilon$, $\chi = \epsilon^3$ et $\chi' = \epsilon^2$ ou $\bar\epsilon^2$. Pour des raisons de caractère, le terme correspondant à $\chi' = \bar\epsilon^2$ disparaît. Il vient

$$r_\epsilon = \frac{13^2\sqrt{13}}{2^4 \cdot 3}(1 + \zeta_6) \cdot \frac{12}{13\pi i} \cdot L(f_\epsilon, 2)L(f_\epsilon, \epsilon^2, 1)$$
$$= \frac{13\sqrt{13}}{4\pi i}(1 + \zeta_6) \cdot L(f_\epsilon, 2)L(f_\epsilon, \epsilon^2, 1).$$

L'identité $r_{\bar\epsilon} = -\overline{r_\epsilon}$ résulte d'un calcul analogue à celui effectué au cours de la démonstration du lemme 83. □

Nous allons maintenant donner une expression de $L'(f_\epsilon, 0)$ en termes d'intégrale de la forme différentielle $\eta(x, y)$ associée aux fonctions $x$ et $y$. Noter que l'équation fonctionnelle de $L(f_\epsilon, s)$ permet d'exprimer $L(f_\epsilon, 2)$ en termes de $L'(f_\epsilon, 0)$. Nous avons une suite exacte

$$0 \to \mathbf{Z} \to \mathrm{Div}\, P_{13} \to H_1(Y_1(13)(\mathbf{C}), \mathbf{Z}) \to H_1(X_1(13)(\mathbf{C}), \mathbf{Z}) \to 0, \tag{3.168}$$

où nous rappelons que $\mathrm{Div}\, P_{13}$ désigne le groupe des diviseurs sur $P_{13}$. L'application $\mathrm{Div}\, P_{13} \to H_1(Y_1(13)(\mathbf{C}), \mathbf{Z})$ envoie une pointe $P \in P_{13}$ vers un lacet $\gamma_P$ autour du point $P$, orienté dans le sens trigonométrique. D'après (3.109) et la proposition 104, nous avons

$$\int_{\gamma_P} \eta(x, y) = 0 \qquad (P \in P_{13}). \tag{3.169}$$

L'intégration de la forme différentielle $\eta(x, y)$ induit donc une application de $H_1(X_1(13)(\mathbf{C}), \mathbf{Z})$ vers $\mathbf{C}$, qui s'étend par linéarité en un morphisme

$$\int \eta(x, y) : H_1(X_1(13)(\mathbf{C}), \mathbf{C}) \to \mathbf{C}. \tag{3.170}$$

Noter que ce morphisme est nul sur $H_1^+(X_1(13)(\mathbf{C}), \mathbf{C})$, d'après la propriété $c^*\eta(x, y) = -\eta(x, y)$.

**Définition 108.** *Pour tout caractère de Dirichlet $\psi$ pair modulo* 13, *posons*

$$H^\pm(\psi) = \{\gamma \in H_1^\pm(X_1(13)(\mathbf{C}), \mathbf{C}),\ \langle d \rangle_* \gamma = \psi(d)\gamma \quad \forall d \in (\mathbf{Z}/13\mathbf{Z})^*\}. \tag{3.171}$$

**Lemme 109.** *Pour tout caractère $\psi \in \{\epsilon, \bar\epsilon\}$, l'espace vectoriel complexe $H^\pm(\psi)$ est de dimension* 1. *Des générateurs de $H^+(\epsilon)$ et $H^-(\epsilon)$ sont donnés respectivement par*

$$\gamma_\epsilon^+ = 2(1 - 2\zeta_6)\xi(1, 2)^\epsilon + \xi(1, 3)^\epsilon + \xi(1, -3)^\epsilon \tag{3.172}$$
$$\gamma_\epsilon^- = \xi(1, 3)^\epsilon - \xi(1, -3)^\epsilon. \tag{3.173}$$

*Des générateurs de $H^+(\bar\epsilon)$ et $H^-(\bar\epsilon)$ sont donnés respectivement par*

$$\gamma_{\bar\epsilon}^+ = \overline{\gamma_\epsilon^+} = 2(2\zeta_6 - 1)\xi(1, 2)^{\bar\epsilon} + \xi(1, 3)^{\bar\epsilon} + \xi(1, -3)^{\bar\epsilon} \tag{3.174}$$
$$\gamma_{\bar\epsilon}^- = \overline{\gamma_\epsilon^-} = \xi(1, 3)^{\bar\epsilon} - \xi(1, -3)^{\bar\epsilon}. \tag{3.175}$$

*Démonstration.* L'intégration induit un accouplement parfait

$$H^\pm(\psi) \times S_2(\Gamma_1(13), \psi) \to \mathbf{C}.$$

L'espace $S_2(\Gamma_1(13),\psi)$ étant de dimension 1 pour $\psi = \epsilon$ ou $\bar\epsilon$, il en va de même pour $H^\pm(\psi)$. Au cours de la démonstration du théorème 5 (p. 98), nous avons vu que l'espace $H^+(\epsilon)$ est engendré par les cycles $\xi^\epsilon_\chi$, $\chi \neq 1, \epsilon$ définis en (3.97). Nous pouvons calculer $\xi^\epsilon_{\epsilon^3}$ grâce à [55, Lemme 5]

$$\xi^\epsilon_{\epsilon^3} = 2(1-2\zeta_6)\xi(1,2)^\epsilon + \xi(1,3)^\epsilon + \xi(1,-3)^\epsilon,$$

ce qui montre en outre $\xi^\epsilon_{\epsilon^3} \neq 0$. Nous pouvons donc choisir $\gamma^+_\epsilon = \xi^\epsilon_{\epsilon^3}$. Passons à l'espace $H^-(\epsilon)$. L'élément $\gamma^-_\epsilon$ est anti-invariant par la conjugaison complexe, et de caractère $\epsilon$. D'après (3.101), son bord est nul, d'où $\gamma^-_\epsilon \in H^-(\epsilon)$. À nouveau [55, Lemme 5] montre que $\gamma^-_\epsilon \neq 0$. Enfin, les générateurs de $H^+(\bar\epsilon)$ et $H^-(\bar\epsilon)$ sont obtenus en remarquant que la conjugaison complexe sur les coefficients $\gamma \mapsto \bar\gamma$ définit un isomorphisme de $H^\pm(\epsilon)$ vers $H^\pm(\bar\epsilon)$. □

La propriété $\overline{\eta(x,y)} = \eta(x,y)$ entraîne que le morphisme (3.170) est caractérisé par sa valeur en $\gamma^-_\epsilon$. Pour résumer, nous avons

$$\int_{\gamma^+_\epsilon} \eta(x,y) = \int_{\gamma^+_{\bar\epsilon}} \eta(x,y) = 0 \quad \text{et} \quad \int_{\gamma^-_{\bar\epsilon}} \eta(x,y) = \overline{\int_{\gamma^-_\epsilon} \eta(x,y)}. \tag{3.176}$$

**Théorème 110.** *Nous avons l'identité*

$$\int_{\gamma^-_\epsilon} \eta(x,y) = 4\pi\zeta_6 L'(f_\epsilon, 0), \tag{3.177}$$

*où $L'(f_\epsilon, 0)$ désigne la dérivée de $L(f_\epsilon, s)$ en $s=0$.*

*Démonstration.* D'après le lemme 107, nous avons

$$\eta(x,y) = \frac{13^2\sqrt{13}}{2^4 \cdot 3}\left((1+\zeta_6)\eta(u_{\epsilon^3}, u_{\widehat{\epsilon^2}}) + (2-\zeta_6)\eta(u_{\epsilon^3}, u_{\widehat{\bar\epsilon^2}})\right).$$

D'après (3.74) et (3.75), la forme différentielle $\eta(u_{\epsilon^3}, u_{\widehat{\epsilon^2}})$ (resp. $\eta(u_{\epsilon^3}, u_{\widehat{\bar\epsilon^2}})$) est de caractère $\bar\epsilon$ (resp. $\epsilon$). Il vient donc

$$\int_{\gamma^-_{\bar\epsilon}} \eta(x,y) = \frac{13^2\sqrt{13}}{2^4 \cdot 3}(1+\zeta_6)\int_{\gamma^-_{\bar\epsilon}} \eta(u_{\epsilon^3}, u_{\widehat{\epsilon^2}}). \tag{3.178}$$

Utilisons le théorème 3 avec la forme primitive $f = f_\epsilon$ et le caractère $\chi = \epsilon^2$. Au cours (3.125) de la démonstration de ce théorème, nous avons obtenu

$$\begin{aligned}
L(f_\epsilon, 2)L(f_\epsilon, \epsilon^2, 1) &= \frac{13}{24\pi i}\int_{c(\epsilon^3, \widehat{\epsilon^2})} \omega_{f_\epsilon} \\
&= \frac{13i}{24}\sum_{x \in E_{13}/\pm 1}\left(\int_{g_x\rho}^{g_x\rho^2} \eta(\epsilon^3, \widehat{\epsilon^2})\right)\xi_{f_\epsilon}(x) \\
&= -\frac{13\pi^2}{24}\sum_{x \in E_{13}/\pm 1}\left(\int_{g_x\rho}^{g_x\rho^2} \eta(u_{\epsilon^3}, u_{\widehat{\epsilon^2}})\right)\xi_{f_\epsilon}(x) \\
&= -\frac{13\pi^2}{24}\int_{a_{\bar\epsilon}} \eta(u_{\epsilon^3}, u_{\widehat{\epsilon^2}}),
\end{aligned}$$

où nous avons noté $a_{\bar\epsilon}$ le cycle relatif suivant de $Y_1(13)(\mathbf{C})$

$$a_{\bar\epsilon} = \sum_{x \in E_{13}/\pm 1} \xi_{f_\epsilon}(x)\{g_x\rho, g_x\rho^2\}. \tag{3.179}$$

Rappelons que la conjugaison complexe envoie $x = (u, v) \in E_{13}$ sur $x^c = (-u, v)$. Un calcul simple montre que

$$c_* a_{\bar\epsilon} = \sum_{x \in E_{13}/\pm 1} \xi_{f_\epsilon}(x) \{g_{x^c}\rho^2, g_{x^c}\rho\}$$

$$= -\sum_{x \in E_{13}/\pm 1} \xi_{f_\epsilon}(x^c) \{g_x\rho, g_x\rho^2\}.$$

Définissant $b_{\bar\epsilon} = (a_{\bar\epsilon} - c_* a_{\bar\epsilon})/2$, nous avons alors

$$b_{\bar\epsilon} = \sum_{x \in E_{13}/\pm 1} \xi^+_{f_\epsilon}(x) \{g_x\rho, g_x\rho^2\}.$$

Puisque $c^* \eta(u_{\epsilon^3}, u_{\widehat{\epsilon^2}}) = -\eta(u_{\epsilon^3}, u_{\widehat{\epsilon^2}})$, nous en tirons

$$L(f_\epsilon, 2)L(f_\epsilon, \epsilon^2, 1) = -\frac{13\pi^2}{24} \int_{-c_* a_{\bar\epsilon}} \eta(u_{\epsilon^3}, u_{\widehat{\epsilon^2}}) = -\frac{13\pi^2}{24} \int_{b_{\bar\epsilon}} \eta(u_{\epsilon^3}, u_{\widehat{\epsilon^2}}). \tag{3.180}$$

Remarquons que le cycle $b_{\bar\epsilon}$ est de caractère $\bar\epsilon$, et anti-invariant par la conjugaison complexe $c_*$. À l'aide du lemme suivant, nous allons maintenant montrer que le cycle $b_{\bar\epsilon}$ est fermé et l'exprimer explicitement en termes de $\gamma^-_{\bar\epsilon}$. Rappelons que les matrices $\sigma, \tau \in \mathrm{SL}_2(\mathbf{Z})$ sont définies par

$$\sigma = \begin{pmatrix} 0 & -1 \\ 1 & 0 \end{pmatrix} \quad \text{et} \quad \tau = \begin{pmatrix} 0 & -1 \\ 1 & -1 \end{pmatrix}.$$

**Lemme 111.** *Soit $\gamma$ un cycle relatif de $Y_1(N)(\mathbf{C})$ de la forme*

$$\gamma = \sum_{x \in E_N} \alpha_x \{g_x\rho, g_x\rho^2\} \tag{3.181}$$

*avec les relations*

$$\alpha_x + \alpha_{x\sigma} = 0 \quad et \quad \alpha_x + \alpha_{x\tau} + \alpha_{x\tau^2} = 0 \qquad (x \in E_N). \tag{3.182}$$

*Alors $\gamma$ est fermé et possède l'expression suivante en termes de symboles de Manin*

$$\gamma = -\frac{1}{3} \sum_{x \in E_N} (\alpha_x + 2\alpha_{x\tau}) \xi(x). \tag{3.183}$$

*Démonstration du lemme.* Le bord de $\gamma$ est donné par

$$\partial \gamma = \sum_{x \in E_N} \alpha_x \left([g_x\rho^2] - [g_x\rho]\right)$$

$$= \sum_{x \in E_N} \alpha_x \left([g_{x\sigma}\rho] - [g_x\rho]\right)$$

$$= \sum_{x \in E_N} \alpha_{x\sigma} [g_x\rho] - \sum_{x \in E_N} \alpha_x [g_x\rho]$$

$$= -2 \sum_{x \in E_N} \alpha_x [g_x\rho].$$

Puisque la matrice $\tau$ fixe $\rho$, nous obtenons

$$\partial \gamma = -2 \sum_{x \in E_N} \alpha_x \cdot \frac{1}{3} \big( [g_x \rho] + [g_{x\tau} \rho] + [g_{x\tau^2} \rho] \big)$$
$$= -\frac{2}{3} \sum_{x \in E_N} (\alpha_x + \alpha_{x\tau} + \alpha_{x\tau^2})[g_x \rho] = 0.$$

D'autre part, nous avons

$$\gamma = \sum_{x \in E_N} \alpha_x \big( \{g_x \rho, g_x \infty\} + \{g_x \infty, g_x \rho^2\} \big)$$
$$= \sum_{x \in E_N} \alpha_x \{g_x \rho, g_x \infty\} + \sum_{x \in E_N} \alpha_x \{g_{x\sigma} 0, g_{x\sigma} \rho\}$$
$$= \sum_{x \in E_N} \alpha_x \{g_x \rho, g_x \infty\} - \sum_{x \in E_N} \alpha_x \{g_x 0, g_x \rho\},$$

d'après $\alpha_x + \alpha_{x\sigma} = 0$. En écrivant $\{g_x 0, g_x \rho\} = \xi(x) + \{g_x \infty, g_x \rho\}$, il suit

$$\gamma = 2 \sum_{x \in E_N} \alpha_x \{g_x \rho, g_x \infty\} - \sum_{x \in E_N} \alpha_x \xi(x).$$

En utilisant la matrice $\tau$, il vient

$$\gamma = \frac{2}{3} \sum_{x \in E_N} \alpha_x \{g_x \rho, g_x \infty\} + \alpha_{x\tau} \{g_{x\tau} \rho, g_{x\tau} \infty\} + \alpha_{x\tau^2} \{g_{x\tau^2} \rho, g_{x\tau^2} \infty\}$$
$$- \sum_{x \in E_N} \alpha_x \xi(x)$$
$$= \frac{2}{3} \sum_{x \in E_N} \alpha_x \{g_x \rho, g_x \infty\} + \alpha_{x\tau} \{g_x \rho, g_x 0\} + \alpha_{x\tau^2} \{g_x \rho, g_x 1\}$$
$$- \sum_{x \in E_N} \alpha_x \xi(x).$$

On fait alors apparaître $(\alpha_x + \alpha_{x\tau} + \alpha_{x\tau^2})\{g_x \rho, g_x \infty\}$ ce qui donne

$$\gamma = \frac{2}{3} \sum_{x \in E_N} \alpha_{x\tau} \{g_x \infty, g_x 0\} + \alpha_{x\tau^2} \{g_x \infty, g_x 1\} - \sum_{x \in E_N} \alpha_x \xi(x).$$

L'identité $\{g_x \infty, g_x 1\} = \{g_{x\tau^2} 0, g_{x\tau^2} \infty\}$ donne finalement

$$\gamma = \frac{2}{3} \sum_{x \in E_N} -\alpha_{x\tau} \xi(x) + \frac{2}{3} \sum_{x \in E_N/\pm 1} \alpha_x \xi(x) - \sum_{x \in E_N} \alpha_x \xi(x)$$
$$= -\frac{1}{3} \sum_{x \in E_N} (\alpha_x + 2\alpha_{x\tau}) \xi(x).$$

□

Les quantités $\xi_{f_\epsilon}^+(x) = \left(\xi_{f_\epsilon}(x) + \overline{\xi_{f_{\overline{\epsilon}}}(x)}\right)$ vérifient les relations de Manin, ce qui permet d'appliquer le lemme précédent au cycle $b_{\overline{\epsilon}}$. Nous obtenons

$$b_{\overline{\epsilon}} = -\frac{1}{3} \sum_{x \in E_{13}/\pm 1} \left(\xi_{f_\epsilon}^+(x) + 2\xi_{f_\epsilon}^+(x\tau)\right)\xi(x).$$

Cette somme comporte 84 termes et il n'est donc pas commode de l'évaluer. Pour réduire le nombre de termes à 14, nous fixons $x_0 \in E_{13}/\pm 1$ et calculons

$$\sum_{\lambda \in (\mathbf{Z}/13\mathbf{Z})^*/\pm 1} \left(\xi_{f_\epsilon}^+(\lambda x_0) + 2\xi_{f_\epsilon}^+(\lambda x_0 \tau)\right)\xi(\lambda x_0)$$
$$= \sum_{\lambda \in (\mathbf{Z}/13\mathbf{Z})^*/\pm 1} \left(\xi_{f_\epsilon}^+(x_0) + 2\xi_{f_\epsilon}^+(x_0 \tau)\right) \cdot \epsilon(\lambda) \langle\lambda\rangle_* \xi(x_0)$$
$$= 6\left(\xi_{f_\epsilon}^+(x_0) + 2\xi_{f_\epsilon}^+(x_0 \tau)\right) \cdot \xi(x_0)^{\overline{\epsilon}}.$$

Notons $\mathcal{E}_0 \subset E_{13}/\pm 1$ l'ensemble des couples $(0,1)$ et $(1,v)$, $v \in \mathbf{Z}/13\mathbf{Z}$. Nous obtenons

$$b_{\overline{\epsilon}} = -2 \sum_{x_0 \in \mathcal{E}_0} \left(\xi_{f_\epsilon}^+(x_0) + 2\xi_{f_\epsilon}^+(x_0 \tau)\right) \cdot \xi(x_0)^{\overline{\epsilon}}.$$

Un calcul simple montre que les termes correspondant à $x_0 = (0,1)$ et $x_0 = (1,0)$ sont opposés l'un de l'autre[1]. D'où

$$b_{\overline{\epsilon}} = -2 \sum_{v \in (\mathbf{Z}/13\mathbf{Z})^*} \left(\xi_{f_\epsilon}^+(1,v) + 2\xi_{f_\epsilon}^+(v,-v-1)\right) \cdot \xi(1,v)^{\overline{\epsilon}}$$
$$= -2 \sum_{v \in (\mathbf{Z}/13\mathbf{Z})^*} \left(\xi_{f_\epsilon}^+(1,v) + 2\epsilon(v)\xi_{f_\epsilon}^+(1,1+\frac{1}{v})\right) \cdot \xi(1,v)^{\overline{\epsilon}}.$$

Or [55, Lemme 5] permet d'exprimer $\xi_{f_\epsilon}^+(1,w)$, $w \neq 0$ en fonction de $\xi_{f_\epsilon}^+(1,2)$ et $\xi_{f_\epsilon}^+(1,3)$. Nous avons la relation $\xi_{f_\epsilon}^+(1,-w) = \xi_{f_\epsilon}^+(1,w)$ et

$$\xi_{f_\epsilon}^+(1,1) = 0 \qquad\qquad \xi_{f_\epsilon}^+(1,4) = (1-\zeta_6)\xi_{f_\epsilon}^+(1,3)$$
$$\xi_{f_\epsilon}^+(1,5) = (\zeta_6-1)\left(\xi_{f_\epsilon}^+(1,2) - \xi_{f_\epsilon}^+(1,3)\right) \qquad \xi_{f_\epsilon}^+(1,6) = (\zeta_6-1)\xi_{f_\epsilon}^+(1,2).$$

D'autre part grâce à [55, Lemme 5], les cycles $\xi(1,v)^{\overline{\epsilon}}$, $v \in (\mathbf{Z}/13\mathbf{Z})^*$ s'expriment comme combinaisons linéaires de $\xi(1,2)^{\overline{\epsilon}}$, $\xi(1,3)^{\overline{\epsilon}}$ et $\xi(1,-3)^{\overline{\epsilon}}$. Il en va donc de même de $b_{\overline{\epsilon}}$. Mais nous savons que $b_{\overline{\epsilon}}$ est multiple de $\gamma_{\overline{\epsilon}}^- = \xi(1,3)^{\overline{\epsilon}} - \xi(1,-3)^{\overline{\epsilon}}$. Il suffit donc de calculer le coefficient devant $\xi(1,3)^{\overline{\epsilon}}$, ce qui conduit à l'identité

$$b_{\overline{\epsilon}} = \left(6\xi_{f_\epsilon}^+(1,2) + (4\zeta_6-2)\xi_{f_\epsilon}^+(1,3)\right)\left(\xi(1,3)^{\overline{\epsilon}} - \xi(1,-3)^{\overline{\epsilon}}\right). \tag{3.184}$$

Dans le premier facteur du produit (3.184), nous reconnaissons l'intégrale de $-\frac{1}{2\pi}\omega_{f_\epsilon}$ le long du cycle $(2\zeta_6-1)\gamma_\epsilon^+$. Nous obtenons finalement

$$b_{\overline{\epsilon}} = (2\zeta_6-1)\left(-\frac{1}{2\pi}\int_{\gamma_\epsilon^+} \omega_{f_\epsilon}\right)\gamma_{\overline{\epsilon}}^-. \tag{3.185}$$

---

[1] Cela résulte aussi de considérations de bord, compte tenu de (3.101).

En remplaçant dans (3.180), il vient

$$L(f_\epsilon, 2)L(f_\epsilon, \epsilon^2, 1) = \frac{13\pi}{48}(2\zeta_6 - 1)\Big(\int_{\gamma_\epsilon^+} \omega_{f_\epsilon}\Big) \int_{\gamma_{\bar\epsilon}^-} \eta(u_{\epsilon^3}, u_{\widehat{\epsilon^2}}). \tag{3.186}$$

Le caractère $\epsilon^2$ étant primitif, nous avons [45, Thm 4.2.b)]

$$L(f_\epsilon, \epsilon^2, 1) = -\frac{\tau(\epsilon^2)}{13} \sum_{a \in (\mathbf{Z}/13\mathbf{Z})^*} \bar\epsilon^2(a) \int_{\widetilde{a}/13}^{\infty} \omega_{f_\epsilon}, \tag{3.187}$$

où $\widetilde{a}$ désigne un représentant de $a$ dans $\mathbf{Z}$. En écrivant le symbole modulaire intervenant dans l'intégrale (3.187) en termes de symboles de Manin, nous obtenons

$$\sum_{a \in (\mathbf{Z}/13\mathbf{Z})^*} \bar\epsilon^2(a)\Big\{\frac{\widetilde{a}}{13}, \infty\Big\} = (\zeta_6 - 1)\gamma_\epsilon^+.$$

En remplaçant dans (3.186) puis en divisant par la période de $\omega_{f_\epsilon}$ (que l'on sait être non nulle), nous obtenons

$$L(f_\epsilon, 2) \cdot -\frac{\tau(\epsilon^2)}{13}(\zeta_6 - 1) = \frac{13\pi}{48}(2\zeta_6 - 1) \int_{\gamma_{\bar\epsilon}^-} \eta(u_{\epsilon^3}, u_{\widehat{\epsilon^2}}),$$

d'où l'égalité

$$\int_{\gamma_{\bar\epsilon}^-} \eta(u_{\epsilon^3}, u_{\widehat{\epsilon^2}}) = -\frac{16}{13^2\pi}(1 + \zeta_6) \cdot \tau(\epsilon^2) \cdot L(f_\epsilon, 2).$$

En remplaçant dans (3.178), nous obtenons finalement

$$\int_{\gamma_{\bar\epsilon}^-} \eta(x, y) = -\frac{\sqrt{13}}{\pi}\zeta_6 \cdot \tau(\epsilon^2) \cdot L(f_\epsilon, 2). \tag{3.188}$$

Cette formule se simplifie si l'on utilise l'équation fonctionnelle de $L(f_\epsilon, s)$. Notons $w(f_\epsilon) \in \mathbf{C}$ la pseudo-valeur propre de $f_\epsilon$ pour $W_{13}$. Elle est définie par l'identité

$$f_\epsilon \mid W_{13} = w(f_\epsilon) f_{\bar\epsilon}. \tag{3.189}$$

Le treizième coefficient de Fourier de $f_\epsilon$ vaut $a_{13}(f_\epsilon) = -3\zeta_6 - 1$. D'après [4, Thm 2.1], nous avons

$$w(f_\epsilon) = \frac{3\zeta_6 - 4}{13} \cdot \tau(\epsilon). \tag{3.190}$$

Les propriétés des sommes de Gauß et de Jacobi permettent d'exprimer $\tau(\epsilon^2)$ en fonction de $\tau(\epsilon)$. Grâce à [35, §8.3, Thm 1 (d)], nous avons

$$\tau(\epsilon)\tau(\epsilon^2) = \Big(\sum_{a \in \mathbf{Z}/13\mathbf{Z}} \epsilon(a)\,\epsilon^2(1 - a)\Big)\tau(\epsilon^3) = (4\zeta_6 - 3)\tau(\epsilon^3).$$

Or d'après [35, §6.4, Thm 1], nous avons $\tau(\epsilon^3) = \sqrt{13}$, ce qui conduit à

$$\tau(\epsilon^2) = \frac{(4\zeta_6 - 3)\sqrt{13}}{\tau(\epsilon)}. \tag{3.191}$$

Posons

$$\Lambda(f,s) = 13^{s/2}(2\pi)^{-s}\Gamma(s)L(f,s) \qquad (f \in S_2(\Gamma_1(13))).$$

L'équation fonctionnelle pour $L(f_\epsilon, s)$ s'écrit [58, Prop 18]

$$\Lambda(f_\epsilon, s) = -w(f_\epsilon)\Lambda(f_{\bar\epsilon}, s).$$

Un calcul classique montre alors

$$L(f_\epsilon, 2) = \frac{4\pi^2}{13^2}(4 - 3\zeta_6)\tau(\epsilon)L'(f_{\bar\epsilon}, 0), \tag{3.192}$$

En reportant dans (3.188) et en utilisant (3.191), nous trouvons au final

$$\int_{\gamma_{\bar\epsilon}^-} \eta(x,y) = 4\pi(1 - \zeta_6)L'(f_{\bar\epsilon}, 0).$$

En appliquant la conjugaison complexe à cette identité, nous trouvons (3.177). □

*Remarque* 112. Il est possible de déduire du théorème 110 le calcul de l'intégrale $\int_\gamma \eta(x,y)$, où $\gamma$ est un élément de $H_1^-(X_1(13)(\mathbf{C}), \mathbf{Z})$. Par un argument galoisien, il est facile de voir que $\int_\gamma \eta(x,y)$ s'exprime en termes de $L'(f,0)$, où $f \in S_2(\Gamma_1(13))$ est une forme parabolique à coefficients *rationnels*.

Nous allons maintenant utiliser la théorie du chapitre 2 pour exprimer le régulateur (3.163). Nous considérons $X_1(13)(\mathbf{C})$ incluse dans $J_1(13)(\mathbf{C})$ au moyen de l'immersion fermée $i$ associée à la pointe infinie de $X_1(13)(\mathbf{C})$. Nous notons $P = i(P_2) \in J_1(13)(\mathbf{Q})$ l'image de la pointe $P_2$ par cette immersion fermée. Le point $P$ est d'ordre 19 et engendre $J_1(13)(\mathbf{Q})$ [46].

Nous pouvons considérer $R_{J_1(13)} = R_{J_1(13)(\mathbf{C})}$ comme une fonction à valeurs dans $\mathbf{T} \otimes \mathbf{C} \cong \mathrm{Hom}_\mathbf{C}(S_2(\Gamma_1(13)), \mathbf{C})$. Pour $\psi = \epsilon$ ou $\bar\epsilon$, posons $R_\psi = \langle R_{J_1(13)}, f_\psi \rangle$. Définissons des périodes réelles pour $f_\epsilon$ et $f_{\bar\epsilon}$ par

$$\Omega_\psi^+ = \int_{\gamma_\psi^+} \omega_{f_\psi} \qquad (\psi \in \{\epsilon, \bar\epsilon\}). \tag{3.193}$$

D'après le lemme 107, $\Omega_\psi^+$ est non nul, et nous vérifions la formule $\Omega_{\bar\epsilon}^+ = \overline{\Omega_\epsilon^+}$. Pour $\psi = \epsilon$ ou $\bar\epsilon$, définissons encore une fonction réduite $R'_\psi = R_\psi/\Omega_\psi^+$.

**Théorème 113.** *Nous avons les identités*

$$L(f_\epsilon, 2) = \frac{8\pi i \tau(\epsilon)}{13^2}(2 + 3\zeta_6)\bigl(R'_\epsilon(2P) - 2R'_\epsilon(P)\bigr) \tag{3.194}$$

$$= \frac{8\pi i \tau(\epsilon)}{13^2}(1 - 4\zeta_6)\bigl(R'_\epsilon(P) + R'_\epsilon(2P) + (2\zeta_6 - 1)R'_\epsilon(4P)\bigr). \tag{3.195}$$

*Remarque* 114. Dans ces identités, il est curieux que les normes des éléments $2 + 3\zeta_6$ et $1 - 4\zeta_6$ soient respectivement 19 et 13.

*Démonstration.* Les opérateurs diamants $\langle d \rangle$, $d \in (\mathbf{Z}/13\mathbf{Z})^*/\pm 1$ induisent des automorphismes de $J_1(13)$, par fonctorialité d'Albanese (2.87). Pour connaître leur action sur $J_1(13)(\mathbf{Q})$, il suffit de déterminer l'action du diamant $\langle 2 \rangle$ sur le point $P$. Posons $\langle 2 \rangle P = aP$ avec $a \in \frac{\mathbf{Z}}{19\mathbf{Z}}$. Noter que $\langle 5 \rangle = \langle 2 \rangle^3$ est l'involution hyperelliptique de $X_1(13)$, puisque le quotient de $X_1(13)$ par le sous-groupe engendré par $\langle 5 \rangle$ est de genre 0. Nous avons donc $a^3 \equiv -1 \pmod{19}$, soit

$$(a+1)(a^2 - a + 1) \equiv 0 \pmod{19}.$$

La congruence $a \equiv -1 \pmod{19}$ entraîne $P_4 - P_2 = -(P_2 - P_1)$, soit $P_4 - P_2 = P_1 - P_2$, ce qui est absurde puisque $P_4 \neq P_1$. Donc $a^2 \equiv a - 1 \pmod{19}$. En appliquant successivement l'opérateur $\langle 2 \rangle$, nous trouvons alors

$$\begin{aligned} P_4 - P_1 &= (a+1)P & P_5 - P_1 &= 2aP \\ P_3 - P_1 &= (2a-1)P & P_6 - P_1 &= (a-1)P. \end{aligned} \qquad (3.196)$$

D'après [53, p. 224], le diviseur 0 0 0 3 1 −4 est principal, ce qui fournit la relation

$$\bigl(3(a+1) + 2a - 4(a-1)\bigr)P = 0.$$

Nous en déduisons $a \equiv -7 \pmod{19}$, soit $\langle 2 \rangle P = -7P$. Les orbites de l'action des diamants sur $J_1(13)(\mathbf{Q})$ sont donc données par

$$\{0\} \qquad \{\pm P, \pm 7P, \pm 8P\} \qquad \{\pm 2P, \pm 3P, \pm 5P\} \qquad \{\pm 4P, \pm 6P, \pm 9P\}.$$

En appliquant la proposition 59 avec la forme différentielle $\omega = \omega_{f_\epsilon}$, nous trouvons

$$R_\epsilon\bigl(\langle d \rangle x\bigr) = \epsilon(d)\, R_\epsilon(x) \qquad (d \in (\mathbf{Z}/13\mathbf{Z})^*/\pm 1,\ x \in J_1(13)(\mathbf{C})). \qquad (3.197)$$

Pour $d = 5$, la formule (3.197) redonne l'imparité de $R_\epsilon$. Pour déterminer la fonction $R_\epsilon$ sur $J_1(13)(\mathbf{Q})$, il suffit donc de connaître $R_\epsilon(P)$, $R_\epsilon(2P)$ et $R_\epsilon(4P)$. Les autres valeurs de $R_\epsilon$ s'obtiennent grâce à l'imparité de $R_\epsilon$ et la table suivante

$$\begin{aligned} R_\epsilon(7P) &= -\zeta_6 R_\epsilon(P) & R_\epsilon(8P) &= (1-\zeta_6)R_\epsilon(P) \\ R_\epsilon(3P) &= (\zeta_6 - 1)R_\epsilon(2P) & R_\epsilon(5P) &= \zeta_6 R_\epsilon(2P) \\ R_\epsilon(6P) &= (\zeta_6 - 1)R_\epsilon(4P) & R_\epsilon(9P) &= -\zeta_6 R_\epsilon(4P). \end{aligned} \qquad (3.198)$$

Soit maintenant $\chi$ un caractère de Dirichlet modulo 13, pair et distinct de 1 et $\epsilon$. Utilisons le théorème 7 avec $N = 13$, $\psi = \epsilon$ et le caractère $\overline{\chi}$, qui est bien primitif et distinct de $\overline{\epsilon}$. Il vient

$$L(\mathbf{T}^\epsilon, 2)L(\mathbf{T}^\epsilon, \overline{\chi}, 1) = C_{\epsilon,\overline{\chi}} \sum_{\lambda,\mu \in (\frac{\mathbf{Z}}{13\mathbf{Z}})^*/\pm 1} \chi\overline{\epsilon}(\lambda)\overline{\chi}(\mu) R_{J_1(13)}(P_\lambda - P_\mu)$$

soit en appliquant à $f_\epsilon$

$$L(f_\epsilon, 2)L(f_\epsilon, \overline{\chi}, 1) = C_{\epsilon,\overline{\chi}} \sum_{\lambda,\mu \in (\frac{\mathbf{Z}}{13\mathbf{Z}})^*/\pm 1} \chi\overline{\epsilon}(\lambda)\overline{\chi}(\mu) R_\epsilon(P_\lambda - P_\mu). \qquad (3.199)$$

Grâce à (3.197) et à l'imparité de $R_\epsilon$, il vient

$$\begin{aligned} &\sum_{\lambda,\mu \in (\frac{\mathbf{Z}}{13\mathbf{Z}})^*/\pm 1} \chi\overline{\epsilon}(\lambda)\overline{\chi}(\mu) R_\epsilon\bigl(\langle \lambda \rangle(P_1 - P_{\mu/\lambda})\bigr) \\ &= \sum_{\lambda,\mu \in (\frac{\mathbf{Z}}{13\mathbf{Z}})^*/\pm 1} \chi(\lambda)\overline{\chi}(\mu) R_\epsilon(P_1 - P_{\mu/\lambda}) \\ &= -6 \sum_{\lambda \in (\frac{\mathbf{Z}}{13\mathbf{Z}})^*/\pm 1} \overline{\chi}(\lambda) R_\epsilon(P_\lambda - P_1). \end{aligned}$$

Pour tout caractère $\chi$ pair modulo 13, posons

$$S_\epsilon(\chi) = \sum_{\lambda \in (\mathbf{Z}/13\mathbf{Z})^*/\pm 1} \overline{\chi}(\lambda) R_\epsilon(P_\lambda - P_1). \tag{3.200}$$

L'égalité (3.199) peut donc s'écrire

$$L(f_\epsilon, 2) L(f_\epsilon, \chi, 1) = -6\, C_{\epsilon, \overline{\chi}}\, S_\epsilon(\chi) \qquad (\chi \neq 1, \epsilon). \tag{3.201}$$

Le changement de variable $\lambda \mapsto 1/\lambda$ dans (3.200) permet de constater que $S_\epsilon(\overline{\chi}\epsilon) = -S_\epsilon(\chi)$. Posons $\zeta = \chi(2) \in \mu_6$. Grâce à (3.196) et (3.198), nous trouvons explicitement

$$S_\epsilon(\chi) = \big(\overline{\zeta} + \zeta(\zeta_6 - 1)\big) R_\epsilon(P) + \zeta^3 \zeta_6 R_\epsilon(2P) + \big(\zeta^2 + \overline{\zeta}^2(1 - \zeta_6)\big) R_\epsilon(4P).$$

En particulier, nous avons

$$\begin{pmatrix} S_\epsilon(1) \\ S_\epsilon(\epsilon^2) \\ S_\epsilon(\epsilon^3) \end{pmatrix} = \begin{pmatrix} \zeta_6 & \zeta_6 & 2 - \zeta_6 \\ -2\zeta_6 & \zeta_6 & 0 \\ -\zeta_6 & -\zeta_6 & 2 - \zeta_6 \end{pmatrix} \begin{pmatrix} R_\epsilon(P) \\ R_\epsilon(2P) \\ R_\epsilon(4P) \end{pmatrix}. \tag{3.202}$$

L'explicitation de $S_\epsilon(\epsilon) = -S_\epsilon(1)$, $S_\epsilon(\overline{\epsilon}) = -S_\epsilon(\epsilon^2)$ et $S_\epsilon(\overline{\epsilon}) = -S_\epsilon(\epsilon^3)$ ne livre pas d'information supplémentaire sur la fonction $R_\epsilon$ : la matrice permettant d'exprimer ces quantités en fonction de $R_\epsilon(P)$, $R_\epsilon(2P)$ et $R_\epsilon(4P)$ est l'opposée de la matrice (3.202).

Calculons maintenant le produit $L(\chi, 2) L(\overline{\chi}\epsilon, 2)/\pi^4$ apparaissant dans la constante $C_{\epsilon, \overline{\chi}}$, en utilisant les formules (3.165) et [35, §8.3, Thm 1 (d)]

$$\begin{aligned} \frac{L(\chi, 2) L(\overline{\chi}\epsilon, 2)}{\pi^4} &= \frac{\tau(\chi) \tau(\overline{\chi}\epsilon)}{13^4} B_{2, \overline{\chi}} B_{2, \chi\overline{\epsilon}} \\ &= \Big( \sum_{a \in \frac{\mathbf{Z}}{13\mathbf{Z}}} \chi(a) \overline{\chi}\epsilon(1 - a) \Big) \frac{\tau(\epsilon)}{13^4} B_{2, \overline{\chi}} B_{2, \chi\overline{\epsilon}} \qquad (\chi \neq 1, \epsilon). \end{aligned}$$

En menant les calculs explicitement, nous trouvons pour $\chi \neq 1, \epsilon$

$$\begin{aligned} \frac{L(\chi, 2) L(\overline{\chi}\epsilon, 2)}{\pi^4} = \frac{4\tau(\epsilon)}{3 \cdot 13^5} \Big( &(504 - 101\zeta - 119\zeta^2 - 19\overline{\zeta}^2 - 157\overline{\zeta}) \\ &+ (-274 + 157\zeta + 138\zeta^2 - 119\overline{\zeta}^2 + 56\overline{\zeta})\zeta_6 \Big). \end{aligned}$$

De la même manière qu'en (3.187), nous avons

$$L(f_\epsilon, \overline{\chi}, 1) = -\frac{\tau(\overline{\chi})}{13} \sum_{a \in (\mathbf{Z}/13\mathbf{Z})^*} \chi(a) \int_{\widetilde{a}/13}^{\infty} \omega_{f_\epsilon}.$$

Notons $\theta(\chi)$ le symbole modulaire suivant sur $X_1(13)(\mathbf{C})$

$$\theta(\chi) = \sum_{a \in (\mathbf{Z}/13\mathbf{Z})^*/\pm 1} \chi(a) \Big\{ \frac{\widetilde{a}}{13}, \infty \Big\}, \tag{3.203}$$

et $\theta(\chi)^\epsilon$ la projection de $\theta(\chi)$ sur la composante $\epsilon$-isotypique. Lorsque $\chi \neq 1, \epsilon$, le cycle $\theta(\chi)^\epsilon$ est fermé, et nous trouvons explicitement

$$\theta(\chi)^\epsilon = (\zeta^2 + \zeta^3 + \overline{\zeta}^2 - \overline{\zeta}^2 \zeta_6) \cdot \gamma_\epsilon^+ \qquad (\chi \neq 1, \epsilon). \tag{3.204}$$

Il vient donc

$$\frac{L(f_\epsilon, \overline{\chi}, 1)}{\tau(\overline{\chi})} = -\frac{1}{13}(\zeta^2 + \zeta^3 + \overline{\zeta}^2 - \overline{\zeta}^2 \zeta_6) \cdot \Omega_\epsilon^+ \qquad (\chi \neq 1, \epsilon). \qquad (3.205)$$

Injectons tous les calculs précédents dans la formule (3.201). Pour $\chi = \epsilon^2$, nous obtenons (3.194), tandis que pour $\chi = \epsilon^3$, nous trouvons (3.195). Notons que les caractères $\chi = \overline{\epsilon}$ et $\chi = \overline{\epsilon}^2$ donnent respectivement les mêmes égalités (3.194) et (3.195), ce qui découle des propriétés de symétrie vues précédemment. □

*Remarque* 115. Il serait bon de vérifier numériquement les identités (3.194) et (3.195), ne serait-ce que pour confirmer ces calculs. Cela nécessite de trouver un moyen de calcul efficace de la fonction $R'_\epsilon$. Il serait alors extrêmement intéressant de chercher numériquement d'autres relations entre la valeur spéciale $L(f_\epsilon, 2)$ et la fonction $R'_\epsilon$. Par exemple, chacune des valeurs $R'_\epsilon(P)$, $R'_\epsilon(2P)$ et $R'_\epsilon(4P)$ est-elle proportionnelle à $L(f_\epsilon, 2)$ ?

*Remarque* 116. Il est intéressant de chercher une version du théorème 113 linéaire en la forme modulaire $f \in S_2(\Gamma_1(13))$. Nous proposons pour cela une méthode. Notons $L(2)$ la forme linéaire $f \mapsto L(f, 2)$ sur $S_2(\Gamma_1(13))$. Soit $R : J_1(13)(\mathbf{C}) \to \mathrm{Hom}_{\mathbf{C}}(S_2(\Gamma_1(13)), \mathbf{C})$ la fonction définie de manière unique par

$$\langle R(x), f_\epsilon \rangle = \tau(\epsilon) R'_\epsilon(x) \quad \text{et} \quad \langle R(x), f_{\overline{\epsilon}} \rangle = -\tau(\overline{\epsilon}) R'_{\overline{\epsilon}}(x) \qquad (x \in J_1(13)(\mathbf{C})).$$

Nous déduisons de (3.194) et (3.198)

$$L(f_\epsilon, 2) = \frac{8\pi i \tau(\epsilon)}{13^2} \big( -4 R'_\epsilon(P) + 5 R'_\epsilon(2P) + 3 R'_\epsilon(3P) + 6 R'_\epsilon(7P) \big).$$

Puisque $L(f_{\overline{\epsilon}}, 2) = \overline{L(f_\epsilon, 2)}$ et $\overline{R'_\epsilon} = R'_{\overline{\epsilon}} \circ c$, nous obtenons alors

$$L(2) = \frac{8\pi i}{13^2} \big( -4 R(P) + 5 R(2P) + 3 R(3P) + 6 R(7P) \big). \qquad (3.206)$$

Remarquons que le diviseur $-4[P] + 5[2P] + 3[3P] + 6[7P]$ est de somme nulle dans $J_1(13)(\mathbf{Q})$. De même, (3.195) mène à

$$L(2) = \frac{8\pi i}{13^2} \big( R(P) - 3R(2P) - 4R(3P) + 5R(4P) - 2R(6P) + 4R(7P) \big), \qquad (3.207)$$

et le diviseur $[P] - 3[2P] - 4[3P] + 5[4P] - 2[6P] + 4[7P]$ est de somme nulle dans $J_1(13)(\mathbf{Q})$. Revenant à l'égalité (3.206), nous remarquons que chacun des diviseurs $-4[P] + 6[7P]$ et $5[2P] + 3[3P]$ est de somme nulle. Nous sommes donc tentés de conjecturer que chacune des quantités $-4R(P) + 6R(7P)$ et $5R(2P) + 3R(3P)$ est proportionnelle à $L(2)$. Voici encore quelques questions. Est-il possible de trouver une définition plus intrinsèque pour la fonction $R$ ? D'autre part, une formule de type (3.206) a-t-elle lieu dans $J_1(N)$, pour tout entier $N$ ?

*Remarque* 117. La courbe $X_1(N)$ est de genre 2 pour deux autres valeurs de $N$, à savoir $N = 16$ et $N = 18$. Le théorème 113 devrait pouvoir s'adapter à ces cas-là.

## 3.9 Vers la mesure de Mahler

Notons $E = X_1(11)$ la courbe elliptique de conducteur 11 considérée dans la section 3.7. D'après les observations et travaux de Smyth, Boyd, Deninger, Rodriguez-Villegas, Bertin…, la valeur spéciale $L(E, 2)$ semble être liée numériquement à la mesure de Mahler de certains polynômes en deux variables, par exemple [16, (1-25)]

$$m\bigl((X+Y+1)(X+1)(Y+1)+XY\bigr) \stackrel{?}{=} \frac{7\cdot 11}{4\pi^2}L(E,2), \tag{3.208}$$

où pour un polynôme $P(X,Y) \in \mathbf{C}[X,Y]$, $P \neq 0$, $m(P)$ désigne la *mesure de Mahler logarithmique de $P$*, définie par

$$m(P) = \int_0^1 \int_0^1 \log|P(e^{2\pi i u}, e^{2\pi i v})| du\, dv. \tag{3.209}$$

Nous démontrons maintenant le résultat suivant.

**Théorème 118.** *L'identité (3.208) est vraie.*

La preuve de ce théorème repose sur des résultats et travaux de Deninger [23], Rodriguez-Villegas [57] et Bertin [10]. L'ingrédient final est la formule explicite pour $L(E,2)$ obtenue à la section 3.7. Nous renvoyons à l'article de Boyd [16] pour une introduction à la mesure de Mahler ainsi qu'à son lien avec les valeurs spéciales de fonctions $L$.

*Démonstration.* Dans l'article [10], Bertin a calculé la mesure de Mahler (3.208) en termes d'un régulateur sur la courbe elliptique $E$. Dans la remarque finale de *op. cit.*, Bertin montre que la relation (3.208) découle de l'égalité conjecturale suivante

$$|r(\{x,y\})| \stackrel{?}{=} \frac{11}{4\pi^2}L(E,2), \tag{3.210}$$

où le régulateur $r(\{x,y\})$ est défini par [56, 57]

$$r(\{x,y\}) = \frac{1}{2\pi} \int_\gamma \eta(x,y) = \frac{1}{2\pi} \int_\gamma \log|x| d\arg y - \log|y| d\arg x, \tag{3.211}$$

et $\gamma$ est un chemin fermé quelconque de $E(\mathbf{C})$ ne passant pas par les zéros et pôles de $x$ et $y$, et engendrant[2] le sous-groupe $H_1^-(E(\mathbf{C}), \mathbf{Z})$ formé des cycles anti-invariants par la conjugaison complexe agissant sur $E(\mathbf{C})$. Grâce à (3.151), la preuve de (3.210) se ramène à celle de l'identité

$$|r(\{x,y\})| \stackrel{?}{=} \frac{5}{2\pi}D_E(P).$$

Or, cette dernière égalité résulte de [10, Th. 6 et Cor. 6.1]. □

*Remarque* 119. En utilisant [10, Th. 5], nous obtenons également l'identité suivante, conjecturée par Boyd [16, (1-24)]

$$m\bigl(Y^2 + (X^2+2X-1)Y + X^3\bigr) = \frac{5\cdot 11}{4\pi^2}L(E,2). \tag{3.212}$$

## 3.10 Questions et perspectives

Soit $A$ une variété abélienne de type $\mathrm{GL}_2$ définie sur $\mathbf{Q}$, et $L(A,s)$ la fonction $L$ associée à $A$. Conjecturalement, il existe un entier $N \geq 1$ tel que la variété abélienne $A$ soit quotient (ou facteur, à isogénie près) de la jacobienne $J_1(N)$ de la courbe modulaire $X_1(N)$. Supposant cela pour $A$, il semble que l'association des chapitres 2 et 3 permette de montrer une version faible de la conjecture de Beĭlinson concernant la valeur spéciale $L(A,2)$. Plus précisément, considérons le régulateur de Beĭlinson

---

[2]Le groupe $H_1^-(E(\mathbf{C}), \mathbf{Z})$ possédant deux générateurs, le nombre $r(\{x,y\})$ est défini au signe près. Puisque l'identité (3.208) relie deux nombres positifs, cette ambiguïté de signe ne pose pas problème.

$$r_A : K_2^{(2)}(A) \to V_{\mathbf{R}} = H^1(A(\mathbf{C}), \mathbf{R}(1))^- \qquad (3.213)$$

(voir la section 2.6). L'espace $V_{\mathbf{Q}} = H^1(A(\mathbf{C}), \mathbf{Q}(1))^-$ est une $\mathbf{Q}$-structure de $V_{\mathbf{R}}$. La forme faible de la conjecture de Beĭlinson pour $L(A, 2)$ prédit l'existence d'un sous-$\mathbf{Q}$-espace vectoriel $\mathcal{P}_A \subset K_2^{(2)}(A)$ tel que l'image $r_A(\mathcal{P}_A)$ soit une $\mathbf{Q}$-structure de $V_{\mathbf{R}}$, et tel que les deux $\mathbf{Q}$-structures ainsi définies soient reliées par

$$\det_{V_{\mathbf{Q}}} r_A(\mathcal{P}_A) \sim_{\mathbf{Q}^*} L^{(g)}(A, 0) \sim_{\mathbf{Q}^*} \pi^{-2g} L(A, 2) \qquad (g = \dim A),$$

où nous avons écrit $\sim_{\mathbf{Q}^*}$ pour indiquer l'égalité à un facteur rationnel non nul près[3]. Les résultats des chapitres 2 et 3, joints à la fonctorialité du régulateur de Beĭlinson pour le morphisme $A \to J_1(N)$, entraînent la conjecture faible énoncée ci-dessus. Nous n'avons pas abordé la question de l'intégralité des éléments construits dans le groupe $K_2^{(2)}(A)$. Pour cette dernière question, il nous semble intéressant d'envisager la condition d'intégralité à l'aide du modèle de Néron de $A$. Dans le cas d'une courbe elliptique, cela a été fait par Schappacher et Scholl [63, 64], voir également [76, §1.4]. Dans le cas général, la définition usuelle du groupe $K_2^{(2)}(A)_{\mathbf{Z}}$ fait intervenir un modèle propre et plat de $A$ sur $\mathbf{Z}$. Un tel modèle a l'avantage d'exister, mais la définition de $K_2^{(2)}(A)_{\mathbf{Z}}$ qui en résulte a le désavantage de dépendre du choix du modèle [36].

Enfin, nous pensons et espérons que les méthodes présentées dans les sections 3.7 à 3.9 auront des applications à la détermination exacte d'autres mesures de Mahler. On peut se référer à [16] pour une liste impressionnante d'identités conjecturales reliant mesures de Mahler et valeurs spéciales de fonctions $L$ de courbes elliptiques.

---

[3] La forme forte de la conjecture prédit que l'on peut prendre $\mathcal{P}_A = K_2^{(2)}(A)_{\mathbf{Z}}$, un groupe qui reste cependant à définir (à ma connaissance).

# Bibliographie

# Appendice :
# Symboles de Manin et valeurs de fonctions L

Loïc Merel

Ce texte est un complément à la thèse de F. Brunault, qui fait usage des théorèmes A, C et D ci-dessous et de leurs corollaires. Nous espérons donner un compte-rendu plus complet, incluant une généralisation du théorème A qui ne se limite pas au poids 2, et qui tire avantage du langage adélique dans une publication ultérieure.

## 1. L'interprétation arithmétique des symboles de Manin

**1.a.** Soit $N$ un entier $> 0$. Soit $f$ une forme modulaire primitive (c'est-à-dire propre pour l'algèbre de Hecke, nouvelle et normalisée) de poids $k = 2$ pour le groupe de congruence $\Gamma_1(N)$.

Pour tout entier $m \geq 1$, on note $\Sigma_m$ le support de $m$ dans l'ensemble des nombres premiers. Soit $\chi$ un caractère de Dirichlet de conducteur à support dans $\Sigma_N$. Notons $f \otimes \chi$ la forme primitive dont le $p$-ième coefficient de Fourier est $a_p(f)\chi(p)$ ($p$ nombre premier ne divisant pas $N$). Notons $N_\chi$ le niveau de $f \otimes \chi$. Notons $L(f \otimes \chi, s)$ la fonction $L$ de $f \otimes \chi$. Elle admet un développement en série de Dirichlet $\sum_{n=1}^\infty a_n(f \otimes \chi)/n^s$ et en produit eulerien $\prod_p L_p(f \otimes \chi, p^{-s})$, où $L_p(f \otimes \chi, X) = 1/(1 - a_p(f \otimes \chi)X + a_{p,p}(f \otimes \chi)p^{k-1}X^2)$ ($p$ nombre premier) ; on complète ce produit pour former $\Lambda(f \otimes \chi, s) = (2\pi)^{-s}\Gamma(s)N_\chi^{s/2}L(f \otimes \chi, s)$. On pose $a_p = a_p(f)$ et $a_{p,p} = a_{p,p}(f)$. Notons $\psi$ le caractère de nebentypus de $f$ ; il vérifie $\psi(p) = a_{p,p}(f)$ ($p$ nombre premier ne divisant pas $N$). On pose $\bar{f} = f \otimes \bar{\psi}$, et on a $a_n(\bar{f}) = \bar{a}_n(f)$ ($n$ entier $\geq 1$).

Lorsque $T^+$ et $T^-$ sont des ensembles finis de nombres premiers, on prive $\Lambda(f \otimes \chi, s)$ de certains facteurs d'Euler en posant

$$\Lambda^{[T^+, T^-]}(f \otimes \chi, s) = \Lambda(f \otimes \chi, s)/(\prod_{p \in T^+} L_p(f \otimes \chi, p^{-s}) \prod_{p \in T^-} L_p(\bar{f} \otimes \bar{\chi}, p^{s-k})).$$

Lorsque $R^+$ et $R^-$ sont des sous-ensembles de $T^+$ et $T^-$ respectivement, on pose

$$\Lambda^{[\frac{T^+}{R^+}, \frac{T^-}{R^-}]}(f \otimes \chi, s) = \Lambda^{[T^+, T^-]}(f \otimes \chi, s) \prod_{p \in R^+} \frac{L_p(f \otimes \chi, p^{-s})}{L_p(f \otimes \chi, p^{-s-1})} \prod_{p \in R^-} \frac{L_p(\bar{f} \otimes \bar{\chi}, p^{s-k})}{L_p(\bar{f} \otimes \bar{\chi}, p^{s-k+1})}.$$

Nous dirons que les nombres premiers $p$ qui vérifient $v_p(N) = 1$ et $\psi$ non ramifié en $p$ sont *spéciaux* pour $f$ (ils correspondent aux représentation spéciales de $\mathrm{GL}_2(\mathbf{Q}_p)$). Notons $\Sigma_f$ l'ensemble des nombres premiers spéciaux pour $f$. Le cas qui nous intéressera est le cas où $R^+$ et $R^-$ sont composés de nombres premiers spéciaux pour $f$.

Soit $S$ un sous-ensemble de $\Sigma_N$. Posons $\bar{S} = \Sigma_N - S$. Pour $M$ nombre entier $\geq 1$ à support dans $\Sigma_N$, posons $M = M_S M_{\bar{S}}$ où $M_S$ et $M_{\bar{S}}$ sont à supports dans $S$ et $\bar{S}$ respectivement. On pose $S(M) = \Sigma_M \cap S$ et $\bar{S}(M) = \Sigma_M \cap \bar{S}$. On note $w_S(\bar{f} \otimes \chi)$ la pseudo-valeur propre de $\bar{f} \otimes \chi$ pour l'opérateur d'Atkin-Lehner associé à $S$ (voir [A-L] ou la mise au point de la section **2.b.**). On note de plus $w(f \otimes \chi) = w_{\Sigma_N}(f \otimes \chi)$ ; on a $\Lambda^{[T^+, T^-]}(f \otimes \chi, s) = -w(f \otimes \chi)\Lambda^{[T^-, T^+]}(\bar{f} \otimes \bar{\chi}, k - s)$. Pour $\alpha$ caractère de niveau à support dans $\Sigma_N$, on convient de décomposer $\alpha$ sous la forme $\alpha = \alpha_S \alpha_{\bar{S}}$, où $\alpha_S$ et $\alpha_{\bar{S}}$ sont des caractères de Dirichlet de niveaux à supports dans $S$ et $\bar{S}$ respectivement.

Soit $(u, v) \in (\mathbf{Z}/N\mathbf{Z})^2$. Notons $N'$ l'ordre de $uv$ dans $\mathbf{Z}/N\mathbf{Z}$. Pour $p \in \Sigma_N$ et $\chi$ caractère de Dirichlet primitif de conducteur $m_\chi$ divisant $N$, notons $Q_{p,f,\chi}(X)$ la fraction rationnelle suivante : $Q_{p,f,\chi}(X) = (\bar{a}_p p^{1-k/2})^{v_p(N'/m_\chi)}$ sauf si $a_p = 0$, $v_p(N') = 1$ et $v_p(m_\chi) = 0$, auquel cas on a $Q_{p,f,\chi}(X) = -\bar{\chi}(p)X^{-1}$. Cet objet désagréable dépend de $p$, $a_p$, $\chi(p)$, $v_p(m_\chi)$, $k$ et $v_p(N')$ ; c'est donc un objet local.

**1.b.** Les notations qui précèdent sont valables même lorsque $k > 2$. Supposons, jusqu'à la fin de cette section, que $k = 2$.

On pose $\xi_f(u,v) = 0$ si $(u,v)$ n'est pas d'ordre $N$ dans le groupe additif $(\mathbf{Z}/N\mathbf{Z})^2$. Sinon on considère une matrice $g = \begin{pmatrix} a & b \\ c & d \end{pmatrix} \in \mathrm{SL}_2(\mathbf{Z})$ telle que $(c,d) \in (u,v)$ et on pose

$$\xi_f(u,v) = -i \int_{g0}^{g\infty} f(z)dz,$$

où l'intégrale est prise le long d'un chemin continu du demi-plan de Poincaré. On dispose donc de

$$\xi_f : (\mathbf{Z}/N\mathbf{Z})^2 \longrightarrow \mathbf{C}.$$

L'application $f \mapsto \xi_f$ est injective. On peut même être plus précis. Si on pose $\xi_f^+(u,v) = (\xi_f(u,v) + \xi_f(-u,v))/2$ et $\xi_f^-(u,v) = (\xi_f(u,v) - \xi_f(-u,v))/2$. Les applications $f \mapsto \xi_f^+$ et $f \mapsto \xi_f^-$ sont injectives [M].

Notre nous proposons de donner un sens à $\xi_f$ en terme des invariants arithmétiques de $f$. Supposons $(u,v)$ d'ordre $N$. Pour cela nous calculons la transformée de Fourier multiplicative de cette fonction. En effet, on identifie $(\mathbf{Z}/N\mathbf{Z})$ à $\cup_{d|N} (\mathbf{Z}/d\mathbf{Z})^*$ (par $w \mapsto wN'_w/N \pmod{N'_w}$, où $N'_w$ est l'ordre de $w$ dans $(\mathbf{Z}/N\mathbf{Z})$). Soit $S$ un sous-ensemble de $\Sigma_N$ contenant le support de $u$ mais disjoint du support de $v$. Posons $\bar{S} = \Sigma_N - S$. Les entiers $N'_S$ et $N'_{\bar{S}}$ expriment les composantes associées à $u$ et $v$ via cette identification ; les entiers $N_S/N'_S$ et $N_{\bar{S}}/N'_{\bar{S}}$ ne dépendent pas du choix de $S$. Toute fonction $\xi : (\mathbf{Z}/N\mathbf{Z})^2 \to \mathbf{C}$ s'écrit sous la forme $\xi(u,v) = \sum_{\alpha,\beta} c_{\alpha,\beta} \alpha(N'_{\bar{S}}v/N_{\bar{S}})\beta(N'_S u/N_S)$, où $c_{\alpha,\beta}$ dépend seulement de $\xi$, $\alpha$, $\beta$ et $N'$ et où $\alpha$ et $\beta$ parcourent les caractères de Dirichlet primitifs de niveau divisant $N'$. L'intérêt du théorème A réside dans une telle écriture pour $\xi = \xi_f$, qui met en évidence le sens arithmétique des coefficients $c_{\alpha,\beta}$.

Pour $\omega$ caractère de Dirichlet, notons $\tau'(\omega)$ la somme de Gauss et $m_\omega$ le conducteur du caractère primitif associé à $\omega$. Notons $\phi$ la fonction indicatrice d'Euler.

THÉORÈME A. — *On a*
$$\xi_f(u,v) =$$
$$\frac{w(f)}{\phi(N')} \sum_\chi \chi_{\bar{S}}(m_{\chi,S})(\bar{\psi}_S \bar{\chi}_S)(m_{\bar{\psi}\bar{\chi},\bar{S}})\chi_S(-1) \frac{\tau'(\chi_S)\tau'(\bar{\psi}_{\bar{S}}\bar{\chi}_{\bar{S}})}{\sqrt{N_\chi}} \Big( \prod_{p \in S(N')} Q_{p,f,\bar{\chi}}(1) \Big) \Big( \prod_{p \in \bar{S}(N')} Q_{p,f,\chi\psi}(1) \Big)$$
$$(\bar{\psi}_{\bar{S}}\bar{\chi}_{\bar{S}}^2)(N_{\chi,S})(\bar{\psi}\bar{\chi})\Big(\frac{N'_{\bar{S}}v}{N_{\bar{S}}}\Big)\chi\Big(\frac{N'_S u}{N_S}\Big)\overline{w_S(f \otimes \chi)}\Lambda^{\Big[\frac{\Sigma_{N'} - S(m_\chi)}{(S(N') - S(m_\chi)) \cap \Sigma_f}, \frac{\Sigma_{N'} - \bar{S}(m_{\bar{\psi}\bar{\chi}})}{(\bar{S}(N') - \bar{S}(m_{\bar{\psi}\bar{\chi}})) \cap \Sigma_f}\Big]}(f \otimes \chi, 1),$$

*où $\chi$ parcourt les caractères de Dirichlet primitifs tels que $m_{\chi,S} m_{\psi\chi,\bar{S}} | N'$.*

On déduit du théorème A une formule analogue pour $\xi_f^+$ (resp. $\xi_f^-$) qui fait disparaître les termes faisant intervenir les expressions $\Lambda(f \otimes \chi, 1)$ pour $\chi$ impair (resp. pair).

COROLLAIRE. — *La forme modulaire $f$ primitive de poids 2 pour $\Gamma_1(N)$ est caractérisée par les données suivantes, où on fait parcourir à $\chi$ les caractères de Dirichlet pairs (resp. impairs) de conducteur divisant $N$ :*
*(i) le caractère de $f$,*
*(ii) les niveaux des formes primitives $f \otimes \chi$,*
*(iii) les pseudo-valeurs propres $w_S(f \otimes \chi)$,*
*(iv) les facteurs d'Euler $L_p(f \otimes \chi, s)$, pour $p \in \Sigma_N$ et*
*(v) les nombres $\Lambda(f \otimes \chi, 1)$.*

On peut voir la fonction $\xi_f$ comme une façon commode de comprimer les données $(i)$, $(ii)$, $(iii)$, $(iv)$ et $(v)$. Il résulte du théorème A un énoncé de théorie analytique des nombres.

COROLLAIRE 2. — *Il existe un caractère de Dirichlet primitif $\chi$, qu'on peut choisir pair ou impair, de conducteur divisant $N$ et tel que $L(f \otimes \chi, 1) \neq 0$.*

## 2. Formulaire préliminaire

Cette section consiste en des mises au points concernant des questions essentiellement déjà connues. Elles concernent en **2.a.** la suppression des facteurs d'Euler des fonctions $L$, en **2.b.** les opérateurs d'Atkin-Lehner, en **2.c.** et **2.d.** la torsion des formes modulaires par des caractères non nécessairement primitifs, en **2.e.** la translation des formes modulaires par des nombres rationnels.

**2.a.** On note $\mathrm{GL}_2(\mathbf{Q})^+$ le sous-groupe de $\mathrm{GL}_2(\mathbf{Q})$ formé par les matrices de déterminant $> 0$. On pose, pour $\begin{pmatrix} A & B \\ C & D \end{pmatrix} \in \mathrm{GL}_2(\mathbf{Q})^+$, et $F$ forme primitive de poids $k$ et de niveau $M$ :

$$F_{|\begin{pmatrix} A & B \\ C & D \end{pmatrix}}(z) = \frac{(AD-BC)^{k/2}}{(Cz+D)^k} F(\frac{Az+B}{Cz+D}).$$

Cette opération s'étend $\mathbf{C}$-linéairement à $\mathbf{C}[\mathrm{GL}_2(\mathbf{Q})^+]$ ; elle se factorise par $\mathbf{C}[\mathrm{PGL}_2(\mathbf{Q})^+]$. Gardons à l'esprit la formule suivante

$$(2\pi)^{-s}\Gamma(s)L(F,s) = \int_0^\infty F(iy) y^s \frac{dy}{y}.$$

On a, pour $h = \begin{pmatrix} A & 0 \\ 0 & D \end{pmatrix} \in \mathrm{GL}_2(\mathbf{Q})^+$,

$$\int_0^\infty F_{|h}(iy) y^s \frac{dy}{y} = (\frac{A}{D})^{k/2-s} \int_0^\infty F(iy) y^s \frac{dy}{y} = (\frac{A}{D})^{k/2-s} (2\pi)^{-s} \Gamma(s) L(F,s).$$

Soient $T^+$ et $T^-$ deux ensembles de nombres premiers. On pose

$$F^{[T^+,T^-]} = F_{|\prod_{p\in T^+} L_p(F, p^{-k/2}\begin{pmatrix}p & 0 \\ 0 & 1\end{pmatrix})^{-1} \prod_{p\in T^-} L_p(\bar{F}, p^{-k/2}\begin{pmatrix}1 & 0 \\ 0 & p\end{pmatrix})^{-1}},$$

de telle sorte que

$$\int_0^\infty F^{[T^+,T^-]}(iy) y^s \frac{dy}{y} = \frac{(2\pi)^{-s}\Gamma(s)L(F,s)}{\prod_{p\in T^+} L_p(F, p^{-s}) \prod_{p\in T^-} L_p(\bar{F}, p^{s-k})} = M^{-s/2} \Lambda^{[T^+,T^-]}(F,s).$$

On pose, lorsque $R^+$ et $R^-$ sont des sous-ensembles de $T^+$ et $T^-$ respectivement,

$$F^{[\frac{T^+}{R^+}, \frac{T^-}{R^-}]} = F^{[T^+-R^+, T^--R^-]}_{|\prod_{p\in R^+} L_p(F, p^{1-k/2}\begin{pmatrix}p & 0 \\ 0 & 1\end{pmatrix})^{-1} \prod_{p\in R^-} L_p(\bar{F}, p^{1-k/2}\begin{pmatrix}1 & 0 \\ 0 & p\end{pmatrix})^{-1}}$$

si bien que

$$\int_0^\infty F^{[\frac{T^+}{R^+}, \frac{T^-}{R^-}]}(iy) y^s \frac{dy}{y} = M^{-s/2} \Lambda^{[\frac{T^+}{R^+}, \frac{T^-}{R^-}]}(F,s).$$

**2.b.** Mettons au point les normalisations pour les opérateurs d'Atkin-Lehner. Notons $\psi'$ le caractère de nebentypus de $F$. Notons $M'$ le conducteur de $\psi'$. Supposons que $M$ soit à support dans $\Sigma_N$. Soit $S$ un sous-ensemble de $\Sigma_M$. Notons $\bar{S} = \Sigma_M - S$. Posons $M = M_S M_{\bar{S}}$ et $M' = M'_S M'_{\bar{S}}$ et $\psi' = \psi'_S \psi'_{\bar{S}}$. Soit $\begin{pmatrix} A & B \\ C & D \end{pmatrix} \in \mathrm{M}_2(\mathbf{Z})$ telle que $M_S | A$, $M_S | D$, $M | C$, $M_{\bar{S}} | B$, $AD - BC = M_S$, $A \equiv M_S \pmod{M'}$ et $B \equiv 1$ (mod $M'_S$); on pose alors, comme Atkin et Li dans [Atkin-Li], $W_S F = F_{|}\begin{pmatrix} A & B \\ C & D \end{pmatrix}$ et il existe un nombre complexe $w_S(F)$ de module 1 tel que $W_S F = w_S(F) F \otimes \bar\psi'_S$. Lorsque $\begin{pmatrix} A & B \\ C & D \end{pmatrix} \in \mathrm{M}_2(\mathbf{Z})$ avec $M_S|A$, $M_S|D$, $M|C$, $M_{\bar S}|B$ et $AD - BC = M_S$, on a de plus [A-L]

$$F_{|}\begin{pmatrix} A & B \\ C & D \end{pmatrix} = \psi'_S(B)\psi'_{\bar S}(A/M_S) W_S F.$$

Lorsque $M|NN'$, $M'|N$ et lorsque $\begin{pmatrix} A & B \\ C & D \end{pmatrix} \in \mathrm{M}_2(\mathbf{Z})$ vérifie les conditions $N_S N'_S|A$, $N_S N'_S|D$, $NN'|C$, $N_{\bar S} N'_{\bar S}|B$ et $AD - BC = N_S N'_S$, on a

$$F_{|}\begin{pmatrix} A & B \\ C & D \end{pmatrix} = w_S(F) \bar\psi'_S(B) \bar\psi'_{\bar S}(A/(N_S N'_S)) F_{|}\begin{pmatrix} N_S N'_S/M_S & 0 \\ 0 & 1 \end{pmatrix}.$$

Lorsque $S$ est égal au support de $M$, on pose $w_S(F) = w(F)$.

On a de plus [A-L, proposition 1.1],

$$w_S(F) w_S(F \otimes \bar\psi'_S) = \psi'_S(-1) \bar\psi'_{\bar S}(M_S).$$

Mentionnons enfin la formule, pour $S_1$ et $S_2$ deux sous-ensembles disjoints de $\Sigma_M$,

$$W_{S_2}(W_{S_1} F) = \psi'_{S_2}(M_{S_1}) W_{S_1 \cup S_2} F.$$

Cela permet de ramener le calcul de $w_S(F)$ aux cas où $S$ est un singleton.

Ajoutons la formule suivante. Soit $p$ un nombre premier tel que $a_p(F) \neq 0$ (c'est le cas si et seulement si $v_p(M) = v_p(m_\psi)$ ou si $v_p(M) \leq 1$). On a

$$w_{\{p\}}(F) = \frac{p^{v_p(N)(k/2-1)} \tau(\psi'_S)}{a_p(F)^{v_p(N)}}$$

où $\tau(\psi'_S)$ est la somme de Gauss du caractère (non nécessairement primitif) $\psi'_S$. Si $p$ est spécial pour $F$, on a $a_p(F)\bar a_p(F) = p^{k-2}$.

**2.c.** Revenons maintenant sur la torsion des formes modulaires par des caractères. Soit $\alpha$ un caractère de Dirichlet de niveau $m$, de caractère de Dirichlet primitif associé $\omega$, lui-même de conducteur $m_\omega$. Notons $\bar f_\alpha$ la forme modulaire (non nécessairement primitive) donnée par le développement

$$\bar f_\alpha(z) = \sum_{n=1}^\infty \bar a_n \alpha(n) q^n.$$

Elle est liée à la forme primitive $\bar f \otimes \omega$ par la formule

$$\bar f_\alpha = (\bar f \otimes \omega)^{[\Sigma_m, \emptyset]}.$$

Posons de plus

$$S_\alpha \bar f = \sum_{a \bmod m} \bar\alpha(a) \bar f_{|}\begin{pmatrix} 1 & a/m \\ 0 & 1 \end{pmatrix}.$$

On a, lorsque $\alpha$ est primitif (et donc égal à $\omega$),

$$S_\omega \bar f = \tau(\bar\omega) \bar f_\omega.$$

Soit $p$ un nombre premier divisant $m/m_\omega$. Notons $\beta$ le caractère de Dirichlet de niveau $m/p$ qui coïncide avec $\alpha$ sur les entiers premiers à $p$. On a

$$S_\alpha \bar{f} = \bar{a}_p p^{1-k/2}(S_\beta \bar{f})_{|\begin{pmatrix} p & 0 \\ 0 & 1 \end{pmatrix}} - \bar{\beta}(p) S_\beta \bar{f}.$$

Posons, dans $\mathbf{C}[X]$, $R_p(X) = (\bar{a}_p p^{1-k/2} X)^{v_p(m/m_\omega)-1}(\bar{a}_p p^{1-k/2} X - \bar{\omega}(p))$. Par une application répétée de la formule ci-dessus, on obtient

$$S_\alpha \bar{f} = \tau(\bar{\omega})(\bar{f}_\omega)_{|\prod_p R_p(\begin{pmatrix} p & 0 \\ 0 & 1 \end{pmatrix})},$$

où le produit porte sur les nombres premiers divisant $m/m_\omega$.

Il est nécessaire maintenant de distinguer plusieurs cas. Si $v_p(m/m_\omega) = 0$, on a $R_p = 1$. Si $v_p(m/m_\omega) = 1$ et $a_p = 0$, on a $R_p = 0$. Si $v_p(m/m_\omega) > 1$ et $a_p = 0$, on a $R_p = -\bar{\omega}(p)$.

Or on a, lorsque $a_p \neq 0$ et $p|N$ non spécial pour $\bar{f}$, $a_p \bar{a}_p = p^{k-1}$ et donc, lorsque de plus $p|m$ on a, dans $\mathbf{C}[\mathrm{PGL}_2(\mathbf{Q})^+]$,

$$R_p(\begin{pmatrix} p & 0 \\ 0 & 1 \end{pmatrix}) = (\bar{a}_p p^{1-k/2}\begin{pmatrix} p & 0 \\ 0 & 1 \end{pmatrix})^{v_p(m/m_\omega)}(1 - \bar{\omega}(p) a_p p^{-k/2}\begin{pmatrix} 1 & 0 \\ 0 & p \end{pmatrix}).$$

Cette dernière formule est encore valable lorsque $a_p = 0$ et $v_p(m) > 1$.

Lorsque $p|(m/m_\omega)$ et $p$ est spécial pour $\bar{f}$, on a $a_p \bar{a}_p = p^{k-2}$ (et donc $a_p \neq 0$). On a donc

$$R_p(\begin{pmatrix} p & 0 \\ 0 & 1 \end{pmatrix}) = (\bar{a}_p p^{1-k/2}\begin{pmatrix} p & 0 \\ 0 & 1 \end{pmatrix})^{v_p(m/m_\omega)}(1 - \bar{\omega}(p) a_p p^{1-k/2}\begin{pmatrix} 1 & 0 \\ 0 & p \end{pmatrix}).$$

Lorsque $a_p = 0$, $v_p(m) = 1$ et $v_p(m_\omega) = 0$, on a

$$R_p(\begin{pmatrix} p & 0 \\ 0 & 1 \end{pmatrix}) = -\bar{\omega}(p).$$

On a donc

$$S_\alpha \bar{f} = \tau(\bar{\omega})(\bar{f} \otimes \omega)^{[\Sigma_m, \frac{\Sigma_{m/m_\omega}}{\Sigma_{m/m_\omega} \cap \Sigma_f}]}_{|\prod_p P_p(\begin{pmatrix} p & 0 \\ 0 & 1 \end{pmatrix})},$$

où le monôme $P_p(X)$ vaut $(\bar{a}_p p^{1-k/2} X)^{v_p(m/m_\omega)}$ sauf si $a_p = 0$, $v_p(m) = 1$ et $v_p(m_\chi) = 0$, auquel cas on a $P_p(\begin{pmatrix} p & 0 \\ 0 & 1 \end{pmatrix}) = -\bar{\omega}(p)$.

**2.d.** Reprenons la situation laissée en **2.c** en nous plaçant dans le cas où $N' = m$ est un diviseur de $N$.

*Lemme .* — *Soit $p$ un nombre premier tel que $p|m_\omega$ et $p|(N'/m_\omega)$. On a*

$$(\bar{f} \otimes \omega)^{[\emptyset, p]}_{|P_p(\begin{pmatrix} p & 0 \\ 0 & 1 \end{pmatrix})} = (\bar{f} \otimes \omega)_{|P_p(\begin{pmatrix} p & 0 \\ 0 & 1 \end{pmatrix})}.$$

*Démonstration.* — Il suffit de montrer que $P_p = 0$ ou que $a_p(\bar{f} \otimes \omega) = 0 = a_{p,p}(\bar{f} \otimes \omega)$. Supposons $P_p \neq 0$. Si $a_p = 0$, on a $v_p(m_\omega) = 0$ et $v_p(N') = 1$, ce qui entraîne $a_p(\bar{f} \otimes \omega) = 0 = a_{p,p}(\bar{f} \otimes \omega)$. Restreignons maintenant notre attention au cas où $a_p \neq 0$. Rappelons d'abord que cela entraîne que le conducteur de $\psi$ a pour valuation $p$-adique $v_p(N)$ (ce qui entraîne $a_{p,p}(\bar{f} \otimes \omega) = 0$) ou que $v_p(N) = 1$. Les hypothèses excluent le cas $v_p(N) = 1$. On a de plus $a_p(\bar{f} \otimes \omega) \neq 0$ si et seulement si $\omega$ est de conducteur premier à $p$ (impossible par hypothèse) ou $\bar{\psi}/\omega$ est de conducteur premier à $p$ ; ce dernier cas est impossible, en effet on a $p|(N'/m_\omega)$, et donc $p|(N/m_\omega)$ et les valuations $p$-adiques des conducteurs de $\psi$ et $\bar\psi/\omega$ sont égales et donc non nulles. On a bien $a_p(\bar f \otimes \omega) = 0$.

On a donc
$$S_\alpha \bar f = \tau(\bar\omega)(\bar f \otimes \omega)_{|\prod_p P_p(\begin{pmatrix} p & 0 \\ 0 & 1 \end{pmatrix})}^{[\Sigma_{N'}, \frac{\Sigma_{N'}-\Sigma_{m_\omega}}{(\Sigma_{N'}-\Sigma_{m_\omega})\cap\Sigma_f}]},$$

**2.e.** Soit $n \in \mathbf{Z}$. Récrivons la forme modulaire $\bar f_{|\begin{pmatrix} 1 & n/N \\ 0 & 1 \end{pmatrix}}$ comme combinaison linéaire de $F_{|\begin{pmatrix} d & 0 \\ 0 & 1 \end{pmatrix}}$ où $d$ parcourt les diviseurs de $N$ et où $F$ parcourt les formes primitives de niveau divisant $N^2/d$. Nous ne savons pas si un pareil calcul a déjà été rédigé. Notons $n'$ le nombre entier et $N'$ le diviseur $>0$ de $N$ qui vérifient $n'/N' = n/N$. On a, par inversion de Fourier,

$$\bar f_{|\begin{pmatrix} 1 & n'/N' \\ 0 & 1 \end{pmatrix}} = \sum_\alpha \frac{\alpha(n')}{\phi(N')} S_\alpha \bar f,$$

où $\alpha$ parcourt les caractères de Dirichlet de niveau $N'$. En combinant avec la formule trouvée en **2.d.**, on obtient

$$\bar f_{|\begin{pmatrix} 1 & n/N \\ 0 & 1 \end{pmatrix}} = \sum_\omega \frac{\omega(n')}{\phi(N')} \tau(\bar\omega)(\bar f \otimes \omega)_{|\prod_p P_p(\begin{pmatrix} p & 0 \\ 0 & 1 \end{pmatrix})}^{[\Sigma_{N'}, \frac{\Sigma_{N'}-\Sigma_{m_\omega}}{(\Sigma_{N'}-\Sigma_{m_\omega})\cap\Sigma_f}]}$$

où $\omega$ parcourt les caractères de Dirichlet primitifs de conducteur $m_\omega$ divisant $N'$, le produit portant sur les nombres premiers divisant $N'/m_\omega$.

### 3. La démonstration du théorème A

Soit $g = \begin{pmatrix} a & b \\ c & d \end{pmatrix} \in \mathrm{SL}_2(\mathbf{Z})$ telle que la classe modulo $N$ de $(c,d)$ soit égale à $(u,v)$.

On peut comprendre notre démarche ainsi. La fonction $f_{|g}$ est une forme modulaire pour le groupe de congruence $\Gamma(N)$, si bien que la fonction $f_{|g\begin{pmatrix} N & 0 \\ 0 & 1 \end{pmatrix}}$ est modulaire pour le groupe de congruence $\Gamma_1(N) \cap \Gamma_0(N^2)$. Cette dernière forme modulaire s'écrit donc comme combinaison linéaire de fonctions du type $F_{|\begin{pmatrix} d & 0 \\ 0 & 1 \end{pmatrix}}$, où $F$ parcourt les formes primitives de niveau $M$ divisant $N^2$ et $d$ les entiers divisant $N^2/M$. Nous allons montrer que les formes primitives qui interviennent dans cette écriture sont de la forme $f \otimes \chi$, où $\chi$ parcourt les caractères de Dirichlet de niveau divisant $N$ et donner explicitement les coefficients de cette combinaison linéaire.

Lorsque $k=2$, on a $\xi_f(u,v) = \int_0^\infty f_{|g}(iy)\,dy$. Lorsque, de plus, $s=1$ et $h$ est une matrice diagonale de $\mathrm{PGL}_2(\mathbf{Q})^+$, et $\chi$ est un caractère de Dirichlet de conducteur divisant $N$, on a

$$\int_0^\infty (f \otimes \chi)_{|h}(iy)\,dy = \frac{1}{2\pi} L(f \otimes \chi, 1) = \frac{1}{\sqrt{N_\chi}} \Lambda(f \otimes \chi, 1).$$

C'est pourquoi le théorème A se déduit de la proposition B suivante, par intégration de chaque membre de l'égalité ci-dessous le long de la géodésique reliant 0 à $\infty$ dans le demi-plan de Poincaré.

Remarquons que la proposition B permet de démontrer des analogues du théorème A pour les formes modulaires de poids $\neq 2$.

Proposition B. — *On a*

$$f_{|g} = \frac{w(f)}{\phi(N')} \sum_\chi \chi_{\bar S}(m_{\chi,S})(\bar\psi_S \bar\chi_S)(m_{\psi\chi,\bar S})(\psi_S \chi_S)(-1)\tau'(\chi_S)\tau'(\bar\psi_{\bar S}\bar\chi_{\bar S})(\bar\psi\bar\chi)(\frac{N'_{\bar S} v}{N_{\bar S}})\chi(\frac{N'_S u}{N_S})(\bar\psi_{\bar S}\bar\chi_{\bar S}^2)(N_{\chi,S})$$

$$\overline{w_S(f\otimes\chi)}(f\otimes\chi)^{[\frac{\Sigma_{N'}-S(m_\chi)}{(S(N')-S(m_\chi))\cap\Sigma_f},\frac{\Sigma_{N'}-\bar S(m_{\bar\psi_{\bar\chi}})}{(\bar S(N')-\bar S(m_{\bar\psi_{\bar\chi}}))\cap\Sigma_f}]}_{|\begin{pmatrix} \frac{N'_{\bar S}}{N_{\chi,S}N_{\bar S}m_{\psi_\chi,\bar S}} & 0 \\ 0 & \frac{N'_S}{N_S m_{\bar\chi,S}} \end{pmatrix}} \prod_{p\in S(N')} Q_{p,f,\bar\chi}(\begin{pmatrix} 1 & 0 \\ 0 & p \end{pmatrix}) \prod_{p\in \bar S(N')} Q_{p,f,\chi\psi}(\begin{pmatrix} p & 0 \\ 0 & 1 \end{pmatrix}),$$

où $\chi$ parcourt les caractères de Dirichlet primitifs tel que $m_{\chi,S}m_{\psi_\chi,\bar S}|N'$.

*Démonstration.* — Considérons $\begin{pmatrix} A & B \\ C & D \end{pmatrix} \in M_2(\mathbf{Z})$ telle que $N_S N'_S|A$, $N_S N'_S|D$, $NN'|C$, $N_{\bar S}N'_{\bar S}|B$, $AD-BC = N_S N'_S$, $A \equiv uN'_S \pmod{N_{\bar S}}$ et $B \equiv v/N_{\bar S} \pmod{N_S}$. Soit $k \in \mathbf{Z}$ tel que $n \equiv uv \pmod{N_{\bar S}}$ et $n \equiv -uv \pmod{N_S}$. Notre point de départ réside dans l'identité

$$\Gamma_1(N)g = \Gamma_1(N)\begin{pmatrix} 0 & -1 \\ N & 0 \end{pmatrix}\begin{pmatrix} 1 & n/N \\ 0 & 1 \end{pmatrix}\begin{pmatrix} A & B \\ C & D \end{pmatrix}\begin{pmatrix} NN'_S & 0 \\ 0 & N_S \end{pmatrix}^{-1},$$

que le lecteur vérifiera grâce au lemme chinois. Comme $w(f)\bar f = f_{|\begin{pmatrix} 0 & 1 \\ -N & 0 \end{pmatrix}}$, on a la formule

$$f_{|g} = w(f)\bar f_{|\begin{pmatrix} 1 & n/N \\ 0 & 1 \end{pmatrix}\begin{pmatrix} A & B \\ C & D \end{pmatrix}\begin{pmatrix} NN'_S & 0 \\ 0 & N_S \end{pmatrix}^{-1}},$$

et donc, d'après la formule trouvée en **2.e.**,

$$f_{|g} = w(f)\sum_\omega \frac{\omega(n')}{\phi(N')}\tau(\bar\omega)(\bar f\otimes\omega)^{[\Sigma_{N'},\frac{\Sigma_{N'}-\Sigma_{m_\omega}}{(\Sigma_{N'}-\Sigma_{m_\omega})\cap\Sigma_f}]}_{|\prod_p P_p(\begin{pmatrix} p & 0 \\ 0 & 1 \end{pmatrix})\begin{pmatrix} A & B \\ C & D \end{pmatrix}\begin{pmatrix} NN'_S & 0 \\ 0 & N_S \end{pmatrix}^{-1}}$$

où $\omega$ parcourt les caractères de Dirichlet primitifs de conducteur $m_\omega$ divisant $N'$, le produit portant sur les nombres premiers divisant $N'$. Appliquons les formules de **2.b.** à $F = \bar f \otimes \omega$ ; on a $M = N_\omega$, $\psi' = \bar\psi\omega^2$ et $F\otimes\bar\psi'_S = \bar f\otimes\omega\psi_S\bar\omega_S^2 = \bar f\otimes\psi_S\bar\omega_S\omega_{\bar S} = \overline{f\otimes\bar\psi_{\bar S}\omega_S\bar\omega_{\bar S}}$.

Soit $p \in \Sigma_N$. Soit $r$ un entier $\geq 0$. On a

$$(\bar f\otimes\omega)_{|\begin{pmatrix} p^r & 0 \\ 0 & 1 \end{pmatrix}\begin{pmatrix} A & B \\ C & D \end{pmatrix}} = (\bar f\otimes\omega)_{|\begin{pmatrix} p^r A & B \\ C & D/p^r \end{pmatrix}\begin{pmatrix} 1 & 0 \\ 0 & p^r \end{pmatrix}}$$

si $p \in S$ et

$$(\bar f\otimes\omega)_{|\begin{pmatrix} p^r & 0 \\ 0 & 1 \end{pmatrix}\begin{pmatrix} A & B \\ C & D \end{pmatrix}} = (\bar f\otimes\omega)_{|\begin{pmatrix} A & p^r B \\ C/p^r & D \end{pmatrix}\begin{pmatrix} p^r & 0 \\ 0 & 1 \end{pmatrix}}$$

si $p \in \bar S$. On a de plus les formules

$$(\bar f\otimes\omega)_{|\begin{pmatrix} p^r A & B \\ C & D/p^r \end{pmatrix}} = (\psi_S\bar\omega_S^2)(B)(\psi_{\bar S}\bar\omega_{\bar S}^2)(p^r A/(N_S N'_S))w_S(\bar f\otimes\omega)(\bar f\otimes\omega_{\bar S}\bar\omega_S\psi_S)_{|\begin{pmatrix} N_S N'_S/N_{\bar\omega,S} & 0 \\ 0 & 1 \end{pmatrix}}$$

lorsque $p^r|(N_S N'_S/N_{\omega,S})$ et

$$(\bar f\otimes\omega)_{|\begin{pmatrix} A & p^r B \\ C/p^r & D \end{pmatrix}} = (\psi_S\bar\omega_S^2)(p^r B)(\psi_{\bar S}\bar\omega_{\bar S}^2)(A/(N_S N'_S))w_S(\bar f\otimes\omega)(\bar f\otimes\omega_{\bar S}\bar\omega_S\psi_S)_{|\begin{pmatrix} N_S N'_S/N_{\bar\omega,S} & 0 \\ 0 & 1 \end{pmatrix}}$$

lorsque $p^r|(N_{\bar S}N'_{\bar S}/N_{\omega,\bar S})$. Soit $P \in \mathbf{C}[X]$. On a alors $(\bar f\otimes\omega)_{|P(\begin{pmatrix} p & 0 \\ 0 & 1 \end{pmatrix})\begin{pmatrix} A & B \\ C & D \end{pmatrix}} =$

$$(\psi_S\bar\omega_S^2)(B)(\psi_{\bar S}\bar\omega_{\bar S}^2)(A/(N_S N'_S))w_S(\bar f\otimes\omega)(\bar f\otimes\omega_{\bar S}\bar\omega_S\psi_S)_{|P((\psi_{\bar S}\bar\omega_{\bar S}^2)(p))\begin{pmatrix} 1 & 0 \\ 0 & p \end{pmatrix})\begin{pmatrix} N_S N'_S/N_{\bar\omega,S} & 0 \\ 0 & 1 \end{pmatrix}}.$$

si $p \in S$ et $P$ de degré $\leq v_p(N_S N'_S/N_{\omega,S})$ et on a $(\bar{f} \otimes \omega)_{|P(\begin{pmatrix} p & 0 \\ 0 & 1 \end{pmatrix})}\begin{pmatrix} A & B \\ C & D \end{pmatrix} =$

$$(\psi_S \bar{\omega}_S^2)(B)(\psi_{\bar{S}} \bar{\omega}_{\bar{S}}^2)(A/(N_S N'_S)) w_S(\bar{f} \otimes \omega)(\bar{f} \otimes \omega_{\bar{S}} \bar{\omega}_S \psi_S)_{|P((\psi_S \bar{\omega}_S^2)(p)\begin{pmatrix} p & 0 \\ 0 & 1 \end{pmatrix})}\begin{pmatrix} N_S N'_S/N_{\bar{\omega},S} & 0 \\ 0 & 1 \end{pmatrix}.$$

si $p \in \bar{S}$ et $P$ de degré $\leq v_p(N_{\bar{S}} N'_{\bar{S}}/N_{\omega,\bar{S}})$.

On en déduit que

$$(\bar{f} \otimes \omega)^{[\{p\},\emptyset]}_{|\begin{pmatrix} A & B \\ C & D \end{pmatrix}} = (\psi_S \bar{\omega}_S^2)(B)(\psi_{\bar{S}} \bar{\omega}_{\bar{S}}^2)(A/(N_S N'_S)) w_S(\bar{f} \otimes \omega)(\bar{f} \otimes \omega_{\bar{S}} \bar{\omega}_S \psi_S)^{[\emptyset,\{p\}]}_{|\begin{pmatrix} N_S N'_S/N_{\bar{\omega},S} & 0 \\ 0 & 1 \end{pmatrix}}$$

si et $p \in S$ et

$$(\bar{f} \otimes \omega)^{[\{p\},\emptyset]}| \begin{pmatrix} A & B \\ C & D \end{pmatrix} = (\psi_S \bar{\omega}_S^2)(B)(\psi_{\bar{S}} \bar{\omega}_{\bar{S}}^2)(A/(N_S N'_S)) w_S(\bar{f} \otimes \omega)(\bar{f} \otimes \omega_{\bar{S}} \bar{\omega}_S \psi_S)^{[\{p\},\emptyset]}_{|\begin{pmatrix} N_S N'_S/N_{\bar{\omega},S} & 0 \\ 0 & 1 \end{pmatrix}}$$

si $p \in \bar{S}$. Un calcul analogue donne les formules

$$(\bar{f} \otimes \omega)^{[\emptyset,\{p\}]}_{|\begin{pmatrix} A & B \\ C & D \end{pmatrix}} = (\psi_S \bar{\omega}_S^2)(B)(\psi_{\bar{S}} \bar{\omega}_{\bar{S}}^2)(A/(N_S N'_S)) w_S(\bar{f} \otimes \omega)(\bar{f} \otimes \omega_{\bar{S}} \bar{\omega}_S \psi_S)^{[\{p\},\emptyset]}_{|\begin{pmatrix} N_S N'_S/N_{\bar{\omega},S} & 0 \\ 0 & 1 \end{pmatrix}}$$

si et $p \in S$ et

$$(\bar{f} \otimes \omega)^{[\emptyset,\{p\}]}| \begin{pmatrix} A & B \\ C & D \end{pmatrix} = (\psi_S \bar{\omega}_S^2)(B)(\psi_{\bar{S}} \bar{\omega}_{\bar{S}}^2)(A/(N_S N'_S)) w_S(\bar{f} \otimes \omega)(\bar{f} \otimes \omega_{\bar{S}} \bar{\omega}_S \psi_S)^{[\emptyset,\{p\}]}_{|\begin{pmatrix} N_S N'_S/N_{\bar{\omega},S} & 0 \\ 0 & 1 \end{pmatrix}}$$

si $p \in \bar{S}$.

Dans les quatre formules qui précèdent, on peut remplacer, partout où il intervient, le symbole $\{p\}$ par $\frac{\{p\}}{\{p\} \cap \Sigma_f}$.

Remarquons qu'on a, dans $\mathbf{C}(X)$, $Q_{p,f,\omega}(X) = X^{-v_p(N'/m_\chi)} P_p(X)$. On a

$$\prod_{p \in S(N')} P_p((\psi_{\bar{S}} \bar{\omega}_{\bar{S}}^2)(p)\begin{pmatrix} 1 & 0 \\ 0 & p \end{pmatrix}) \prod_{p \in \bar{S}(N')} P_p((\psi_S \bar{\omega}_S^2)(p)\begin{pmatrix} p & 0 \\ 0 & 1 \end{pmatrix}) = (\psi_{\bar{S}} \bar{\omega}_{\bar{S}}^2)(N'_S/m_{\omega,S})(\psi_S \bar{\omega}_S^2)(N'_{\bar{S}}/m_{\omega,\bar{S}})$$

$$\begin{pmatrix} N'_{\bar{S}}/m_{\omega,S} & 0 \\ 0 & N'_S/m_{\omega,\bar{S}} \end{pmatrix} \prod_{p \in S(N')} Q_{p,f,\omega}((\psi_{\bar{S}} \bar{\omega}_{\bar{S}}^2)(p)\begin{pmatrix} 1 & 0 \\ 0 & p \end{pmatrix}) \prod_{p \in \bar{S}(N')} Q_{p,f,\omega}((\psi_S \bar{\omega}_S^2)(p)\begin{pmatrix} p & 0 \\ 0 & 1 \end{pmatrix}).$$

Revenons à notre calcul principal. On a

$$(\bar{f} \otimes \omega)^{[\Sigma_{N'}, \frac{\Sigma_{N'} - \Sigma_{m_\omega}}{(\Sigma_{N'} - \Sigma_{m_\omega}) \cap \Sigma_f}]}_{|\prod_p P_p(\begin{pmatrix} p & 0 \\ 0 & 1 \end{pmatrix})\begin{pmatrix} A & B \\ C & D \end{pmatrix}\begin{pmatrix} NN'_S & 0 \\ 0 & N_S \end{pmatrix}^{-1}} = (\psi_S \bar{\omega}_S^2)(N'_{\bar{S}} B/m_{\omega,\bar{S}})(\psi_{\bar{S}} \bar{\omega}_{\bar{S}}^2)(A/(N_S m_{\omega,S})) w_S(\bar{f} \otimes \omega)$$

$$(\bar{f} \otimes \omega_{\bar{S}} \bar{\omega}_S \psi_S)^{[\frac{\Sigma_{N'} - S(m_\chi)}{(S(N') - S(m_\chi)) \cap \Sigma_f}, \frac{\Sigma_{N'} - \bar{S}(m_{\bar{\psi}\bar{\chi}})}{(\bar{S}(N') - \bar{S}(m_{\bar{\psi}\bar{\chi}})) \cap \Sigma_f}]}_{|\begin{pmatrix} \frac{N'_{\bar{S}}}{N_{\bar{\omega},S} N_{\bar{S}} m_{\omega,\bar{S}}} & 0 \\ 0 & \frac{N'_S}{N_S m_{\omega,S}} \end{pmatrix} \prod_{p \in S(N')} Q_{p,f,\omega}((\psi_{\bar{S}} \bar{\omega}_{\bar{S}}^2)(p)\begin{pmatrix} 1 & 0 \\ 0 & p \end{pmatrix}) \prod_{p \in \bar{S}(N')} Q_{p,f,\omega}((\psi_S \bar{\omega}_S^2)(p)\begin{pmatrix} p & 0 \\ 0 & 1 \end{pmatrix})}.$$

Par ailleurs, on a $(\psi_S\bar{\omega}_S^2)(N'_{\bar{S}}B) = (\psi_S\bar{\omega}_S^2)(N'_{\bar{S}}v/N_{\bar{S}})$, $(\psi_{\bar{S}}\bar{\omega}_{\bar{S}}^2)(A/N_S) = (\psi_{\bar{S}}\bar{\omega}_{\bar{S}}^2)(uN'_S/N_S)$ et $\omega(n') = \omega(nN'/N) = \omega_S(nN'/N)\omega_{\bar{S}}(nN'/N) = \omega_S(-uvN'/N)\omega_{\bar{S}}(uvN'/N)$ et donc

$$\omega(n') = \omega_S(-1)\omega(uN'_S/N_S)\omega(vN'_{\bar{S}}/N_{\bar{S}}).$$

On a donc la simplification

$$\omega(n')(\psi_S\bar{\omega}_S^2)(N'_{\bar{S}}B)(\psi_{\bar{S}}\bar{\omega}_{\bar{S}}^2)(A/N_S) = \omega_S(-1)\psi_S\bar{\omega}_S\omega_{\bar{S}}(N'_{\bar{S}}v/N_{\bar{S}})\psi_{\bar{S}}\bar{\omega}_{\bar{S}}\omega_S(N'_Su/N_S).$$

De plus on a

$$Q_{p,f,\omega}((\psi_{\bar{S}}\bar{\omega}_{\bar{S}}^2)(p)\begin{pmatrix} 1 & 0 \\ 0 & p \end{pmatrix}) = Q_{p,f,\bar{\omega}_{\bar{S}}\omega_S\psi_S}(\begin{pmatrix} 1 & 0 \\ 0 & p \end{pmatrix})$$

si $p \in S$ et

$$Q_{p,f,\omega}((\psi_S\bar{\omega}_S^2)(p)\begin{pmatrix} p & 0 \\ 0 & 1 \end{pmatrix}) = Q_{p,f,\omega_{\bar{S}}\bar{\omega}_S\psi_S}(\begin{pmatrix} p & 0 \\ 0 & 1 \end{pmatrix})$$

si $p \in \bar{S}$.

En combinant ces formules, on obtient

$$f_{|g} = \frac{w(f)}{\phi(N')}\sum_\omega \tau(\bar{\omega})\omega_S(-1)(\psi_S\bar{\omega}_S\omega_{\bar{S}})(\frac{N'_{\bar{S}}v}{N_{\bar{S}}})(\psi_{\bar{S}}\bar{\omega}_{\bar{S}}\omega_S)(\frac{N'_Su}{N_S})w_S(\bar{f}\otimes\omega)(\bar{\psi}_{\bar{S}}\omega_{\bar{S}}^2)(m_{\omega,S})(\bar{\psi}_S\omega_S^2)(m_{\omega,\bar{S}})$$

$$(\bar{f}\otimes\omega_{\bar{S}}\bar{\omega}_S\psi_S)_{|\begin{pmatrix} \frac{N'_{\bar{S}}}{N_{\bar{\omega},S}N_{\bar{S}}m_{\omega,\bar{S}}} & 0 \\ 0 & \frac{N'_S}{N_Sm_{\omega,S}} \end{pmatrix}}^{[\frac{\Sigma_{N'}-S(m_\chi)}{(S(N')-S(m_\chi))\cap\Sigma_f},\frac{\Sigma_{N'}-\bar{S}(m_{\bar{\psi}\bar{\chi}})}{(\bar{S}(N')-\bar{S}(m_{\bar{\psi}\bar{\chi}}))\cap\Sigma_f}]} \prod_{p\in S(N')} Q_{p,f,\bar{\omega}_{\bar{S}}\omega_S\psi_S}(\begin{pmatrix} 1 & 0 \\ 0 & p \end{pmatrix})\prod_{p\in\bar{S}(N')} Q_{p,f,\omega_{\bar{S}}\bar{\omega}_S\psi_S}(\begin{pmatrix} p & 0 \\ 0 & 1 \end{pmatrix}),$$

où $\omega$ parcourt les caractères de Dirichlet primitifs de conducteur divisant $N'$. Simplifions encore cette formule.

On a la relation entre sommes de Gauss

$$\tau(\bar{\omega}) = \bar{\omega}_S(m_{\bar{\omega},\bar{S}})\bar{\omega}_{\bar{S}}(m_{\bar{\omega},S})\tau(\bar{\omega}_S)\tau(\bar{\omega}_{\bar{S}}).$$

Cela donne

$$\tau(\bar{\omega})(\bar{\psi}_{\bar{S}}\omega_{\bar{S}}^2)(m_{\omega,S})(\bar{\psi}_S\omega_S^2)(m_{\omega,\bar{S}}) = \tau(\bar{\omega}_S)\tau(\bar{\omega}_{\bar{S}})(\bar{\psi}_{\bar{S}}\omega_{\bar{S}})(m_{\omega,S})(\bar{\psi}_S\omega_S)(m_{\omega,\bar{S}}).$$

Récrivons notre formule en notant $\chi$ le caractère de Dirichlet primitif associé à $\omega_{\bar{S}}\bar{\omega}_S\bar{\psi}_S$. On a donc $\chi_S = \bar{\omega}_S$ et $\chi_{\bar{S}} = \omega_{\bar{S}}\bar{\psi}_{\bar{S}}$, $S(m_\omega) = S(m_\chi)$, $\bar{S}(m_\omega) = \bar{S}(m_{\psi\chi})$, $N_{\bar{\omega},S} = N_{\chi,S}$ et $\omega_S(-1) = \chi_S(-1)$.

On obtient

$$f_{|g} = \frac{w(f)}{\phi(N')}\sum_\chi \tau'(\chi_S)\tau'(\bar{\psi}_{\bar{S}}\bar{\chi}_{\bar{S}})\chi_S(-1)(\bar{\psi}\bar{\chi})(\frac{N'_{\bar{S}}v}{N_{\bar{S}}})\chi(\frac{N'_Su}{N_S})w_S(f\otimes\bar{\psi}_S\bar{\chi}_S\chi_{\bar{S}})\chi_{\bar{S}}(m_{\chi,S})(\bar{\psi}_S\bar{\chi}_S)(m_{\psi\chi,\bar{S}})$$

$$(f\otimes\chi)_{|\begin{pmatrix} \frac{N'_{\bar{S}}}{N_{\chi,S}N_{\bar{S}}m_{\psi\chi,\bar{S}}} & 0 \\ 0 & \frac{N'_S}{N_Sm_{\bar{\chi},S}} \end{pmatrix}}^{[\frac{\Sigma_{N'}-S(m_\chi)}{(S(N')-S(m_\chi))\cap\Sigma_f},\frac{\Sigma_{N'}-\bar{S}(m_{\bar{\psi}\bar{\chi}})}{(\bar{S}(N')-\bar{S}(m_{\bar{\psi}\bar{\chi}}))\cap\Sigma_f}]} \prod_{p\in S(N')} Q_{p,f,\bar{\chi}}(\begin{pmatrix} 1 & 0 \\ 0 & p \end{pmatrix})\prod_{p\in\bar{S}(N')} Q_{p,f,\chi\psi}(\begin{pmatrix} p & 0 \\ 0 & 1 \end{pmatrix}),$$

où $\chi$ parcourt les caractères de Dirichlet primitifs de conducteur divisant $N'$. Appliquons la relation reliant $w_S(F)$ et $w_S(F\otimes\bar{\psi}'_S)$ (voir **2.b**), pour $F = f\otimes\chi$ (et donc $\psi' = \psi\chi^2$). On obtient

$$w_S(f\otimes\chi)w_S(f\otimes\bar{\psi}_S\bar{\chi}_S\chi_{\bar{S}}) = \psi_S(-1)(\bar{\psi}_{\bar{S}}\bar{\chi}_{\bar{S}}^2)(N_{\chi,S}).$$

Cela permet de substituer $w_S(f\otimes\bar{\psi}_S\bar{\chi}_S\chi_{\bar{S}})$ pour obtenir la proposition B.

**4. Le produit scalaire de Petersson**

Soit $\Gamma$ un sous-groupe d'indice fini de $\mathrm{SL}_2(\mathbf{Z})$ contenant la matrice $\begin{pmatrix} -1 & -0 \\ 0 & -1 \end{pmatrix}$. Soient $f_1$ et $f_2$ deux formes modulaires paraboliques de poids 2 pour $\Gamma$. Rappelons que le produit scalaire de Petersson de $f_1$ et $f_2$ est donné par la formule :

$$<f_1, f_2> = \frac{1}{[\mathrm{SL}_2(\mathbf{Z}):\Gamma]} \int_{D_\Gamma} f_1(z)\overline{f}_2(z) dx\, dy,$$

où $D_\Gamma$ est un domaine fondamental pour $\Gamma$ dans le demi-plan de Poincaré $\mathbf{H}$. Posons $\tau = \begin{pmatrix} 0 & -1 \\ 1 & -1 \end{pmatrix}$ et $\sigma = \begin{pmatrix} 0 & 0 \\ 1 & -1 \end{pmatrix}$. Posons de plus $\rho = e^{2i\pi/3} \in \mathbf{H}$. Soit $R$ un système de représentants de $\Gamma \backslash \mathrm{SL}_2(\mathbf{Z})$.

THÉORÈME C. — *On a*

$$<f_1, f_2> = \frac{1}{2i[\mathrm{SL}_2(\mathbf{Z}):\Gamma]} \sum_{g \in R} \int_{g0}^{g\infty} f_1(z)\, dz \int_{gi}^{g\rho} \overline{f_2(z)\, dz},$$

*et*

$$<f_1, f_2> = \frac{-i}{12[\mathrm{SL}_2(\mathbf{Z}):\Gamma]} \sum_{g \in R} \int_{g\tau 0}^{g\tau \infty} f_1(z)\, dz \int_{g0}^{g\infty} \overline{f_2(z)\, dz} - \int_{g0}^{g\infty} f_1(z)\, dz \int_{g\tau 0}^{g\tau \infty} \overline{f_2(z)\, dz}.$$

*Démonstration.* — Posons $\omega_1 = f_1(z)\, dz$ et $\omega_2 = f_2(z)\, dz$. Pour $g$ dans $\mathrm{SL}_2(\mathbf{Z})$, posons $\omega_{i|g} = f_{i|g}\, dz$ ($i \in \{1, 2\}$). Considérons le domaine fondamental $D_0$ pour $\mathrm{SL}_2(\mathbf{Z})$ constitué par le triangle hyperbolique de sommets $\infty$, $0$ et $\rho$. On a

$$<f_1, f_2> = \frac{1}{2i[\mathrm{SL}_2(\mathbf{Z}):\Gamma]} \int_{D_\Gamma} \omega_1 \wedge \overline{\omega_2} = \frac{1}{2i[\mathrm{SL}_2(\mathbf{Z}):\Gamma]} \sum_{g \in R} \int_{D_0} \omega_{1|g} \wedge \overline{\omega_{2|g}}.$$

Posons $F_g(z) = \int_\rho^z \overline{f_{2|g}(u)\, du}$. On a $df_{1|g}F_g(z)\, dz = \omega_1 \wedge \overline{\omega_2}$. Cela donne, par la formule de Stokes,

$$2i[\mathrm{SL}_2(\mathbf{Z}):\Gamma]<f_1, f_2> = \sum_{g \in R} \int_{\partial D_0} f_{1|g} F_g(z)\, dz$$

$$= \sum_{g \in R} \int_\infty^0 f_{1|g} F_g(z)\, dz + \int_0^\rho f_{1|g} F_g(z) + \int_0^\rho f_{1|g} F_g(z).$$

Utilisons que $\sigma$ est d'ordre 2 dans $\mathrm{PSL}_2(\mathbf{Z})$ et qu'on a $\tau\rho = \rho$ et $\tau\infty = 0$. Cela donne

$$2i[\mathrm{SL}_2(\mathbf{Z}):\Gamma]<f_1, f_2> =$$

$$\frac{1}{2}\sum_{g \in R} \int_\infty^0 f_{1|g} F_g(z)\, dz + \int_\infty^0 f_{1|g\sigma} F_{g\sigma}(z)\, dz + \sum_{g \in R} \int_\infty^\rho f_{1|g} F_g(z)\, dz + \int_\rho^i nfty f_{1|g\tau} F_{g\tau}(z)\, dz.$$

Utilisons la relation $F_{gh}(hz) = \int_{h^{-1}\rho}^z \overline{f_{2|g}(u)\, du}$. On obtient

$$2i[\mathrm{SL}_2(\mathbf{Z}):\Gamma]<f_1, f_2> = \frac{1}{2}\sum_{g \in R} \int_\infty^0 \omega_{1|g} \int_{\sigma\rho}^\rho \overline{\omega_{2|g}} + \int_\rho^\infty \omega_{1|g\tau} \int_{\tau^2\rho}^\rho \overline{\omega_{2|g}}.$$

Le dernier terme est nul. Décomposons le deuxième facteur du premier terme. On a

$$2i[\mathrm{SL}_2(\mathbf{Z}):\Gamma]<f_1,f_2>=\frac{1}{2}\sum_{g\in R}\int_0^\infty \omega_{1|g}(\int_{\sigma\rho}^{\sigma i}\overline{\omega_{2|g}}+\int_{\sigma i}^\rho \overline{\omega_{2|g}}).$$

Comme $\sigma i = i$, et comme $\int_0^\infty \omega_{1|g}(\int_{\sigma\rho}^{\sigma i}\overline{\omega_{2|g}}) = \int_0^\infty \omega_{1|g\sigma}(\int_i^\rho \overline{\omega_{2|g\sigma}})$, on a la première formule du théorème.

Démontrons maintenant la deuxième formule. On a

$$2i[\mathrm{SL}_2(\mathbf{Z}):\Gamma]<f_1,f_2>=\sum_{g\in R}\int_0^\infty \omega_{1|g}\int_i^\infty\overline{\omega_{2|g}}-\int_0^\infty \omega_{1|g}\int_\rho^\infty\overline{\omega_{2|g}}.$$

Calculons séparément les deux séries de termes. On a

$$\sum_{g\in R}\int_0^\infty \omega_{1|g}\int_i^\infty\overline{\omega_{2|g}}=\frac{1}{2}\sum_{g\in R}\int_0^\infty \omega_{1|g}\int_i^\infty\overline{\omega_{2|g}}+\int_0^\infty \omega_{1|g\sigma}\int_i^\infty\overline{\omega_{2|g\sigma}},$$

et comme $\sigma\infty = 0$,

$$\sum_{g\in R}\int_0^\infty \omega_{1|g}\int_i^\infty\overline{\omega_{2|g}}=\frac{1}{2}\sum_{g\in R}\int_0^\infty \omega_{1|g}\int_0^\infty\overline{\omega_{2|g}}.$$

Par ailleurs, on a

$$\int_0^\infty \omega_{1|g}\int_\rho^\infty\overline{\omega_{2|g}}=\frac{1}{3}\sum_{g\in R}\int_0^\infty \omega_{1|g}\int_\rho^\infty\overline{\omega_{2|g}}+\int_0^\infty \omega_{1|g\tau}\int_\rho^\infty\overline{\omega_{2|g\tau}}+\int_0^\infty \omega_{1|g\tau^2}\int_\rho^\infty\overline{\omega_{2|g\tau^2}}.$$

Remarquons qu'on a $\int_0^\infty \omega_{1|g}+\int_0^\infty \omega_{1|g\tau}+\int_0^\infty \omega_{1|g\tau^2}=0$. C'est pourquoi on a

$$\int_0^\infty \omega_{1|g}\int_\rho^\infty\overline{\omega_{2|g}}=\frac{1}{3}\sum_{g\in R}\int_0^\infty \omega_{1|g\tau}(\int_\rho^\infty\overline{\omega_{2|g\tau}}-\int_\rho^\infty\overline{\omega_{2|g}})+\int_0^\infty \omega_{1|g\tau^2}(\int_\rho^\infty\overline{\omega_{2|g\tau^2}}-\int_\rho^\infty\overline{\omega_{2|g}}).$$

Comme $(\int_\rho^\infty\overline{\omega_{2|g\tau}}-\int_\rho^\infty\overline{\omega_{2|g}}) = -\int_0^\infty\overline{\omega_{2|g}}$ et comme $\int_\rho^\infty\overline{\omega_{2|g\tau^2}}-\int_\rho^\infty\overline{\omega_{2|g}}=\int_0^\infty\overline{\omega_{2|g\tau^2}}$, on obtient, en posant $\lambda_i(g)=\int_0^\infty \omega_i$ ($i\in\{1,2\}$),

$$\int_0^\infty \omega_{1|g}\int_i^\rho\overline{\omega_{2|g}}=$$

$$\frac{1}{2}\sum_{g\in R}\lambda_1(g)\overline{\lambda_2(g)}+\frac{1}{3}\sum_{g\in R}\lambda_1(g\tau)\overline{\lambda_2(g)}-\frac{1}{3}\sum_{g\in R}\lambda_1(g)\overline{\lambda_2(g)}=\frac{1}{6}\sum_{g\in R}\lambda_1(g)\overline{\lambda_2(g)}+\frac{1}{3}\sum_{g\in R}\lambda_1(g\tau)\overline{\lambda_2(g)}.$$

En utilisant la relation $\lambda_1(g)+\lambda_1(g\tau)+\lambda_1(g\tau^2)=0$, on obtient finalement

$$\int_0^\infty \omega_{1|g}\int_i^\rho\overline{\omega_{2|g}}=\frac{1}{6}\sum_{g\in R}\lambda_1(g\tau)\overline{\lambda_2(g)}-\frac{1}{6}\sum_{g\in R}\lambda_1(g\tau)\overline{\lambda_2(g\tau^2)}.$$

Cela achève la démonstration.

COROLLAIRE . — *Supposons qu'on ait $\Gamma=\Gamma_1(N)$, on a*

$$<f_1,f_2>=\frac{i}{12[\mathrm{SL}_2(\mathbf{Z}):\Gamma_1(N)]}\sum_{(u,v)\in(\mathbf{Z}/N\mathbf{Z})^2}\xi_{f_1}(v,-u-v)\overline{\xi_{f_2}(u,v)}-\xi_{f_1}(u,v)\overline{\xi_{f_2}(v,-u-v)}.$$

*Démonstration.* — C'est une application directe de la deuxième formule du théorème C, en tenant compte de la formule ($i \in \{1,2\}$)

$$\xi_{f_i}(u,v) = -i \int_{g0}^{g\infty} f_i(z)\,dz,$$

où $g = \begin{pmatrix} a & b \\ c & d \end{pmatrix} \in \mathrm{SL}_2(\mathbf{Z})$ vérifie $(c,d) \in (u,v)$.

## 5. La fonction $L$ du carré tensoriel de $f$

Supposons la forme modulaire $f$ de caractère de nebentypus trivial et de niveau $N$ premier. Considérons la série de Dirichlet

$$L(f \otimes f, s) = \sum_{n=1}^{\infty} \frac{a_n^2}{n^s}.$$

Elle admet un prolongement méromorphe au plan complexe, avec un pôle simple en $s=2$ et on a [P,S] (notre normalisation pour le produit scalaire de Petersson est en rapport $\pi/3 = \mathrm{vol}(\mathrm{SL}_2(\mathbf{Z})\backslash \mathbf{H})$ avec celle de [S])

$$\mathrm{Res}_{s=2} L(f \otimes f, s) = 12\pi <f,f>.$$

THÉORÈME D. — *On a*

$$\mathrm{Res}_{s=2} L(f \otimes f, s) = \frac{2\pi i}{(N+1)(N-1)^2} \sum_{\chi,\chi',\chi\chi'(-1)=-1} \frac{\Lambda(f \otimes \chi', 1)\Lambda(f \otimes \chi, 1)}{\tau(\chi\chi')},$$

*où $\chi$ et $\chi'$ parcourent les caractères de Dirichlet primitifs de conducteur $N$ tels que $\chi\chi'$ soit impair.*

*Démonstration.* — Partons de la formule donnée dans le corollaire du théorème C. On a, puisque $N$ est premier, $[\mathrm{SL}_2(\mathbf{Z}) : \Gamma_1(N)] = N^2 - 1$ et donc

$$<f,f> = \frac{i}{12(N^2-1)} \sum_{(u,v) \in (\mathbf{Z}/N\mathbf{Z})^2} \xi_f(v,-u-v)\overline{\xi_f(u,v)} - \xi_f(u,v)\overline{\xi_f(v,-u-v)}.$$

Venons-en à la fonction $\xi_f$. Puisque le nebentypus de $f$ est trivial, la fonction $\xi_f$ est homogène. C'est pourquoi on pose, pour $u/v \in \mathbf{P}_1(\mathbf{Z}/N\mathbf{Z}) = \mathbf{Z}/N\mathbf{Z} \cup \{\infty\}$,

$$\xi_f(u/v) = \xi_f(u,v).$$

Cela permet d'écrire

$$<f,f> = \frac{i}{12(N+1)} \sum_{x \in \mathbf{P}_1(\mathbf{Z}/N\mathbf{Z})} \xi_f(-\frac{1}{x+1})\overline{\xi_f(x)} - \xi_f(x)\overline{\xi_f(-\frac{1}{x+1})}.$$

Utilisons la relation de Manin $\xi_f(x) + \xi_f(-1/x) = 0$ ($x \in \mathbf{P}_1(\mathbf{Z}/N\mathbf{Z})$). On obtient

$$<f,f> = \frac{i}{12(N+1)} \sum_{x \in \mathbf{P}_1(\mathbf{Z}/N\mathbf{Z})} \xi_f(x)\overline{\xi_f(x+1)} - \xi_f(x+1)\overline{\xi_f(x)}.$$

Le terme correspondant à $x = \infty$ est nul dans la somme qui précède. Par application des relations de Manin, on a les relations $\xi_f(1/0) + \xi_f(0/1) = 0$ et $\xi_f(1/0) + \xi_f(0/1) + \xi_f(-1/1) = 0$. On en déduit $\xi_f(-1) = 0$ et $\xi_f(1) = 0$. C'est pourquoi on a également la nullité des termes correspondant à $x = 0$ et $x = -1$. On a donc, si on ne conserve que les termes correspondant à $x \neq 0$, $x \neq 1$ et $x \neq \infty$ et si on change de plus $x$ en $-x$,

$$<f,f> = \frac{i}{12(N+1)} \sum_{x \in (\mathbf{Z}/N\mathbf{Z})^* - \{1\}} \xi_f(-x)\overline{\xi_f(1-x)} - \xi_f(1-x)\overline{\xi_f(-x)}.$$

Comme $N$ est premier, rendons plus explicite le théorème A. On a, pour $x \in (\mathbf{Z}/N\mathbf{Z})^*$, $N' = N$ et on peut choisir $S = \emptyset$ et $\bar{S} = \{N\}$. Le terme associé au caractère $\chi$ dans la formule du théorème A est nul lorsque $\chi = 1$, car $N$ est spécial pour $f$ et on a que $\Lambda^{[\frac{\{N\}}{\emptyset}, \frac{\{N\}}{\{N\}}]}(f,1)$ est multiple de $(1 - a_N)\Lambda(f,1)$ qui est nul (en effet $a_N = -w(f)$ et $\Lambda(f,1) = 0$ si $w(f) \neq -1$). Lorsque $\chi$ est non trivial, le terme associé à $\chi$ se simplifie car on a $\chi_S = 1$, $Q_{p,f,\chi} = 1$, $N_\chi = N^2$, $L_p(f \otimes \chi, X) = 1$. On obtient, pour $x \in \mathbf{Z}/N\mathbf{Z}^*$,

$$\xi_f(x) = \frac{w(f)}{N-1} \sum_{\chi \neq 1} \frac{\tau(\bar{\chi})}{N} \chi(x) \Lambda(f \otimes \chi, 1)$$

où la somme porte sur les caractères de Dirichlet primitifs de conducteur $N$.

Comme $f$ est de nebentypus trivial, on a

$$\overline{\Lambda(f \otimes \chi, 1)} = \Lambda(f \otimes \bar{\chi}, 1).$$

On en déduit, en utilisant que $|w(f)| = 1$,

$$<f,f> = \frac{i}{12(N+1)(N-1)^2 N^2} \sum_{x \in (\mathbf{Z}/N\mathbf{Z})^* - \{1\}} \sum_{\chi, \chi'} \tau(\bar{\chi}) \overline{\tau(\bar{\chi}')} \Lambda(f \otimes \chi, 1) \Lambda(f \otimes \bar{\chi}', 1)) F(\chi, \chi'),$$

où

$$F(\chi, \chi') = \sum_{x \in (\mathbf{Z}/N\mathbf{Z})^* - \{1\}} (\chi(-x)\bar{\chi}'(1-x) - \chi(1-x)\bar{\chi}'(-x)).$$

Changeons $\chi'$ en $\bar{\chi}'$ dans la somme. On obtient, en tenant compte de l'égalité $\overline{\tau(\chi')} = \chi'(-1)\tau(\bar{\chi}')$,

$$<f,f> = \frac{i}{12(N+1)(N-1)^2 N^2} \sum_{\chi, \chi'} \tau(\bar{\chi}) \tau(\bar{\chi}') \chi'(-1) \Lambda(f \otimes \chi', 1) \Lambda(f \otimes \chi, 1) F(\chi, \bar{\chi}')$$

Cette somme est antisymétrique (resp. symétrique) en $\chi$ et $\chi'$ lorsque $\chi'(-1) = \chi(-1)$ (resp. $\chi'(-1) = -\chi(-1)$). C'est pourquoi on a

$$<f,f> = \frac{i}{6(N+1)(N-1)^2 N^2} \sum_{\chi, \chi', \chi\chi'(-1)=-1} \tau(\bar{\chi}) \tau(\bar{\chi}') \chi'(-1) \Lambda(f \otimes \chi', 1)) \Lambda(f \otimes \chi, 1)) E(\chi, \chi').$$

où $E(\chi, \chi') = \sum_{x \in (\mathbf{Z}/N\mathbf{Z})^* - \{1\}} \chi(-x)\chi'(1-x)$. On reconnaît dans cette dernière expression une somme de Jacobi donnée par la formule

$$E(\chi, \chi') = \chi(-1) \frac{\tau(\chi)\tau(\chi')}{\tau(\chi\chi')}.$$

On obtient donc, en tenant compte de l'imparité de $\chi\chi'$,

$$<f,f> = \frac{i}{6(N+1)(N-1)^2 N^2} \sum_{\chi, \chi', \chi\chi'(-1)=-1} \frac{\tau(\chi)\tau(\chi')\tau(\bar{\chi})\tau(\bar{\chi}')}{\tau(\chi\chi')} \Lambda(f \otimes \chi', 1)) \Lambda(f \otimes \chi, 1))$$

Utilisons les identités $\tau(\chi)\tau(\bar{\chi}) = \chi(-1)N = -\chi'(-1)N = -\tau(\chi')\tau(\bar{\chi}')$. On obtient

$$<f,f> = \frac{i}{6(N+1)(N-1)^2} \sum_{\chi, \chi', \chi\chi'(-1)=-1} \frac{\Lambda(f \otimes \chi', 1)\Lambda(f \otimes \chi, 1)}{\tau(\chi\chi')}.$$

Cela donne la formule annoncée lorsque l'on combine avec l'identité rappelée avant l'énoncé du théorème D.

## Bibliographie